\documentclass[AMA,STIX1COL,doublespace]{WileyNJD-v2}
\articletype{Research Article}

\usepackage{amsmath}
\usepackage{amsthm}
\usepackage{amsfonts}
\usepackage{graphicx}
\usepackage{float}
\usepackage{bm}
\usepackage{color,colordvi}

\usepackage{caption}
\usepackage{subcaption}
\usepackage{tikz}
\usepackage{comment}

\newcommand{\fref}[1]{Figure~\ref{#1}}

\newcommand{\sref}[1]{Section~\ref{#1}}
\newcommand{\vm}[1]{\bm{#1}}
\newcommand{\vx}{\vm{x}}

\newcommand{\proj}[1]{\Pi^{\varepsilon} #1}

\newcommand{\diverge}[1]{\nabla \cdot #1}

\newcommand{\symgrad}[1]{\nabla_s #1}

\newcommand{\Nmatrix}{\bm{N}^{\partial E}}
\newcommand{\Projbeta}{\Pi_{\beta}\vm{\sigma}}
\newcommand{\Projbetavect}{\overline{\Projbeta}}
\newcommand{\symP}{\mathbb{P}}
\newcommand{\symPspace}{[\mathbb{P}_1(E)]^{2\times 2}_{\text{sym}}}

\newcommand{\veps}{\vm{\varepsilon}}
\newcommand{\vsigma}{\vm{\sigma}}
\newcommand{\vd}{\vm{d}}
\newcommand{\Erot}{E^{\prime}}

\newcommand{\ac}[1]{{#1}}  
\newcommand{\acrev}[1]{{#1}}

\graphicspath{ {./figs/} }

\begin{document}

\title{Stress-hybrid virtual element method on  
       quadrilateral meshes for compressible and 
       nearly-incompressible linear elasticity}
\author[1]{Alvin Chen}
\author[2]{N. Sukumar$\mbox{}^{*,}$}
	
\authormark{ALVIN CHEN \sc{and} N. SUKUMAR}
	
\address[1]{\orgdiv{Department of Mathematics},
	\orgname{University of California},
	\orgaddress{Davis, CA 95616}, \country{USA}}

\address[2]{\orgdiv{Department of Civil and Environmental Engineering},
	\orgname{University of California},
	\orgaddress{Davis, CA 95616}, \country{USA}}

\corres{$\mbox{}^*$N. Sukumar, Department of Civil and 
	Environmental Engineering,
	University of California, One Shields Avenue, Davis, 
	CA 95616, USA \\
	\email{nsukumar@ucdavis.edu}}

\abstract{
In this article, we propose a robust 
low-order stabilization-free virtual element method on
quadrilateral meshes for linear elasticity
that is based on
the stress-hybrid principle. We refer to this approach as the Stress-Hybrid Virtual Element
Method (SH-VEM).  In this method, the Hellinger--Reissner variational principle is adopted, wherein both the equilibrium equations and the 
strain-displacement relations are variationally
enforced.  We consider small-strain deformations of linear 
elastic solids in the compressible and near-incompressible regimes over
quadrilateral (convex and nonconvex) meshes. 
Within an element, the displacement field is approximated as a linear
combination of canonical shape functions that are \emph{virtual}.  The stress field, similar
to the stress-hybrid finite element method of Pian and Sumihara, is represented using a
linear combination of symmetric tensor polynomials.  A 5-parameter expansion of the stress field is used 
in each element, with
stress transformation equations applied on
distorted quadrilaterals. 
In the variational statement of the strain-displacement relations, 
the divergence theorem is invoked to express the stress coefficients in terms of the nodal displacements. 
This results in a formulation with solely the 
nodal displacements as unknowns.
Numerical results are presented for several benchmark
problems from linear elasticity. 
We show that SH-VEM is free of volumetric and shear locking, and it converges optimally in the $L^2$ norm and energy
seminorm of the displacement field, and in the $L^2$ norm of the hydrostatic stress.
}


\maketitle

\section{Introduction}\label{sec:intro}
Finite element formulations that are robust (do not suffer from volumetric and shear locking) for solid
continua remain a long-standing problem in computational mechanics.~\cite{Ainsworth:2022:cmame} Low-order, fully integrated displacement-based finite elements are prone to volumetric locking as the Poisson's ratio $\nu \to 0.5$, and for bending-dominated problems, spurious shear strains (element tends to be overly stiff) lead to shear locking phenomenon. 
The advent of the virtual element method~\cite{basicprinciple,elasticdaveiga} has provided new routes to potentially alleviate locking for nearly-incompressible materials.  Initially, mixed variational principles, hybrid formulations, B-bar
and selective reduced integration strategies that are prominent in finite element formulations for constrained problems have been adopted in the virtual element method.~\cite{Artioli:2018:cmame,Caceres:2019:anm,Park:2020:meccanica,Artioli:2020:M2AN,Dassi:2021:m3as} 
\acrev{B{\"o}hm et al.~\cite{Böhm:2023:MVE} has provided a study of different virtual element methods for incompressible problems and compared the results to classical finite element techniques.}
More recently, in the spirit of assumed-strain methods,~\cite{Simo:1986:jam,Simo:1990:ijnme} projections onto higher order strains have been pursued in the VEM to devise stabilization-free schemes.\cite{berrone2021lowest,enhanced:VEM,Chen:2023:SFV,Chen:2023:SFS,Lamperti:2023:cmech} 
\acrev{
On a quadrilateral, the
stabilization-free virtual element method (SF-VEM)\cite{Chen:2023:SFV} with projection onto an affine
strain field (nine parameters)
suffers from volumetric locking in the near-incompressible limit.} 
In this article, as a point of departure,
we appeal to the Hellinger--Reissner variational principle and assumed stress (referred to as hybrid stress or stress hybrid) techniques~\cite{Pian:1984:ijnme,Pian:1988:ijnme} to devise a virtual element formulation that is robust for compressible and nearly-incompressible linear elasticity over quadrilateral meshes.  We refer to this approach as the Stress-Hybrid Virtual Element Method (SH-VEM).

Stress-hybrid finite element methods can be traced to 
the seminal work of Pian and Sumihara.\cite{Pian:1984:ijnme}
The relationship of this approach to the method of
incompatible modes~\cite{Wilson:1973:IDM} has been established.\cite{Pian:1986:RBI,Yeo:1996:ijnme} Furthermore, limit
principles exist that reveal the connections of finite
elements that are based on the
Hellinger--Reissner (two-field) and the Hu--Washizu (three-field) 
variational principles.\cite{Stolarski:1987:LPM} 
In the stress-hybrid finite element approach, 
the displacement field 
and the stress field are independent. 
The displacement field is conforming and the
stress field (discontinuous) is approximated by a linear combination
of symmetric tensor polynomials on the biunit square (parent element). Since a quadrilateral element has eight displacement degrees of freedom and three rigid-body modes, a minimum of five parameters are needed for the stress approximation.\cite{Pian:1984:ijnme} On the parent element, the stress field is approximated using a five-term 
expansion (referred to as $5\beta$), and the isoparametric map is used to obtain the expressions for the stress field in the physical $xy$-coordinate system. The
 equilibrium equations and strain-displacement relations are variationally enforced, and in so doing, the stress coefficients can be expressed in terms of the nodal displacements. This leads to a displacement formulation that solely contains nodal displacements as unknowns. 
 Belytschko and Bindeman~\cite{Belytschko:1991:cmame} have shown that for linear problems, the stress-hybrid 
 finite element method produces more accurate results than many assumed-strain and hourglass stabilized methods. Barlow~\cite{Barlow:1986:ijnme} provides a rationale for the improved accuracy that the stress-hybrid method delivers by proving via a minimization principle that it is an optimal compromise between the displacement-based (stiffness) and stress-based (flexibility) formulations. For nearly-incompressible problems, the finite element method based on the former is stiff, whereas use of the latter tends to be overly flexible. There have been many different approaches within the stress-hybrid finite element
 method
 to find the best choice of stress expansion to ensure good accuracy and robustness. Several studies~\cite{Pian:1988:ijnme,Punch:1984:cmame,Lee:1986:ijnme,Spilker:1981:ijnme} select a judicious choice of the stress expansion so that
 spurious zero-energy modes are excluded. Alternative approaches by Xie and Zhou~\cite{Xie:2008:cnme,Xie:2004:ijnme} and Cen et al.~\cite{Cen:2011:cmame} have used energy compatibility arguments to derive the stress expansion. Jog~\cite{Jog:2010:mom,Jog:2005:fead} 
 has provided rules for choosing the stress functions and has 
 extended the stress-hybrid approach to higher order elements.
 Xue et al.\cite{Xue:1985:ijss} provide a 
 necessary inf-sup condition to ensure the stability of stress-hybrid 
 methods, and Zhou and Xie~\cite{Zhou:2002:cmame} and 
 Yu et al.\cite{Yu:2011:cmame} have analyzed the convergence of these
 methods. 
 Simo et al.\cite{Simo:1989:cmame} showed the excellent performance of the
 stress formulation for linear shell theory and
 Simo et al.\cite{Simo:1989:CMF} have applied the approach for small-strain
 elastoplasticity.  However, very few studies have explored 
 the stress-hybrid method for nonlinear computations. Notable among them are the contributions of Jog and 
coworkers,\cite{Jog:2006:ijnme,Jog:2009:jmms,Jog:2009:SSA,Jog:2017:MHF,Agrawal:2019:HFE}
who have successfully 
extended the stress-hybrid finite element formulation to nonlinear and 
inelastic analysis of solid continua, plates and shells. 
 
In this article, we adopt the stress-hybrid principle using a five-parameter expansion of the stress field to construct a robust
stabilization-free 
virtual element method for
plane elasticity on quadrilateral (convex and nonconvex) meshes that is free of volumetric and shear locking. 
\acrev{To form the stress-hybrid element as a virtual element method, we start from the Hellinger--Reissner variational principle 
with the weak statement of the equilibrium equations and
the strain-displacement relations. In the virtual element method, there is no explicit construction of the \ac{(polygonal)} basis functions, and therefore a standard procedure is to use projection operators to approximate the different fields. 
The  weak enforcement of the strain-displacement equations in tandem with an assumed stress basis is used to define a projection operator for the stress, while the
bilinear form from the equilibrium equations define the energy
projection for the displacements. The use of the Hellinger--Reissner
variational principle and a separately assumed stress 
% distribution 
ansatz makes the method less prone to shear and volumetric locking.
% Note that on a quadrilateral,
% a linear symmetric tensor polynomial for a higher
% order strain
% field (nine parameters) is used in the
% stabilization-free virtual element method (SF-VEM),\cite{Chen:2023:SFV}
% which proves to be overly stiff and suffers from volumetric locking for
% nearly-incompressible materials.
As noted earlier, projection onto a linear strain space over a
quadrilateral in
the SF-VEM~\cite{Chen:2023:SFV}  proves to be overly stiff and suffers from volumetric locking for nearly-incompressible materials.}
The SH-VEM proposed herein overcomes this limitation of the
SF-VEM.

The remainder of this article is structured as follows. In~\sref{sec:formulation}, we introduce the model problem of linear elasticity and construct a weak formulation from the Hellinger--Reissner variational principle. In \sref{sec:VEM}, we present the virtual element space and construct the required projection operators for the SH-VEM. As a point of distinction from the 
stress-hybrid finite element method, in the SH-VEM we use a $5\beta$ expansion of the stress field in a rotated coordinate system for distorted quadrilaterals,\cite{Cook:1974:jsd,Cook:1975:APS} and then the stress transformation equations are used to transform them to Cartesian 
coordinates.  Furthermore, we then invoke the divergence theorem over the distorted element (since the basis functions in the VEM are \emph{virtual}) to express the
stress coefficients in terms of the nodal displacements.
In~\sref{sec:implementation}, the numerical implementation of the projection operators and the stiffness matrix is given. 
In~\sref{sec:numericalresults}, we first present numerical results for standard VEM,\cite{elasticdaveiga} stabilization-free VEM,\cite{Chen:2023:SFV}
B-bar VEM~\cite{Park:2020:meccanica} and SH-VEM for a manufactured problem with
Poisson's ratio $\nu \in [0.3,0.5)$.  Subsequently, results are presented for
B-bar VEM and SH-VEM on a series of benchmark problems for $\nu$ close
to $0.5$: thin cantilever beam (to demonstrate absence of shear locking), Cook's membrane under shear end load, infinite plate with a hole under uniaxial tension, hollow pressurized cylinder and flat punch.  Our main finding from this study are summarized in~\sref{sec:conclusions}.

 \section{Hellinger--Reissner Variational Formulation}\label{sec:formulation}
Consider an elastic body that occupies the region $\Omega \subset \mathbb{R}^2$ with boundary $\partial \Omega$. We assume that the boundary can be written as a disjoint union $\partial \Omega = \Gamma_u \cup \Gamma_t$ 
with $\Gamma_u \cap \Gamma_t =\emptyset$. Displacement boundary conditions
$\vm{u} = \vm{u}_0$ are prescribed on 
$\Gamma_u$ and tractions $\vm{t} = \bar{\vm{t}}$ are imposed on 
$\Gamma_t$.
Let $\vm{u}$ denote the displacement field, 
$\veps = \nabla_s \vm{u}$ ($\nabla_s$ is the symmetric gradient
operator) is the small-strain
tensor, $\vsigma$ is the Cauchy stress tensor,
and $\vm{b} \in [L^2(\Omega)]^2$ is the body force per volume.

To construct the weak form, we start from the Hellinger--Reissner functional for linear elasticity, which is given by:
\begin{align*}
    {\Pi}_{\textrm{HR}}[\vm{u},\vsigma] = -\frac{1}{2}\int_{\Omega}{\vsigma: \mathbb{C}^{-1}:\vsigma \, d\vm{x}}  + \int_{\Omega}{\vsigma:\symgrad{\vm{u}}\, d\vm{x}} - \int_{\Omega}{\vm{b}\cdot \vm{u} \, d\vm{x}} - \int_{\Gamma_t}{\Bar{\vm{t}}\cdot \vm{u} \, ds} .
\end{align*}
On taking the first variation of $\Pi_{\textrm{HR}}(\cdot,\cdot)$ and 
requiring it to be stationary, we obtain
\begin{align*}
   \delta {\Pi}_{\textrm{HR}}[\vm{u},\vsigma;\delta \vm{u},\delta \vsigma]
   =  \int_{\Omega} \! \delta\vsigma : \left(\symgrad{\vm{u}}-\mathbb{C}^{-1}:\vsigma \right) \, d\vm{x} +\int_{\Omega} \! \vsigma : \symgrad{(\delta\vm{u}) \, d\vm{x}} - \int_{\Omega} \! \vm{b} \cdot \delta \vm{u} \, d\vm{x} - \int_{\Gamma_t} \! \Bar{\vm{t}}\cdot \delta\vm{u} \, ds =0 \ \ \forall 
    \delta \vm{u} \in {\cal V}_u, \ \delta \vsigma \in {\cal V}_\sigma, 
\end{align*}
where ${\cal V}_u$ contains vector-valued functions in the Hilbert space
$[H^1(\Omega)]^2$ that also vanish on $\Gamma_u$, whereas
${\cal V}_\sigma$ contains functions in $(L^2)_{\textrm{sym}}^{2 \times 2}$.
This gives us the weak statement of the equilibrium equations and strain-displacement relations:
\begin{subequations}\label{eq:weak_equations}
\begin{align}
    &\int_{\Omega}{\vsigma : \symgrad{(\delta\vm{u})} \, d\vm{x}} - \int_{\Omega}{\vm{b} \cdot \delta \vm{u} \, d\vm{x}} - \int_{\Gamma_t}{\Bar{\vm{t}}\cdot \delta\vm{u} \, ds} =0 \ \ \forall 
    \delta \vm{u} \in {\cal V}_u, \label{eq:weak_equilibrium} \\ 
     &\int_{\Omega}{\delta\vsigma : \left(\nabla_s\vm{u}-\mathbb{C}^{-1}:\vsigma\right) \, d\vm{x}} = 0 \ \ \forall \delta \vsigma \in {\cal V}_\sigma . \label{eq:weak_compatibility}
\end{align}
\end{subequations}

 \section{Virtual Element Discretization}\label{sec:VEM}
 Let $\mathcal{T}^h
 $ be a decomposition of $\Omega$ into nonoverlapping quadrilaterals. For each quadrilateral $E\in \mathcal{T}^h$, let $h_E$
 denote its diameter, $\vx_E$ its centroid, and $\vx_i = (x_i,y_i)$ 
the coordinate of the $i$-th vertex.

 \subsection{Polynomial basis }
In the virtual element method, we need a polynomial 
approximation space for the displacement field. For each 
element $E$, define $[\mathbb{P}_1(E)]^2$ as the space of 
two-dimensional vector-valued polynomials of degree less 
than or equal to one. For this space, we choose a basis of 
scaled vector monomials given by:
\begin{subequations}
    \begin{align}
        \widehat{\vm{M}}(E) &= \begin{bmatrix}
    \begin{Bmatrix}
    1 \\ 0 
    \end{Bmatrix}, 
    \begin{Bmatrix}
    0 \\ 1
    \end{Bmatrix}, 
    \begin{Bmatrix}
    -\eta \\ \xi
    \end{Bmatrix}, 
    \begin{Bmatrix}
    \eta \\ \xi
    \end{Bmatrix}, 
    \begin{Bmatrix}
    \xi \\ 0
    \end{Bmatrix}, 
    \begin{Bmatrix}
    0 \\ \eta
    \end{Bmatrix}
        \end{bmatrix}\label{eq:disp_poly_basis},\\
        \intertext{where}
        \xi &= \frac{x-x_E}{h_E}, \quad \eta = \frac{y-y_E}{h_E}.
    \end{align}
\end{subequations}
We denote the $\alpha$-th vector of $\widehat{\vm{M}}(E)$ by $\vm{m}_\alpha$.
We also need to define an approximation space for the stress field. 
Using Voigt representation, we define a basis for the space $[\mathbb{P}_1(E)]^{2\times 2}_{\text{sym}}$ of $2\times 2$ \acrev{polynomial symmetric tensors} of degree less than or equal to one: 
\acrev{\begin{align}\label{eq:complete_linear_stress}
    \widehat{\vm{M}}^{2\times 2}(E) =  \begin{bmatrix}
    \begin{Bmatrix}
    1 \\ 0 \\ 0
    \end{Bmatrix}, 
    \begin{Bmatrix}
    0 \\ 1 \\ 0
    \end{Bmatrix}, 
    \begin{Bmatrix}
    0 \\ 0 \\ 1
    \end{Bmatrix}, 
    \begin{Bmatrix}
    \xi \\ 0 \\ 0
    \end{Bmatrix}, 
    \begin{Bmatrix}
    0 \\ \xi \\ 0
    \end{Bmatrix}, 
    \begin{Bmatrix}
    0 \\ 0 \\ \xi
    \end{Bmatrix}, 
    \begin{Bmatrix}
    \eta \\ 0 \\ 0
    \end{Bmatrix},
    \begin{Bmatrix}
    0 \\ \eta \\ 0
    \end{Bmatrix},
    \begin{Bmatrix}
    0 \\ 0 \\ \eta
    \end{Bmatrix}
    \end{bmatrix}.
\end{align}}

\subsection{Energy projection of the displacement field}
We first define the energy projection of the displacement field from standard virtual element formulations. For each element $E$, \acrev{denote $H^1(E)$ as the Sobolev space of functions defined on $E$ with square integrable derivatives up to order 1}. \acrev{Now define the operator} $\proj: [H^1(E)]^2 \to [\symP_1(E)]^2$ \acrev{such that for every $\vm{v}\in [H^1(E)]^2 $, $\proj{\vm{v}}$} satisfies the orthogonality relation:
\begin{align}
    \int_{E}{\veps(\vm{m}_\alpha):\mathbb{C}:\veps(\vm{v}-\proj{\vm{v}}) \, d\vx} =0 \quad \forall \vm{m}_\alpha \in \widehat{\vm{M}}(E).
\end{align}
We note that $\veps(\vm{m}_\alpha) =\vm{0}$ for $\alpha =1,2,3$, since these $\vm{m}_{\alpha}$ correspond to rigid-body modes. Therefore, we obtain three trivial equations $0 =0$. In order to define a unique projection, we choose a projection operator $P_0: [H^1(E)]^2 \times [H^1(E)]^2 \to \mathbb{R}$ to replace these three trivial equations. In particular, we choose a discrete $L^2$ inner product:
\begin{subequations}
  \begin{align}
    P_0(\vm{u},\vm{v}) &= \frac{1}{4} \sum_{k=1}^{4}{\vm{u}(\vx_k)\cdot \vm{v}(\vx_k)}, \\
    \intertext{and require three additional conditions:}
    P_0(\vm{m}_{\alpha},\vm{v}-\proj{\vm{v}}) &=\frac{1}{4}\sum_{k=1}^{4}{\vm{m}_{\alpha}}(\vx_k)\cdot (\vm{v}-\proj{\vm{v}})(\vx_k) = 0 \quad (\alpha=1,2,3).
\end{align}  
\end{subequations}
Now we define the energy projection $\proj{\vm{v}}$ as the unique function that satisfies the following two relations:
\begin{subequations}\label{eq:energy_proj_def}
    \begin{align}
    &\frac{1}{4}\sum_{k=1}^{4}{\vm{m}_{\alpha}}(\vx_k)\cdot (\vm{v}-\proj{\vm{v}})(\vx_k) =0 \quad (\alpha = 1,2,3), \label{eq:proj_P0_condtion} \\
    &\int_{E}{\veps(\vm{m}_\alpha):\mathbb{C}:\veps(\vm{v}-\proj{\vm{v}}) \, d\vx} =0 \quad (\alpha = 4,5,6) \label{eq:proj_integral_condition}.
    \end{align}
\end{subequations}
Following Chen and Sukumar,~\cite{Chen:2023:SFV} we write~\eqref{eq:proj_integral_condition} in matrix-vector form as 
\begin{subequations}
\begin{align}
\int_{E}{\left(\bm{S}\proj{\bm{v}}\right)^T\left(\bm{CS}\bm{m}_{\alpha}\right)\,d\bm{x}} &= \int_{E}{\left(\bm{S}\bm{v}\right)^T\left(\bm{CS}\bm{m}_{\alpha}\right)\,d\bm{x}},\label{eq:proj_matrix_condition}\\
    \vm{S} =\begin{bmatrix}
        \frac{\partial}{\partial x} & 0 \\ 0 & \frac{\partial}{\partial y} \\ \frac{\partial}{\partial y} & \frac{\partial}{\partial x}
    \end{bmatrix}, \quad
    \vm{C} &= \frac{E_Y}{(1+\nu)(1-2\nu)} \begin{bmatrix}
         1-\nu &  \nu    & 0 \\
         \nu   & 1 - \nu & 0 \\
         0     &   0     & \frac{1 -2\nu}{2}
        \end{bmatrix}, \label{eq:S_matrix}
        \end{align}
\end{subequations}
where $\vm{C}$ is the matrix representation (plane strain) of the material moduli tensor, and $E_Y$ and $\nu$ are
the Young's modulus and Poisson's ratio, respectively,
of the material.

 \subsection{Stress-hybrid projection operator}
We now define the projection operator for the stress-hybrid formulation. 
On choosing $\delta \vsigma = \symP \in \symPspace$ and $\vsigma \in \symPspace$ in~\eqref{eq:weak_compatibility}, we have
\begin{align*}
    \int_{\Omega}{\mathbb{P} : \left(\nabla_s\vm{u}-\mathbb{C}^{-1}:\vsigma\right) \, d\vm{x}} = 0.
\end{align*}
This condition is true for all $\symP \in \symPspace$, 
so we can view 
$\mathbb{C}^{-1}:\vsigma$ as a projection of $\nabla_s\vm{u}$ with respect to the space $\symPspace$.  
Now, let the assumed stress field be taken as 
$\sigma := \Pi_\beta \sigma$, where $\Pi_\beta$ is the stress projection operator. Then, the orthogonality condition becomes
\begin{subequations}
    \begin{align*}
    &\int_{E}{\symP : \left(\nabla_s\vm{u}-\mathbb{C}^{-1}:\Projbeta\right) \, d\vm{x}} = 0 \quad \forall \symP \in \symPspace, \\
    \intertext{or equivalently}
    &\int_{E}{\symP:\mathbb{C}^{-1}:\Projbeta \, d\vx} = \int_{E}{ \symP : \nabla_s\vm{u}\, d\vx} \quad \forall \symP \in \symPspace.
\end{align*}
\end{subequations}
After applying the divergence theorem and simplifying, we obtain
\begin{align}
    \int_{E}{\symP:\mathbb{C}^{-1}:\Projbeta \, d\vx }= \int_{\partial E}{\left(\symP\cdot \vm{n}\right)\cdot\vm{u} \,ds} - \int_{E}{\left(\diverge{\symP}\right)\cdot \vm{u} \, d\vx}, \label{eq:proj_tensor_equation}
\end{align}
where $\vm{n} = (n_x,n_y)^T$ is the outward unit normal along $\partial E$. For later implementation, we convert this expression into the associated matrix-vector form. Let $\overline{\symP}$, $\overline{\Projbeta}$ be the Voigt representation of $\symP$ and 
$\Projbeta$, respectively. Then, \eqref{eq:proj_tensor_equation} can be written as 

\begin{subequations}
\begin{align}
    \int_{E}{\overline{\symP}^T\vm{C}^{-1}\overline{\Projbeta} \, d\vx} &= \int_{\partial E}{\overline{\symP}^T\Nmatrix\vm{u} \, ds} - \int_{E}{ \left(\vm{\partial}\overline{\symP}\right)^T\vm{u} \, d\vx}\label{eq:matrix_vect_proj} \quad \forall \symP \in \symPspace,\\
    \intertext{where $\Nmatrix$ is the representation of the outward normals and $\vm{\partial}$ is the matrix divergence operator that are given by}
    &\Nmatrix := \begin{bmatrix}
        n_x & 0 \\ 0 & n_y \\ n_y & n_x
    \end{bmatrix}, \quad \vm{\partial} := \begin{bmatrix}
        \frac{\partial}{\partial x} & 0 & \frac{\partial}{\partial y} \\ 0 & \frac{\partial}{\partial y} & \frac{\partial}{\partial x}
    \end{bmatrix}.
\end{align}

\end{subequations}

 \section{Numerical implementation}\label{sec:implementation}
 
 \subsection{Virtual element space}
We define the first order vector-valued virtual element space on an element $E$ as:~\cite{basicprinciple,AHMAD2013376}
\begin{equation}
    \begin{split}
    \vm{V}(E) = \Biggl\{\vm{v}_h \in [H^1(E)]^2 : \Delta \vm{v}_h \in [\mathbb{P}_1(E)]^2, \ \vm{v}_h |_{e} \in [\mathbb{P}_1(e)]^2 \ \forall e \in \partial E, \ \vm{v}_h |_{\partial E} \in [C^0(\partial E)]^2, \\
    \int_{E}{\vm{v}_h \cdot \vm{p} \ d\vx} = \int_{E}{ \proj{\vm{v}_h} \cdot \vm{p} \ d\vx} \quad \forall \vm{p} \in [\mathbb{P}_1(E)]^2   \Biggr\},
    \end{split}
\end{equation}
\acrev{where $\Delta$ is the vector Laplacian operator and $\proj$ is the energy projection operator defined in~\eqref{eq:energy_proj_def}.} On this space, the \acrev{vector-valued functions are continuous on the boundary and affine along each edge $e$}, so we can choose the degrees of freedom (DOFs) to be the function values at the vertices of E. Each element has a total of eight displacement DOFs. For each element $E$ we also assign a basis for the local space $\vm{V}(E)$. Let $\{\phi_i\}$ be the standard scalar basis functions in standard VEM~\cite{basicprinciple} that satisfy the property $\phi_i(\vx_j) = \delta_{ij}$. Using the scalar basis, we define the matrix of vector-valued basis functions by
\begin{subequations}
   \begin{align}
        \vm{\varphi} &= \begin{bmatrix}
    \phi_1 & \phi_2 & \phi_3 &\phi_4 & 0 & 0 & 0 & 0 \\ 
    0 & 0 & 0 & 0 & \phi_1 & \phi_2 & \phi_3 & \phi_4
    \end{bmatrix} := \begin{bmatrix}
        \vm{\varphi_1} & \vm{\varphi_2} & \dots \vm{\varphi_8}
    \end{bmatrix},\\
    \intertext{then any function $\vm{v}_h \in \vm{V}(E)$ can be represented as:}
    \vm{v}_h(\vx) &= \sum_{i=1}^{8}{\vm{\varphi}_{i}(\vx)v_i} = \vm{\varphi}\vd,  
\end{align}
\end{subequations}
where $v_i$ is the $i$-th degree of freedom of $\vm{v}_h$.

 \subsection{Matrix representation of the energy projection}
 To construct the energy projection of the displacement on the virtual element space $\vm{V}(E)$, we follow the construction in Chen and Sukumar.\cite{Chen:2023:SFV} Let $\vm{v}_h = \vm{\varphi}_i$ be the $i$-th basis function $(i=1,2,\dots,8)$, then substituting into~\eqref{eq:proj_P0_condtion} and~\eqref{eq:proj_matrix_condition}, we obtain:
 \begin{subequations}\label{eq:EProj}
    \begin{align}
         \frac{1}{4}\sum_{k=1}^{4}{\vm{m}_{\alpha}(\vx_k)\cdot\proj{\vm{\varphi}_i}(\vx_k)} &= \frac{1}{4}\sum_{k=1}^{4}{\vm{m}_{\alpha}}(\vx_k)\cdot \vm{\varphi}_i(\vx_k) \quad (\alpha=1,2,3), \\
         \int_{E}{\left(\bm{S}\proj{\vm{\varphi}_i}\right)^T\left(\vm{C}\vm{S}\vm{m}_{\alpha}\right)\,d\vx} &= \int_{E}{\left(\vm{S}\vm{\varphi}_i\right)^T\left(\vm{C}\vm{S}\vm{m}_{\alpha}\right)\,d \vx} \quad (\alpha=4,5,6).
 \end{align}  
 \end{subequations}
The projection $\proj{\vm{\varphi}_i}$ is a linear vector polynomial, therefore we expand it in terms of the basis in~\eqref{eq:disp_poly_basis}. That is:
\begin{align*}
    \proj{\vm{\varphi}_i} = \sum_{\mu=1}^{6}{s_{\mu}^{i}\vm{m}_{\mu}} := \sum_{\mu=1}^{6}{(\vm{\Pi}_{*}^{\veps})_{\mu i}\vm{m}_{\mu}} \quad (i=1,2,\dots,8).
\end{align*}
After substituting in~\eqref{eq:EProj} and simplifying, we obtain the energy projection coefficients $\vm{\Pi}_{*}^{\veps}$ as~\cite{Chen:2023:SFV} 
\begin{subequations}
    \begin{align}
    \vm{\Pi}_{*}^{\veps} &= \vm{G}^{-1}\vm{B},\label{eq:Pi_eps_matrix}
    \intertext{where for $\mu=1,2,\dots,6$}
    \vm{G}_{\alpha \mu} &=\begin{cases}
    \frac{1}{4}\sum_{j=1}^{4}{\vm{m}_\mu(\vx_j)\cdot \vm{m}_\alpha(\vx_j)} \quad (\alpha =1,2,3) \\
    \int_{E}{\left(\vm{S}\vm{m}_\mu\right)^T(\vm{C}\vm{S}\vm{m}_{\alpha}) \, d\vx} \quad 
    (\alpha= 4,5,6),
    \end{cases} \\
    \intertext{and for $i=1,2,\dots,8$}
    \vm{B}_{\alpha i} &=
    \begin{cases}
        \frac{1}{4}\sum_{j=1}^{4}{\vm{\vm{\varphi}_i}(\vx_j)\cdot \vm{m}_\alpha(\vx_j)} \quad (\alpha =1,2,3) \\ 
        \int_{E}{\left(\vm{S}\vm{\varphi}_i\right)^T(\vm{C}\vm{S}\vm{m}_{\alpha}) \, d\vx} \quad (\alpha=4,5,6).
    \end{cases}
    \end{align}
\end{subequations}
The matrix $\vm{\Pi}_{*}^{\veps}$ is the coefficients of the energy projection with respect to the polynomial basis $\widehat{\vm{M}}(E)$. For later implementation of the B-bar VEM, it is convenient to define the coefficients with respect to the basis functions $\vm{\varphi}$. Let $\vm{D}$ be the matrix of DOFs of vectors in $\widehat{\vm{M}}(E)$, which is given by
\begin{align}
    \vm{D} = \begin{bmatrix}
        \vm{m}_1(\vx_1) & \vm{m}_2(\vx_1) &\dots & \vm{m}_6(\vx_1)\\ 
        \vm{m}_1(\vx_2) & \vm{m}_2(\vx_2) &\dots & \vm{m}_6(\vx_2)\\
        \dots & \dots & \dots & \dots \\
        \vm{m}_1(\vx_4) & \vm{m}_2(\vx_4) &\dots & \vm{m}_6(\vx_4)
    \end{bmatrix}.
\end{align}
Then the representation of the energy projection with respect to the basis functions $\vm{\varphi}$ is:~\cite{Veiga2014TheHG}
\begin{align}\label{eq:Pi_eps_canonical}
    \vm{\Pi}^{\veps} := \vm{D}\vm{\Pi}_{*}^{\veps}.
\end{align}

\subsection{Matrix representation of the 
stress-hybrid projection}
We present the matrix representation of the
SH-VEM projection, which relies on
the rotated coordinates introduced by Cook~\cite{Cook:1974:jsd} (see Figure~\ref{fig:local_coordinate_construction}).
In the SH-VEM, an assumed stress ansatz is first defined on
$E^\prime$ (rotated element) and then transformed to $E$ using the stress transformation equations. The computation of the element stiffness matrix is carried out on $E$. 

\subsubsection{Rotated coordinates }
\begin{figure}
   \centering
\begin{tikzpicture}[scale=1.5]

            %draws the main shape
			\filldraw[line width=1pt,fill=gray!10] (0,0) -- (2.8978,0.7765) -- (2.1213,3.6742) -- (-.5176, 1.9319) -- cycle;
                %x-y axis
                \draw[line width=1, ->] (1.2176-.075,1.6376-.05) --++ (0,3.15);
                \draw[line width=1,->] (1.2176-.075,1.6376-.05) --++ (3,0);

                %connection of the midpoints
                \draw[line width= 1 , dashed] (1.4489,.3882) -- (0.4,4.3027);
                \draw[line width= 1 , dashed] (-.2588,.9659) -- (2.509+1.5,2.225+.6824);
                
                %theta_2
                \draw[line width=1, <->] (1.2176-.075,1.6376-.05+2) arc (90:90+15:2.0);
                %theta_1 
                 \draw[line width=1, <->] (3.7901,1.6376-.05) arc (0:25:2.60);

                %theta
                \draw[line width=1, <->] (1.2176-.075,1.6376-.05+3) arc (90:90+20:3.0);
                \draw[line width=1, <->] (3.7901-1.5,1.6376-.05) arc (0:20:1);

                %  x' y' axis
                \draw[line width =1.5, ->] (1.2176-.075,1.6376-.05) -- (1.2176-.075+1.8771,1.6376-.050+.6903);
                \draw[line width=1.5, ->](1.2176-.075,1.6376-.05) -- (1.2176-.075-1.0774,1.6376-.05+2.9600);

                %midpoints
                \filldraw [gray] (1.4489,.3882) circle (2pt) ;
                \filldraw[gray] (0.801,2.803) circle (2pt) ; 
                \filldraw[gray](2.509,2.225) circle (2pt) ;
                \filldraw[gray] (-.2588,.9659) circle (2pt) ;

                %label of midpoints
                \draw[] (-.2588-.1,.9659) node[anchor = east]{P} ;
                \draw[] (1.4489,.3882-.1) node[anchor = north]{R}  ;  
                \draw[] (2.6390-.5,1.7424+.7) node[anchor = west]{Q} ;
                \draw[] (.9313-.3+.4,2.3201+.25) node[anchor = south]{S} ;

                %label the angles
                \draw[]  (3.7901,1.3542+.75) node[anchor = west]{$\theta_1$} ;
                \draw[] (1.1901-.3,3.6542) node[anchor = south]{$\theta_2$}  ;  
                \draw[] (1.2176-.075-.75+.2,1.6376-.05+3) node[anchor=south]{$\theta$};
                \draw[] (3.7901-1.5,1.6376+.15)  node[anchor = west]{$\theta$};
                %label the axis
                \draw[] (1.1901,1.3542+3.15)  node[anchor = west]{$y$};
                \draw[] (1.1901+3,1.3542+.15) node[anchor = north]{$x$};
                \draw[] (1.2176-.075+1.8771,1.6376-.050+.6903) node[anchor=north]{$x^{\prime}$};
                \draw[] (1.2176-.075-1.0774-.15,1.6376-.05+2.9600) node[anchor = north]{$y^{\prime}$};
\end{tikzpicture}
    \caption{ Construction of the local coordinate system for a distorted quadrilateral. }
    \label{fig:local_coordinate_construction}
\end{figure}
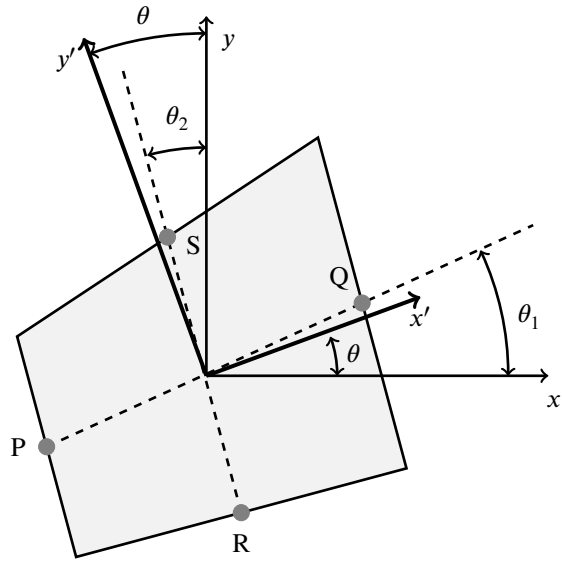

It is known that using global Cartesian coordinates to construct a $5$-term expansion of the stress field leads to an incomplete stress approximation and the resulting element stiffness matrix is not rotationally invariant.~\cite{Pian:1984:ijnme,Cook:1974:jsd,Pilter:1995:ijnme}  We follow the modification proposed by Cook~\cite{Cook:1974:jsd} to construct a local coordinate system for each element $E$. Let $\vx_P, \vx_Q, \vx_R,\vx_S$ be the midpoints of the edges of element $E$ (see Figure~\ref{fig:local_coordinate_construction}).  Define $L_1$ and $L_2$ as the length of the line segments $PQ$ and $RS$, 
respectively. Then, we compute the angles
     \begin{align}\label{eq:thetas}
     \theta_{1} = \arctan{\left(\frac{y_{Q}-y_{P}}{x_{Q}-x_{P}}\right)}, \quad 
     \theta_{2} = \arctan{\left(\frac{x_{R}-x_{S}}{y_{S}-y_{R}}\right)}, \quad
     \theta =\frac{L_{1}\theta_{1}+ L_{2}\theta_{2}}{L_{1}+L_{2}} .
 \end{align} 
 On using the angle $\theta$, we define the rotated coordinates $(x^{\prime},y^{\prime})$ via the transformation
\begin{align}\label{eq:rotated_coords}
      \vx^\prime := 
     \begin{Bmatrix}
         x^\prime \\ y^\prime 
     \end{Bmatrix}= \begin{bmatrix}
            c & s \\ -s & c
         \end{bmatrix}\begin{Bmatrix}
             x \\ y
         \end{Bmatrix}:= \vm{Q}\vx, \ \ c = \cos \theta, \ 
         s = \sin \theta.
 \end{align}

\subsubsection{Stress-hybrid virtual element formulation}
We now present the stress-hybrid approach, which constructs the stress basis functions and stiffness matrix over the
distorted element $E$. \acrev{For a quadrilateral element, using the complete linear stress basis given in~\eqref{eq:complete_linear_stress} results in an overly stiff element; therefore, we only use a five-dimensional subset of this basis.} On a square, the selection
of the stress basis in the stress-hybrid finite element method as~\cite{Pian:1984:ijnme}
\begin{equation*}
\vm{M}_{5\beta}=
    \begin{bmatrix}
        \begin{Bmatrix}
            1 \\ 0 \\ 0
        \end{Bmatrix},
        \begin{Bmatrix}
            0 \\ 1 \\ 0
        \end{Bmatrix},
        \begin{Bmatrix}
            0 \\ 0 \\ 1
        \end{Bmatrix},
        \begin{Bmatrix}
            y \\ 0 \\ 0
        \end{Bmatrix},
        \begin{Bmatrix}
            0 \\ x \\0 
        \end{Bmatrix}
    \end{bmatrix}
\end{equation*}
ensures that uniform stress states as well as pure bending can be exactly represented. To tailor this approach to VEM, we use this $5\beta$ stress expansion in a local coordinate system ($E$ is rotated) and then apply the stress transformation equations to obtain the stress ansatz on $E$. 

Let $\Erot$ be the rotated element with vertices $(x^\prime_i,y^\prime_i)$, centroid $\vx^{\prime}_{\Erot}$ and diameter $h^{\prime}_{\Erot}$. In a rotated element $E^{\prime}$, we assume the stress expansion $\overline{\Projbeta^{\prime}} = \vm{P}^{\prime}\vm{\beta}^{\prime}$, where $\vm{P}^{\prime}$ is given by
 \begin{subequations}
    \begin{align}
    \vm{P}^{\prime} = \begin{bmatrix}
        1 & 0 & 0 & {\eta}^{\prime} & 0 \\ 0 & 1 & 0 & 0 & \xi^{\prime} \\ 0 & 0 & 1 & 0 & 0 
    \end{bmatrix}&=\begin{bmatrix}
        \vm{P}^{\prime}_1 & \vm{P}^{\prime}_2 & \vm{P}^{\prime}_3 & \vm{P}^{\prime}_4 & \vm{P}^{\prime}_5
    \end{bmatrix},\label{eq:Pprime_matrix}\\
    \intertext{and the rotated scaled monomials are}
    {\xi}^{\prime} =\frac{x^{\prime}-x^{\prime}_{\Erot}}{h^{\prime}_{\Erot}} &,\quad \eta^{\prime} =\frac{y^{\prime}-y^{\prime}_{\Erot}}{h^{\prime}_{\Erot}}. 
\end{align}
\end{subequations}
 On viewing each column $\vm{P}_{i}^{\prime}$ of the matrix $\vm{P}^{\prime}$ as an equivalent tensor $\vm{\mathcal{P}}^{\prime}_{i}$, we apply the rotation matrix $\vm{Q}$ given in~\eqref{eq:rotated_coords} to obtain a transformed tensor $\vm{\mathcal{P}}_i$:
 \begin{align}\label{eq:rotated_P_tensor}
     \vm{\mathcal{P}}_i = \vm{Q}^T\vm{\mathcal{P}}_{i}^{\prime}\vm{Q}.
 \end{align}
 After computing each tensor $\vm{\mathcal{P}}_i$ and rewriting them in terms of $3\times 1$ vectors $\vm{P}_i$, we define the matrix $\vm{P}^*$ by~\cite{Cook:1974:jsd}
 \begin{align}\label{eq:P_matrix_full}
     \vm{P}^* = \begin{bmatrix}
         \vm{P}_1 & \vm{P}_2 & \vm{P}_3 & \vm{P}_4 & \vm{P}_5
     \end{bmatrix} = 
     \begin{bmatrix}
         c^2 & s^2 & -2cs & c^2(c\eta-s\xi)  & s^2(c\xi+s\eta) \\
         s^2 & c^2 & 2cs & s^2(c\eta-s\xi) & c^2(c\xi+s\eta) \\
         cs & -cs & c^2-s^2 & cs(c\eta-s\xi) & -cs(c\xi +s\eta)
     \end{bmatrix},
 \end{align}
 where $c$ and $s$ are given in~\eqref{eq:rotated_coords}. Without loss of generality, we choose an orthogonal basis for terms representing constant stresses, which
 results in the matrix
  \begin{align}
     \vm{P} = 
     \begin{bmatrix}
         1 & 0 & 0 & c^2(c\eta-s\xi)  & s^2(c\xi+s\eta) \\
         0 & 1 & 0 & s^2(c\eta-s\xi) & c^2(c\xi+s\eta) \\
         0 & 0 & 1 & cs(c\eta-s\xi) & -cs(c\xi +s\eta)
     \end{bmatrix}.
 \end{align}
 We now construct the stress-hybrid projection operator on the space $\vm{V}(E)$ over the original element $E$ with respect to the basis $\vm{P}$. From~\eqref{eq:matrix_vect_proj},
 we have the relation:
 \begin{align*}
     \int_{E}{\overline{\symP}^T\vm{C}^{-1}\overline{\Projbeta} \, d\vx} &= \int_{\partial E}{\overline{\symP}^T\Nmatrix\vm{u}_h \, ds} - \int_{E}{ \left(\vm{\partial}\overline{\symP}\right)^T\vm{u}_h \, d\vx}.
 \end{align*}
 Expanding $\vm{u}_h$ in terms of the basis in $\vm{V}(E)$, we have $\vm{u}_h = \vm{\varphi}\vd$, where $\vd$ is the displacement vector. We also expand $\overline{\Projbeta}$ in terms of $\vm{P}$: $\overline{\Projbeta} = \vm{P}\vm{\beta}$, 
 and since $\symP$ is arbitrary we take $\overline{\symP} = \vm{P}_i \ \ (i=1,2,\dots,5)$. After substituting in~\eqref{eq:matrix_vect_proj} for each $i=1,2,\dots,5$ and simplifying, we obtain the system:
\begin{equation}
     \left(\int_{E}{\vm{P}^T \vm{C}^{-1}\vm{P} \, d\vx}\right)\vm{\beta}  = \left(\int_{\partial E}{\vm{P}^T\Nmatrix\vm{\varphi} \, ds} - \int_{E}{\left(\vm{\partial} \vm{P}\right)^T\vm{\varphi} \, d\vx}\right)\vd.
\end{equation}
For this choice of $\vm{P}$, we have $\vm{\partial} \vm{P}=\vm{0}$ (divergence-free), so we obtain
\begin{equation}
      \left(\int_{E}{\vm{P}^T \vm{C}^{-1}\vm{P} \, d\vx}\right)\vm{\beta}  =\left(\int_{\partial E}{\vm{P}^T\Nmatrix\vm{\varphi} \, ds} \right)\vd.
 \end{equation}
 Now define the corresponding matrices $\vm{H}$ and $\vm{L}$ by
 \begin{subequations}\label{eq:H_L_beta}
      \begin{align}
      \vm{H} &= \int_{E}{\vm{P}^T \vm{C}^{-1}\vm{P} \, d\vx}, \quad
     \vm{L} = \int_{\partial E}{\vm{P}^T\Nmatrix\vm{\varphi} \, ds}, \label{eq:HandL_matrix}\\
     \intertext{and then the stress coefficients are given by }
     \vm{\beta} &= \vm{H}^{-1}L\vd := \vm{\Pi}_{\beta}\vd , \label{eq:beta_displacement}
 \end{align}
 \end{subequations}
 where $\vm{\Pi}_{\beta}$ is the matrix representation
 of the stress-hybrid projection operator with respect to
 the symmetric tensor polynomial basis $\vm{P}$.

 \subsection{Element stiffness matrix and element force vector}
Following the structure of~\eqref{eq:weak_equilibrium}, we define the discrete system using the projection operator $\Projbeta$ by:
\begin{subequations}
    \begin{align*}
        &a_h^E(\vm{u}_h,\delta\vm{u}_h) = \ell^E_h(\delta \vm{u}_h),\\
        \intertext{where }
        a_h^E(\vm{u}_h,\delta\vm{u}_h)&:=\int_{E}{\overline{\Projbeta(\delta \vm{u}_h)}^T \vm{C}^{-1} \overline{\Projbeta( \vm{u}_h)} \, d\vx}, \\ 
        \ell^E_h(\delta \vm{u}_h) &:= \int_{E}{(\delta \vm{u}_h)^T\vm{b} \, d\vx} + \int_{\Gamma_t\cap \partial E}{(\delta\vm{u}_h)^T\Bar{t} \, ds}.
    \end{align*}
\end{subequations}
 Expanding $\overline{\Projbeta}$ in terms of $\beta$ and applying~\eqref{eq:beta_displacement}, we obtain
 \begin{align}
      a_h^{E}(\vm{u}_h,\delta\vm{u}_h) = (\delta \vd)^T(\vm{\Pi}_\beta)^T\left(\int_{E}{ \vm{P}^T\vm{C}^{-1}\vm{P}  \, d\vx}\right)\vm{\Pi}_\beta \vd := (\delta \vd)^t \vm{K}_E \vd,
 \end{align}
where we identify the element stiffness matrix for SH-VEM as 
\begin{align}\label{eq:SH-VEM_stiffness}
    \vm{K}_E = (\vm{\Pi}_\beta)^T\left(\int_{E}{ \vm{P}^T\vm{C}^{-1}\vm{P}  \, d\vx}\right)\vm{\Pi}_\beta = \vm{\Pi}_\beta^T \vm{H} \vm{\Pi}_\beta.
\end{align}
In the Appendix, we give an alternate stress-hybrid virtual
element formulation based on Cook's approach,~\cite{Cook:1974:jsd} and show that it is
identical to the stress-hybrid element stiffness matrix that is 
obtained using $\vm{P}^*$ in~\eqref{eq:P_matrix_full}.

Now for every element $E$, the element force vector is given by
 \begin{align}\label{eq:force_vect}
    \vm{f}_{E} :=  \int_{E}{\vm{\varphi}^T\vm{b}  \, d\vm{x}} + \int_{\Gamma_t \cap \partial E}{\vm{\varphi}^T\Bar{\vm{t}} \, ds}.
 \end{align}
 For a low-order method, the first term in~\eqref{eq:force_vect} is
 approximated by taking the nodal average of the basis functions $\vm{\varphi}$ and then using a single-point quadrature to compute the integral.\cite{Chen:2023:SFV} The second term is computed using Gauss quadrature over the element edges. 

\subsection{B-bar VEM}
 For later numerical tests, we compare the results of the SH-VEM to a B-bar VEM. Following Park et at.,~\cite{Park:2020:meccanica} we first decompose the material moduli matrix $\vm{C}$ in terms of its eigenvectors:
 \begin{align}\label{eq:C:decomposed}
     \vm{C}  = \sum_{i=1}^{3}\lambda_i \vm{p}_i\vm{p}_i^T,
 \end{align}
 where $(\lambda_i , \vm{p}_i)$ is the $i$-th eigenpair of $\vm{C}$.
 It is known for plane elasticity that $\vm{p}_1 = \frac{1}{\sqrt{2}}[1,1,0]^T$ and $\lambda_1 = 2\kappa + \frac{2\mu}{3}$, where $\kappa$ is the bulk modulus and $\mu$ is the shear modulus. We express~\eqref{eq:C:decomposed} as 
 \begin{align}
     \vm{C} = \lambda_1\vm{p}_1\vm{p}_1^T + \sum_{i=2}^{3}\lambda_i \vm{p}_i\vm{p}_i^T := \vm{C}_{\textrm{dil}} + \vm{C}_{\textrm{dev}}.
 \end{align}
 The element stiffness matrix in the B-bar formulation is the sum of a consistency 
 matrix and a stabilization matrix. For the consistency matrix $\vm{K}^c$, we have after simplification:
 \begin{align}
     \vm{K}^c &= (\vm{\Pi}_{*}^{\veps})^T\left(\int_{E}{ (\vm{S}\vm{M})^T\vm{C}(\vm{S}\vm{M})\,d\vx}\right)\vm{\Pi}_{*}^{\veps} \nonumber \\ &= (\vm{\Pi}_{*}^{\veps} )^T\left(\int_{E}{(\vm{S}\vm{M})^T\vm{C}_{\textrm{dil}}(\vm{S}\vm{M})\,d\vx}\right)\vm{\Pi}_{*}^{\veps}  + (\vm{\Pi}_{*}^{\veps})^T\left(\int_{E}{(\vm{S}\vm{M})^T\vm{C}_{\textrm{dev}} (\vm{S}\vm{M})\,d\vx}\right)\vm{\Pi}_{*}^{\veps}  \nonumber
     \\ &:= \vm{K}^c_{\textrm{dil}} + \vm{K}^c_{\textrm{dev}},
 \end{align}
 where $\vm{S}$ is defined in~\eqref{eq:S_matrix}, $\vm{\Pi}_{*}^{\veps}$ is given in~\eqref{eq:Pi_eps_matrix} and $\vm{M}$ is:
 \begin{align*}
     \vm{M} = \begin{bmatrix}
         1 & 0 & -\eta & \eta & \xi & 0 \\ 
         0 & 1 & \xi & \xi & 0 & \eta
     \end{bmatrix}.
 \end{align*}
 The expression for the stabilization matrix is:
 \begin{subequations}
      \begin{align}
     \vm{K}^s = (\vm{I}-\vm{\Pi}^{\veps})^T\vm{\Lambda}(\vm{I}-\vm{\Pi}^{\veps}), \\ 
     \intertext{where $\vm{\Pi}^{\veps}$ is defined in~\eqref{eq:Pi_eps_canonical} and $\vm{\Lambda}$ is a diagonal matrix with components}
     \vm{\Lambda}_{ii} = \max\left((\vm{K}^c_{\textrm{dev}})_{ii},\frac{\mu}{2}\right).
 \end{align}
 \end{subequations}

 \section{Numerical results}\label{sec:numericalresults}
 We present a collection of two-dimensional numerical examples for linear elasticity under plane strain conditions. We examine the errors of the displacements in the $L^2$ norm and energy seminorm, and the $L^2$ error of the hydrostatic
 stress.  The exact hydrostatic stress (denoted by $\tilde{p}$) and its numerical
 approximation are computed as:
\begin{align}
    \tilde{p} = \frac{\textrm{trace} \, (\vm{\sigma})}{3}, \quad
    \tilde{p}_h = \frac{1+\nu}{3} \left( \left(\overline{\Projbeta}\right)_1 + \left(\overline{\Projbeta}\right)_2\right),
\end{align}
where $\left(\overline{\Projbeta}\right)_i$ is the $i$-th component of $\overline{\Projbeta}$.

The convergence rates of B-bar VEM
and SH-VEM are computed using 
the following discrete error measures:
 \begin{subequations}\label{eq:error_norms}
\begin{align}
\|\vm{u}-\vm{u}_h\|_{\vm{L}^2(\Omega)} &=\sqrt{\sum_{E}{\int_{E}{|\vm{u}-\proj{\vm{u}_h}}|^2 \,d\vm{x}}},  \\ 
\|\tilde{p}-\tilde{p}_h\|_{\vm{L}^2(\Omega)} &=\sqrt{\sum_{E}{\int_{E}{|\tilde{p}-\tilde{p}_h|^2 \,d\vm{x}}}},\\
\|\vm{u} - \vm{u}_h\|_{a} &= \sqrt{\sum_{E}{\int_{E}{(\overline{\vsigma}-\Projbetavect)^T\,\vm{C}^{-1}(\overline{\vsigma}-\Projbetavect)\, d\vm{x}}}}.
\end{align}
\end{subequations}
The matrix $H$ in~\eqref{eq:HandL_matrix} and the integrals that appear in~\eqref{eq:error_norms} are computed using the scaled boundary cubature scheme;\cite{chin:2021:cmame}
see also Chen and Sukumar~\cite{Chen:2023:SFV,Chen:2023:SFS} for its use in the stabilization-free virtual element method.

\subsection{Eigenanalysis of the element stiffness matrix}
We first examine the stability of the SH-VEM for rotated elements through an eigenanalysis. From Cook,~\cite{Cook:1974:jsd} it is known that for a noninvariant method, a rectangular element rotated by $\frac{\pi}{4}$ will contain spurious zero-energy modes. For this test, we take a unit square and rotate it by angle $\ac{\gamma} = 0, \frac{\pi}{6}, \frac{\pi}{4}, \frac{\pi}{3}$, and then compute the eigenvalues of the element stiffness matrix. The material has Young's modulus $E_Y = 1$ psi and Poisson's ratio $\nu=0.4999999$. 
Three formulations are considered: an unrotated $5\beta$, rotated $5\beta$ and an unrotated $7\beta$. In the unrotated formulations, the projection matrix is computed on the original element $E$ without applying the rotated coordinate transformation given in~\eqref{eq:rotated_coords}, and for the $7\beta$ formulation a $7$-term expansion is used for$\vm{P}$:
\begin{align*}
    \vm{P} = \begin{bmatrix}
        1 & 0 & 0 & \eta & 0 & 0 & \xi \\
        0 & 1 & 0 & 0 & \xi & \eta & 0 \\
        0 & 0 & 1 & 0 & 0 & -\xi & -\eta
    \end{bmatrix}.
\end{align*}
Every method has three physical zero eigenvalues that correspond to the zero-energy modes. For a method to be stable the next lowest eigenvalue must be positive and not close to zero. In Table~\ref{tab:eig_analysis}, we indicate the fourth smallest eigenvalue of the element stiffness matrix for the three formulations. The table shows that both the rotated $5\beta$ and the unrotated $7\beta$ have their eigenvalues unaffected for any \ac{angle} $\ac{\gamma}$. However, as $\ac{\gamma}$ is increased to $\ac{\gamma}=\frac{\pi}{4}$, the next lowest eigenvalue of the unrotated $5\beta$ formulation becomes zero. This shows that the $5\beta$ SH-VEM in global coordinates is not rotationally invariant. Numerical tests also reveal that the $7\beta$ formulation ameliorates volumetric locking but is much stiffer for pure bending problems, and therefore both unrotated $5\beta$ and $7\beta$ formulations are not considered in the remainder of this article. 
\begin{table}[!h]
    \centering
    \begin{tabular}{|c|c|c|c|c|}
    \hline
      Method   & $\ac{\gamma} =0$ & $\ac{\gamma} = \frac{\pi}{6}$ & $\ac{\gamma}=\frac{\pi}{4}$ & $\ac{\gamma}=\frac{\pi}{3}$ \\
             \hline
       \textrm{Unrotated} $5\beta$  &0.444 &0.111& 0.000 &0.111 
       \\
       \textrm{Rotated} $5\beta$  &0.444 & 0.444 & 0.444 &0.444\\
       \textrm{Unrotated} $7\beta$ &0.444 & 0.444 & 0.444 & 0.444\\
    \hline
    \end{tabular}
    \caption{Comparison of the fourth-lowest eigenvalue on a square that is rotated by angle $\ac{\gamma}$ for three
    stress-hybrid VEMs. }
    \label{tab:eig_analysis}
\end{table}

\begin{figure}
\begin{minipage}{.48\textwidth}
     \begin{subfigure}{.48\textwidth}
         \centering
         \includegraphics[width=\textwidth]{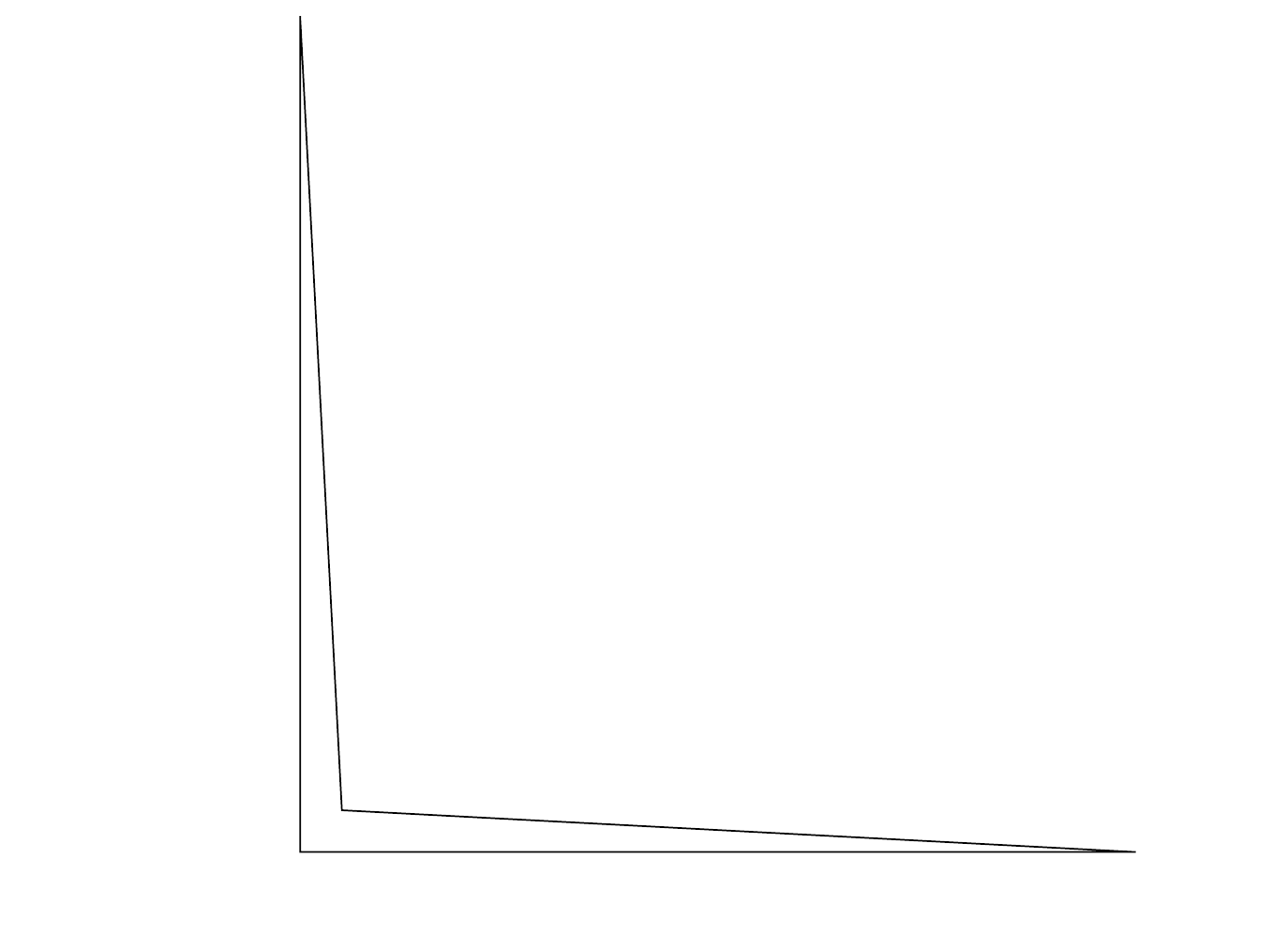}
         \caption{}
     \end{subfigure}
     \hfill
     \begin{subfigure}{.48\textwidth}
         \centering
         \includegraphics[width=\textwidth]{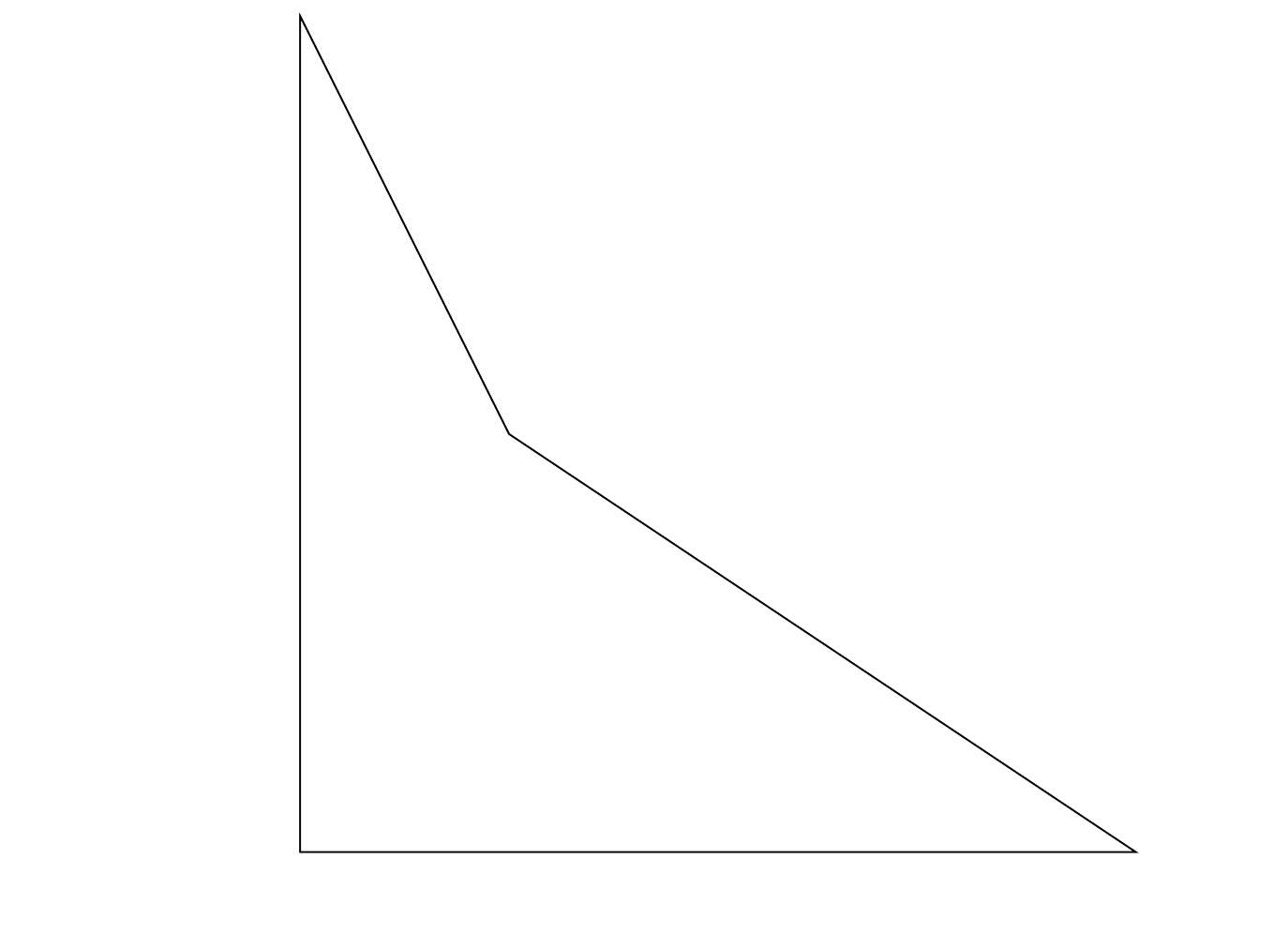}
         \caption{}
     \end{subfigure}
     \vfill
     \begin{subfigure}{.48\textwidth}
         \centering
         \includegraphics[width=\textwidth]{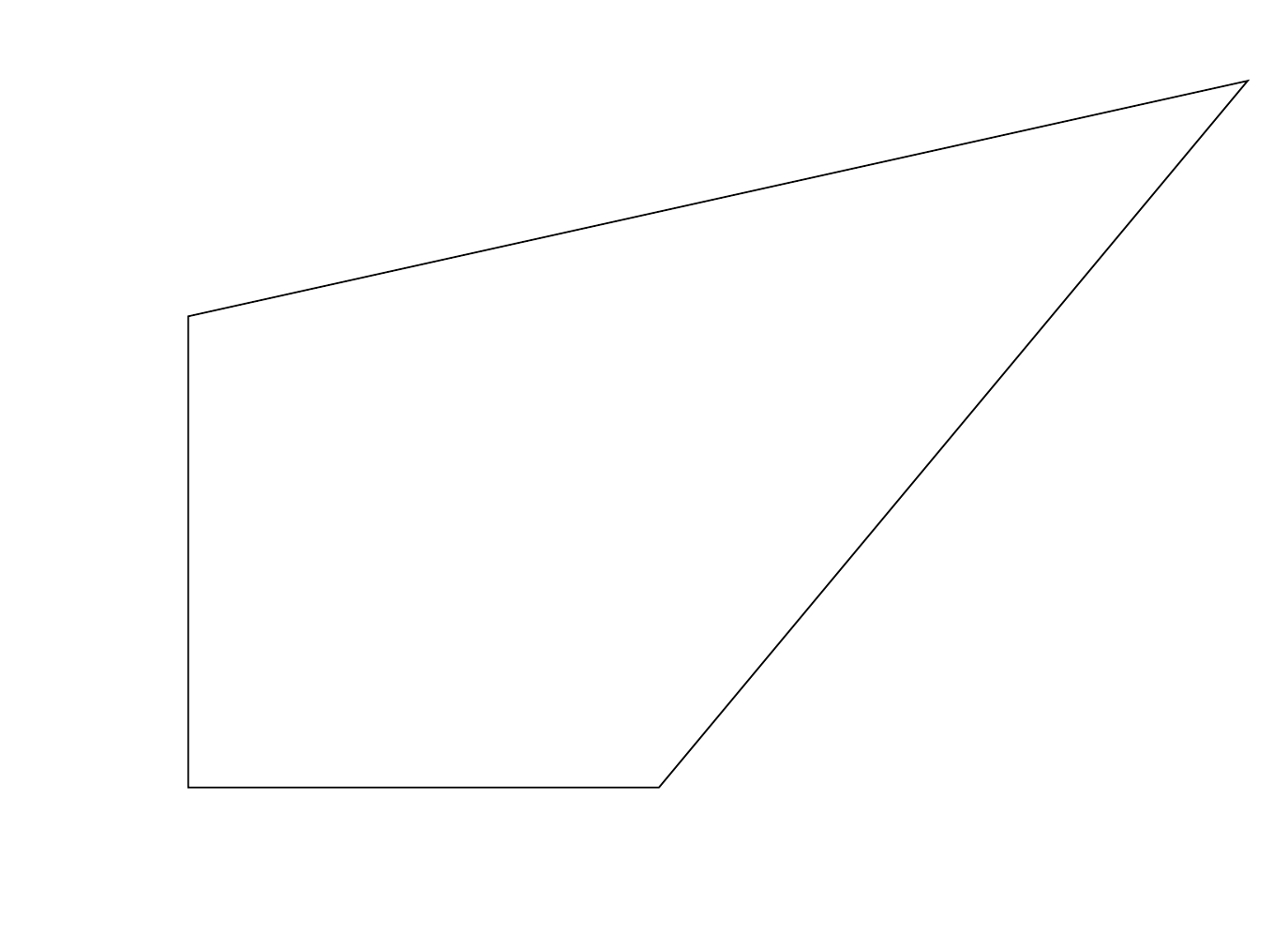}
         \caption{}
     \end{subfigure}
        \hfill
     \begin{subfigure}{.48\textwidth}
         \centering
         \includegraphics[width=\textwidth]{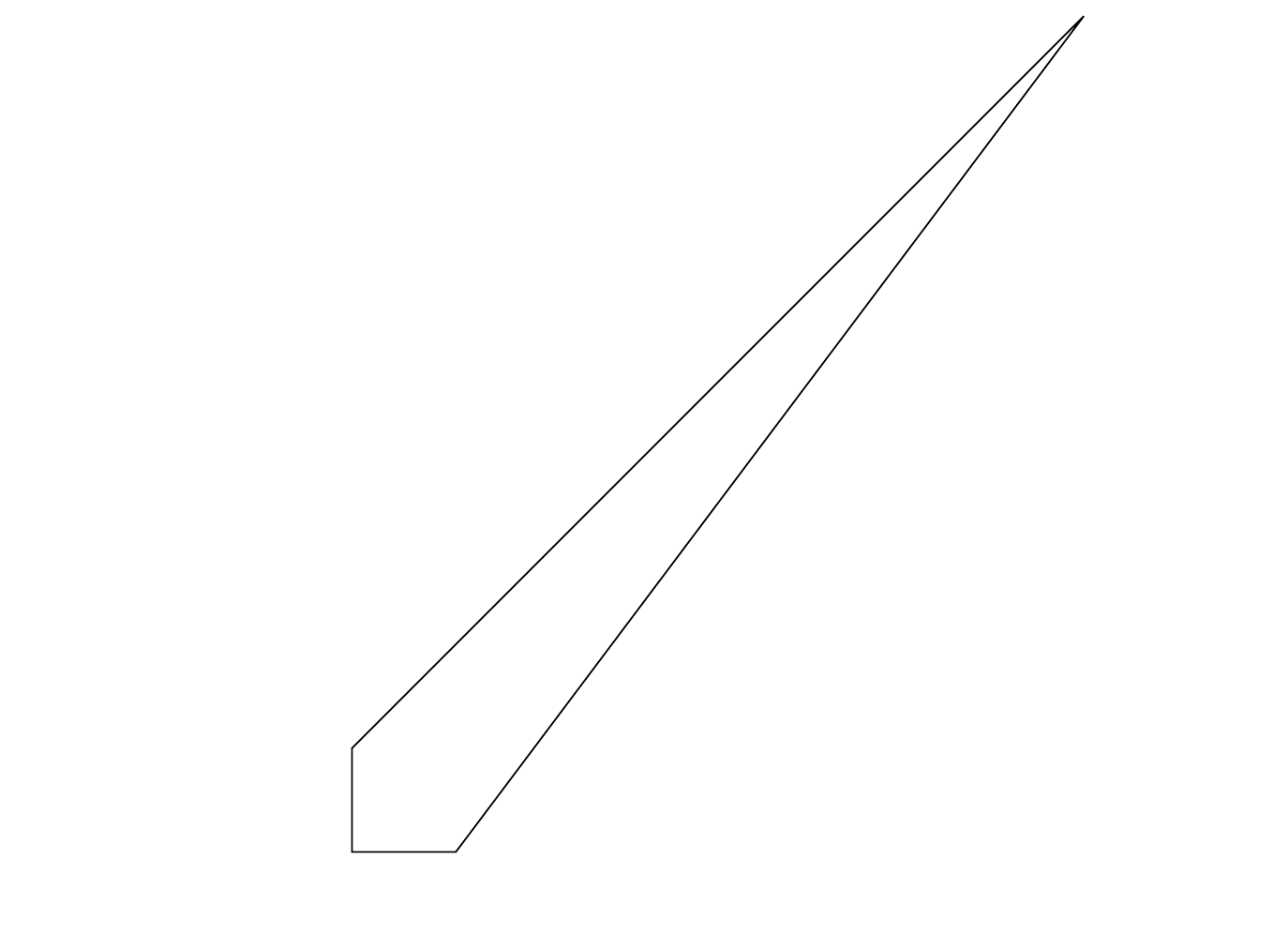}
         \caption{}
     \end{subfigure}
     \end{minipage}
     \hfill
     \begin{minipage}{.48\textwidth}
     \begin{subfigure}{\textwidth}
         \centering
         \includegraphics[width=\textwidth]{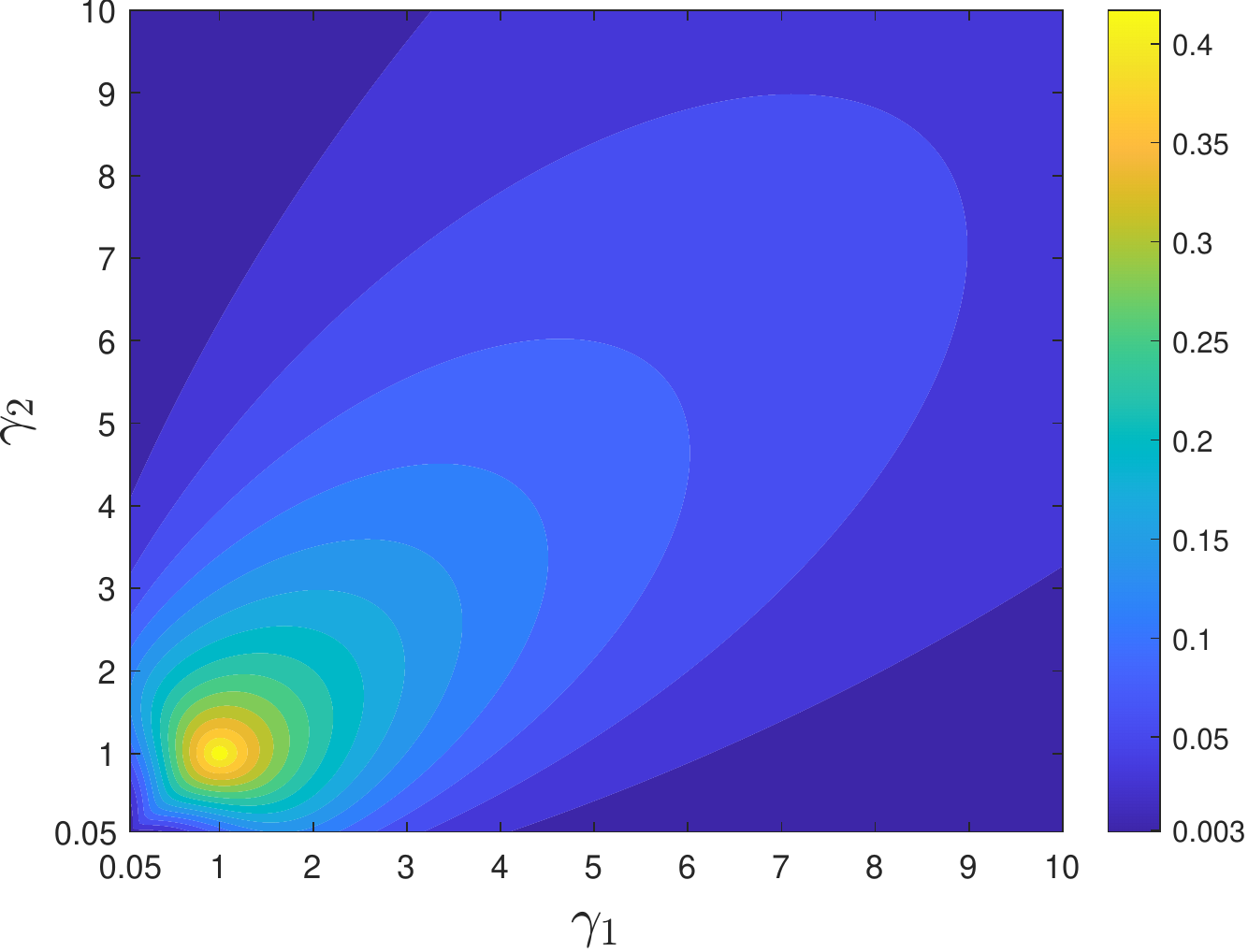}
         \caption{}
     \end{subfigure}
     \end{minipage}
        \caption{(a)--(d) First sequence of distorted quadrilaterals, where the fourth vertex is located at $(\ac{\gamma_1},\ac{\gamma_2})$, and (e) contour plot of 
        the fourth-lowest eigenvalue as a function of
        $\ac{\gamma_1}$ and $\ac{\gamma_2}$.}
        \label{fig:eig_meshes_1}
\end{figure}

To further test the stability of the SH-VEM on different convex and nonconvex element types, we consider two additional tests. For the second test, we study the effects of perturbing a vertex of a unit square. We construct quadrilaterals with coordinates $\{(0,0), \, (1,0), \, (\ac{\gamma_1},\ac{\gamma_2}), \, (0,1)  \},$ where $\ac{\gamma_1},\ac{\gamma_2} \in (0.05,10)$. For every combination of $\ac{\gamma_1}$ and $\ac{\gamma_2}$ we 
compute the element stiffness matrix on this quadrilateral and then determine its fourth smallest eigenvalue. A few 
representative elements and a contour plot of the eigenvalues are shown in Figure~\ref{fig:eig_meshes_1}. The contour plot reveals that deviations from the unit square decreases the value of the fourth smallest eigenvalue; however, the eigenvalue remains positive and away from zero
(greater than $0.003$) in all cases. This test shows that no spurious zero eigenvalues appear even for large perturbations of the unit square.  

\begin{figure}
\begin{minipage}{.48\textwidth}
     \begin{subfigure}{.48\textwidth}
         \centering
         \includegraphics[width=\textwidth]{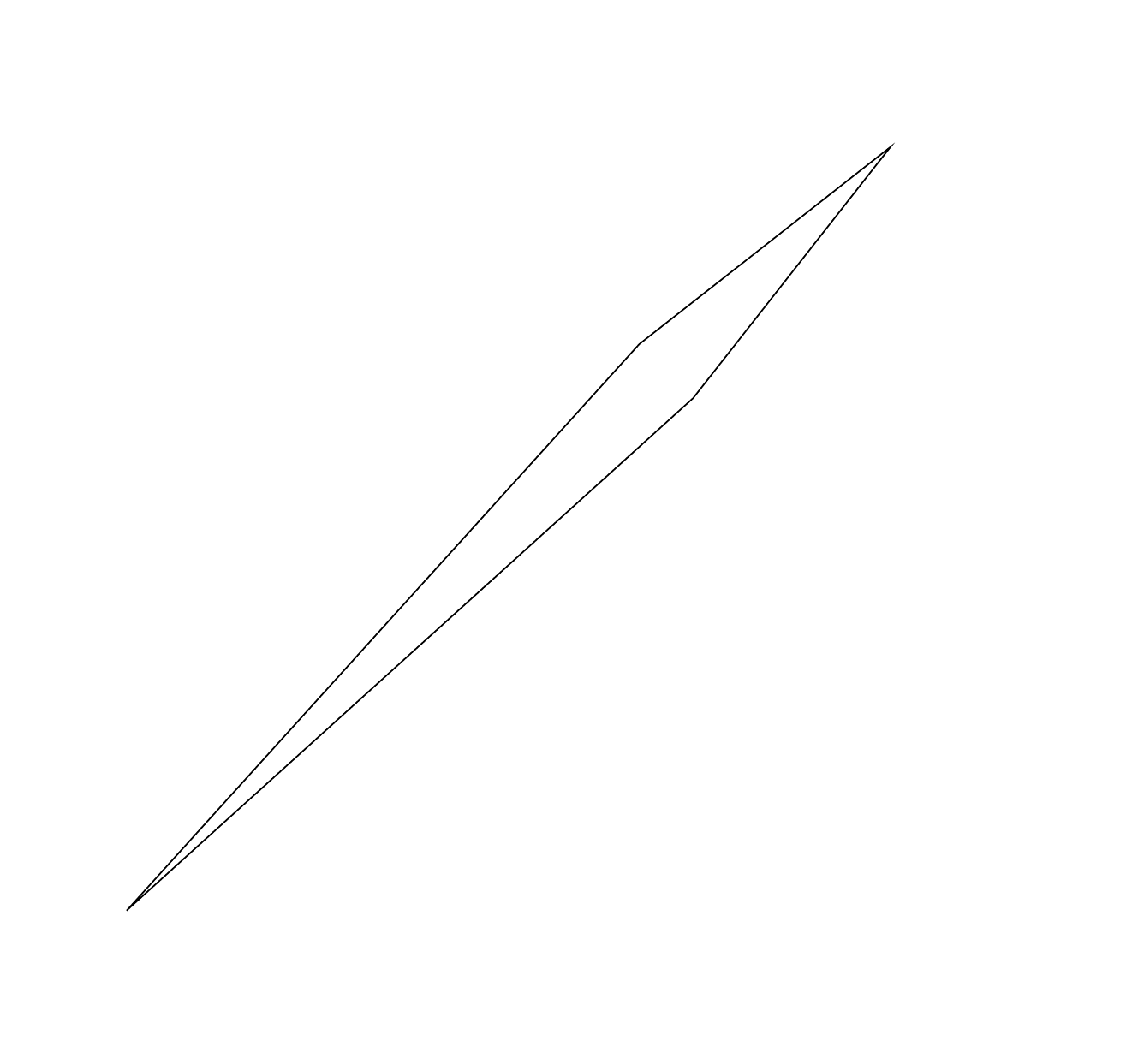}
         \caption{}
     \end{subfigure}
     \hfill
     \begin{subfigure}{.48\textwidth}
         \centering
         \includegraphics[width=\textwidth]{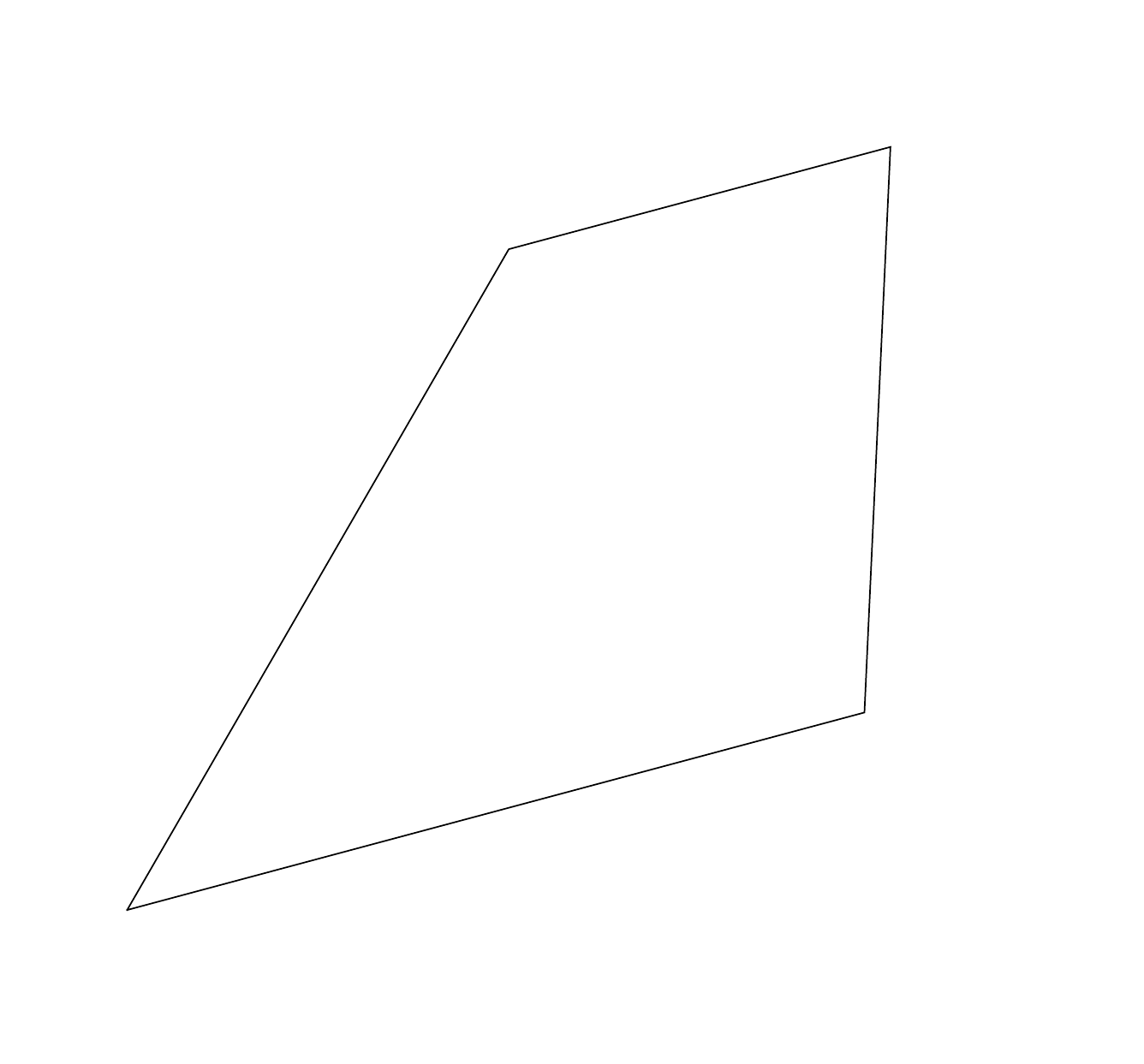}
         \caption{}
     \end{subfigure}
     \vfill
     \begin{subfigure}{.48\textwidth}
         \centering
         \includegraphics[width=\textwidth]{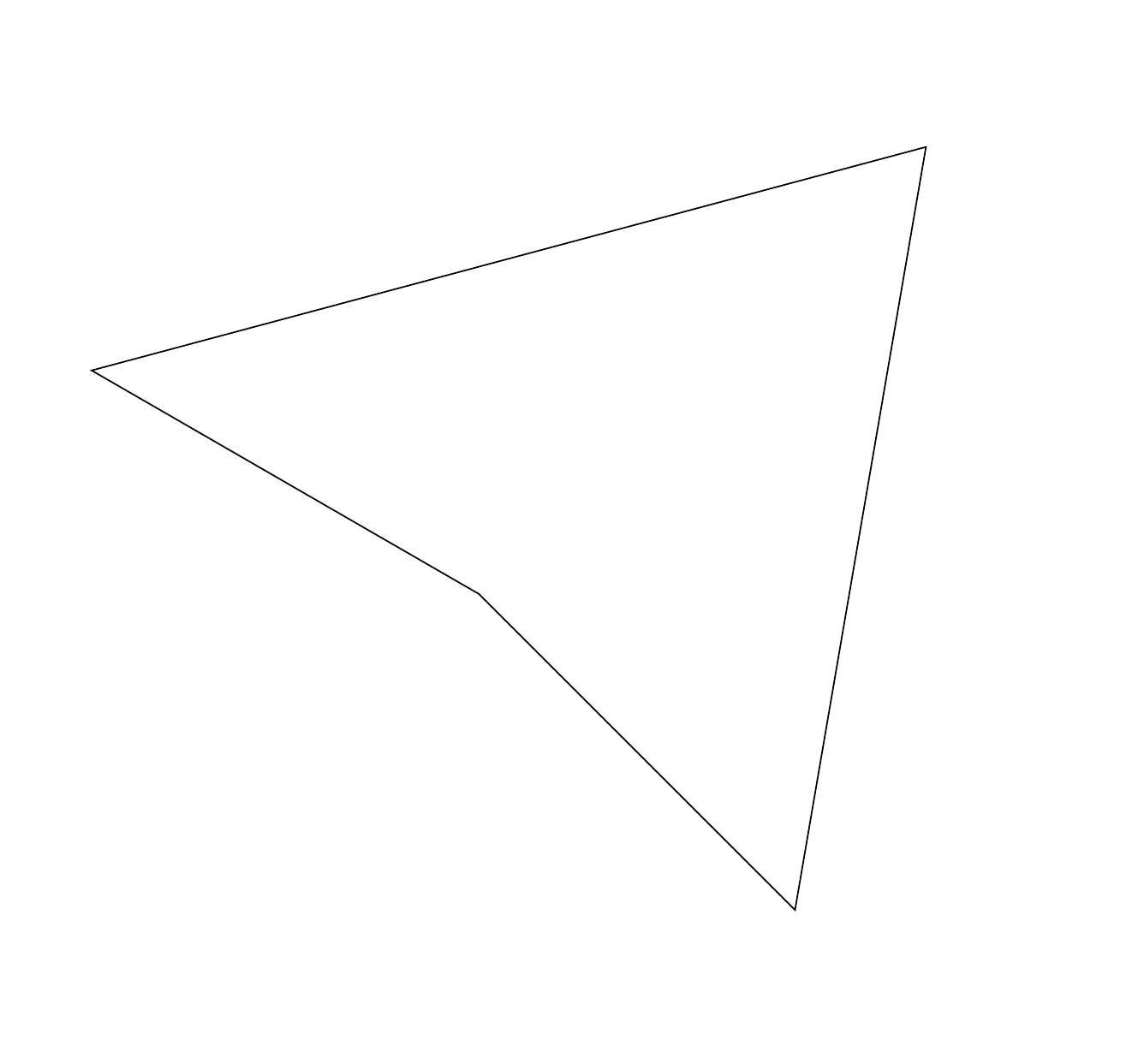}
         \caption{}
     \end{subfigure}
        \hfill
     \begin{subfigure}{.48\textwidth}
         \centering
         \includegraphics[width=\textwidth]{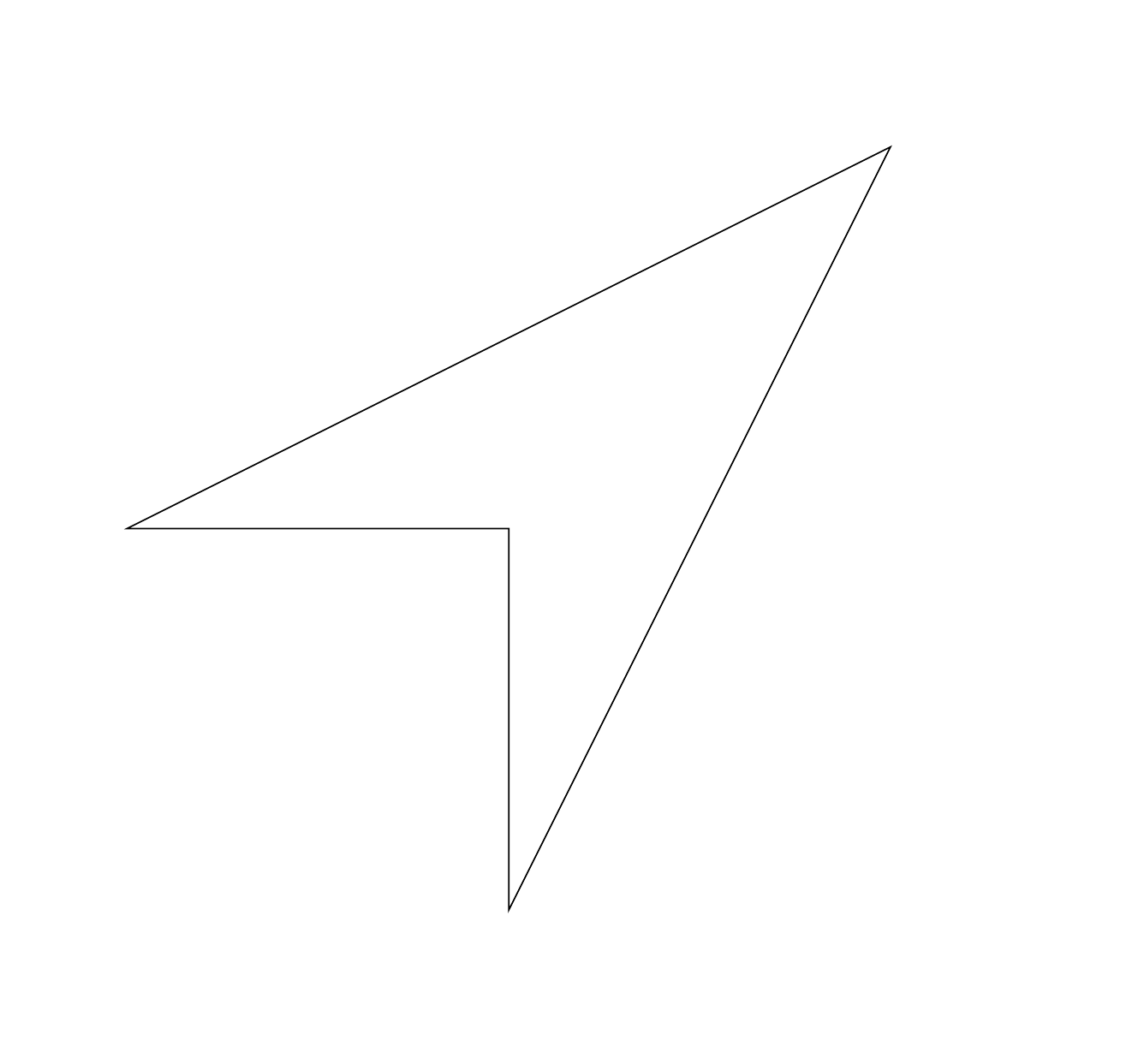}
         \caption{}
     \end{subfigure}
     \end{minipage}
     \hfill
     \begin{minipage}{.48\textwidth}
     \begin{subfigure}{\textwidth}
         \centering
         \includegraphics[width=\textwidth]{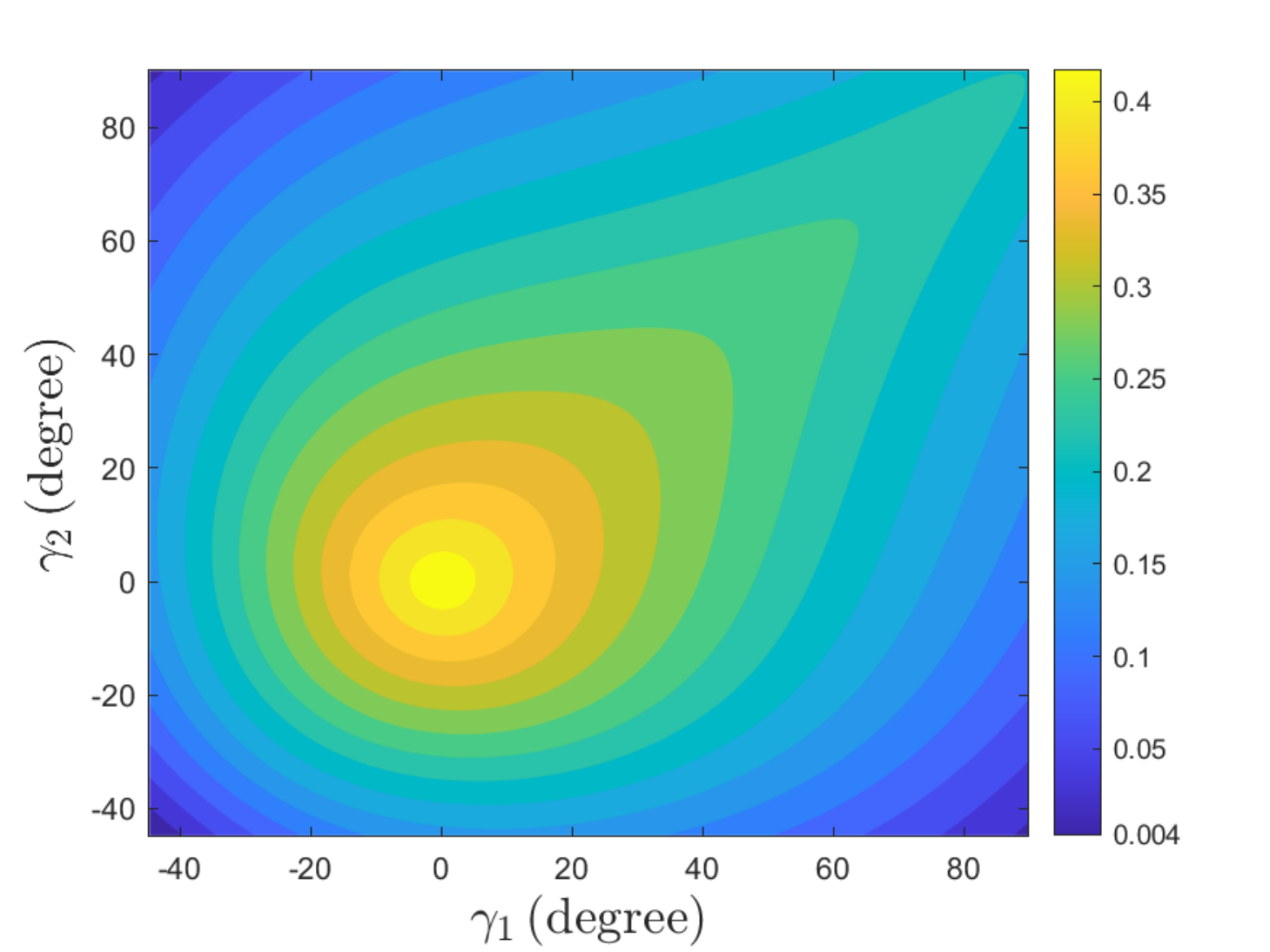}
         \caption{}
     \end{subfigure}
     \end{minipage}
        \caption{(a)--(d) Second sequence of distorted quadrilaterals, where the two vertices that are varied are located at $(\cos\ac{\gamma_1}, -\!\sin\ac{\gamma_1})$ and $(-\!\sin \ac{\gamma_2}, \cos\ac{\gamma_2})$, and
        (e) contour plot of 
        the fourth-lowest eigenvalue as a function of
        $\ac{\gamma_1}$ and $\ac{\gamma_2}$.}
        \label{fig:eig_meshes_2}
\end{figure}

In the third test, we examine the effects of varying the angles of a unit square by constructing quadrilaterals with coordinates $\{(0,0), \, (\cos\ac{\gamma_1}, -\sin\ac{\gamma_1}), \, (1,1), \, (-\sin \ac{\gamma_2}, \cos\ac{\gamma_2})\}$, where $\ac{\gamma_1}, \ac{\gamma_2} \in \bigl[-\frac{\pi}{4},\frac{\pi}{2}\bigr]$. We again compute the eigenvalues of the element stiffness matrix for different combinations of $\ac{\gamma_1}$ and $\ac{\gamma_2}$. Figure~\ref{fig:eig_meshes_2} shows a few representative elements and a contour plot of the fourth smallest eigenvalue. The contour plot shows that the smallest nonzero eigenvalue remains positive and away from zero 
(greater than $0.004$) for any combination of $\ac{\gamma_1},\ac{\gamma_2} \in \bigl[-\frac{\pi}{4},\frac{\pi}{2}\bigr]$, and hence demonstrates that
distorting a quadrilateral by varying its angle does not affect the 
stability of SH-VEM. 

\subsection{Manufactured problem: Convergence in the incompressible limit}
\begin{figure}[!h]
     \centering
     \begin{subfigure}{.32\textwidth}
         \centering
         \includegraphics[width=\textwidth]{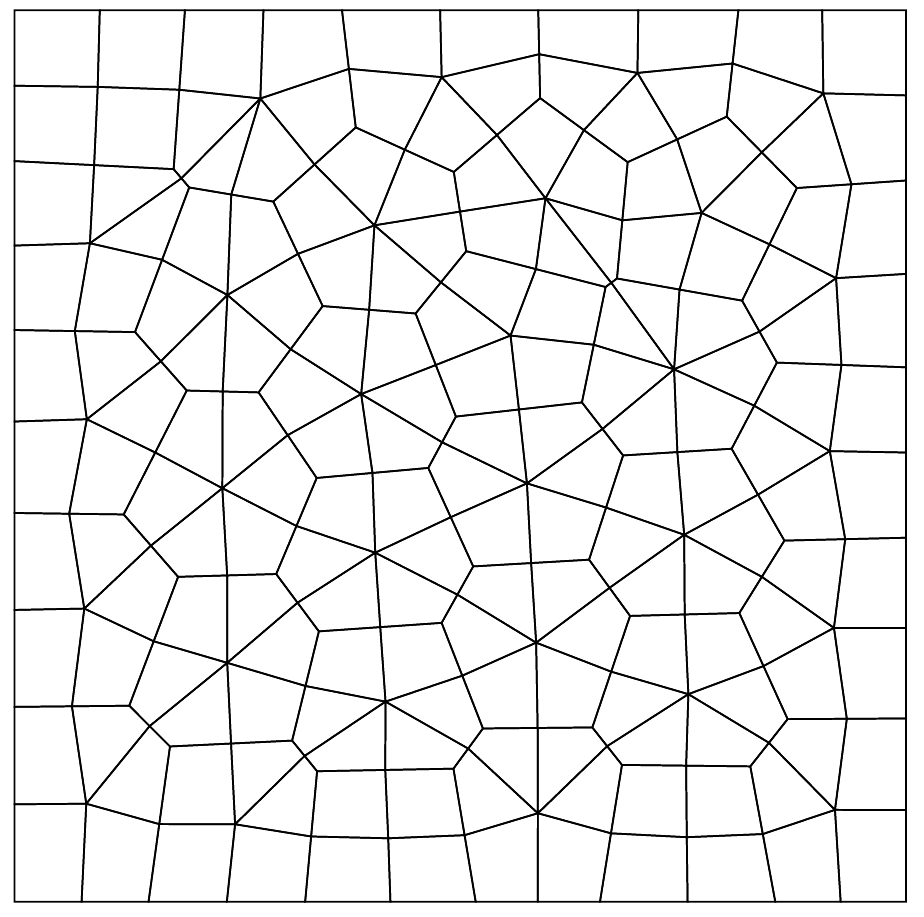}
         \caption{}
     \end{subfigure}
     \hfill
     \begin{subfigure}{.32\textwidth}
         \centering
         \includegraphics[width=\textwidth]{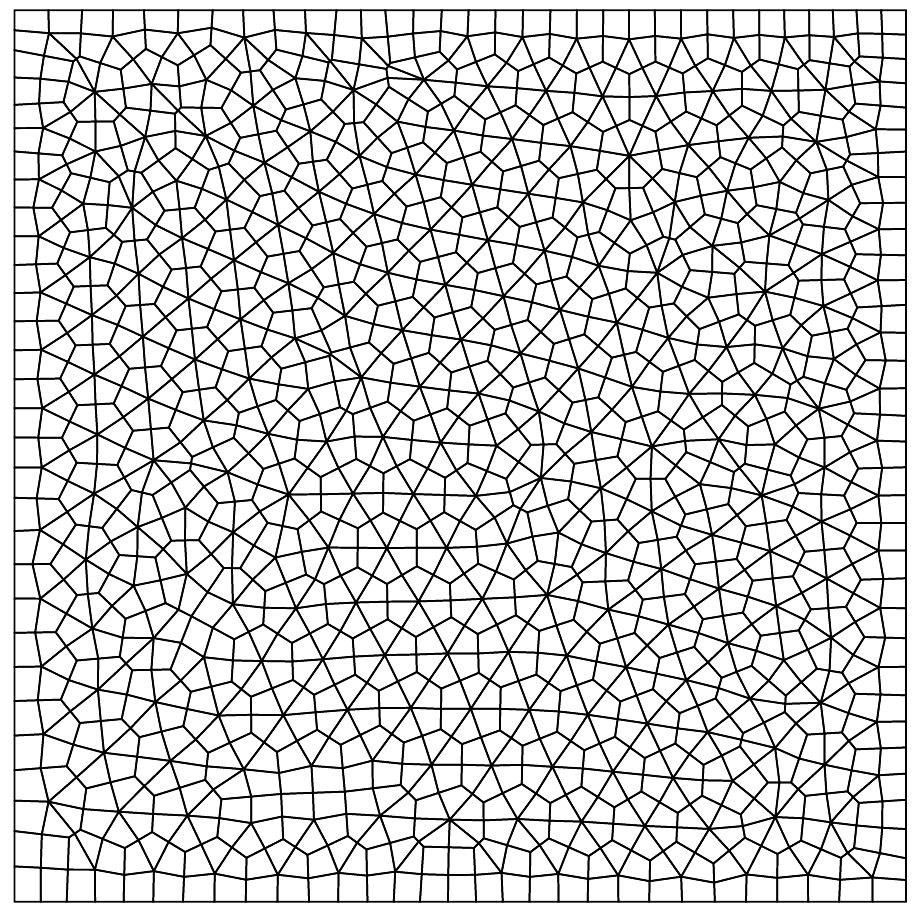}
         \caption{}
     \end{subfigure}
     \hfill
     \begin{subfigure}{.32\textwidth}
         \centering
         \includegraphics[width=\textwidth]{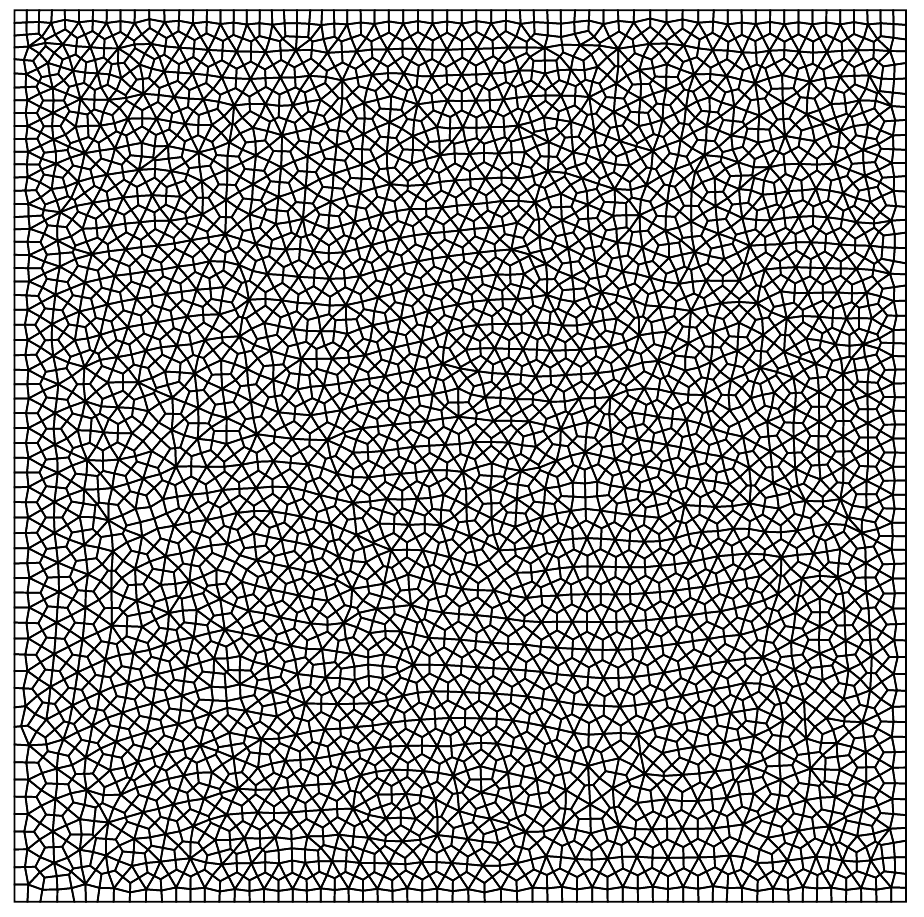}
         \caption{}
     \end{subfigure}
        \caption{Quadrilateral meshes for the manufactured problem. (a) 150 elements, (b) 1500 elements and (c) 6000 elements.}
        \label{fig:squaremesh}
\end{figure}

\begin{figure}[!h]
     \centering
     \begin{subfigure}{.24\textwidth}
         \centering
         \includegraphics[width=\textwidth]{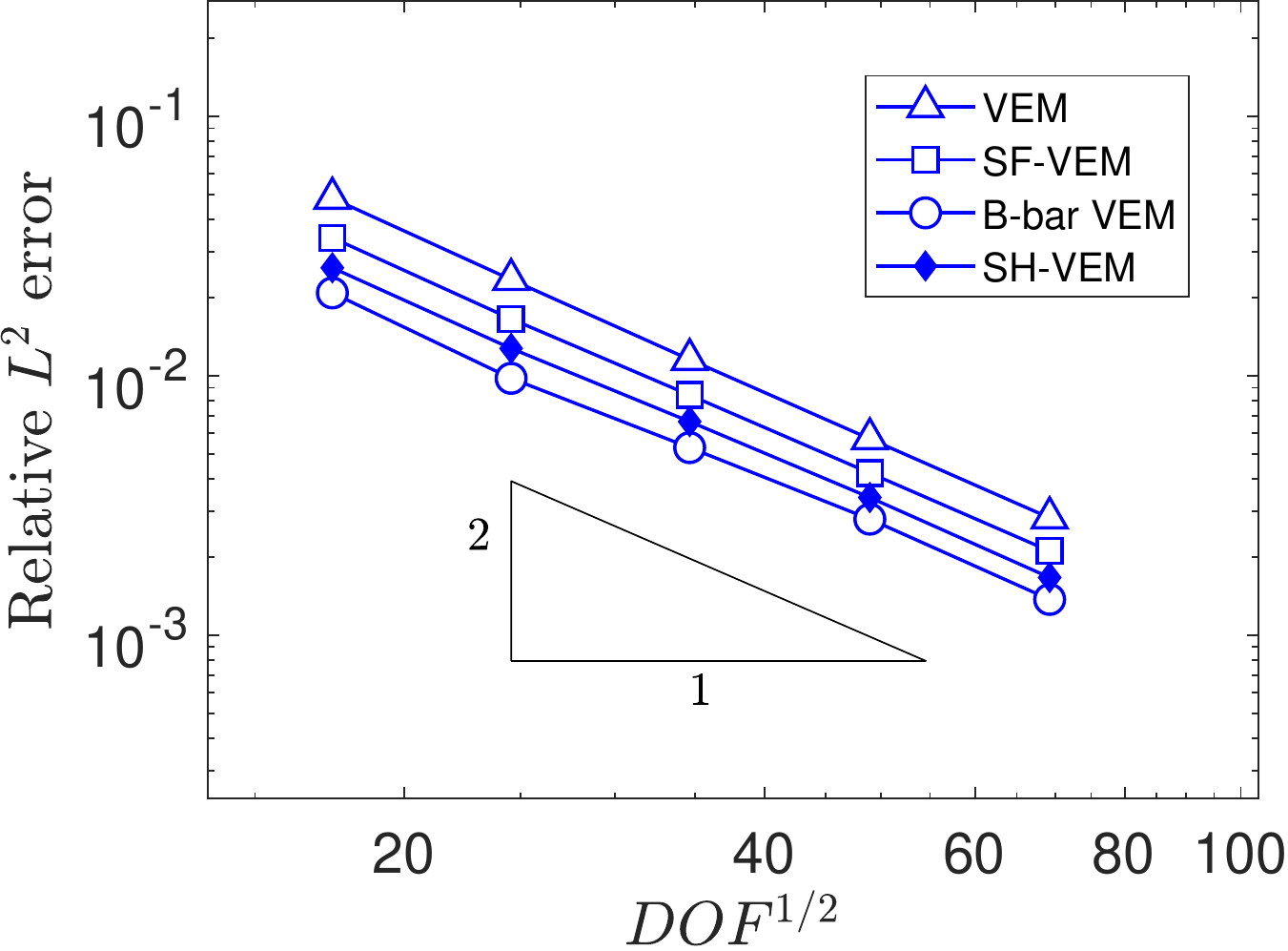}
     \end{subfigure}
     \hfill
     \begin{subfigure}{.24\textwidth}
         \centering
         \includegraphics[width=\textwidth]{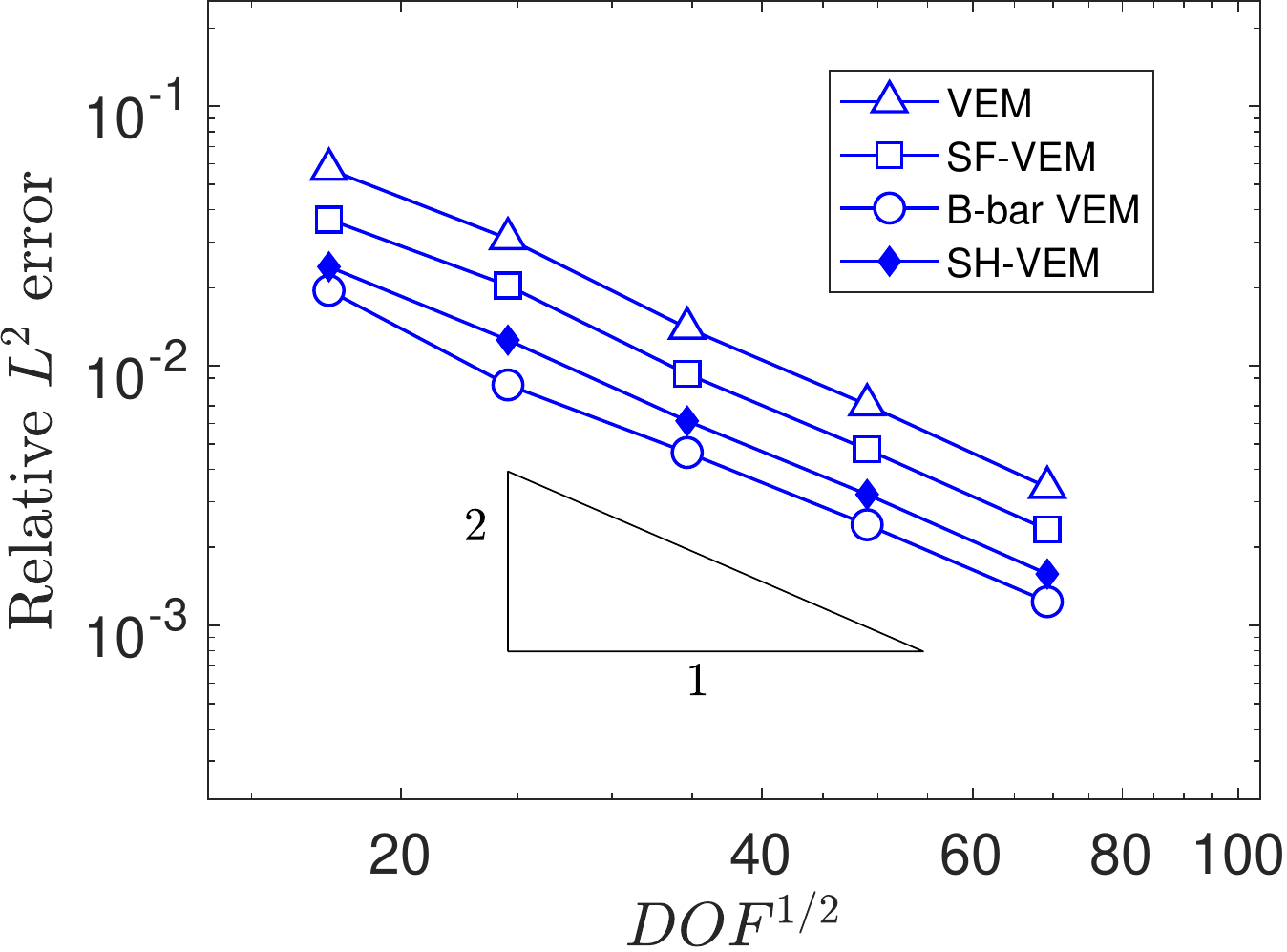}
     \end{subfigure}
     \hfill
     \begin{subfigure}{.24\textwidth}
         \centering
         \includegraphics[width=\textwidth]{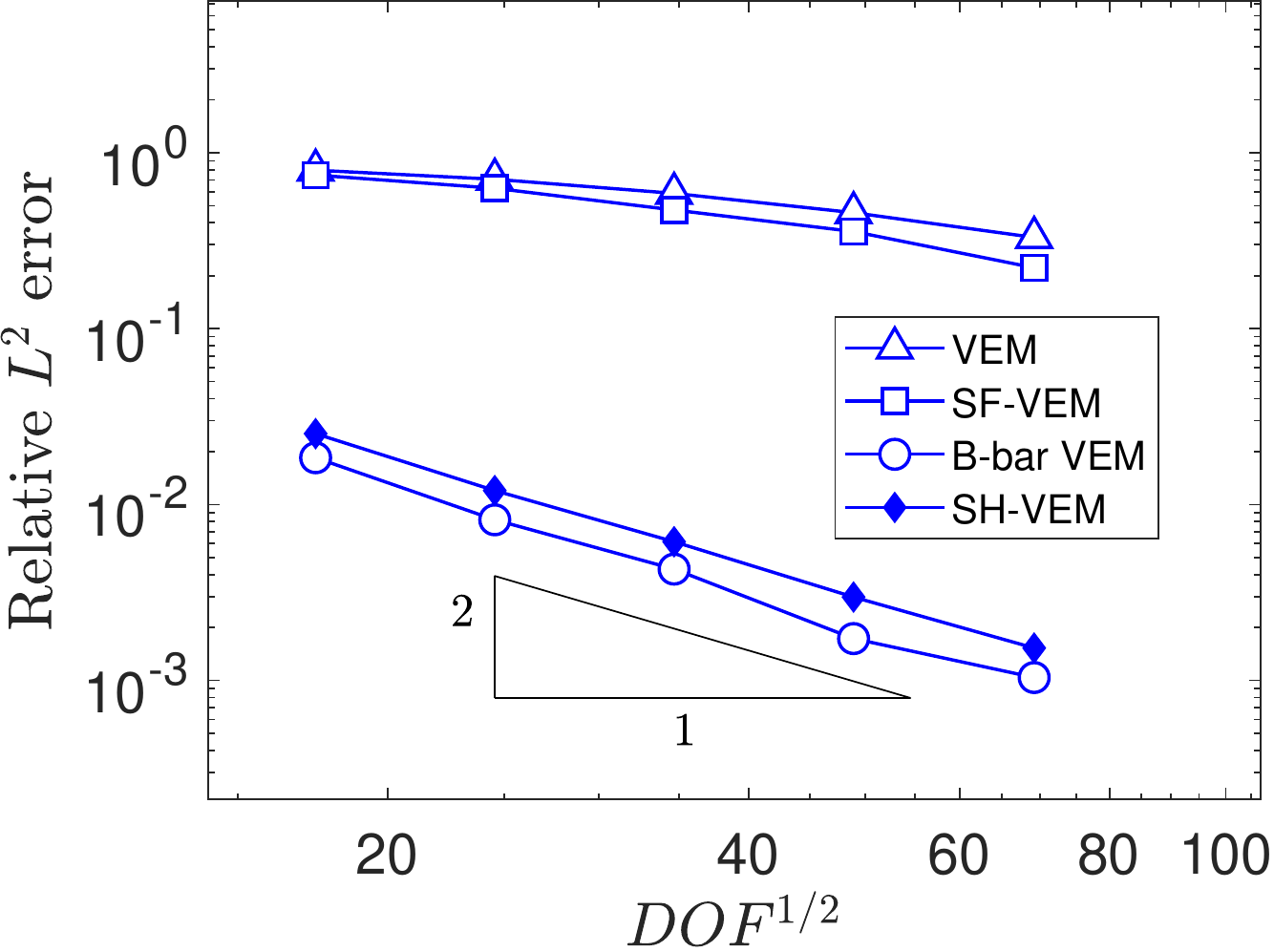}
     \end{subfigure}
          \hfill
     \begin{subfigure}{.24\textwidth}
         \centering
         \includegraphics[width=\textwidth]{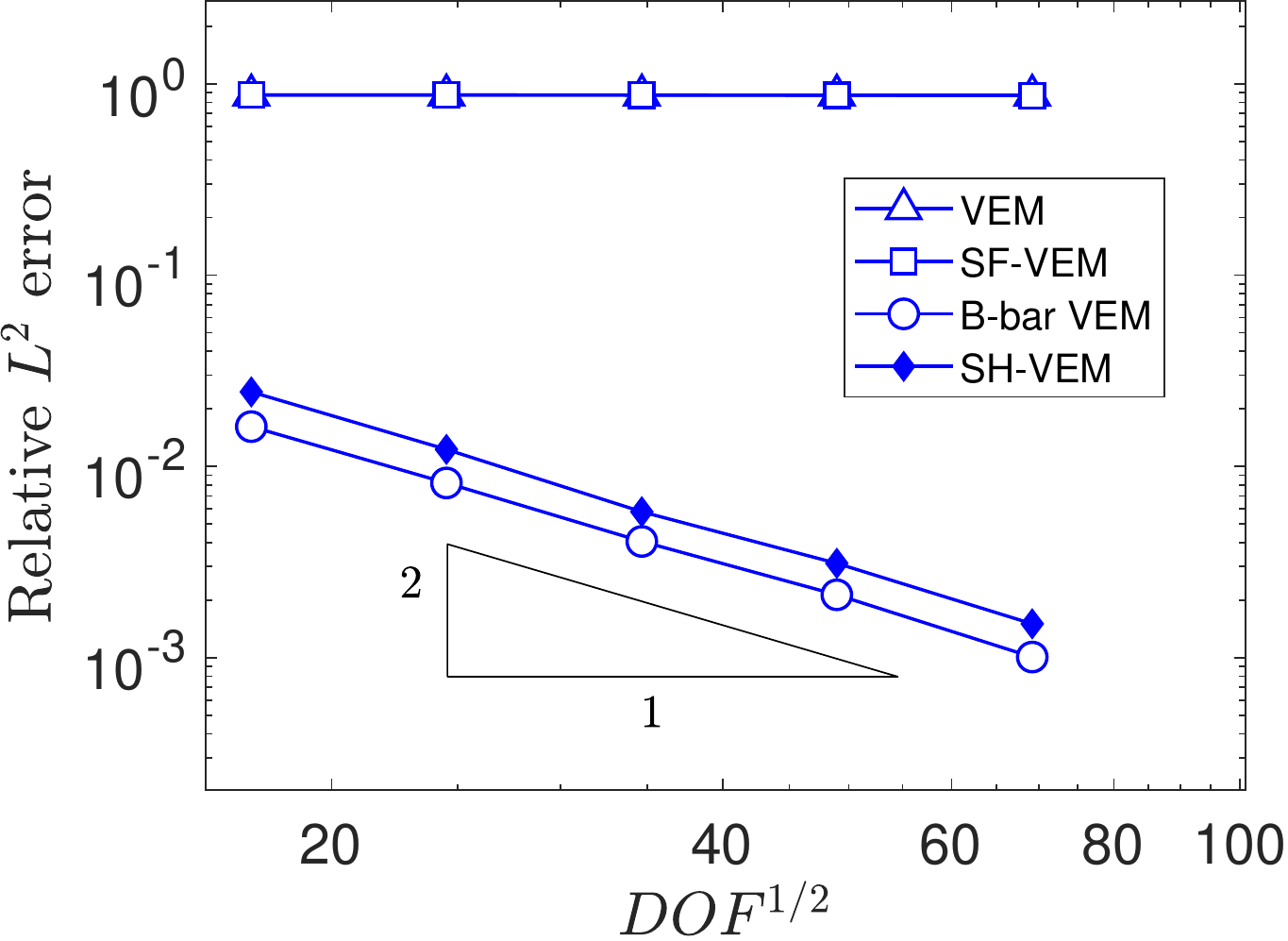}
     \end{subfigure}
     \vfill
          \begin{subfigure}{.24\textwidth}
         \centering
         \includegraphics[width=\textwidth]{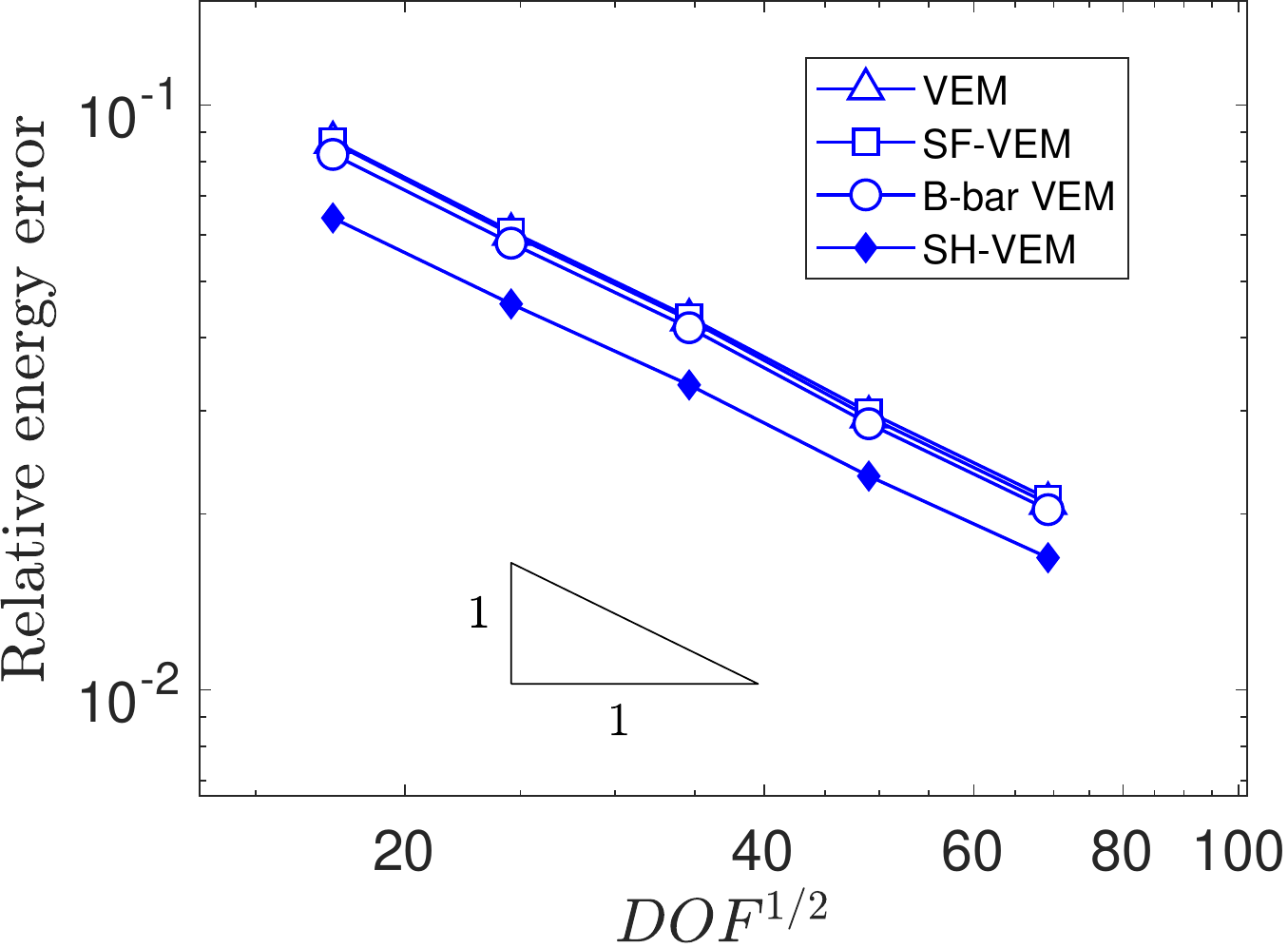}
         \caption{}
     \end{subfigure}
     \hfill
     \begin{subfigure}{.24\textwidth}
         \centering
         \includegraphics[width=\textwidth]{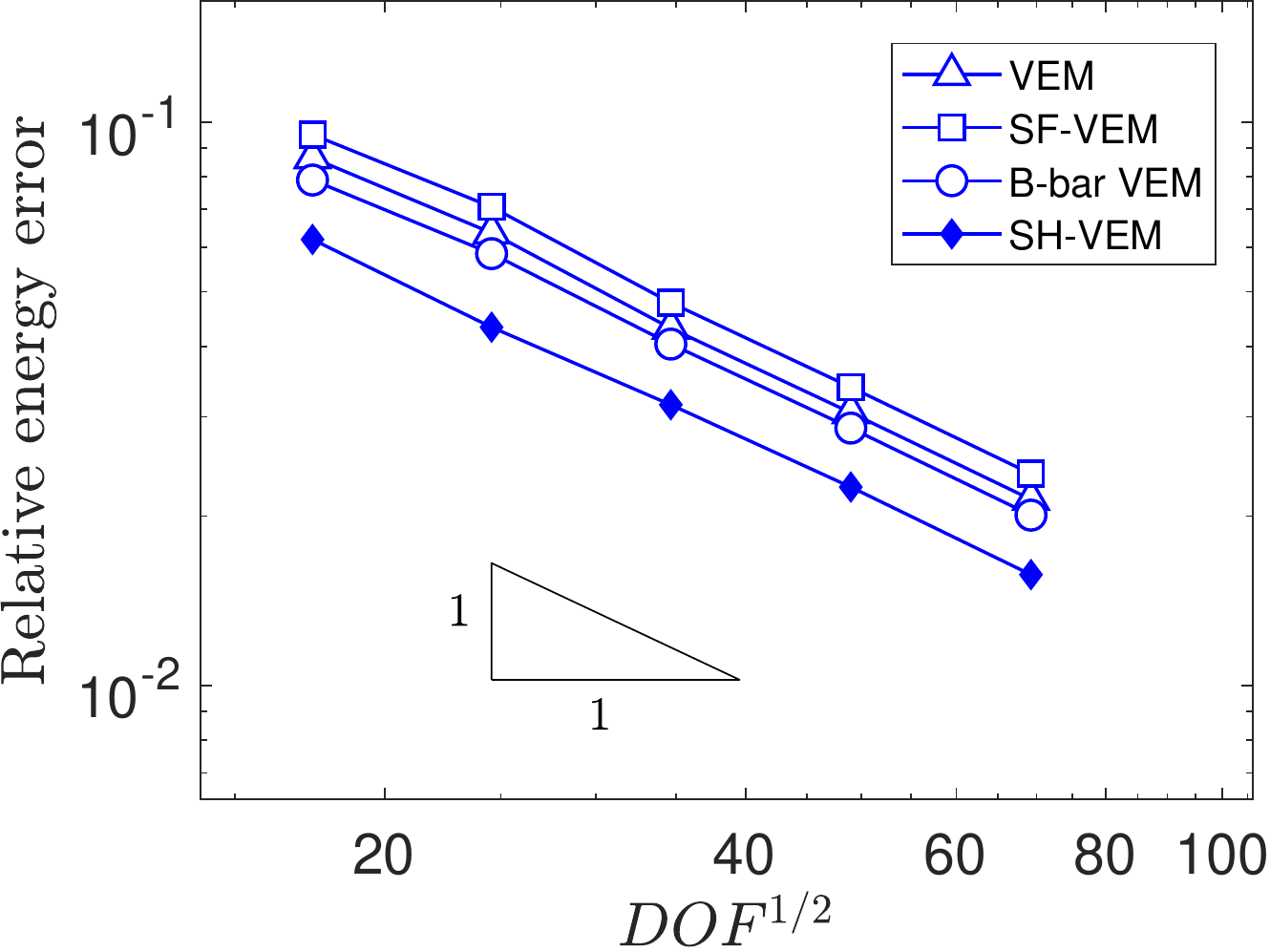}
         \caption{}
     \end{subfigure}
     \hfill
     \begin{subfigure}{.24\textwidth}
         \centering
         \includegraphics[width=\textwidth]{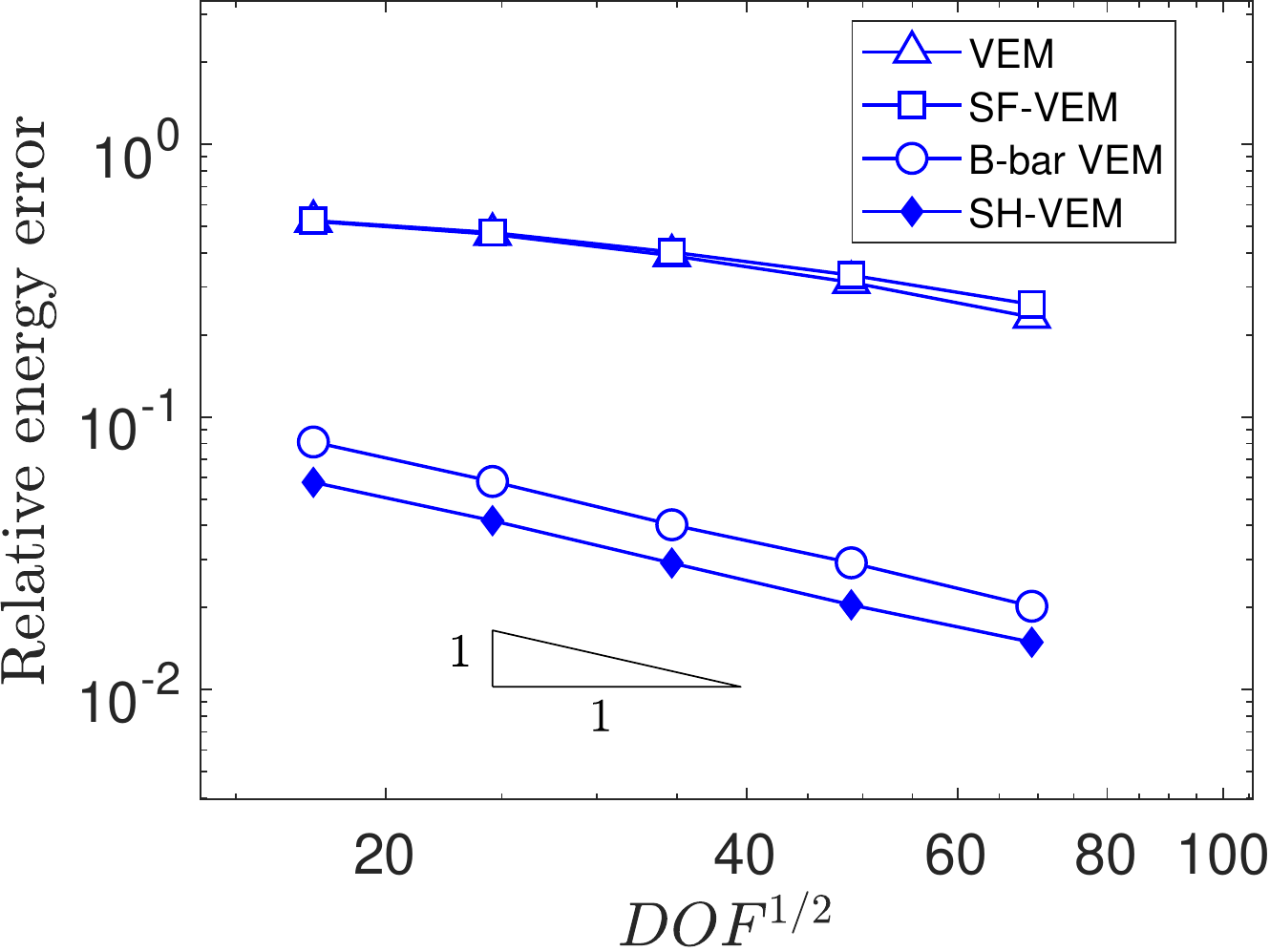}
         \caption{}
     \end{subfigure}
          \hfill
     \begin{subfigure}{.24\textwidth}
         \centering
         \includegraphics[width=\textwidth]{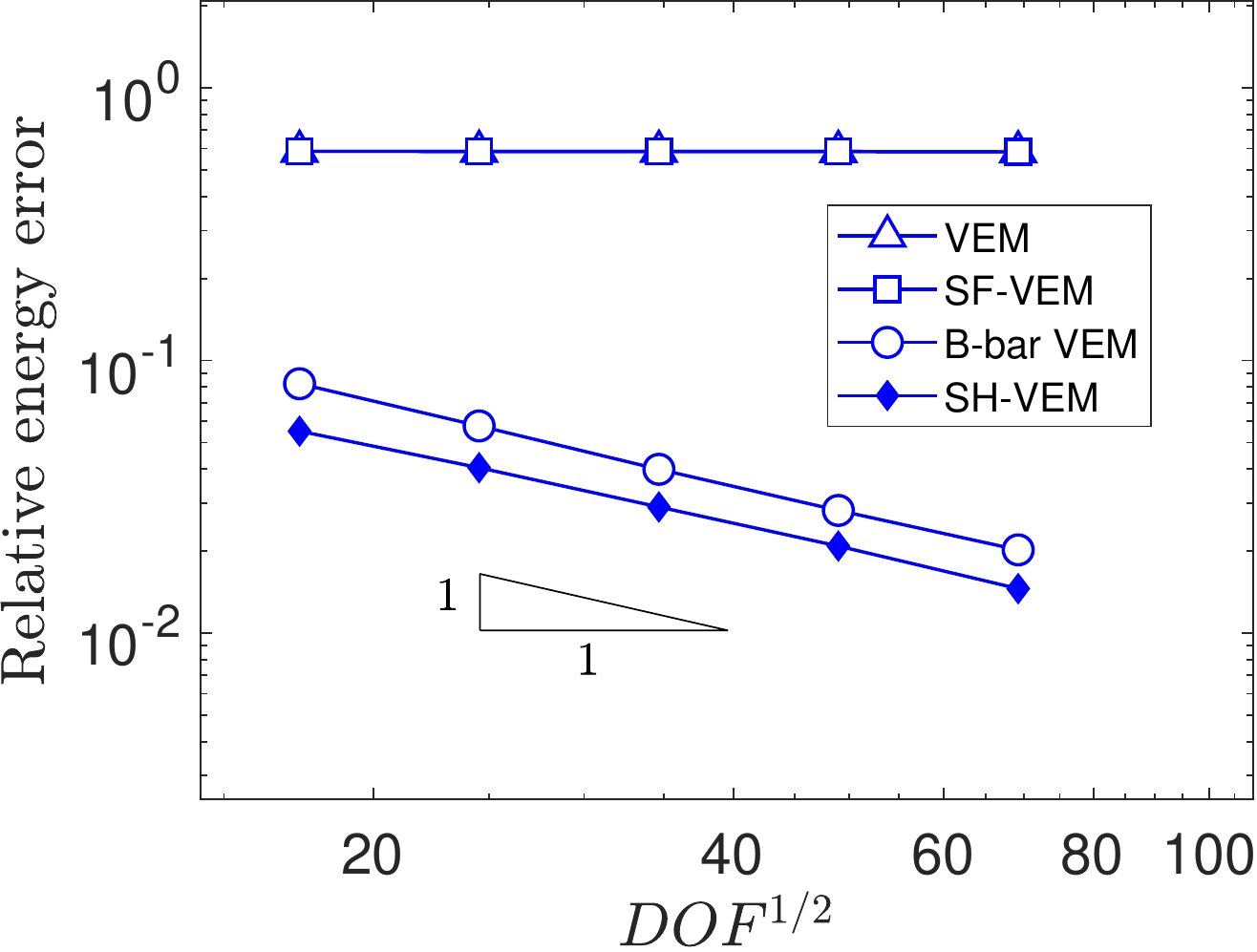}
         \caption{}
     \end{subfigure}
        \caption{Comparison of the convergence of standard VEM, stabilization-free VEM, B-bar VEM and SH-VEM  
        for the manufactured problem \acrev{on unstructured meshes (see~\fref{fig:squaremesh})}. Each column represents a different value of $\nu$. (a) $\nu=0.3$, (b) $\nu=0.4$, (c) $\nu=0.4999$ and (d) $\nu=0.4999999$. }
        \label{fig:manufactured}
\end{figure}

We examine the effects of increasing the Poisson ratio to the incompressible limit $(\nu \to 0.5)$ on a manufactured problem with a known solution.\cite{Ainsworth:2022:cmame} The problem domain is the unit square and the Young's modulus $E_Y = 1$ psi and the Poisson's ratio $\nu \in  \{0.3,\, 0.4,\, 0.4999,\,0.4999999\}$. The exact solution 
with associated loading is given by: 
\begin{equation*}
    u(\vm{x}) = -\cos(\pi x)\sin(\pi y)  \ \ \textrm{and} \ \ v(\vm{x}) =\sin(\pi x)\cos(\pi y), \quad \vm{b}(\vm{x}) = \frac{\pi^2}{1+\nu}\begin{Bmatrix}
         \cos(\pi x)\sin(\pi y) \\
         -\sin(\pi x)\cos(\pi y)
    \end{Bmatrix}.
\end{equation*}
In Figure~\ref{fig:squaremesh}, we show a few sample meshes for the unit square, and in Figure~\ref{fig:manufactured} we show the convergence rates in $L^2$ error of displacement and the energy seminorm as $\nu$ is varied. We assess four
formulations: standard VEM,\cite{elasticdaveiga}
SF-VEM,\cite{Chen:2023:SFV} 
B-bar VEM,\cite{Park:2020:meccanica} and SH-VEM. 
From these plots, we observe that as $\nu$ is
increased, the standard VEM and SF-VEM fail to converge, while both B-bar VEM
and SH-VEM converge with rates that are in agreement with
theory. 

\begin{figure}[!h]
     \centering
     \begin{subfigure}{.32\textwidth}
         \centering
         \includegraphics[width=\textwidth]{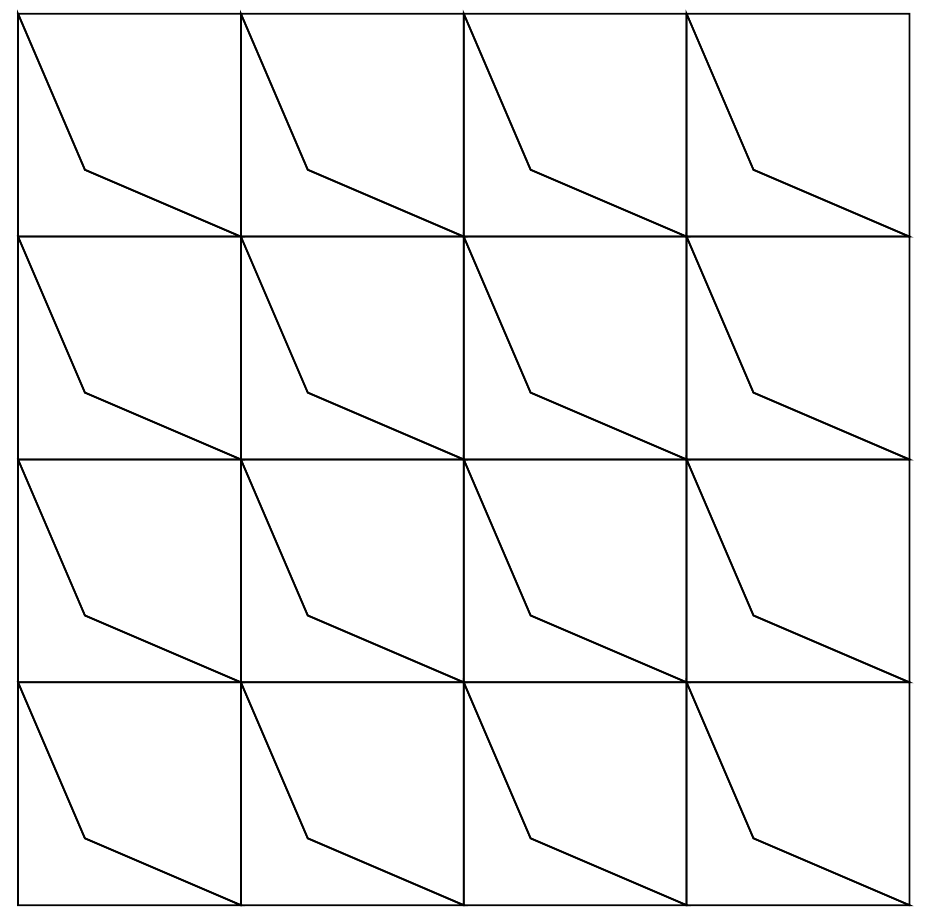}
         \caption{}
     \end{subfigure}
     \hfill
     \begin{subfigure}{.32\textwidth}
         \centering
         \includegraphics[width=\textwidth]{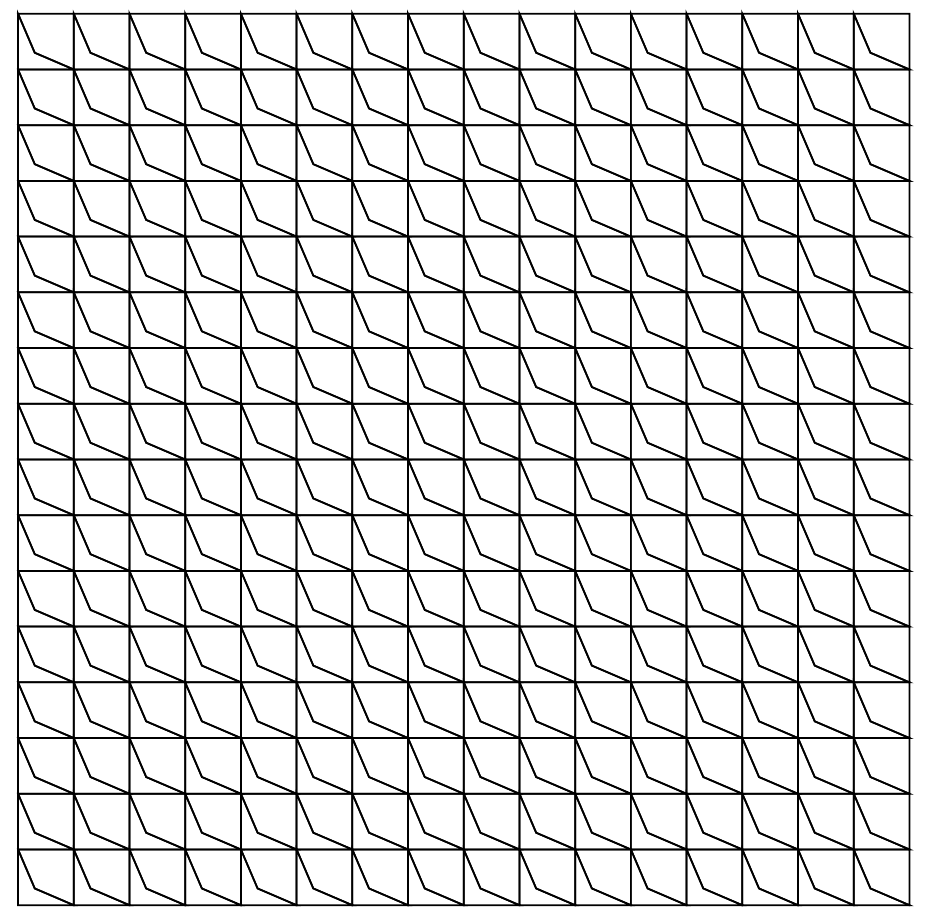}
         \caption{}
     \end{subfigure}
     \hfill
     \begin{subfigure}{.32\textwidth}
         \centering
         \includegraphics[width=\textwidth]{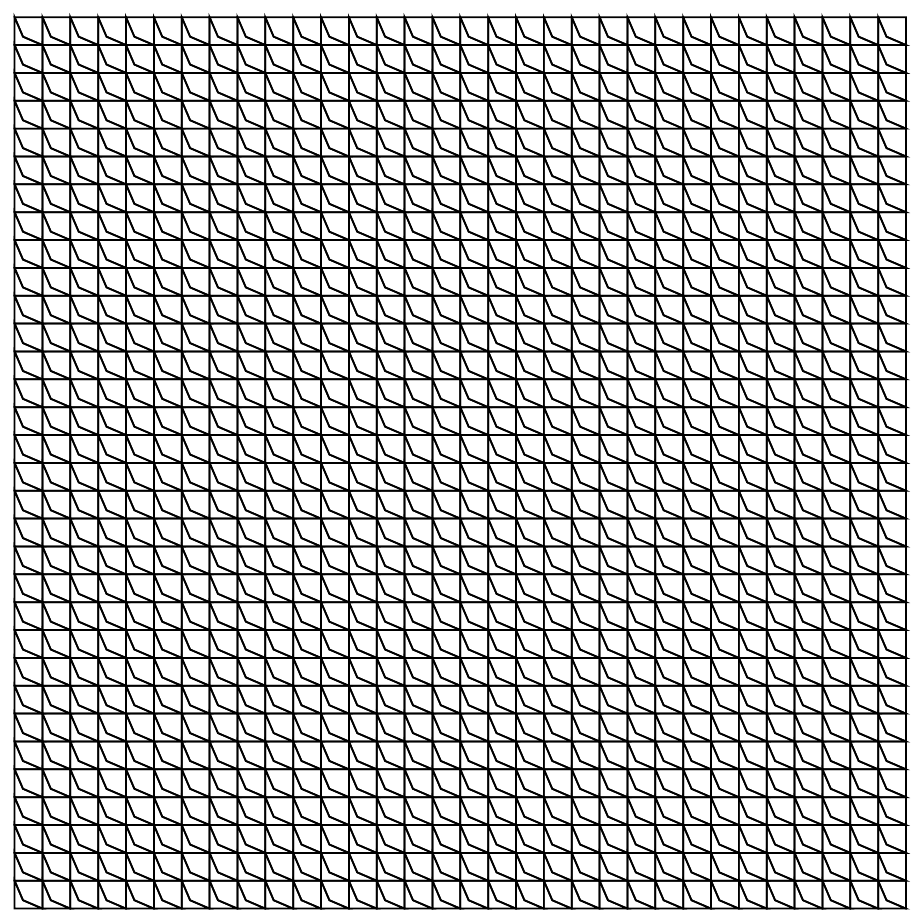}
         \caption{}
     \end{subfigure}
        \caption{Nonconvex quadrilateral meshes for the manufactured problem. (a) 32 elements, (b) 512 elements and (c) 2048 elements.  }
        \label{fig:squaremesh_nonconvex}
\end{figure}

\begin{figure}[!h]
     \centering
     \begin{subfigure}{.24\textwidth}
         \centering
         \includegraphics[width=\textwidth]{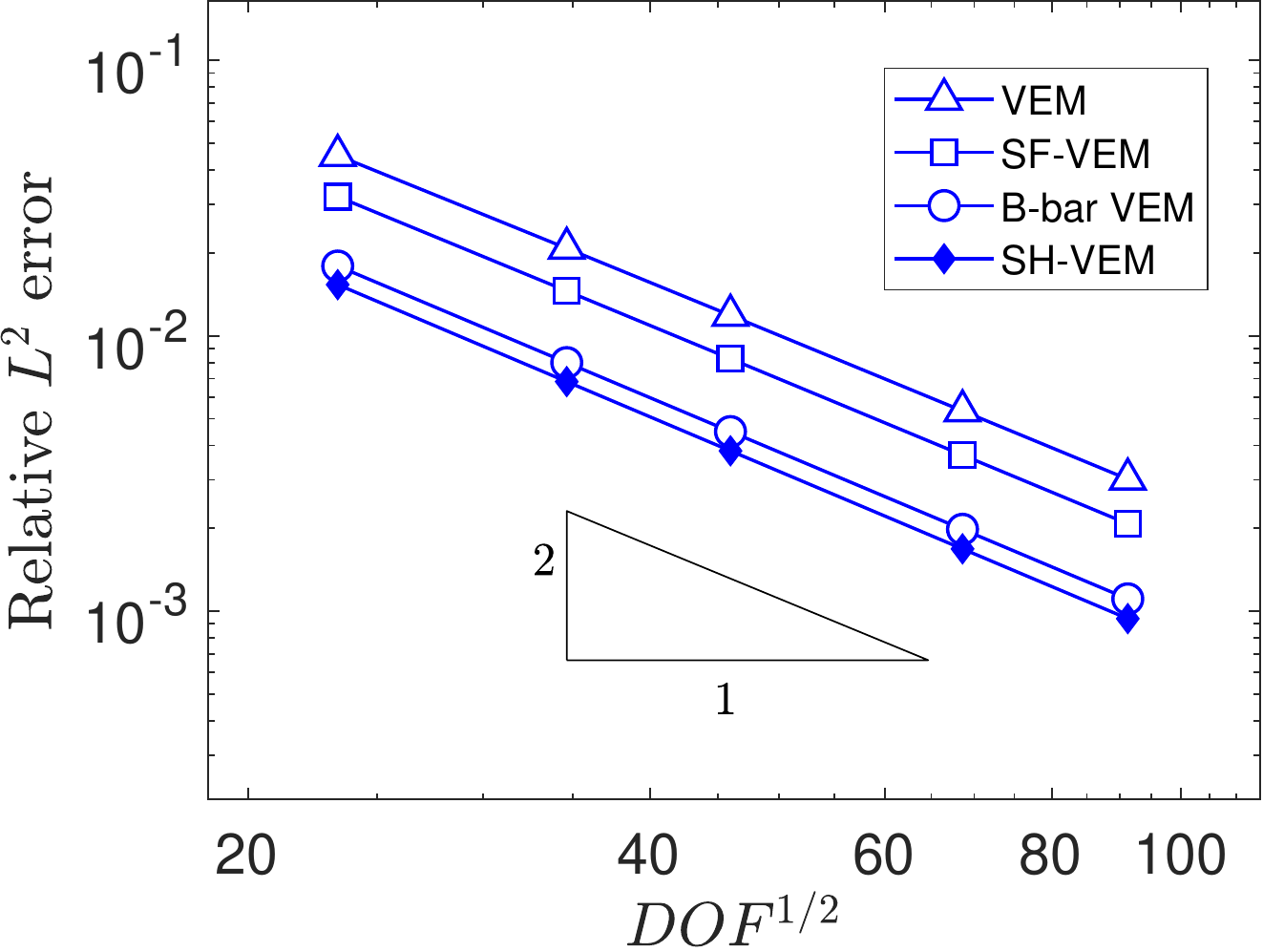}
     \end{subfigure}
     \hfill
     \begin{subfigure}{.24\textwidth}
         \centering
         \includegraphics[width=\textwidth]{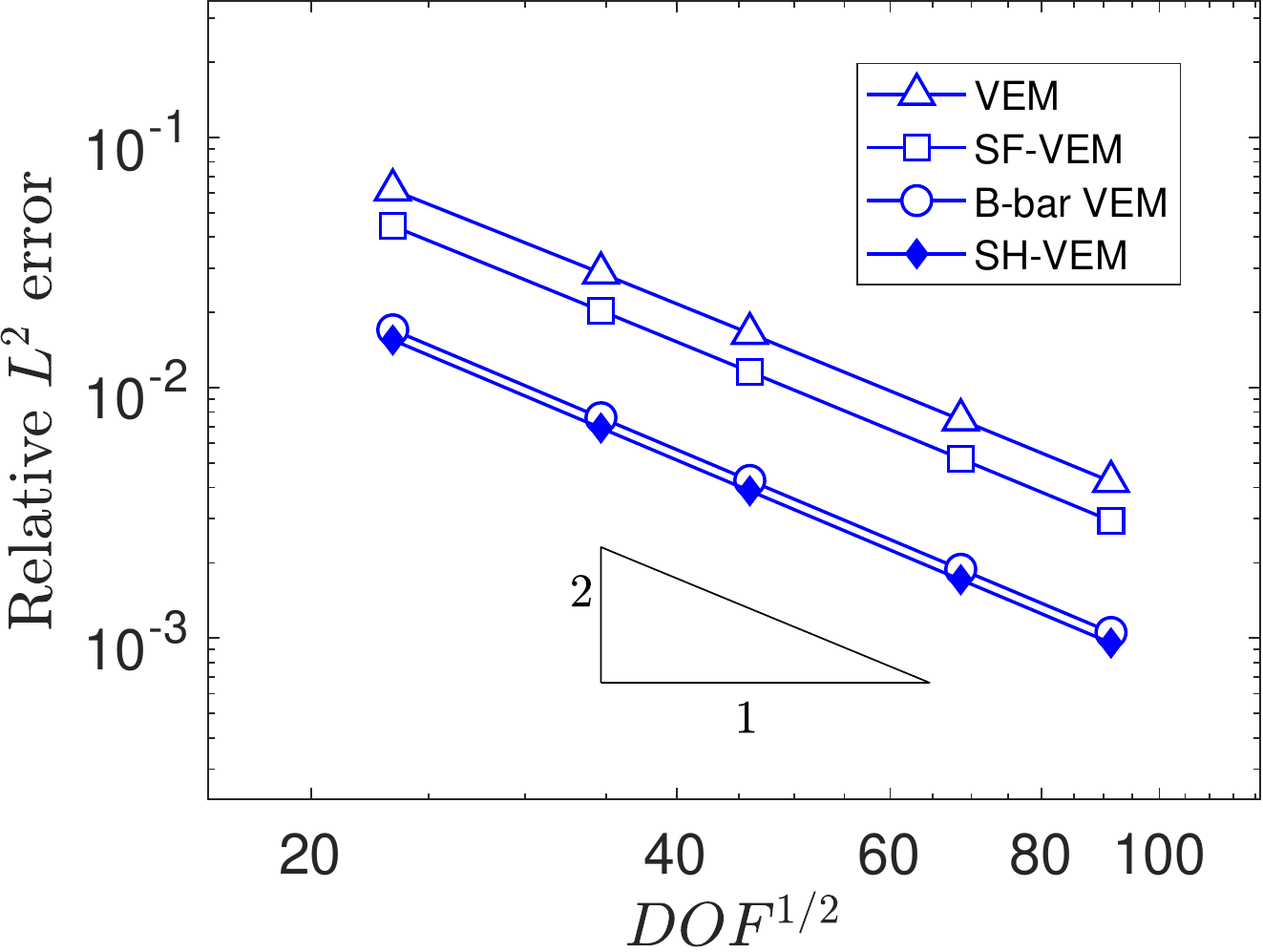}
     \end{subfigure}
     \hfill
     \begin{subfigure}{.24\textwidth}
         \centering
         \includegraphics[width=\textwidth]{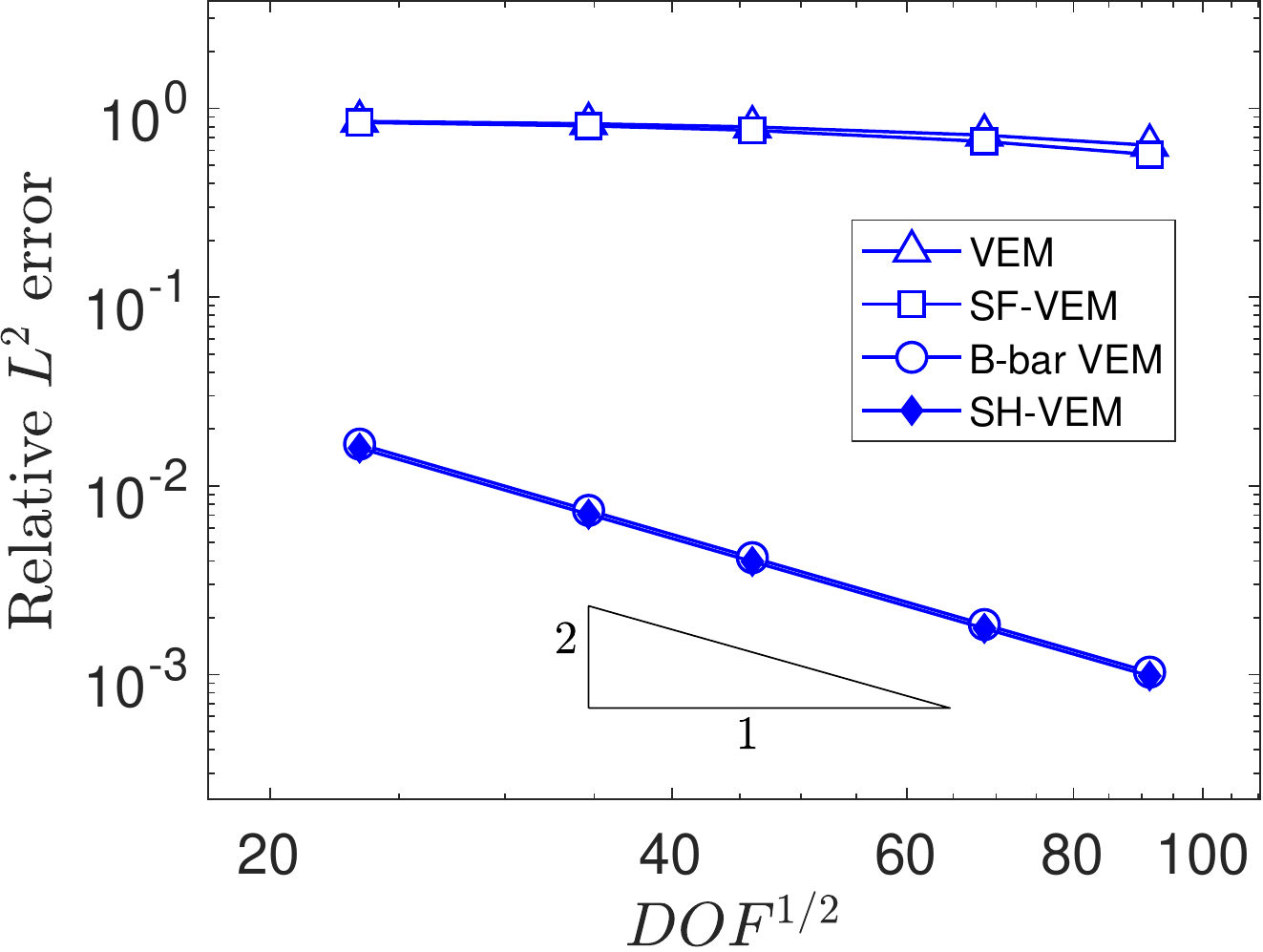}
     \end{subfigure}
          \hfill
     \begin{subfigure}{.24\textwidth}
         \centering
         \includegraphics[width=\textwidth]{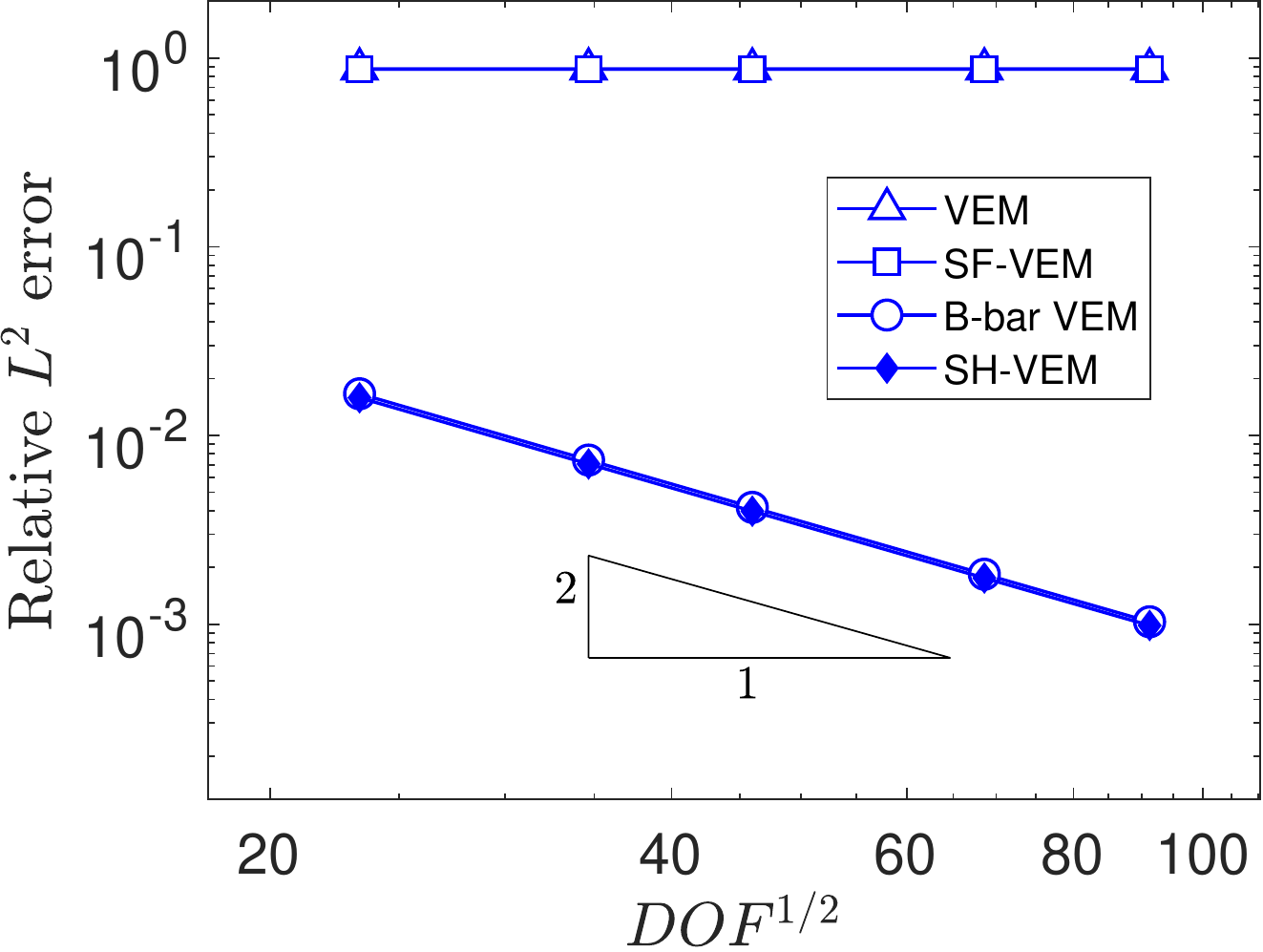}
     \end{subfigure}
     \vfill
          \begin{subfigure}{.24\textwidth}
         \centering
         \includegraphics[width=\textwidth]{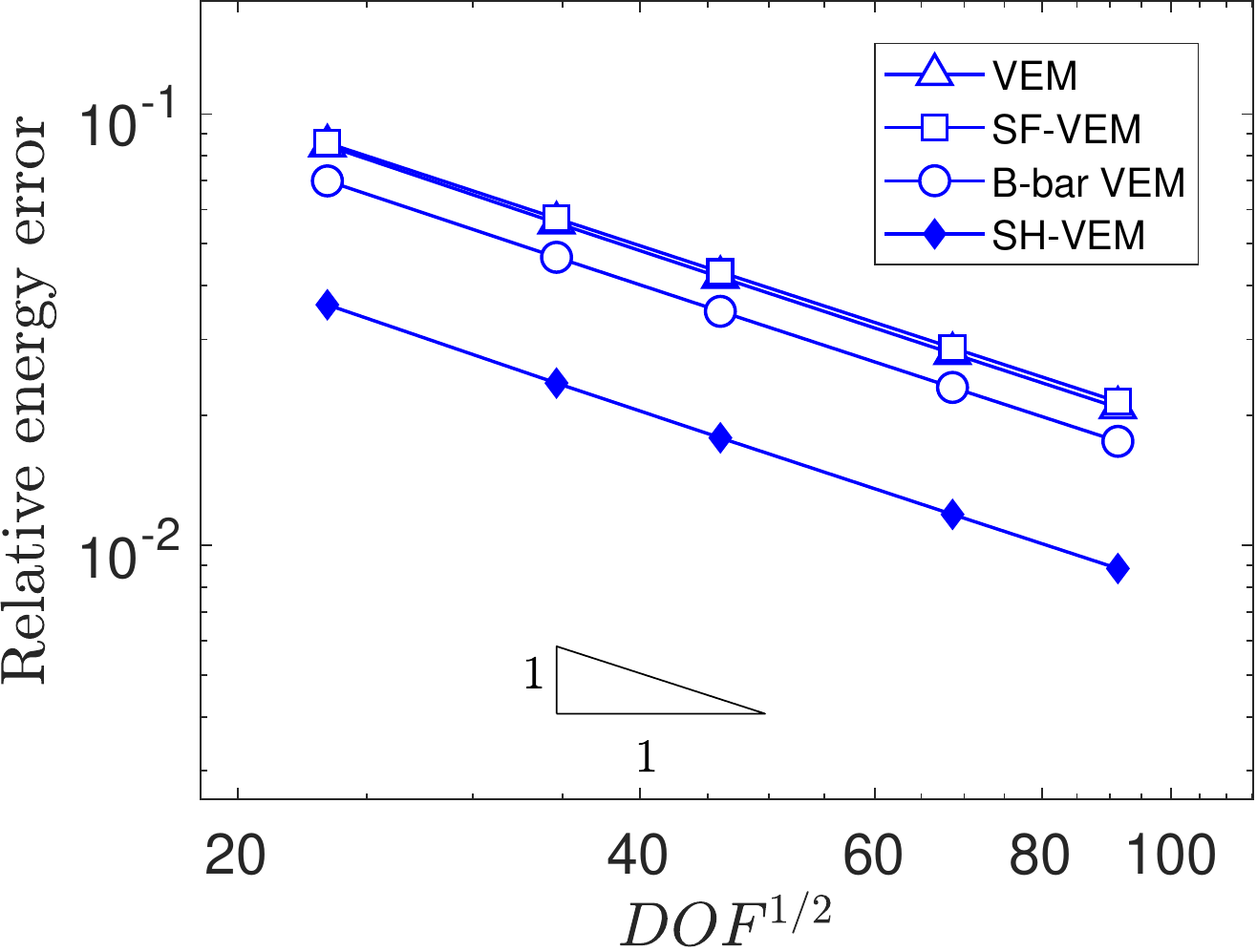}
         \caption{}
     \end{subfigure}
     \hfill
     \begin{subfigure}{.24\textwidth}
         \centering
         \includegraphics[width=\textwidth]{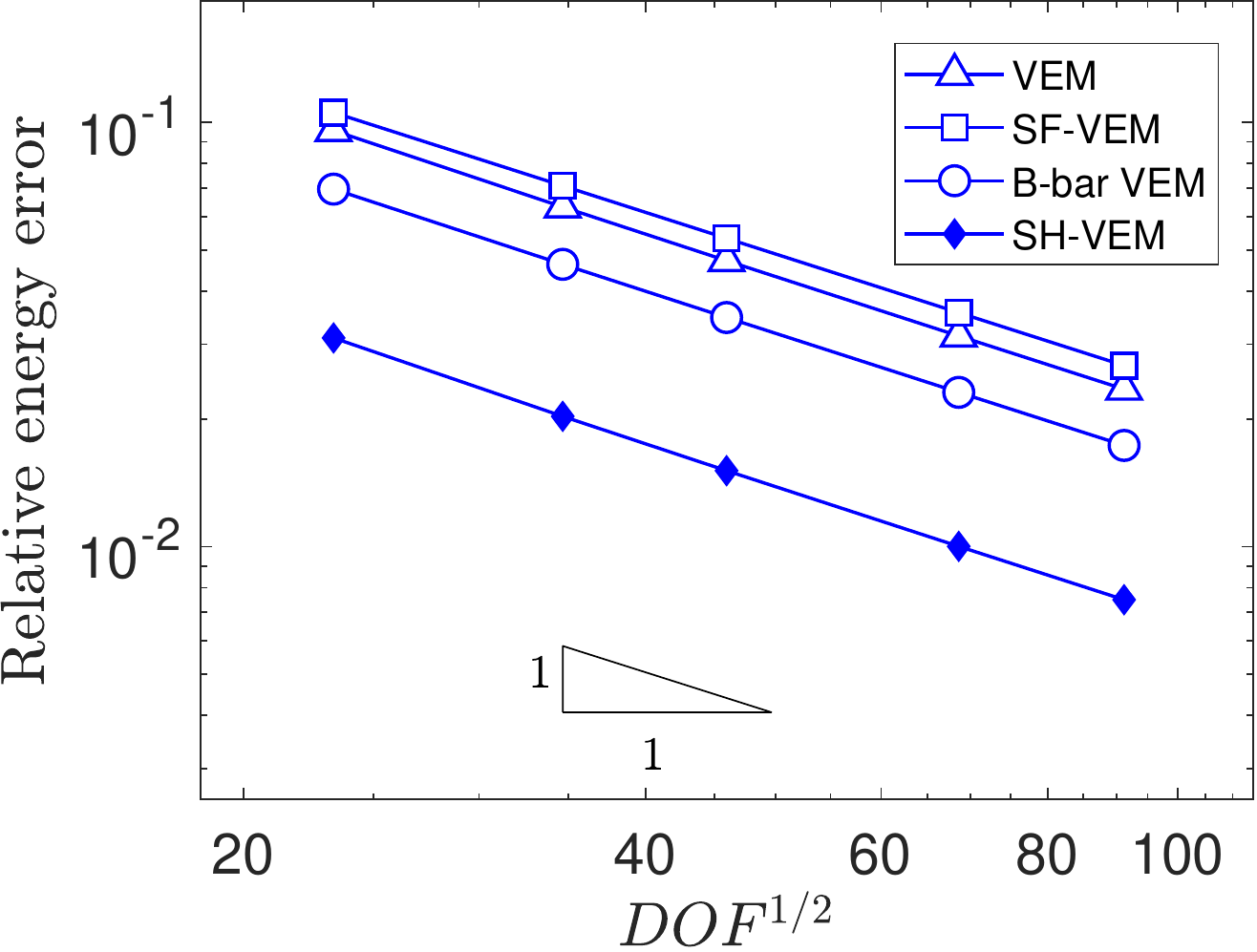}
         \caption{}
     \end{subfigure}
     \hfill
     \begin{subfigure}{.24\textwidth}
         \centering
         \includegraphics[width=\textwidth]{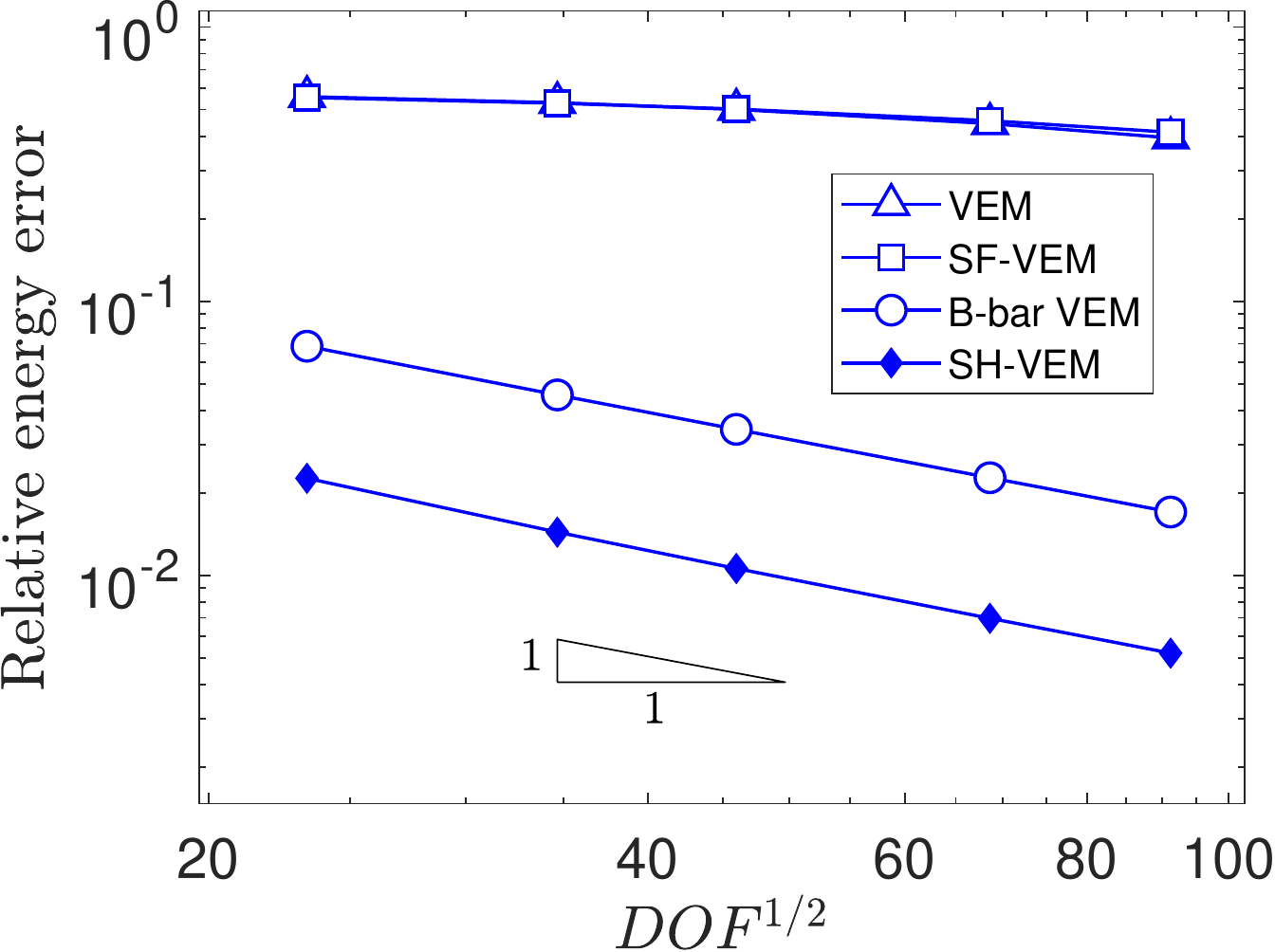}
         \caption{}
     \end{subfigure}
          \hfill
     \begin{subfigure}{.24\textwidth}
         \centering
         \includegraphics[width=\textwidth]{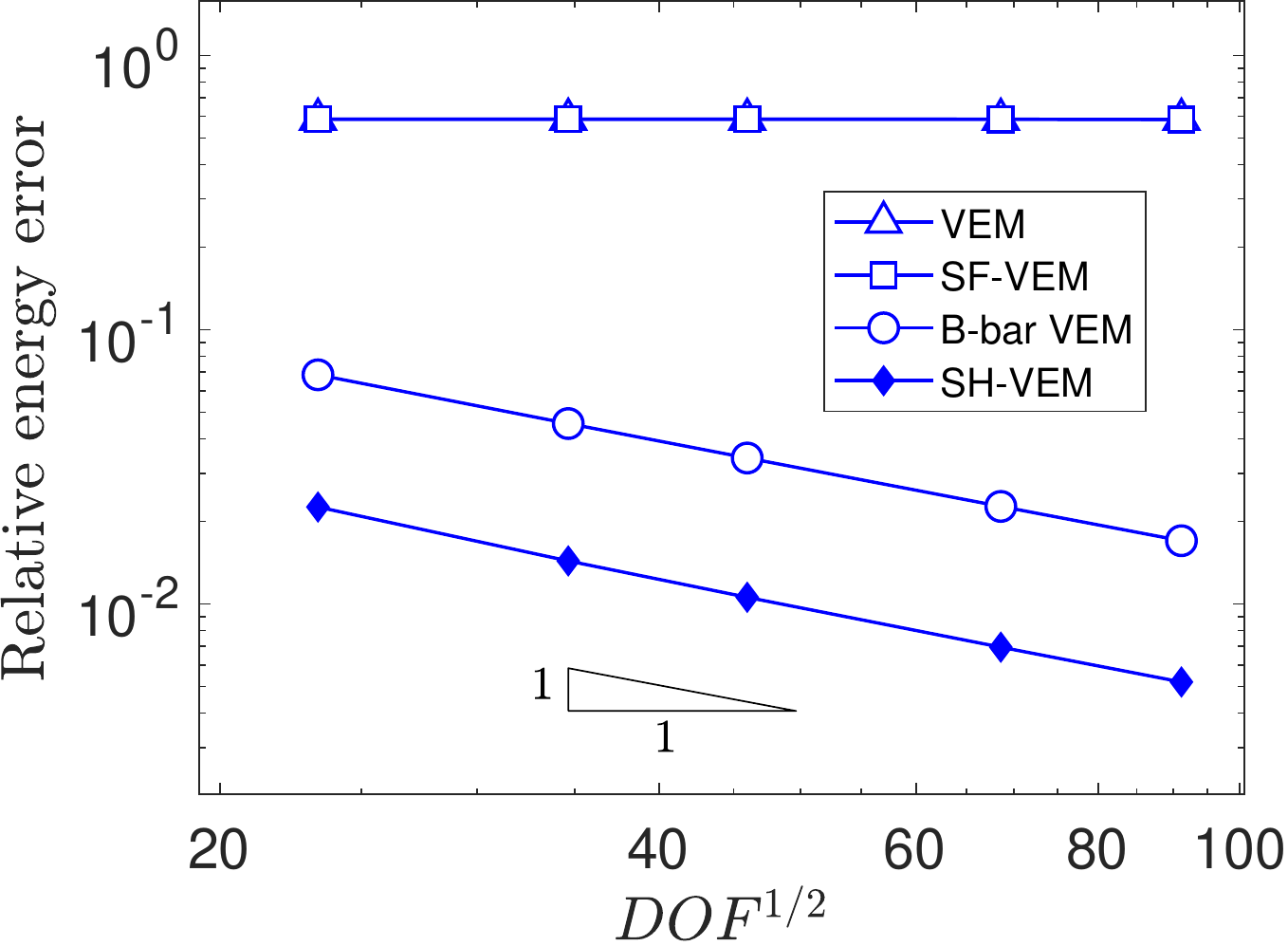}
         \caption{}
     \end{subfigure}
        \caption{Comparison of the convergence of standard VEM, stabilization-free VEM, B-bar VEM and SH-VEM 
        for the manufactured problem on nonconvex meshes \acrev{(see~\fref{fig:squaremesh_nonconvex})}. Each column represents a different value of $\nu$. (a) $\nu=0.3$, (b) $\nu=0.4$, (c) $\nu=0.4999$ and (d) $\nu=0.4999999$. }
        \label{fig:manufactured_nonconvex}
\end{figure}

We also test this problem on noncovex meshes. We begin with a uniform rectangular mesh and then split each element into a convex and a nonconvex quadrilateral. A few sample 
meshes are
shown in Figure~\ref{fig:squaremesh_nonconvex}. Numerical results are presented in Figure~\ref{fig:manufactured_nonconvex}, which reveal that
even on nonconvex meshes B-bar VEM and SH-VEM retain optimal rates of convergence.

\begin{figure}[!h]
     \centering
     \begin{subfigure}{.24\textwidth}
         \centering
         \includegraphics[width=\textwidth]{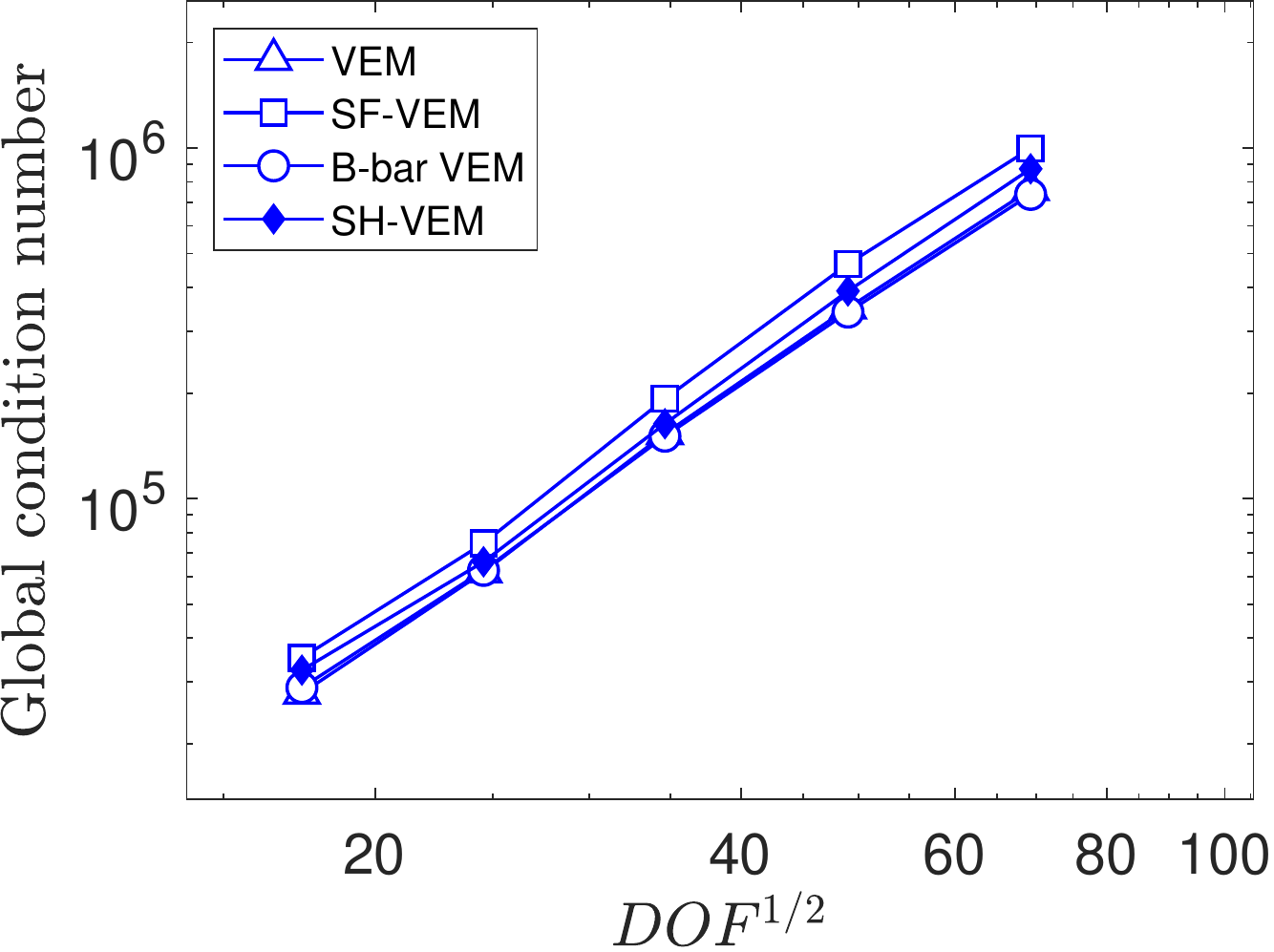}
     \end{subfigure}
     \hfill
     \begin{subfigure}{.24\textwidth}
         \centering
         \includegraphics[width=\textwidth]{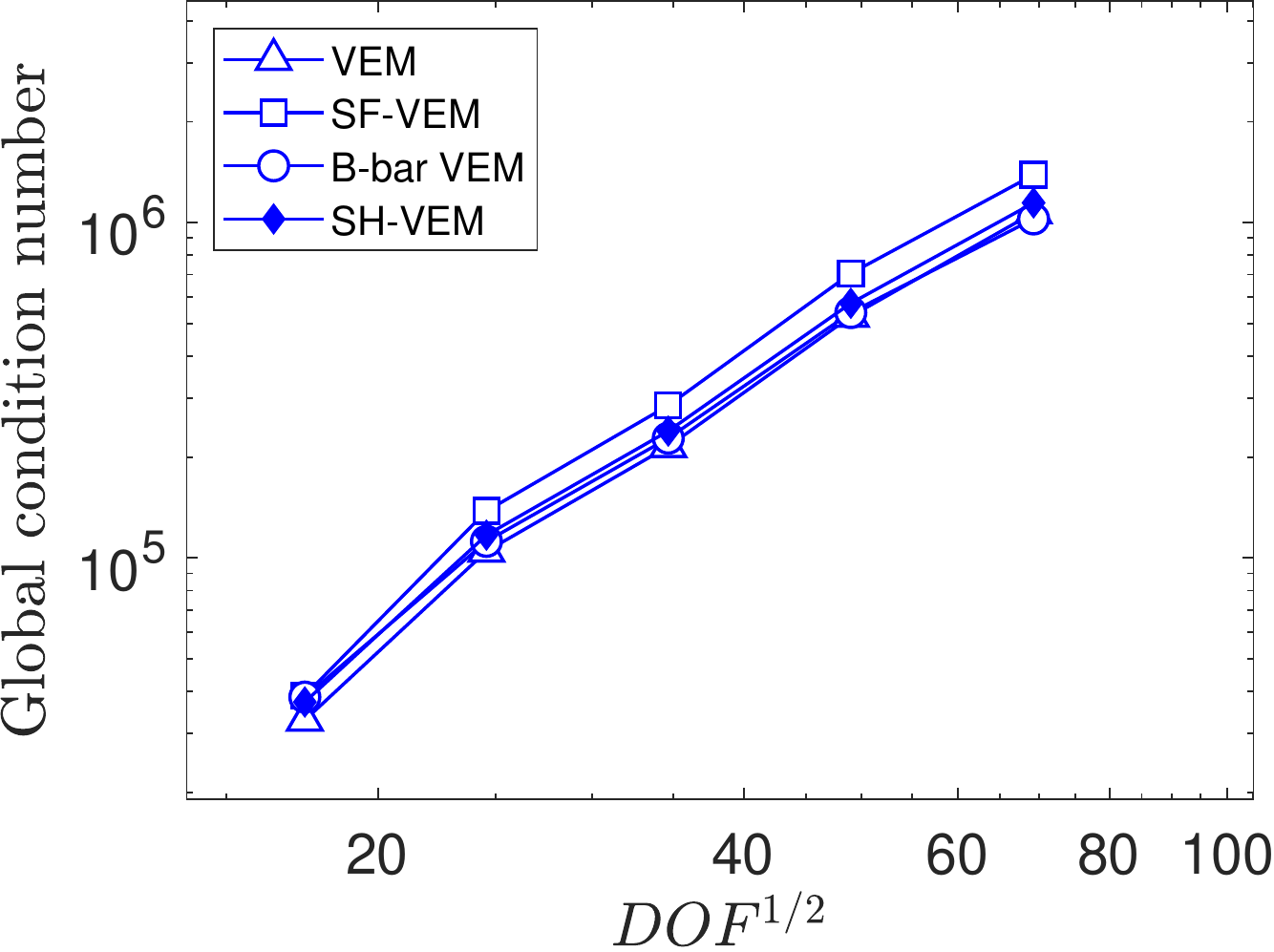}
     \end{subfigure}
     \hfill
     \begin{subfigure}{.24\textwidth}
         \centering
         \includegraphics[width=\textwidth]{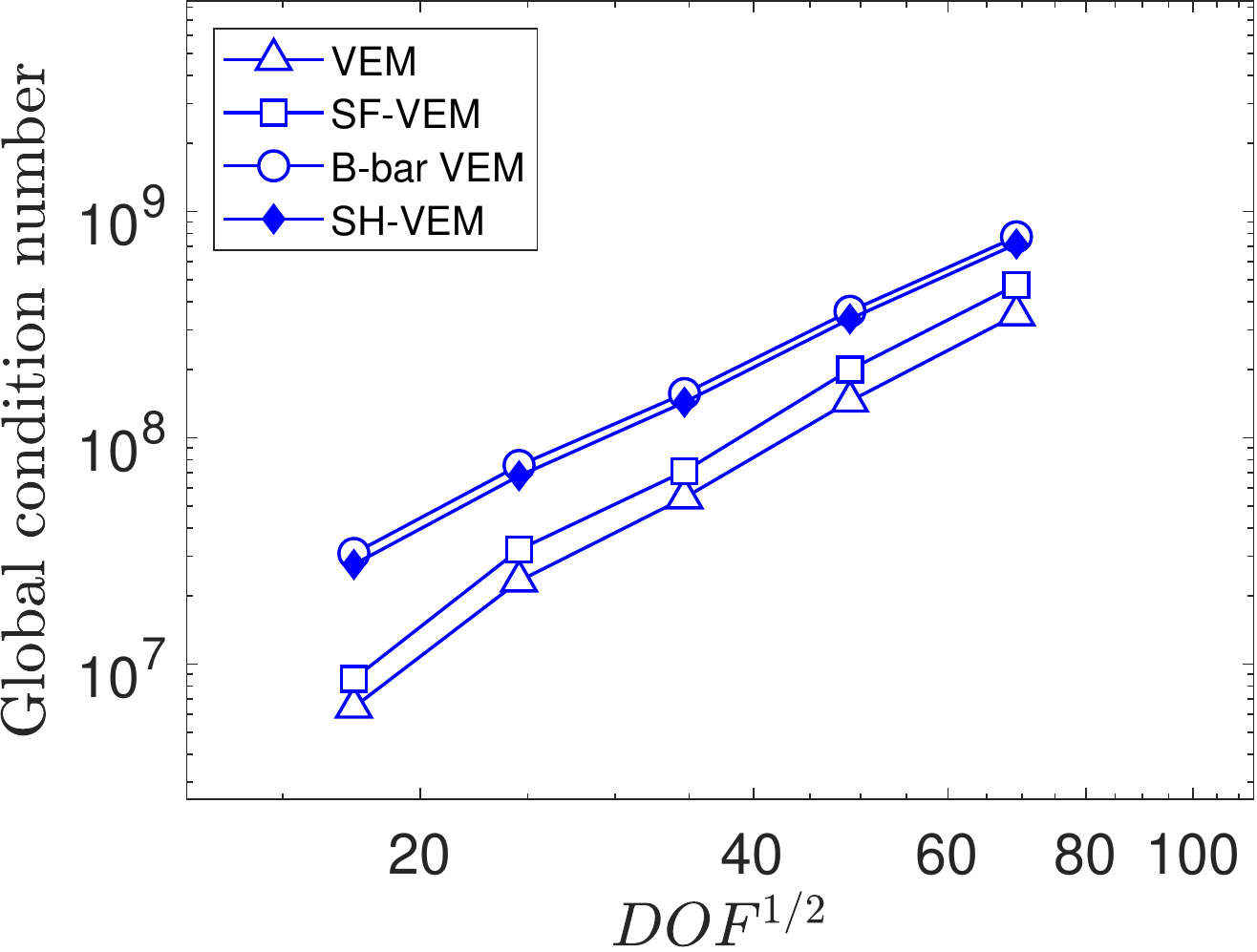}
     \end{subfigure}
          \hfill
     \begin{subfigure}{.24\textwidth}
         \centering
         \includegraphics[width=\textwidth]{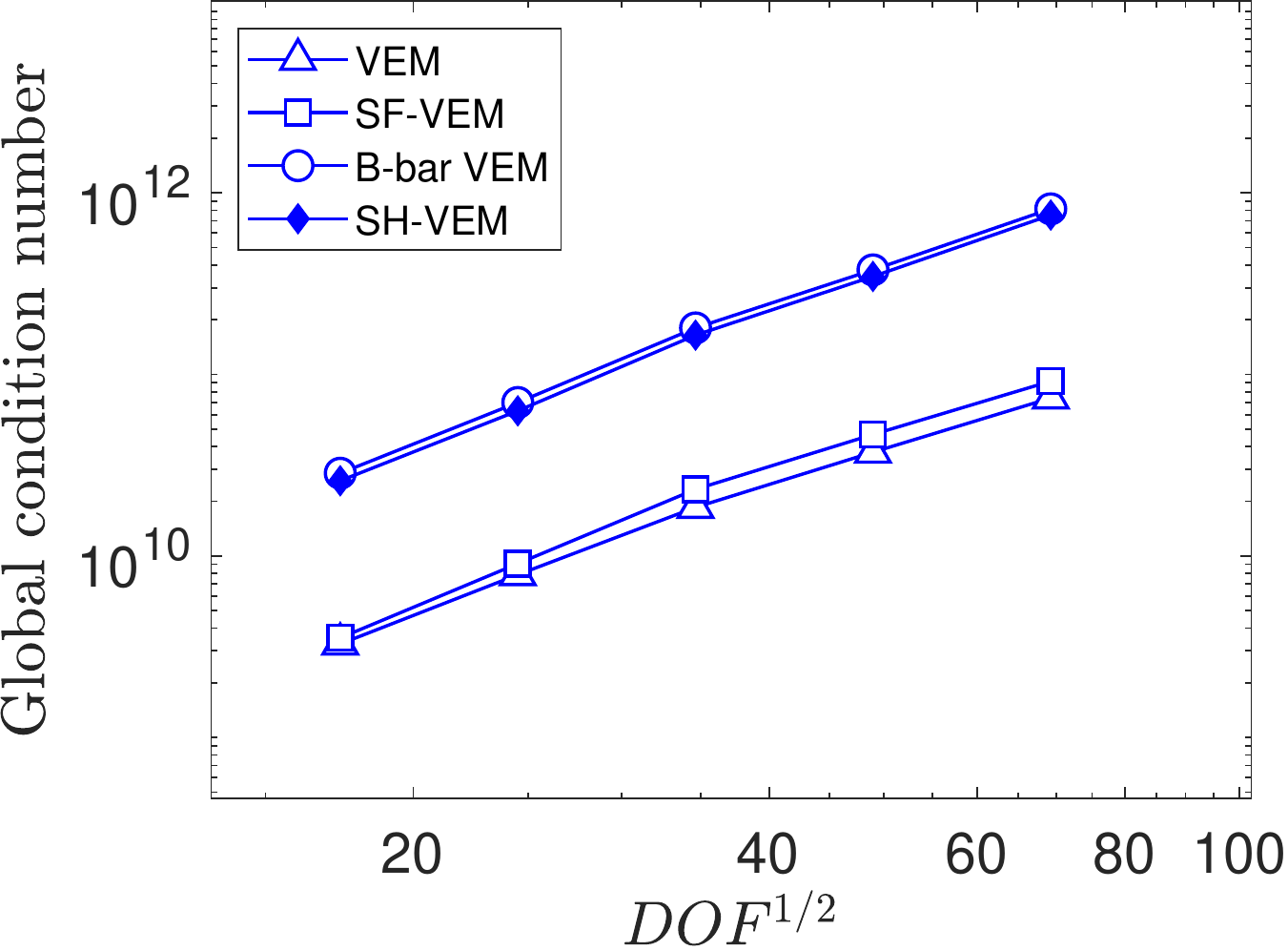}
     \end{subfigure}
     \vfill
          \begin{subfigure}{.24\textwidth}
         \centering
         \includegraphics[width=\textwidth]{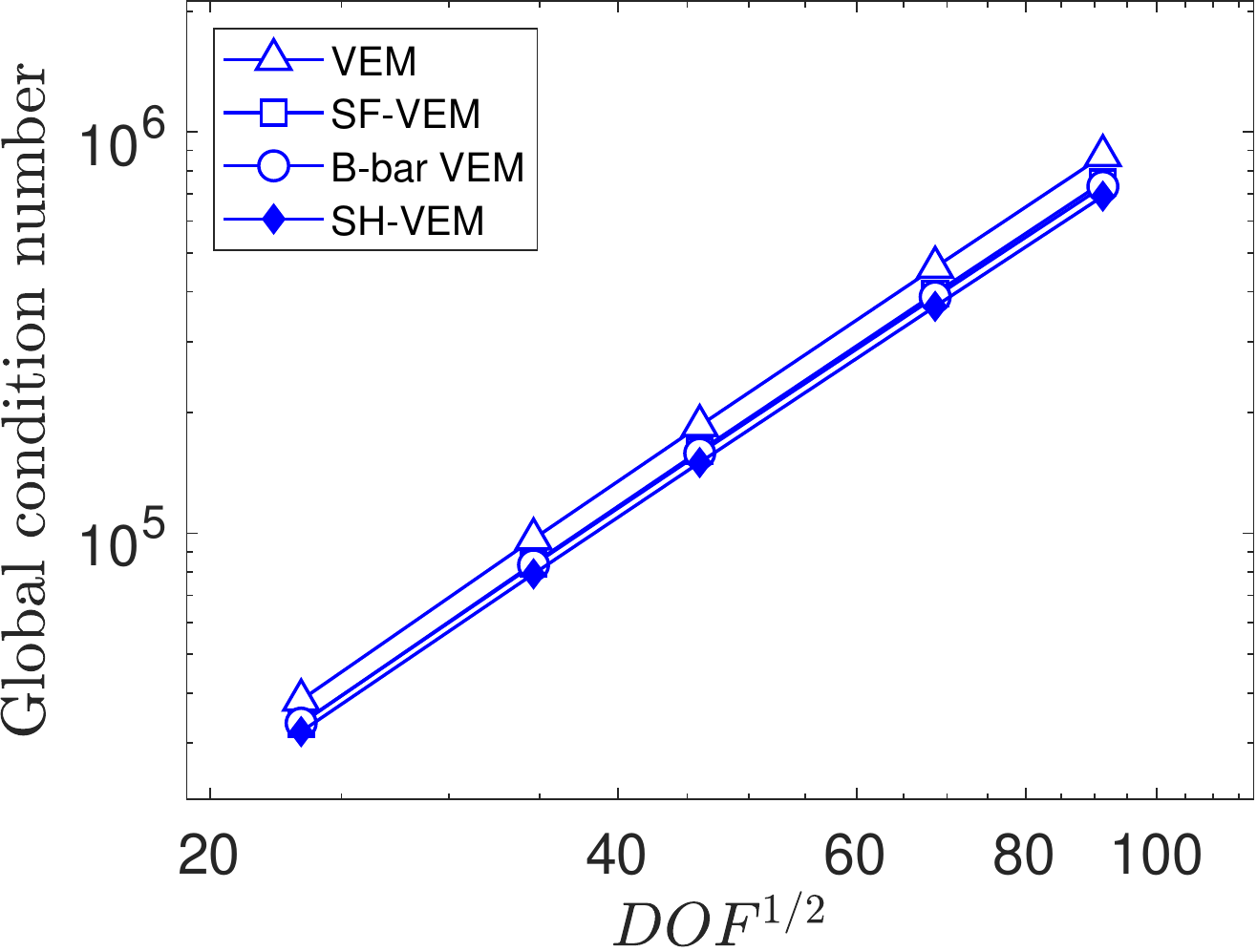}
         \caption{}
     \end{subfigure}
     \hfill
     \begin{subfigure}{.24\textwidth}
         \centering
         \includegraphics[width=\textwidth]{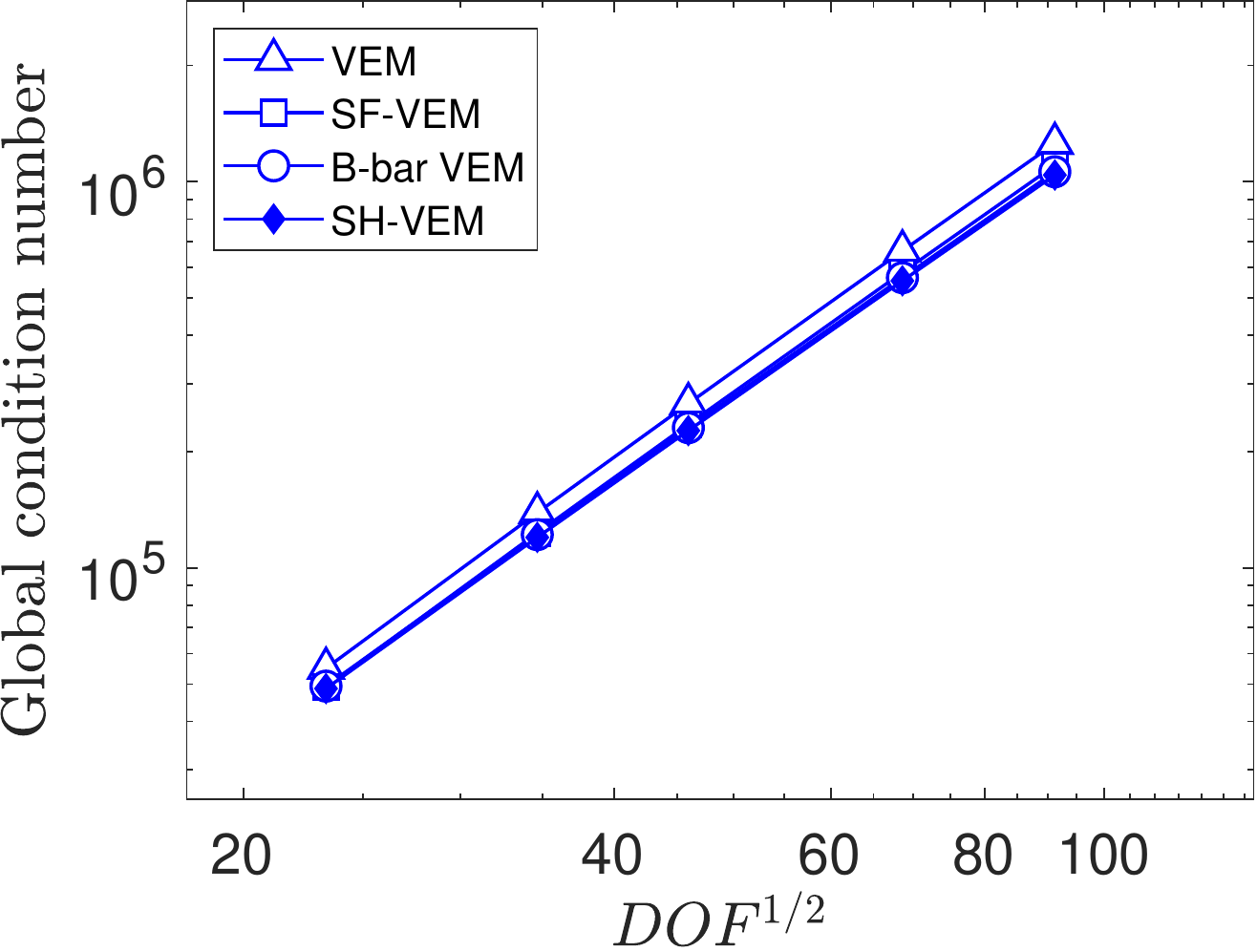}
         \caption{}
     \end{subfigure}
     \hfill
     \begin{subfigure}{.24\textwidth}
         \centering
         \includegraphics[width=\textwidth]{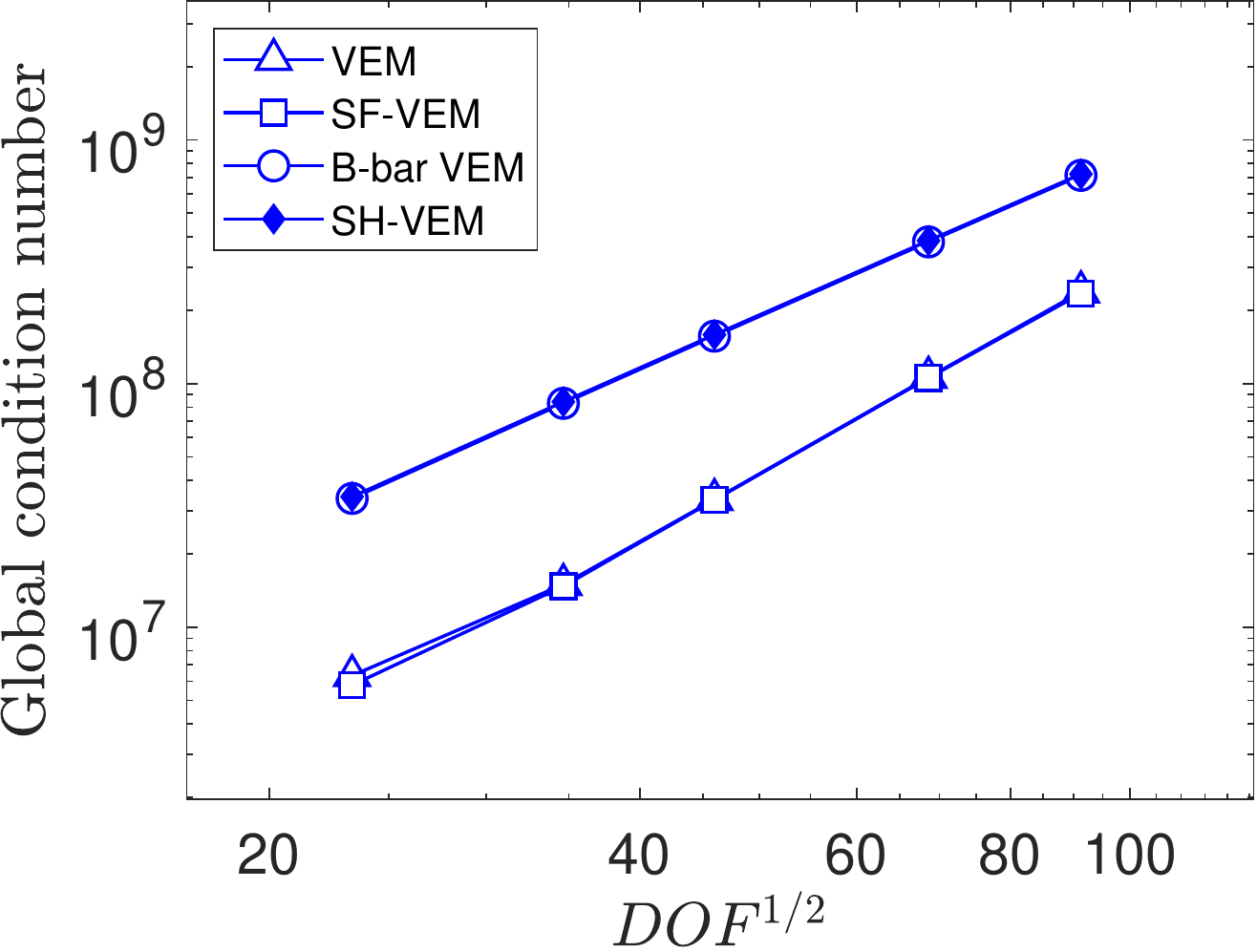}
         \caption{}
     \end{subfigure}
          \hfill
     \begin{subfigure}{.24\textwidth}
         \centering
         \includegraphics[width=\textwidth]{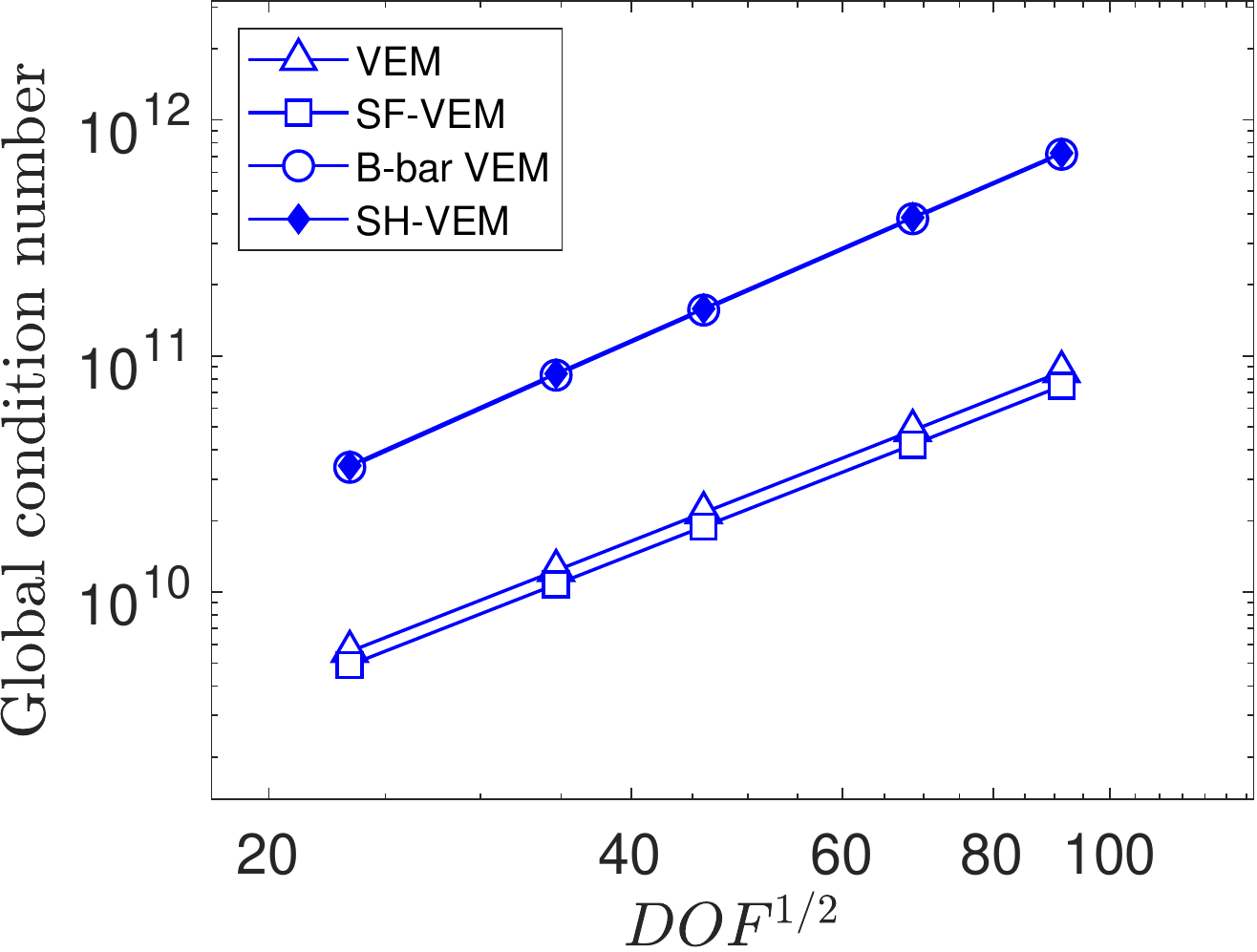}
         \caption{}
     \end{subfigure}
        \caption{Comparison of the conditioning of the global stiffness matrix  of the 
        standard VEM, stabilization-free VEM, B-bar VEM and SH-VEM for the manufactured problem. The first row is for the unstructured quadrilateral mesh \acrev{(see~\fref{fig:squaremesh})}, the second row is for the nonconvex mesh \acrev{(see~\fref{fig:squaremesh_nonconvex})}. Each column represents a different value of $\nu$. (a) $\nu=0.3$, (b) $\nu=0.4$, (c) $\nu=0.4999$ and (d) $\nu=0.4999999$. }
        \label{fig:manufactured_conditioning}
\end{figure}

Lastly, we examine the conditioning of the global stiffness matrix  to ensure that increasing the Poisson's ratio and varying 
the element shapes and refinement does not lead to ill-conditioning. In Figure~\ref{fig:manufactured_conditioning}, we show the condition number of the four methods as $\nu \to 0.5$ on both unstructured and nonconvex meshes. From the plot, we observe that the condition number of SH-VEM is comparable to the other three methods for compressible materials. As the material becomes nearly incompressible, the condition number increases for all methods on the coarsest mesh. However, we observe that the 
growth of the condition number 
with refinement for B-bar VEM and SH-VEM is similar to the
case when $\nu = 0.3$, which is in agreement with the 
${\cal O}(h^{-2})$ increase of the stiffness matrix condition number in the finite element method. 

\subsection{Thin cantilever beam}
\begin{figure}
     \centering
     \begin{subfigure}{\textwidth}
         \centering
         \includegraphics[width=\textwidth]{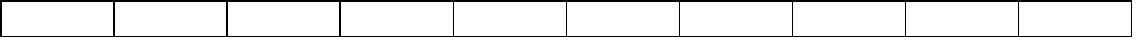}
         \caption{}
     \end{subfigure}
     \vfill
     \begin{subfigure}{\textwidth}
         \centering
         \includegraphics[width=\textwidth]{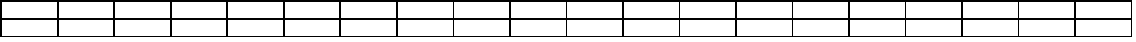}
         \caption{}
     \end{subfigure}
     \vfill
     \begin{subfigure}{\textwidth}
         \centering
         \includegraphics[width=\textwidth]{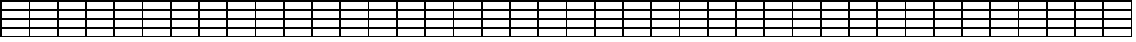}
         \caption{}
     \end{subfigure}
        \caption{Rectangular meshes for the cantilever beam problem. (a) 10 elements, (b)  40 eleme nts and (c) 160 elements.  }
        \label{fig:beammesh}
 \end{figure}

\begin{figure}[!h]
     \centering
     \begin{subfigure}{0.32\textwidth}
         \centering
         \includegraphics[width=\textwidth]{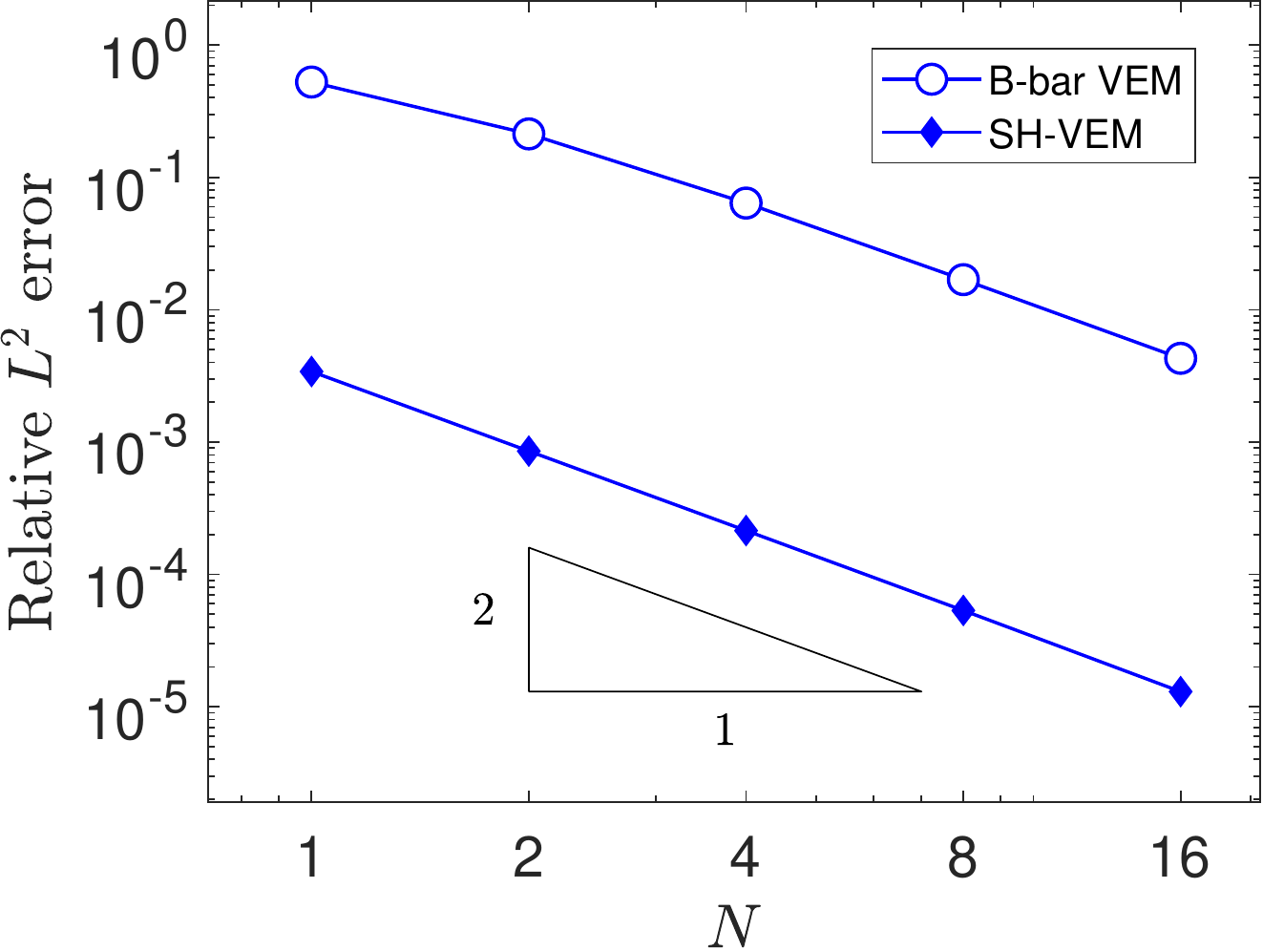}
         \caption{}
     \end{subfigure}
     \hfill
     \begin{subfigure}{0.32\textwidth}
         \centering
         \includegraphics[width=\textwidth]{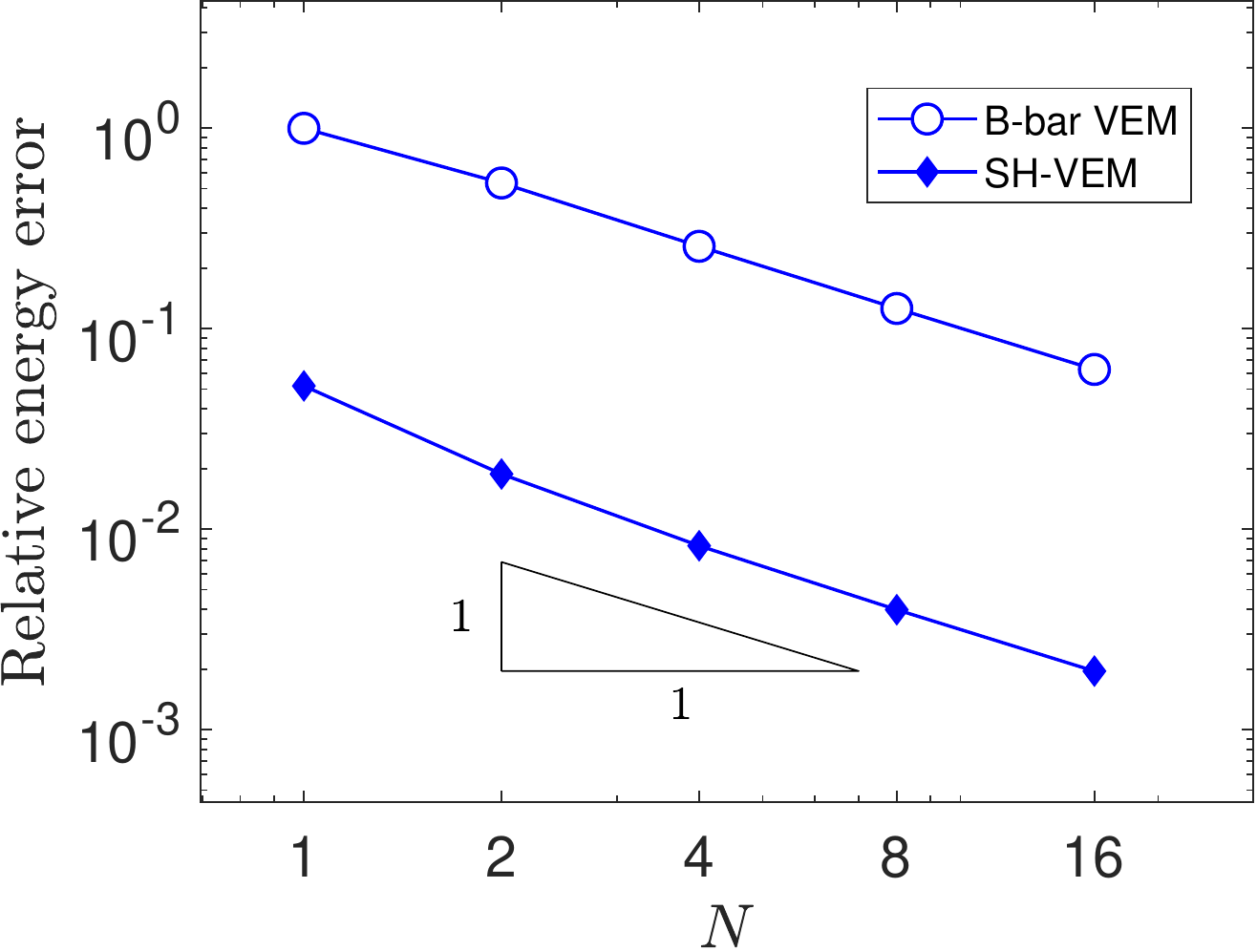}
         \caption{}
     \end{subfigure}
     \hfill
     \begin{subfigure}{0.32\textwidth}
         \centering
         \includegraphics[width=\textwidth]{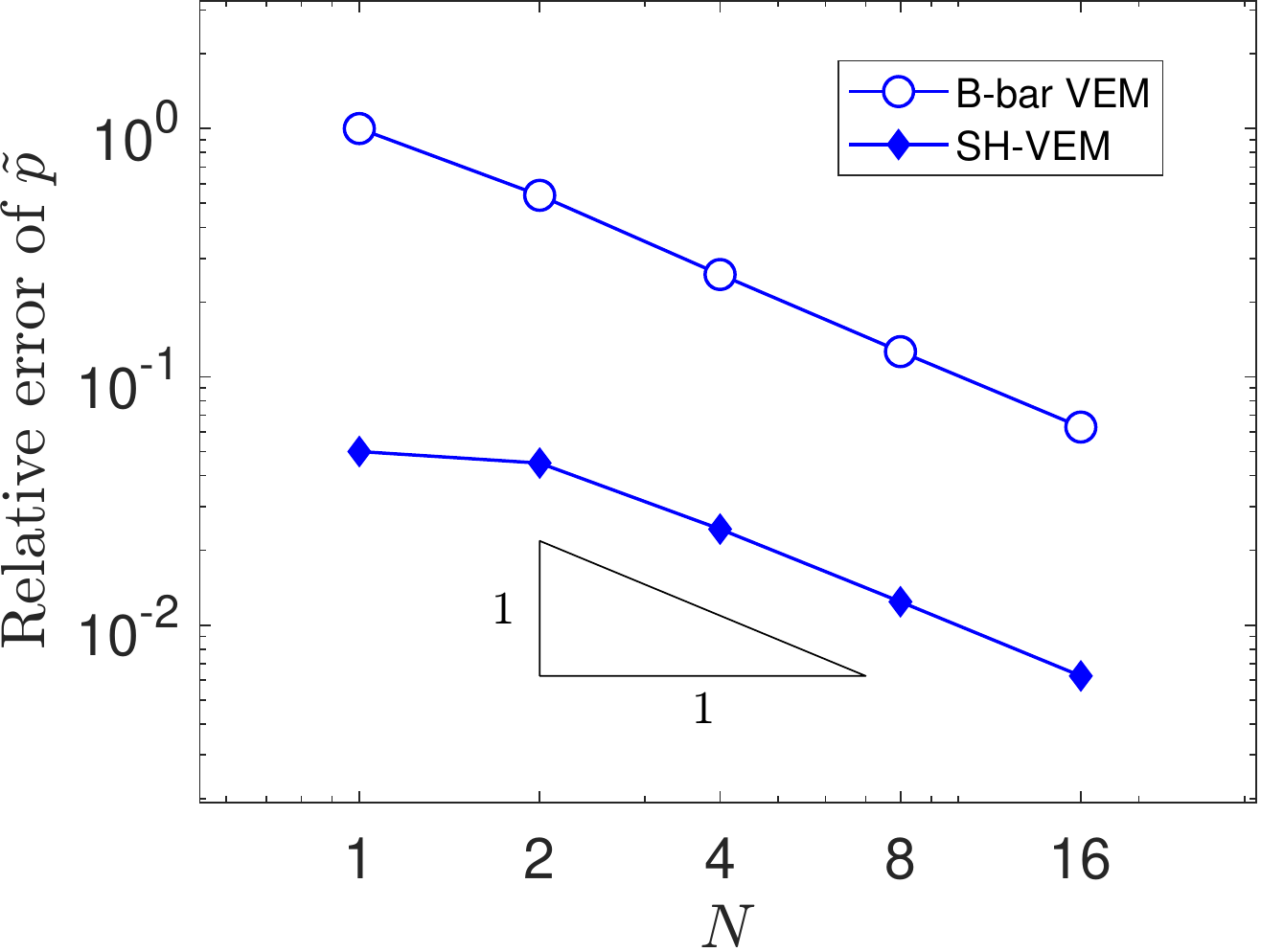}
         \caption{}
     \end{subfigure}
        \caption{Comparison of B-bar VEM and SH-VEM for the thin cantilever beam problem on structured meshes \acrev{(see~\fref{fig:beammesh})}. (a) $L^2$ error of displacement, (b) energy error and (c) $L^2$ error of hydrostatic stress, where $N$ is the number of elements along the height of the beam. }
        \label{fig:beam_convergence}
\end{figure}

\begin{figure}[!h]
    \centering
    \begin{subfigure}{0.48\textwidth}
    \centering
    \includegraphics[width=\textwidth]{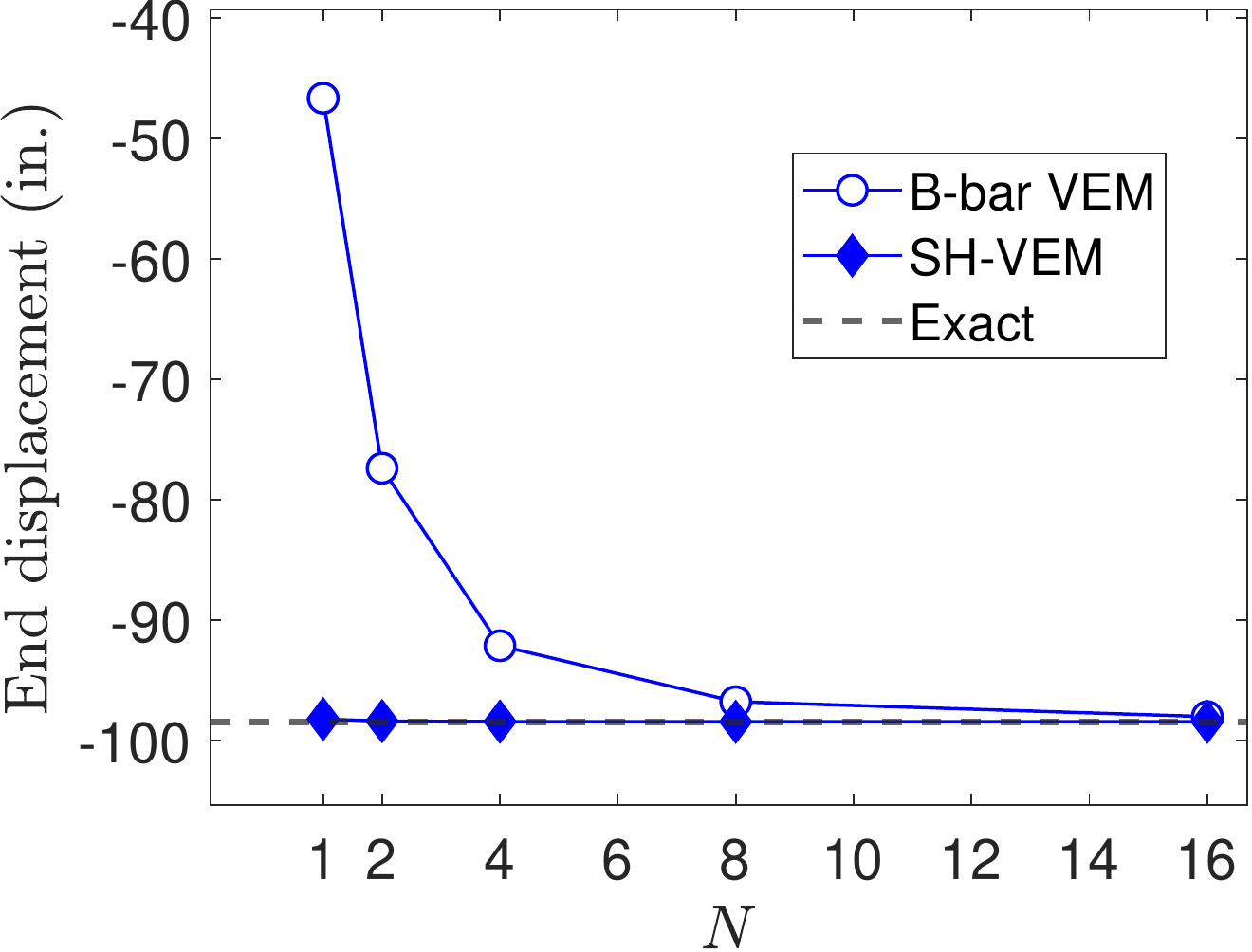}
    \caption{}
    \end{subfigure}
     \hfill
     \begin{subfigure}{0.48\textwidth}
         \centering
         \includegraphics[width=\textwidth]{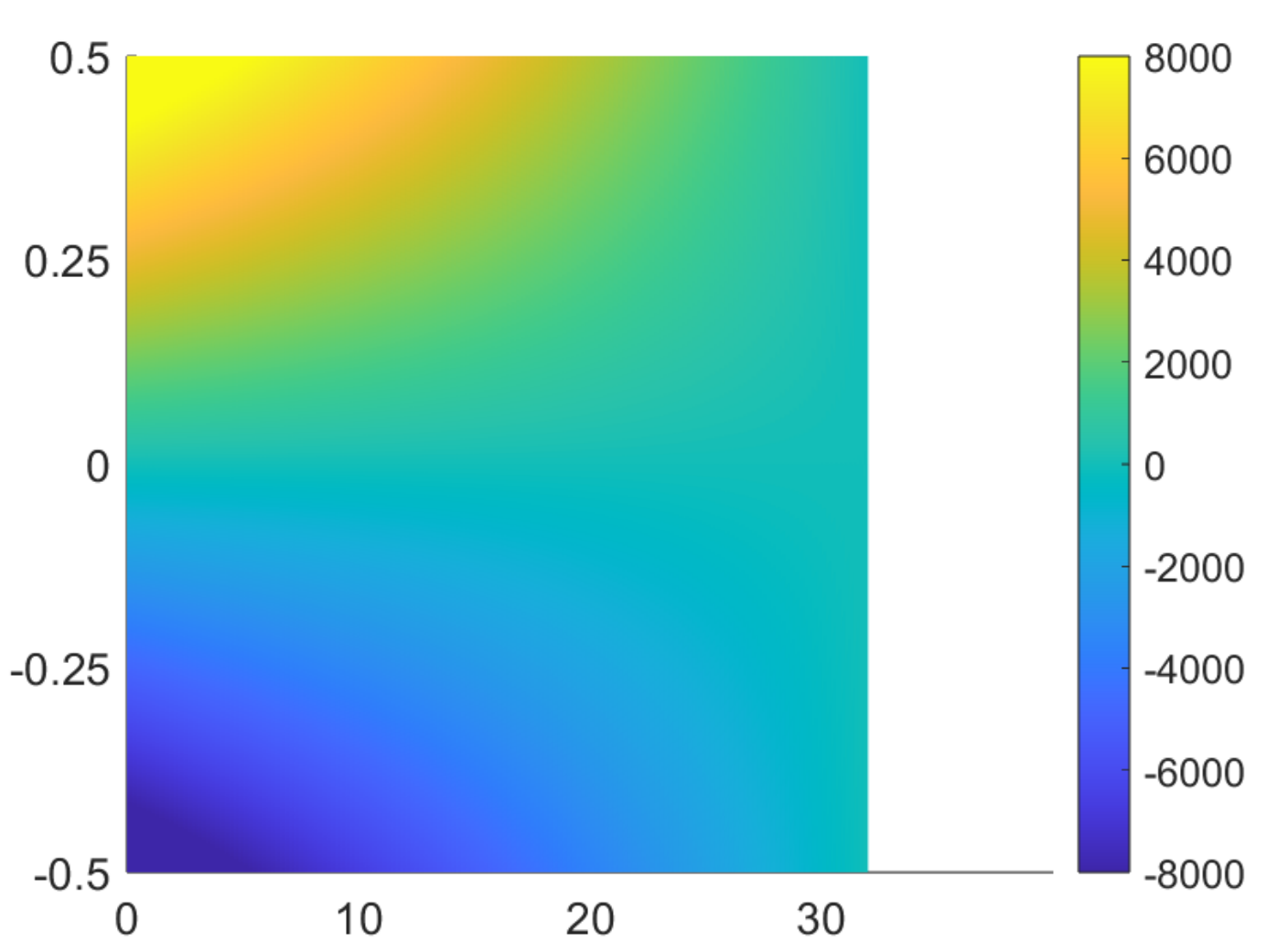}
         \caption{}
     \end{subfigure}
    \caption{(a) Convergence of the end displacement for the cantilever beam problem. The mesh consists of $10N\times N$ rectangular
    elements, where $N$ is the number of elements along the height of the beam \acrev{(see~\fref{fig:beammesh})}, (b) contour plot
    of hydrostatic stress for SH-VEM. }
    \label{fig:tip_dispacement}
\end{figure}

We consider the benchmark 
problem of a thin cantilever beam under a shear end load.\cite{timoshenko1951theory} 
The material has Young's modulus $E_Y = 1\times 10^5$ psi and $\nu = 0.49995$. The beam has length $L = 32$ inch, height $D = 1$ inch and unit thickness. The left boundary is fixed and
a shear end load of $P = -100$ lbf is applied on the right boundary. We use a regular rectangular mesh with $N \in \{1,\, 2,\, 4,\, 8,\, 16\}$ elements along the height and $10N$ elements along the length. In Figure~\ref{fig:beammesh},
we show a few representative meshes and in Figure~\ref{fig:beam_convergence} we compare the rates of convergence of B-bar VEM to SH-VEM in the three error norms. In Figure~\ref{fig:tip_dispacement}, we plot the end displacement of the three methods and contours 
of the hydrostatic stress for SH-VEM. From these results, we observe
that the accuracy of SH-VEM is far superior to B-bar VEM and the displacements in the
SH-VEM display superconvergence (close to the
exact solution) on coarse rectangular meshes. 

\begin{figure}
\begin{minipage}{.48\textwidth}
     \centering
     \begin{subfigure}{\textwidth}
         \centering
         \includegraphics[width=\textwidth]{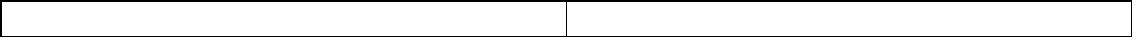}
         \caption{}
     \end{subfigure}
     \vfill
     \begin{subfigure}{\textwidth}
         \centering
         \includegraphics[width=\textwidth]{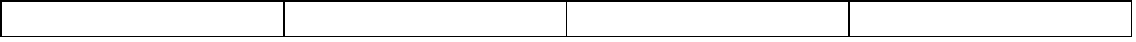}
         \caption{}
     \end{subfigure}
     \vfill
     \begin{subfigure}{\textwidth}
         \centering
         \includegraphics[width=\textwidth]{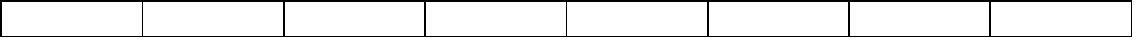}
         \caption{}
     \end{subfigure}
          \vfill
     \begin{subfigure}{\textwidth}
         \centering
         \includegraphics[width=\textwidth]{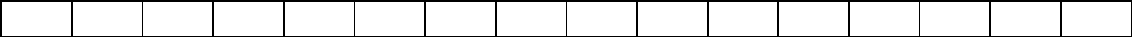}
         \caption{}
     \end{subfigure}
     \end{minipage}
     \hfill
     \begin{minipage}{.48\textwidth}
     \centering
     \begin{subfigure}{\textwidth}
         \centering
         \includegraphics[width=\textwidth]{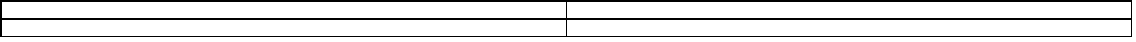}
         \caption{}
     \end{subfigure}
     \vfill
     \begin{subfigure}{\textwidth}
         \centering
         \includegraphics[width=\textwidth]{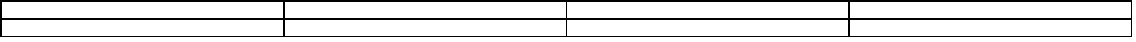}
         \caption{}
     \end{subfigure}
     \vfill
     \begin{subfigure}{\textwidth}
         \centering
         \includegraphics[width=\textwidth]{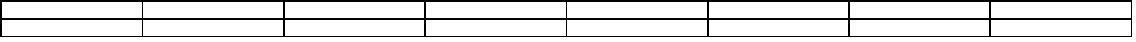}
         \caption{}
     \end{subfigure}
          \vfill
     \begin{subfigure}{\textwidth}
         \centering
         \includegraphics[width=\textwidth]{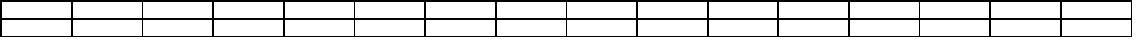}
         \caption{}
     \end{subfigure}
     \end{minipage}
        \caption{Rectangular meshes for the cantilever beam
        problem with fixed length to height ratio for each element.
        (a), (b), (c), (d) 
        16:1, 8:1, 4:1 and 2:1; and
        (e), (f), (g), (h) 32:1, 16:1, 8:1 and 4:1.}
        \label{fig:beammesh_fixed}
\end{figure}

\begin{figure}[!h]
    \centering
    \begin{subfigure}{0.48\textwidth}
    \centering
    \includegraphics[width=\textwidth]{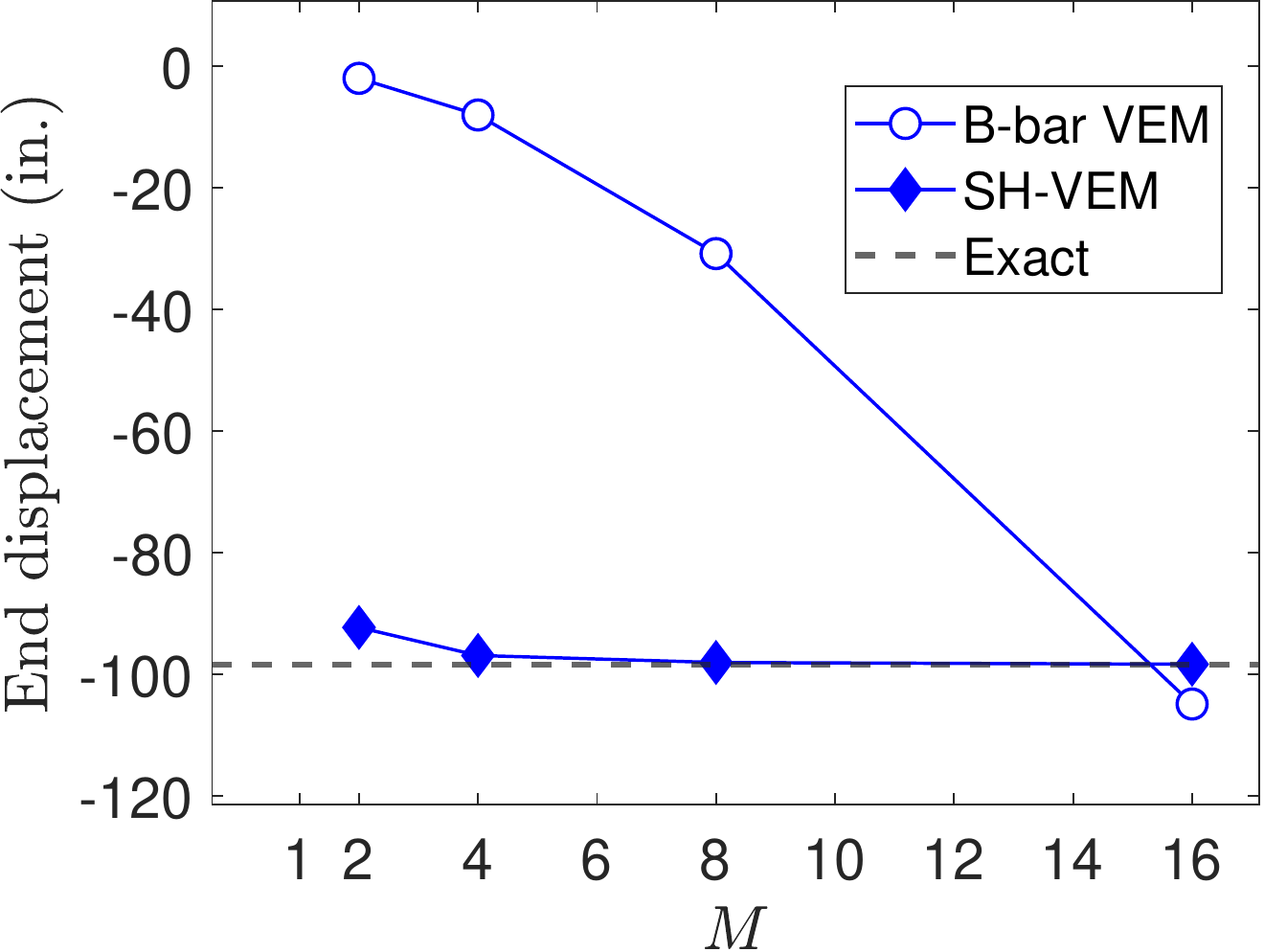}
    \caption{}\label{fig:tip_dispacement_fixed-a}
    \end{subfigure}
    \begin{subfigure}{0.48\textwidth}
    \centering
    \includegraphics[width=\textwidth]{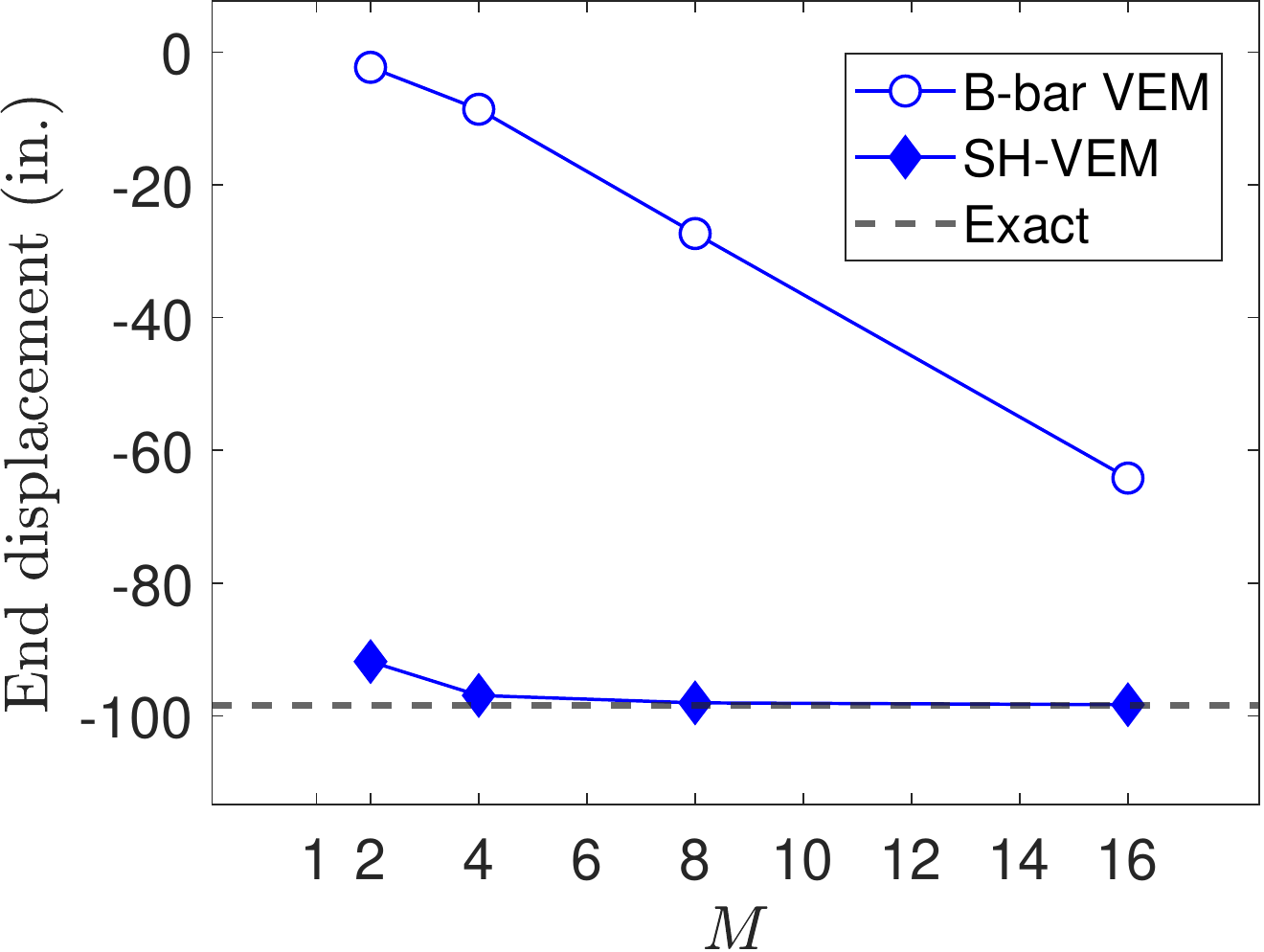}
    \caption{}\label{fig:tip_dispacement_fixed-b}
    \end{subfigure}
    \caption{Convergence of the end displacement for the
    cantilever beam problem. The mesh consists of $M\times N$ quadrilaterals, where $M$ is the number of elements along the length of the beam \acrev{(see~\fref{fig:beammesh_fixed})}. (a) $N = 1$ and (b)
    $N = 2$.}\label{fig:tip_dispacement_fixed}
\end{figure}

We consider another test for the cantilever beam problem using a mesh with 
either one or two elements along the height and $M$ elements along the length. We choose the number of elements $M \in \{2,\, 4,\, 8,\, 16\}$. The
meshes are depicted in Figure~\ref{fig:beammesh_fixed} and the convergence of the tip displacement for the two cases is presented in Figure~\ref{fig:tip_dispacement_fixed}. The plots reveal that SH-VEM is accurate even for high aspect ratio elements and is free of shear locking. However, for one
element along the height and with refinement along the length, we
observe that
B-bar VEM converges to a value below the exact value
(see Figure~\ref{fig:tip_dispacement_fixed-a}) and for the case of
two elements along the height, the end displacement is not accurate
(see Figure~\ref{fig:tip_dispacement_fixed-b}).

\begin{figure}
\begin{minipage}{.48\textwidth}
     \centering
     \begin{subfigure}{\textwidth}
         \centering
         \includegraphics[width=\textwidth]{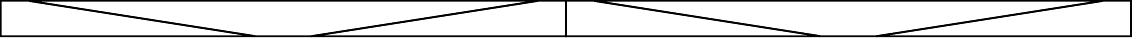}
         \caption{}
     \end{subfigure}
     \vfill
     \begin{subfigure}{\textwidth}
         \centering
         \includegraphics[width=\textwidth]{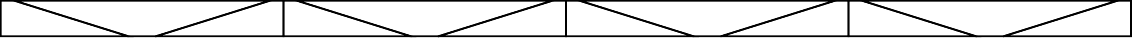}
         \caption{}
     \end{subfigure}
     \vfill
     \begin{subfigure}{\textwidth}
         \centering
         \includegraphics[width=\textwidth]{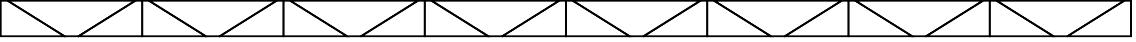}
         \caption{}
     \end{subfigure}
          \vfill
     \begin{subfigure}{\textwidth}
         \centering
         \includegraphics[width=\textwidth]{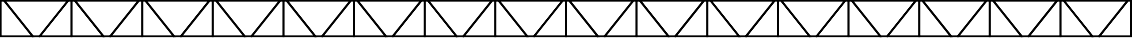}
         \caption{}
     \end{subfigure}
     \end{minipage}
     \hfill
     \begin{minipage}{.48\textwidth}
     \centering
     \begin{subfigure}{\textwidth}
         \centering
         \includegraphics[width=\textwidth]{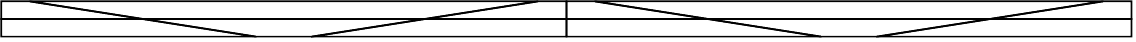}
         \caption{}
     \end{subfigure}
     \vfill
     \begin{subfigure}{\textwidth}
         \centering
         \includegraphics[width=\textwidth]{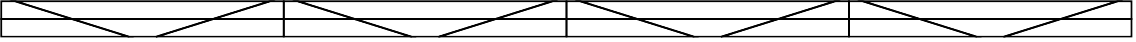}
         \caption{}
     \end{subfigure}
     \vfill
     \begin{subfigure}{\textwidth}
         \centering
         \includegraphics[width=\textwidth]{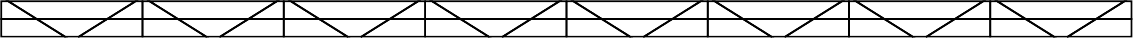}
         \caption{}
     \end{subfigure}
          \vfill
     \begin{subfigure}{\textwidth}
         \centering
         \includegraphics[width=\textwidth]{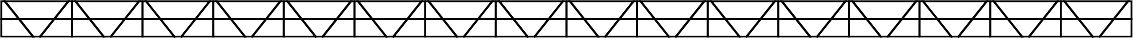}
         \caption{}
     \end{subfigure}
     \end{minipage}
        \caption{Trapezoidal meshes for the cantilever beam problem.
        Mesh is refined along the length with (a) 1 element along the
        height and (b) 2 elements along the height.}
        \label{fig:beammesh_trapezoid}
\end{figure}

\begin{figure}[!h]
    \centering
    \begin{subfigure}{0.48\textwidth}
    \centering
    \includegraphics[width=\textwidth]{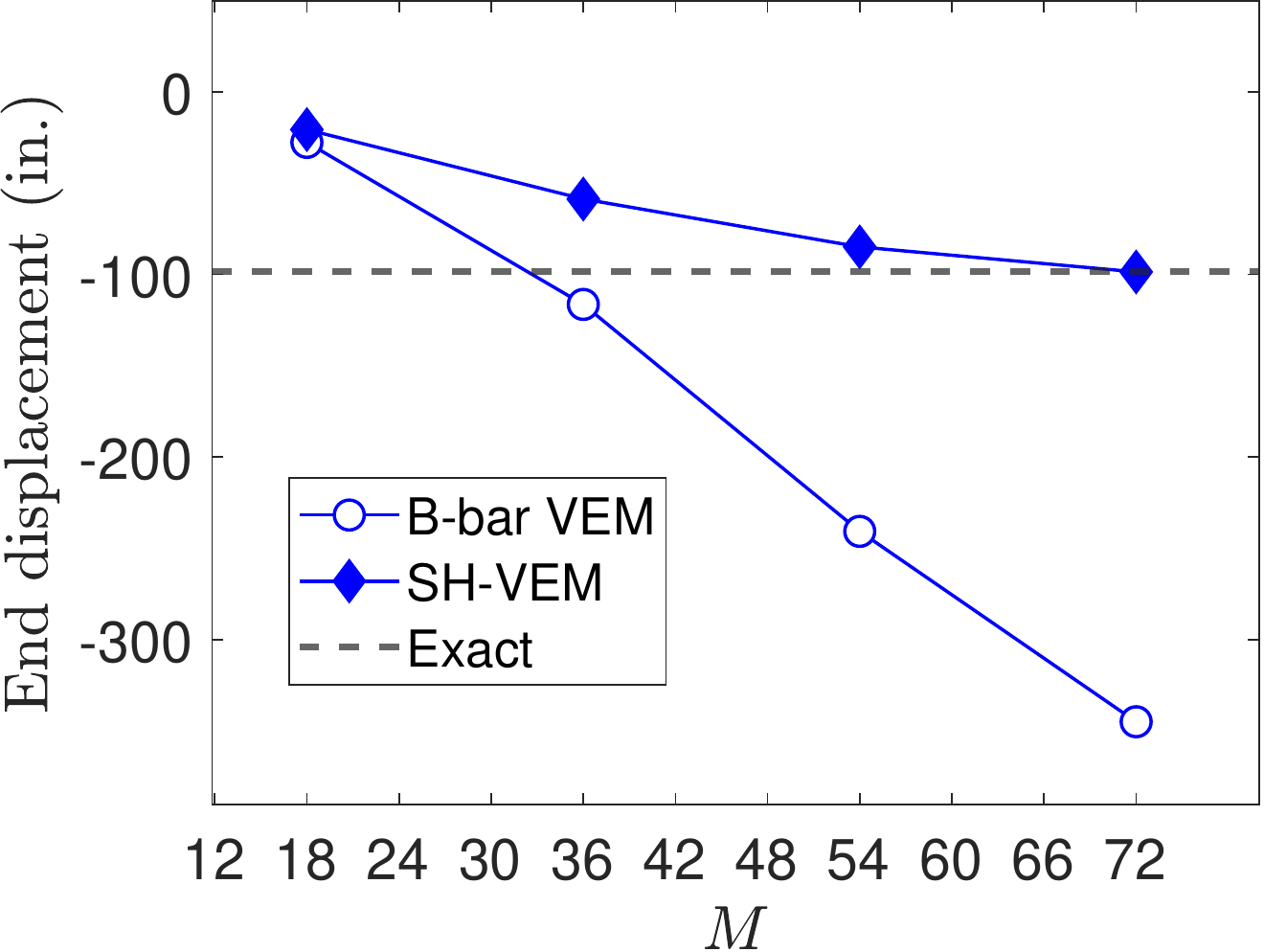}
    \caption{}\label{fig:tip_dispacement_trapezoid-a}
    \end{subfigure}
    \hfill
     \begin{subfigure}{0.48\textwidth}
    \centering
    \includegraphics[width=\textwidth]{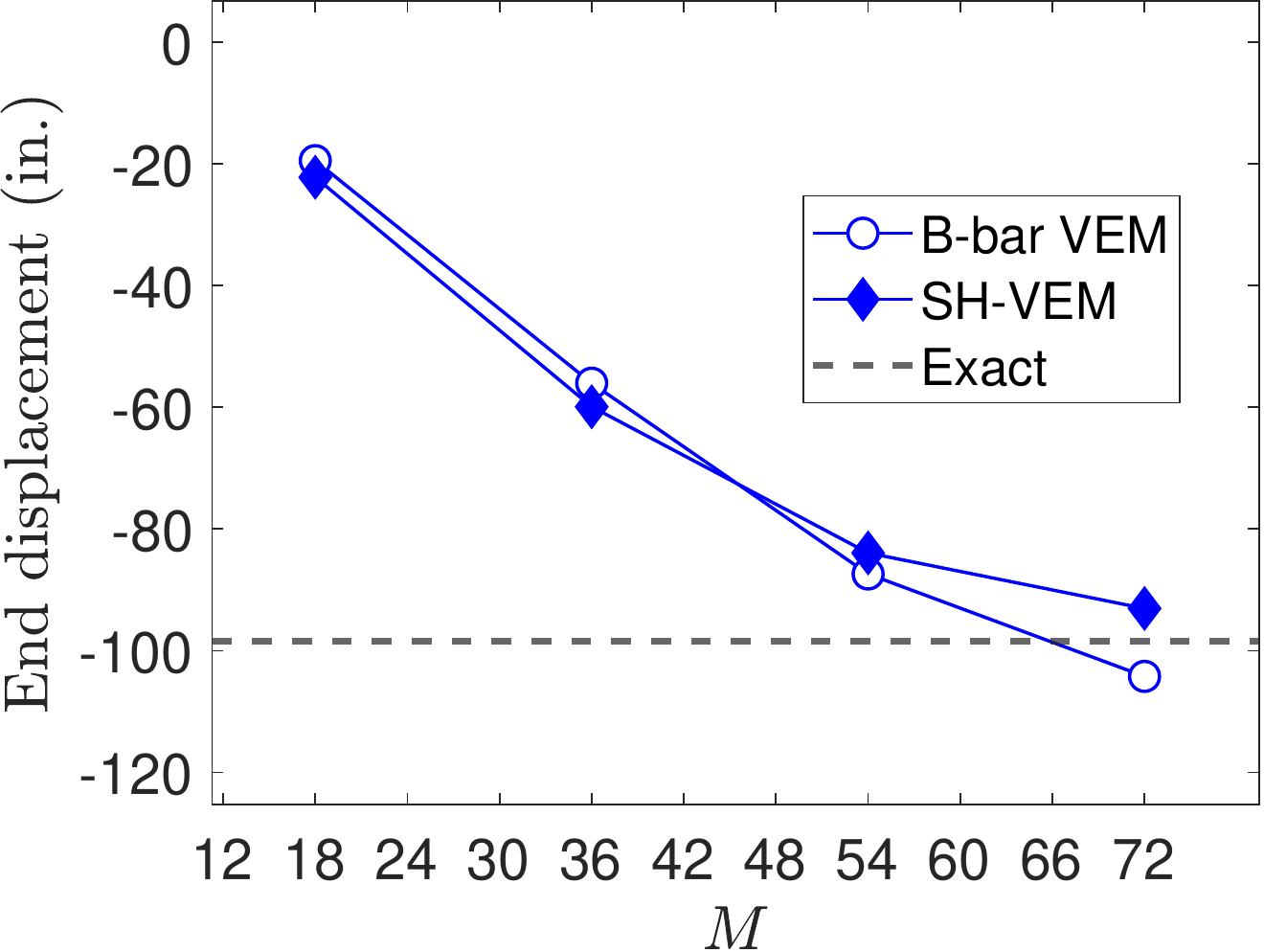}
    \caption{}\label{fig:tip_dispacement_trapezoid-b}
    \end{subfigure}
    \caption{ Convergence of the end displacement for the
    cantilever beam problem. The mesh consists of
    $M\times N$ trapezoids, where $M$ is the number of elements along the length of the beam \acrev{(see~\fref{fig:beammesh_trapezoid})}. (a) $N = 1$ and $N = 2$.}
    \label{fig:tip_dispacement_trapezoid} 
\end{figure}

It is known that distortions of a rectangular mesh can lead to shear locking in the thin beam problem.~\cite{Macneal:1985:fead} We study this 
issue on perturbed trapezoidal meshes that are shown in Figure~\ref{fig:beammesh_trapezoid}. In Figure~\ref{fig:tip_dispacement_trapezoid}, we present the convergence of the end displacement. The plot shows that on such meshes SH-VEM
is convergent but with reduced accuracy; 
however, note that the B-bar formulation fails to converge to the exact 
end displacement for the case $N = 1$ (see Figure~\ref{fig:tip_dispacement_trapezoid-a}).  

\begin{figure}
     \begin{minipage}{.48\textwidth}
     \begin{subfigure}{\textwidth}
         \centering
         \includegraphics[width=\textwidth]{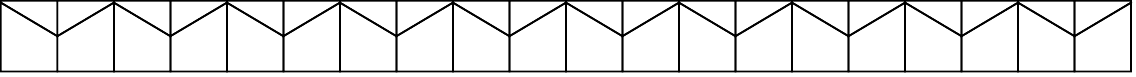}
         \caption{}
     \end{subfigure}
     \vfill
     \begin{subfigure}{\textwidth}
         \centering
         \includegraphics[width=\textwidth]{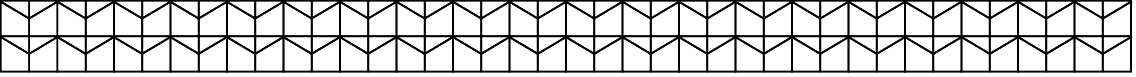}
         \caption{}
     \end{subfigure}
     \vfill
     \begin{subfigure}{\textwidth}
         \centering
         \includegraphics[width=\textwidth]{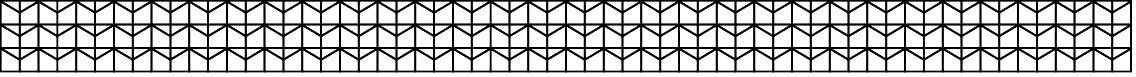}
         \caption{}
     \end{subfigure}
     \end{minipage}
     \hfill 
     \begin{minipage}{.48\textwidth}
        \begin{subfigure}{\textwidth}
         \centering
         \includegraphics[width=\textwidth]{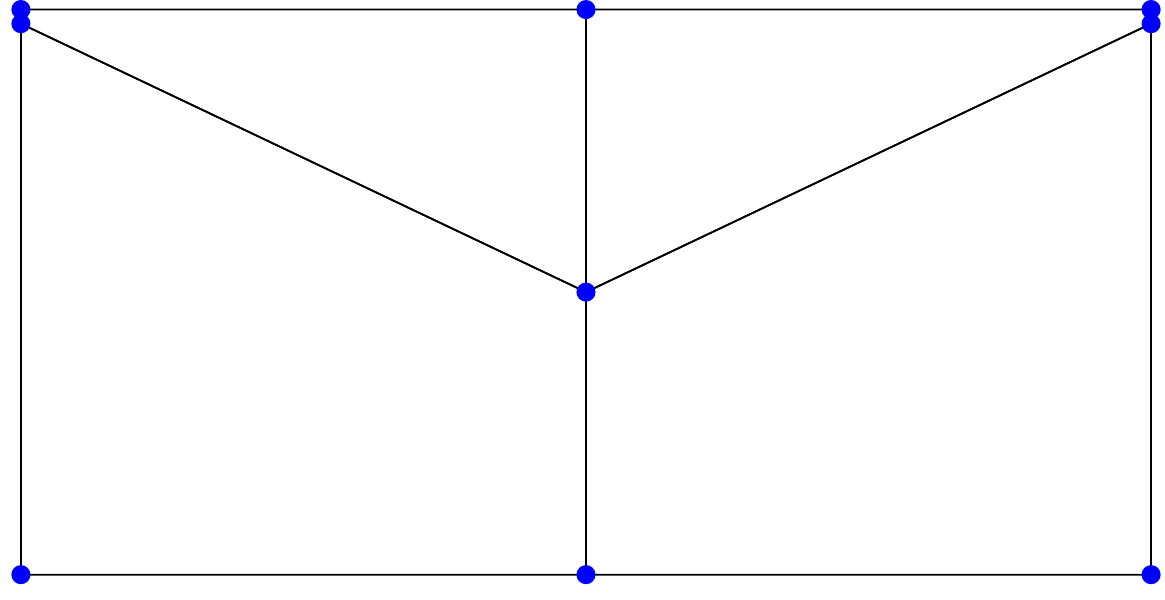}
         \caption{}
     \end{subfigure}
     \end{minipage}
        \caption{Nearly degenerate meshes used for the cantilever beam problem. (a) 40 elements, (b) 160 elements (c) 360 elements, and (d) magnification of a single element split into four quadrilaterals.} 
        \label{fig:beammesh_degen}
\end{figure}

\begin{figure}[!h]
     \centering
     \begin{subfigure}{0.32\textwidth}
         \centering
         \includegraphics[width=\textwidth]{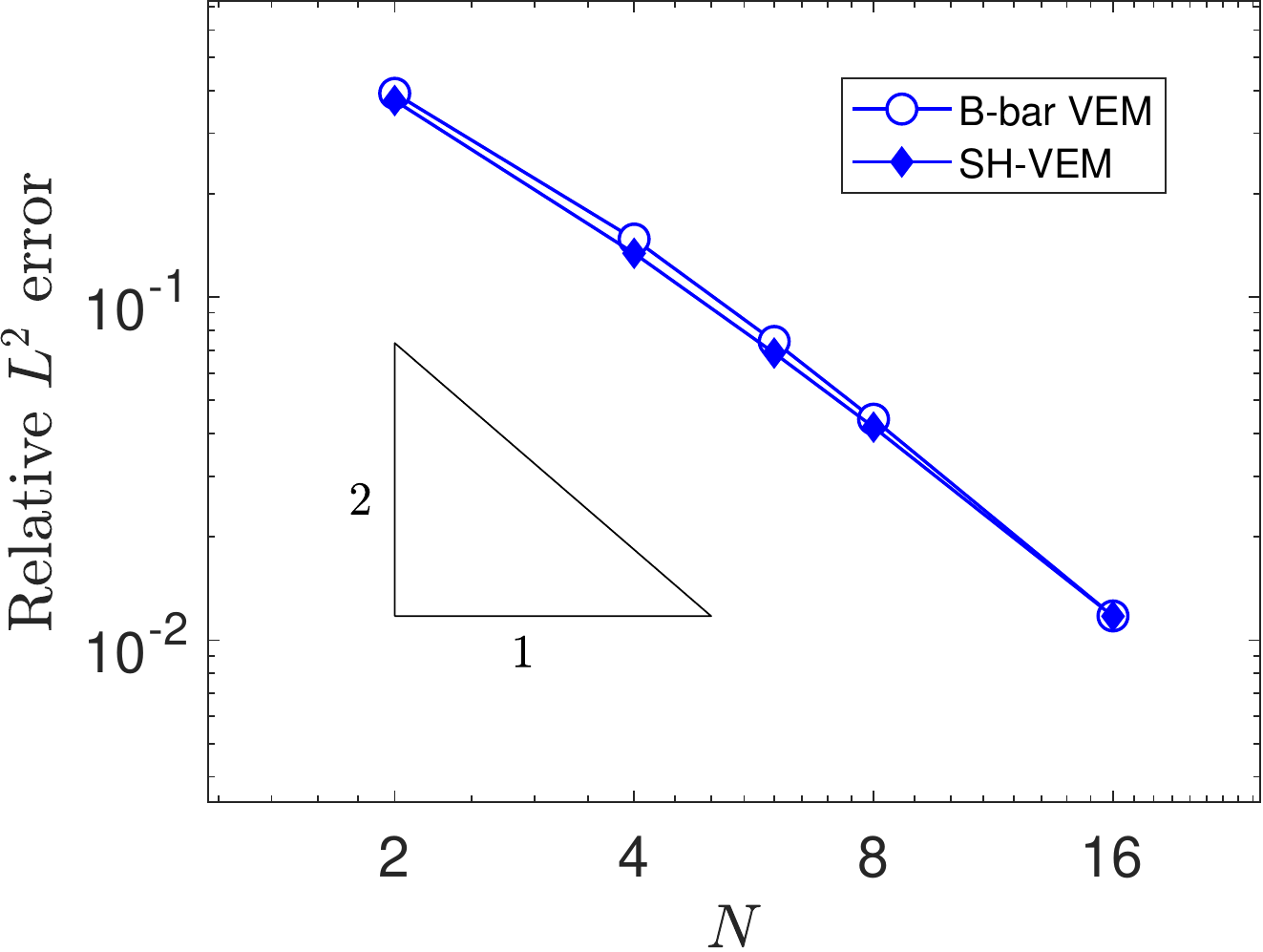}
         \caption{}
     \end{subfigure}
     \hfill
     \begin{subfigure}{0.32\textwidth}
         \centering
         \includegraphics[width=\textwidth]{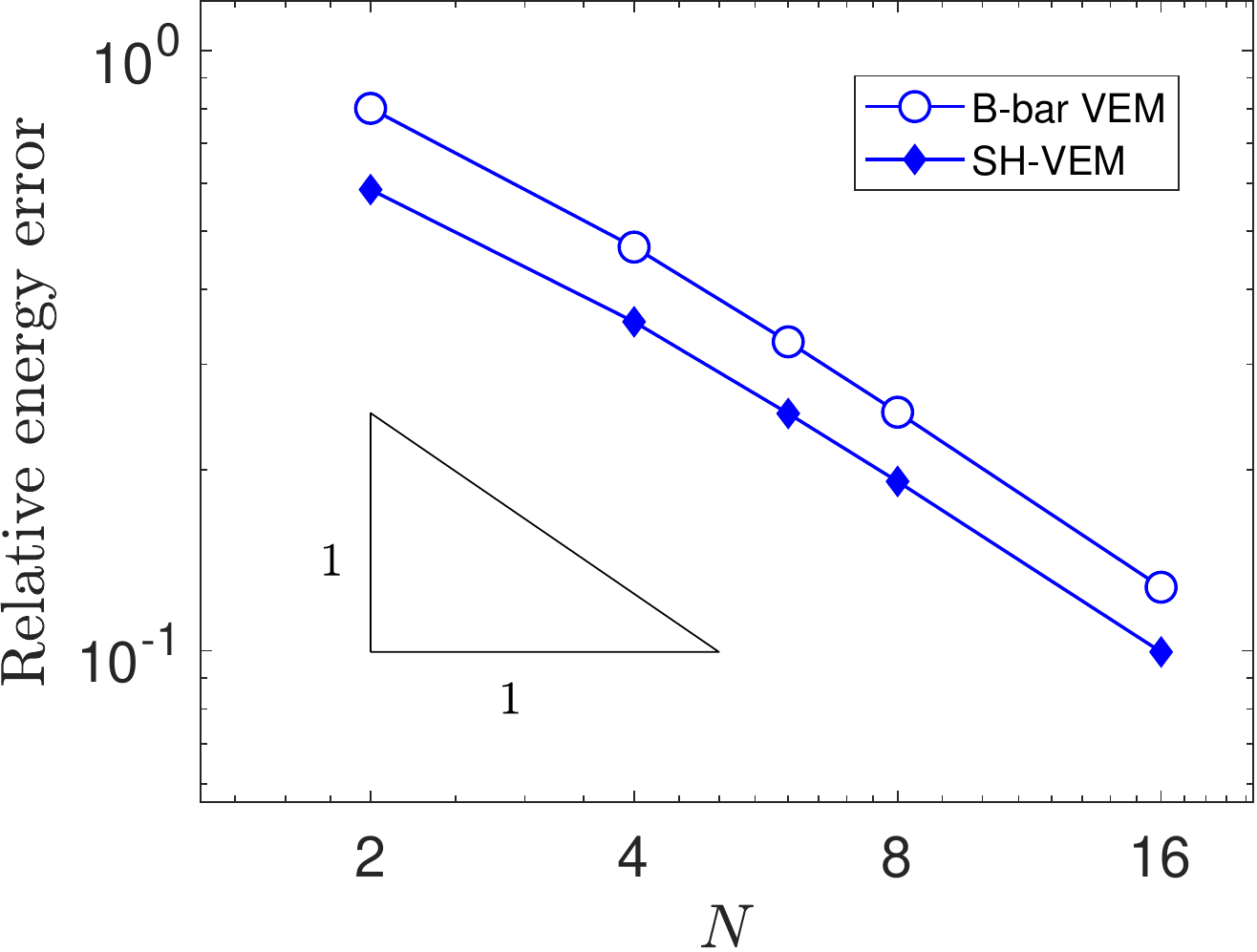}
         \caption{}
     \end{subfigure}
     \hfill
     \begin{subfigure}{0.32\textwidth}
         \centering
         \includegraphics[width=\textwidth]{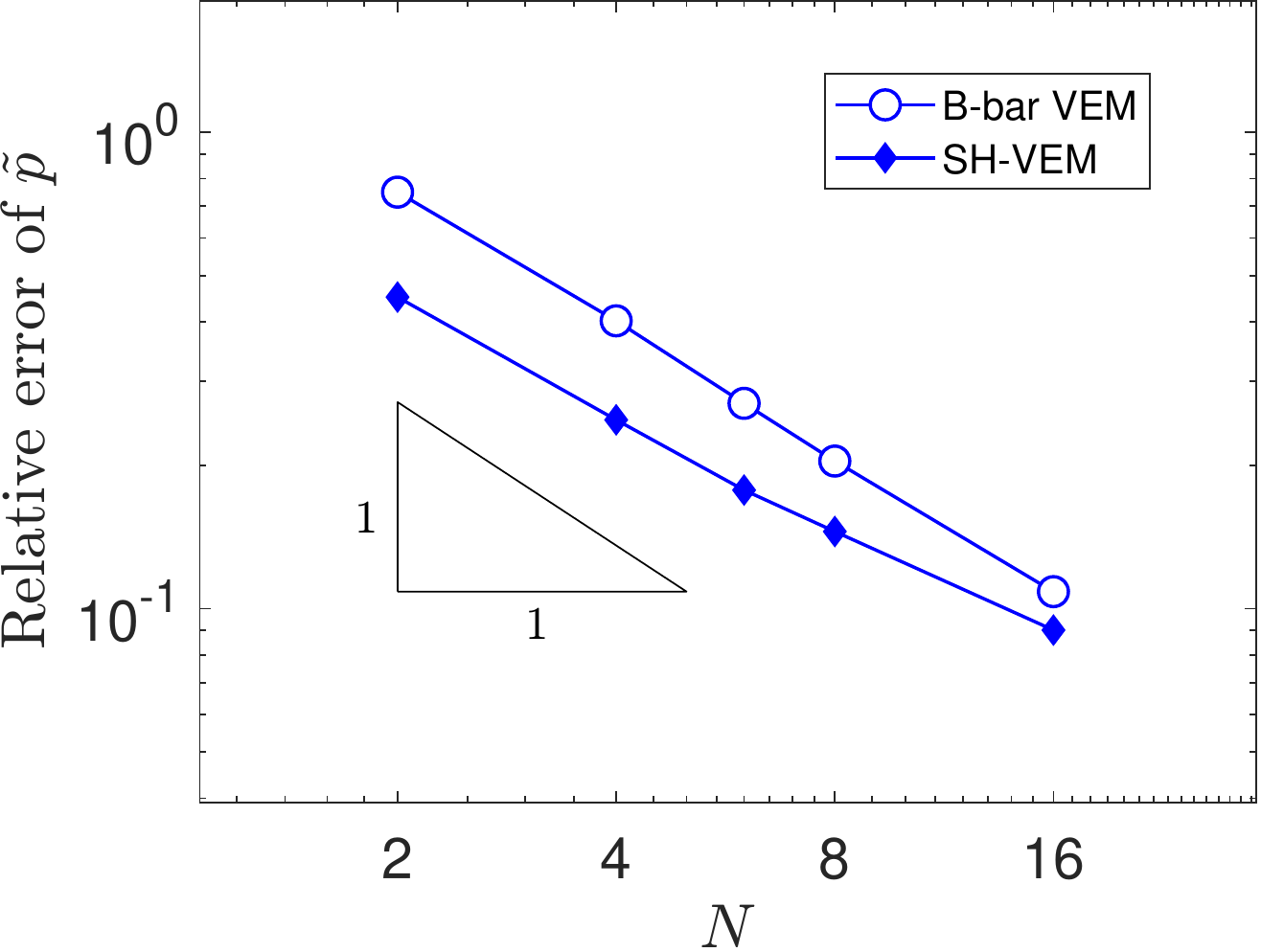}
         \caption{}
     \end{subfigure}
        \caption{Comparison of B-bar VEM and SH-VEM for the thin cantilever beam problem on nearly degenerate meshes \acrev{(see~\fref{fig:beammesh_degen})}. (a) $L^2$ error of displacement, (b) energy error and (c) $L^2$ error of hydrostatic stress, where $N$ is the number of elements along the height of the beam. }
        \label{fig:beam_convergence_degen}
\end{figure}

\begin{figure}[!h]
    \centering
    \begin{subfigure}{0.48\textwidth}
    \centering
    \includegraphics[width=\textwidth]{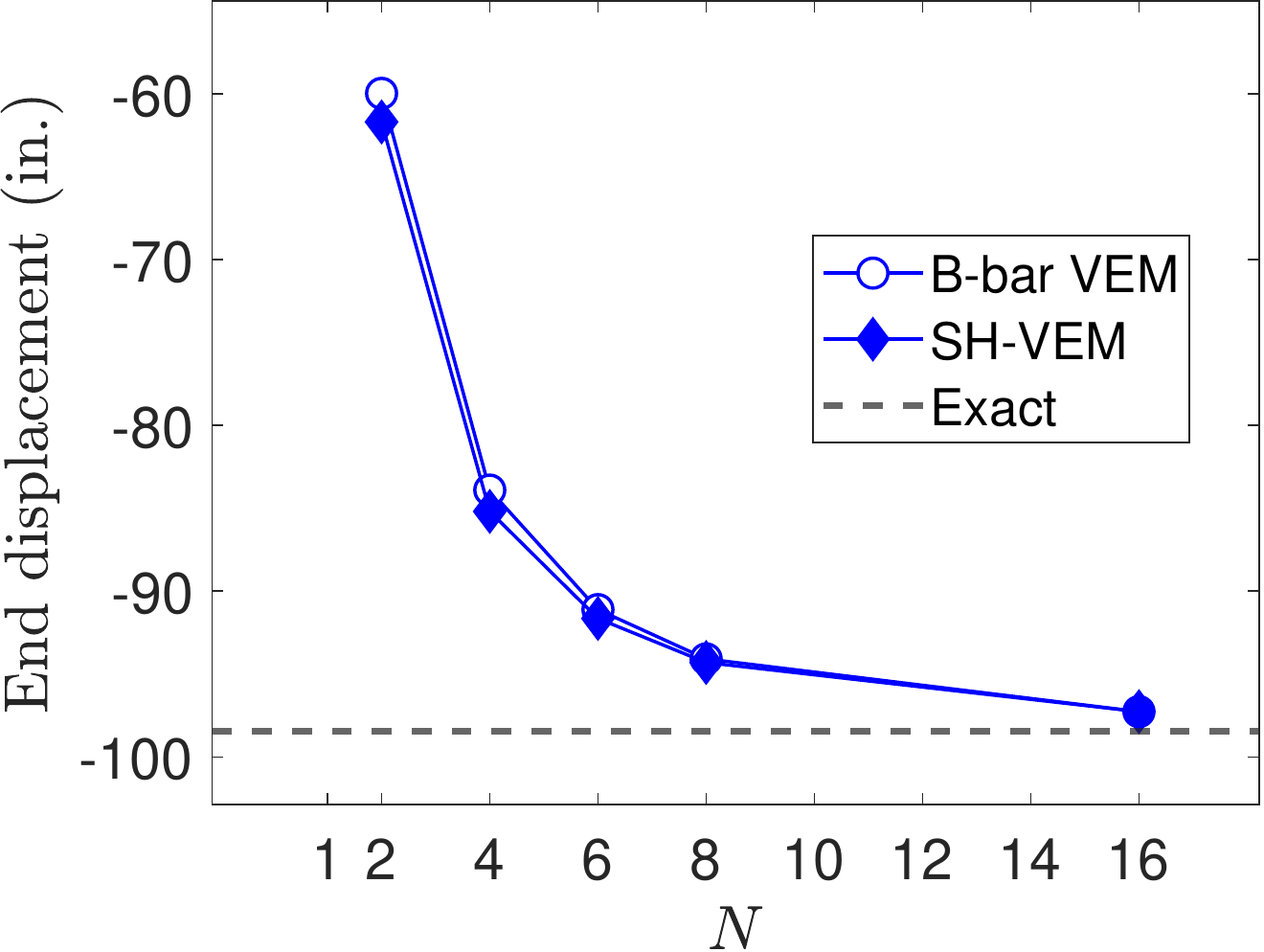}
    \caption{}
    \end{subfigure}
     \hfill
     \begin{subfigure}{0.48\textwidth}
         \centering
         \includegraphics[width=\textwidth]{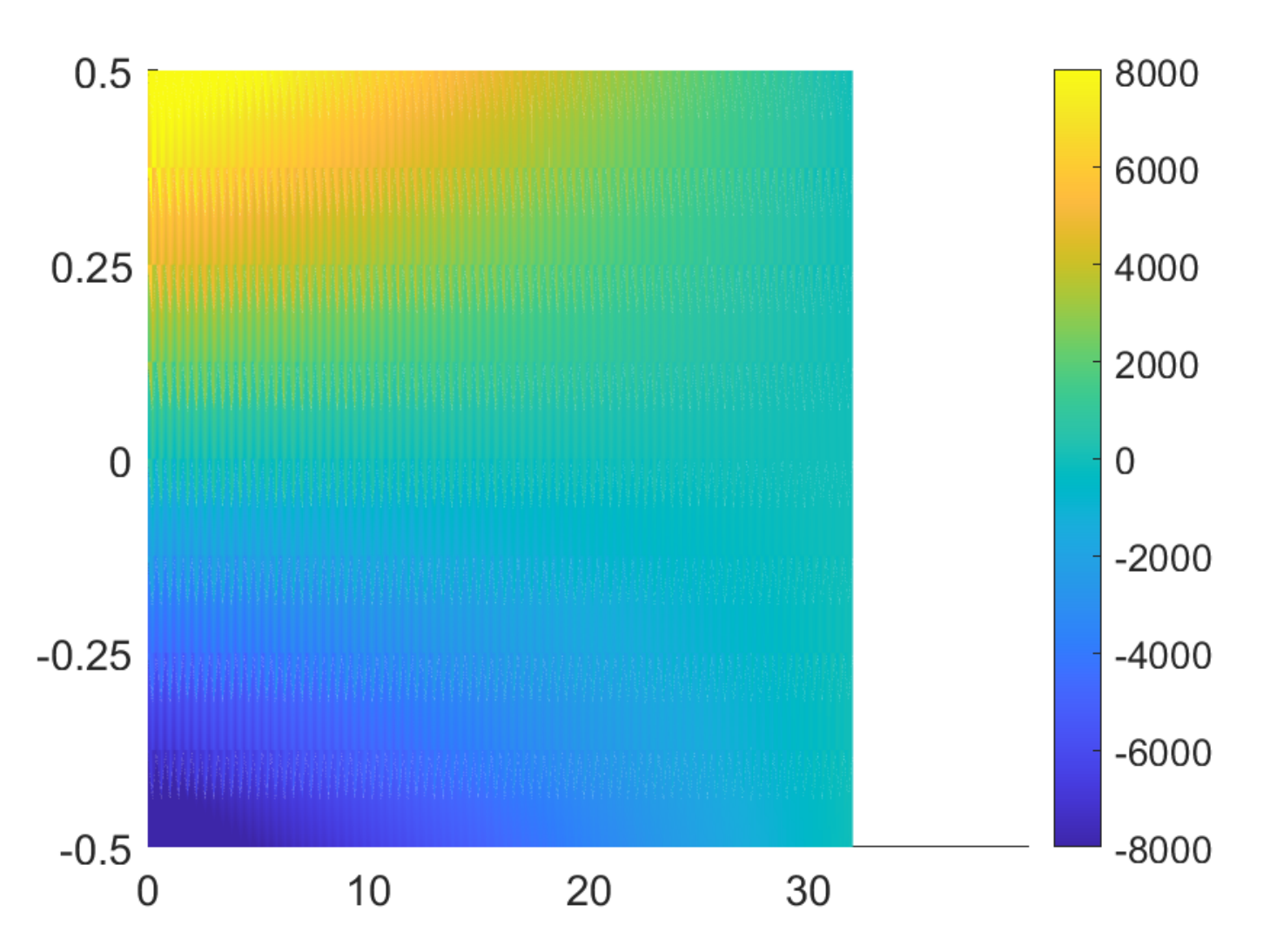}
         \caption{}
     \end{subfigure}
    \caption{(a) Convergence of the end displacement for the cantilever beam problem. The mesh consists of $10N\times N$ \acrev{nearly degenerate} quadrilaterals, where $N$ is the number of elements along the height of the beam \acrev{(see~\fref{fig:beammesh_degen})}, (b) contour plot of hydrostatic stress for SH-VEM. }
    \label{fig:tip_dispacement_degen}
\end{figure}

We also solve the cantilever beam problem on nearly degenerate quadrilateral meshes. We start with a regular rectangular mesh, and then split each element into four quadrilaterals with two of the elements have collapsing edges. A
few sample meshes are shown in Figure~\ref{fig:beammesh_degen}. In Figure~\ref{fig:beam_convergence_degen},
we compare the convergence rates of B-bar VEM and SH-VEM in the three norms, and in Figure~\ref{fig:tip_dispacement_degen} we present the convergence of the tip displacement as well as the contour plot of the hydrostatic stress 
using SH-VEM. The plots reveal
that B-bar VEM and SH-VEM retain optimal convergence rates. Furthermore, the convergence of SH-VEM is monotonic; however its accuracy is worse when compared to the uniform mesh case. This decrease in accuracy
can be attributed to the poor shape (near-degeneracy) quality of the elements.
\subsection{Cook's membrane}
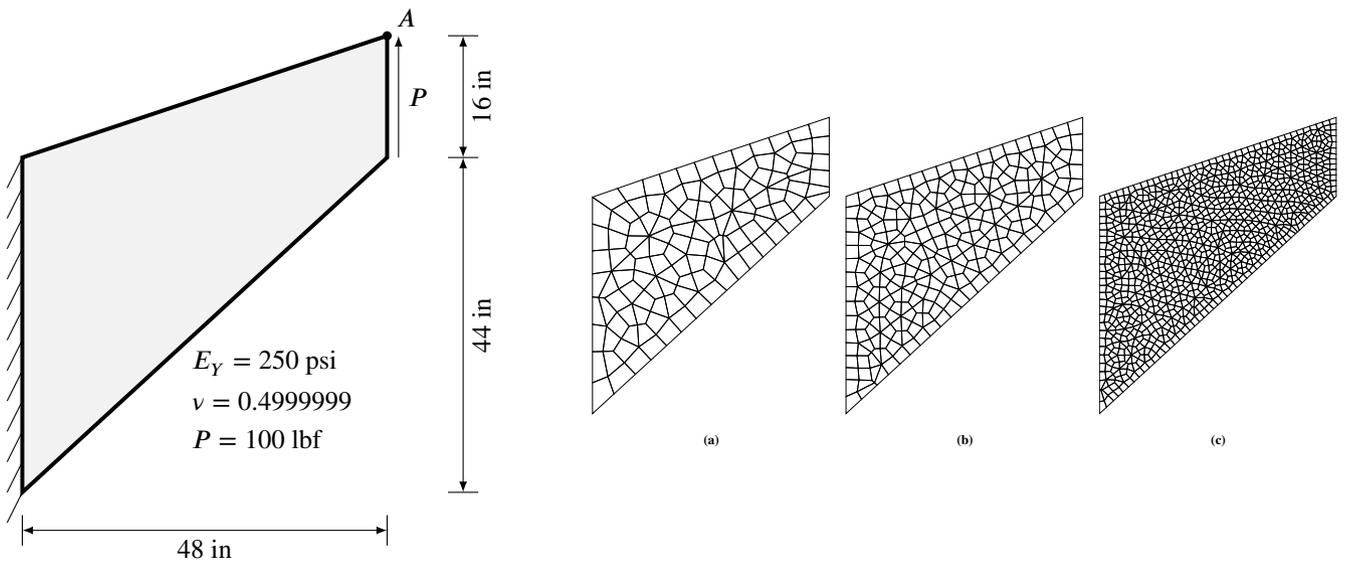
\begin{figure}[!h]
    \begin{minipage}{.4\textwidth}
        \begin{subfigure}{\textwidth}
            \centering
        \begin{tikzpicture}[scale=0.1]
			%draws the cooks membrane shape
			\filldraw[line width=1.5pt,fill=gray!10] (0,0) -- ++(0,44) -- ++(48,16) -- ++(0, -16) -- cycle;

            %draws the arrows on right side
			%\draw[-latex] (48+1.5,44) -- ++(0,16/5);
			%\draw[-latex] (48+1.5,44+16/5) -- ++(0,16/5);
			%\draw[-latex] (48+1.5,44+2*16/5) -- ++(0,16/5);
			%\draw[-latex] (48+1.5,44+3*16/5) -- ++(0,16/5);
			%\draw[-latex] (48+1.5,44+4*16/5) -- ++(0,16/5);
            \draw[-latex] (48+1.5,44) --++(0,16);
			\draw(48+1.5,44+16/2) node[right]{$P$};
			\filldraw[black] (48,44+16) circle (15pt) node[anchor=south west]{$A$};

            %draws the wall on the left hand side
			\draw[] (0,0)--++(-2,-4) ++(2,4+4) -- ++(-2,-4)
			++(2,4+4) -- ++(-2,-4)  ++(2,4+4) -- ++(-2,-4)  ++(2,4+4) -- ++(-2,-4)
			++(2,4+4) -- ++(-2,-4)  ++(2,4+4) -- ++(-2,-4)  ++(2,4+4) -- ++(-2,-4)
			++(2,4+4) -- ++(-2,-4)  ++(2,4+4) -- ++(-2,-4)  ++(2,4+4) -- ++(-2,-4)
			++(2,4+4) -- ++(-2,-4)  ++(2,4+4)
			; 

            %draws the x-axis with label on bottom
			\draw[latex-latex] (0,0-5) -- (48,-5);
            \draw[] (0,-7) --++ (0,4);
            \draw[] (48,-7) --++(0,4);
   
            %(0,0-5)
			%-- ++(-1,-1) -- ++(2,2) ++(-1,-1)
			%-- ++(48,0)
			%-- ++(-1,-1) -- ++(2,2) ++(-1,-1)
			\draw[](48/2,-5) node[anchor = north]{48 in};

            %draws the y axis with label
			\draw[latex-latex] (58,0)--++ (0,44); 
            \draw[latex-latex] (58,44) --++ (0,16);
            \draw[] (56,60) --++(4,0);
            \draw[](56,44) --++(4,0);
            \draw[](56,0) --++(4,0);
            %(44+12,0)
			%-- ++(-1,-1) -- ++(2,2) ++(-1,-1)
			%-- ++(0,44)
			%-- ++(-1,-1) -- ++(2,2) ++(-1,-1)
			%-- ++(0,16)
			%-- ++(-1,-1) -- ++(2,2) ++(-1,-1)
			\draw[] (58,104/2) node[anchor = north, rotate=90]{16 in};
			\draw[] (58,44/2) node[anchor = north, rotate=90]{44 in};
			%;    

			%\draw[](10,40) node[]{$\times$~1 mm};

            %draws the material properties
			\draw[] (21,17) node[anchor=west]{$E_Y = 250$ psi}
			(21,17-5) node[anchor=west]{$\nu = 0.4999999$}
			(21,17-10) node[anchor=west]{$P = 100$ lbf}
			;
        \end{tikzpicture}
        \end{subfigure}
    \end{minipage}
    \hfill
    \begin{minipage}{.55\textwidth}
        \begin{subfigure}{.32\textwidth}
         \centering
         \includegraphics[width=\textwidth]{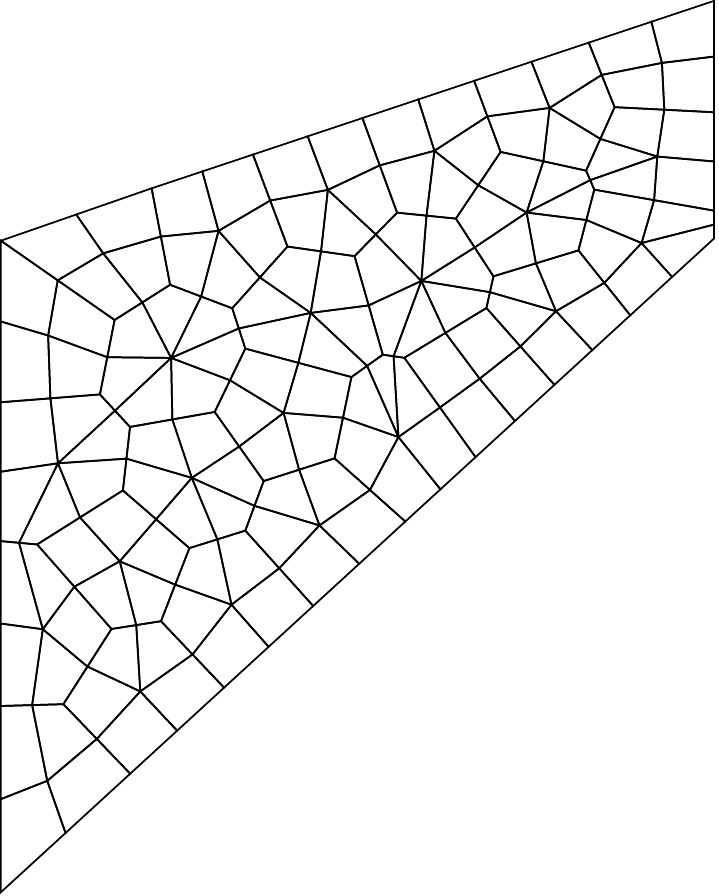}
         \caption{}
     \end{subfigure}
     \hfill
     \begin{subfigure}{.32\textwidth}
         \centering
         \includegraphics[width=\textwidth]{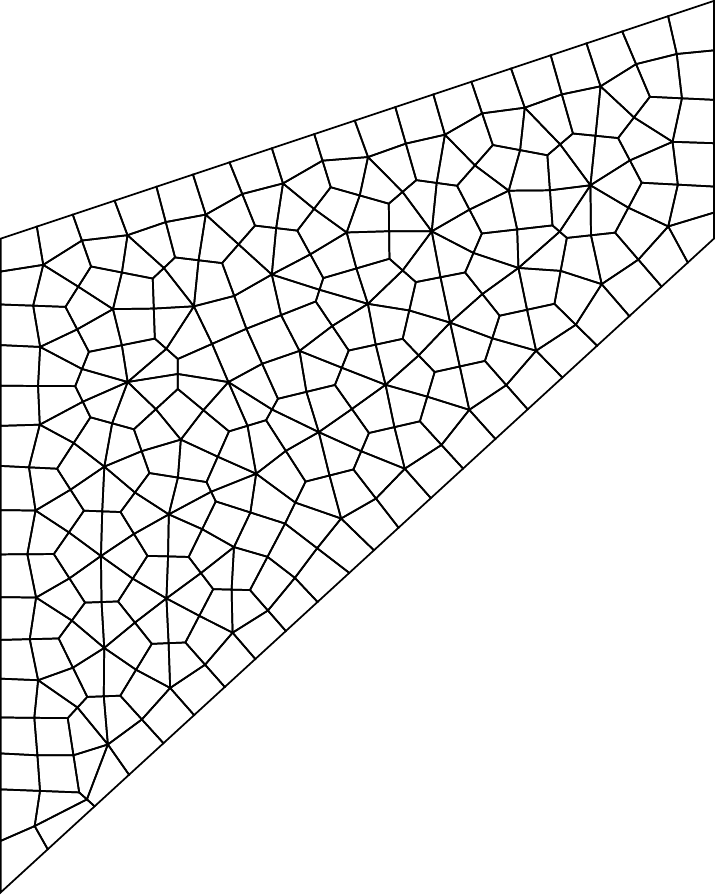}
         \caption{}
     \end{subfigure}
     \hfill
     \begin{subfigure}{.32\textwidth}
         \centering
         \includegraphics[width=\textwidth]{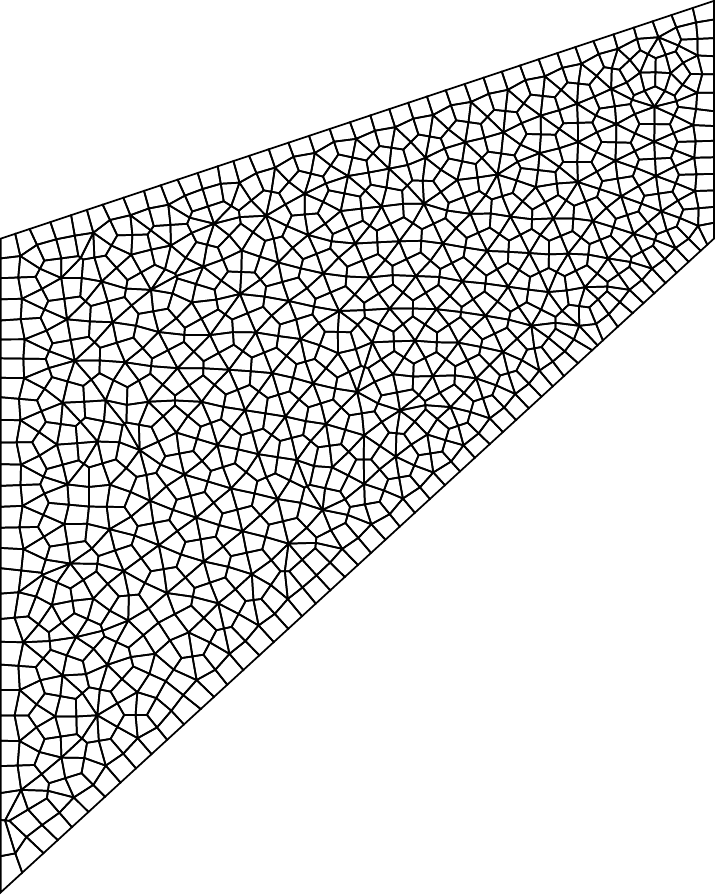}
         \caption{}
     \end{subfigure}
    \end{minipage}
        \caption{(a) Cook's membrane problem. (b), (c), (d)
        Unstructured quadrilateral meshes with 100 elements, 300 elements and 1000 elements,
        respectively.}
        \label{fig:cooksmesh_unstructured}
\end{figure}

\begin{figure}[!h]
    \centering
    \begin{subfigure}{0.48\textwidth}
    \centering
    \includegraphics[width=\textwidth]{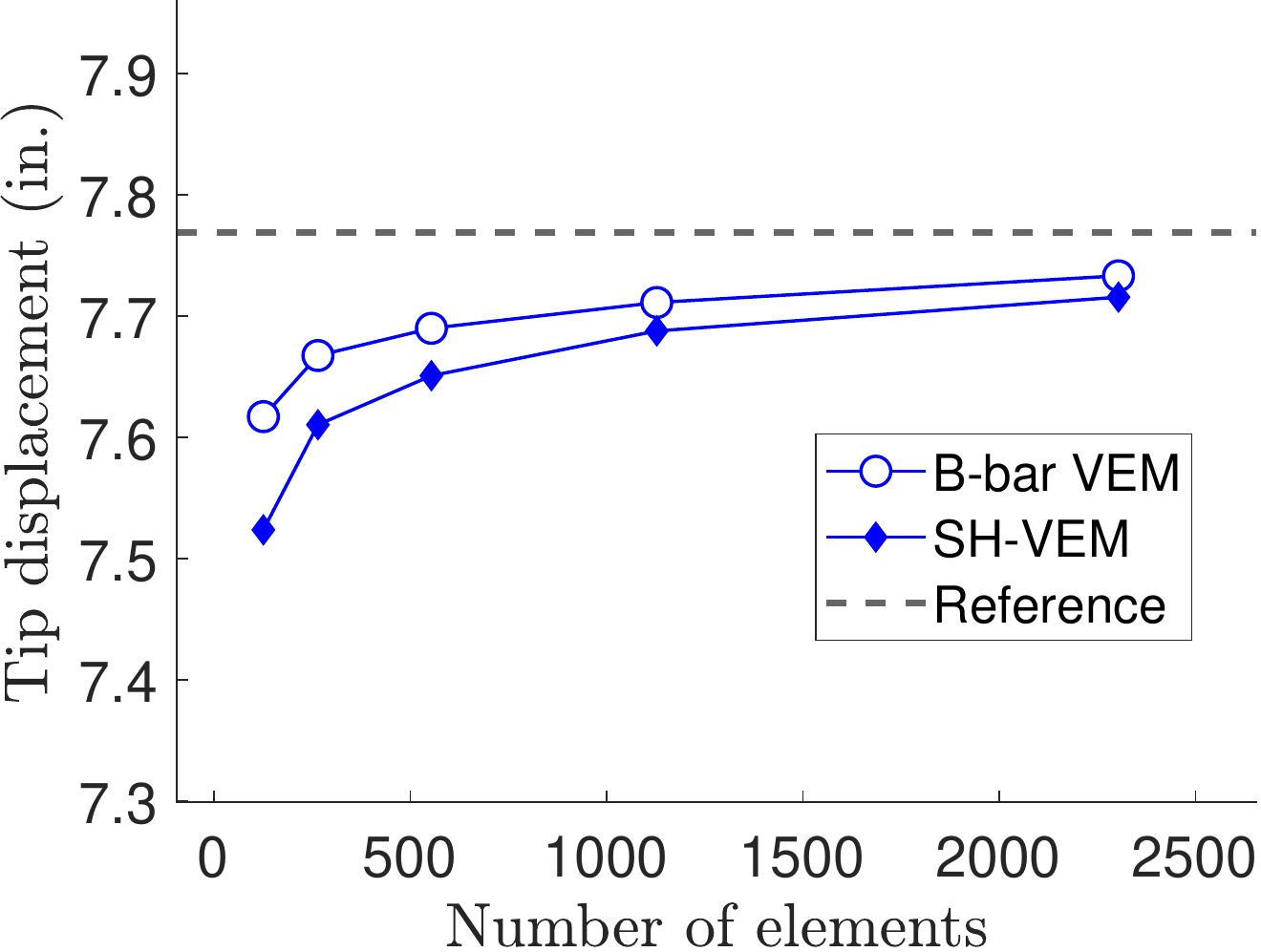}
    \caption{}
    \end{subfigure}
     \hfill
     \begin{subfigure}{0.48\textwidth}
         \centering
         \includegraphics[width=\textwidth]{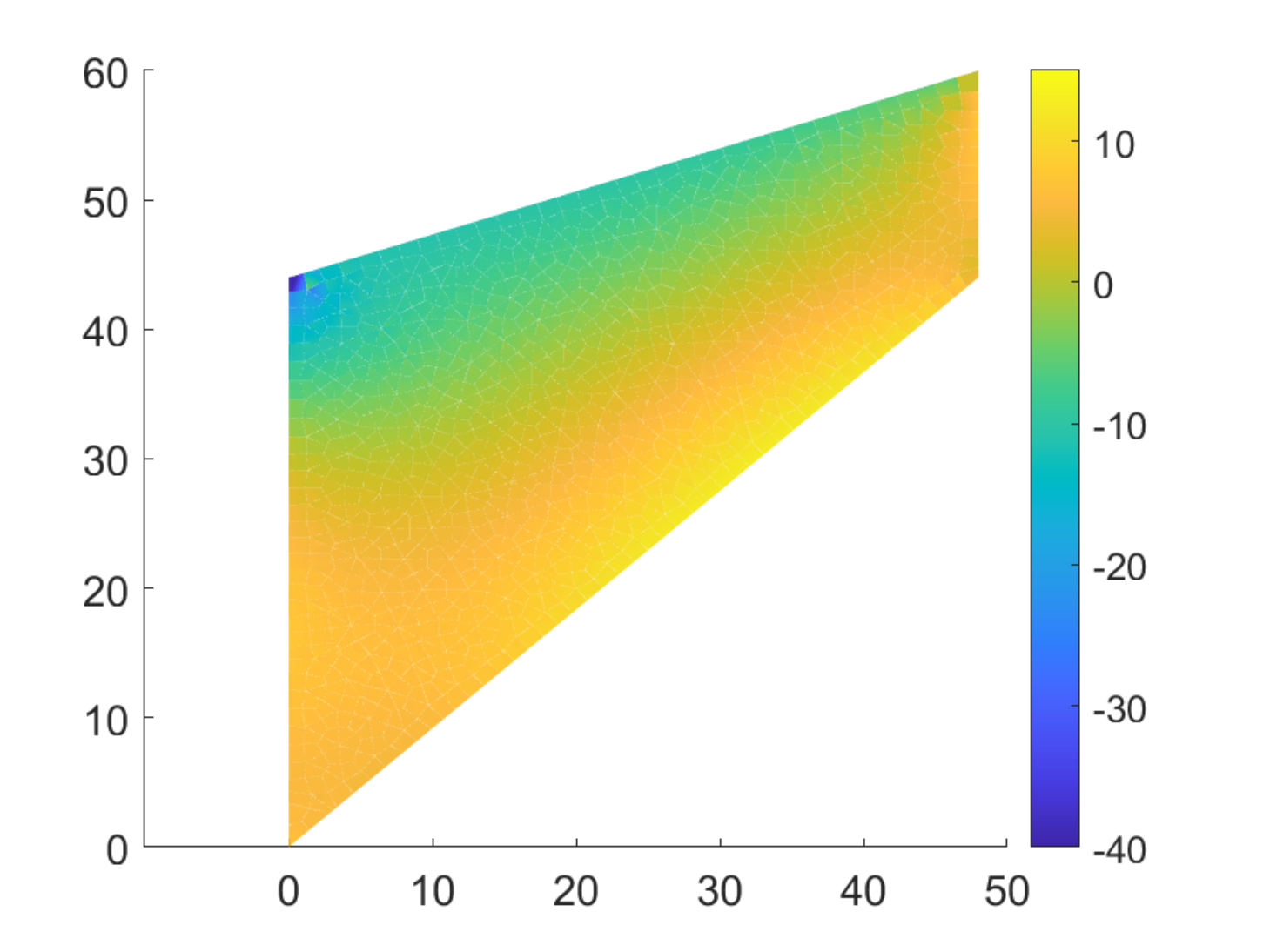}
         \caption{}
     \end{subfigure}
    \caption{(a) Convergence of the tip displacement for Cook's membrane problem. The mesh consists of unstructured quadrilaterals \acrev{(see~\fref{fig:cooksmesh_unstructured})}. (b) Contour plot of the hydrostatic stress for SH-VEM.}
    \label{fig:cook_tip_dispacement_unstructured}
\end{figure}

Here we consider the Cook's membrane problem under shear load~\cite{Cook:1974:jsd} (see Figure~\ref{fig:cooksmesh_unstructured}). This problem is commonly used to test a combination of bending and shear for 
nearly-incompressible materials.
The material has Young's modulus $E_Y=250$ psi and Poisson's ratio $\nu=0.4999999$. The left edge of the membrane is fixed and the right edge has an applied shear load of $F = 6.25$ lbf per unit length. \acrev{This problem does not have an exact solution; a reference solution for the vertical displacement at the tip of the membrane is $v=7.769$ inch.~\cite{Park:2020:meccanica}} We first test this problem on an unstructured quadrilateral mesh. A few sample meshes are shown in Figure~\ref{fig:cooksmesh_unstructured}. In Figure~\ref{fig:cook_tip_dispacement_unstructured},
the convergence of the tip displacement and that of the hydrostatic stress 
are presented. The plot shows that the B-bar VEM and SH-VEM have comparable accuracy and convergence for the tip displacement. In
addition, SH-VEM is able to produce a relatively smooth hydrostatic stress field on an unstructured mesh. 

\begin{figure}[!h]
     \centering
     \begin{subfigure}{.32\textwidth}
         \centering
         \includegraphics[width=\textwidth]{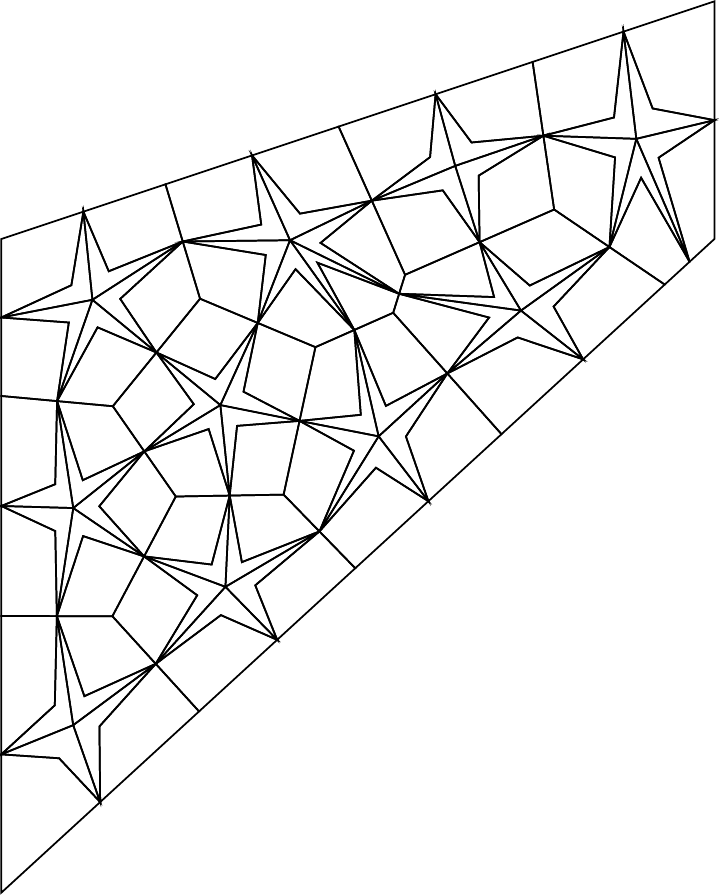}
         \caption{}
     \end{subfigure}
     \hfill
     \begin{subfigure}{.32\textwidth}
         \centering
         \includegraphics[width=\textwidth]{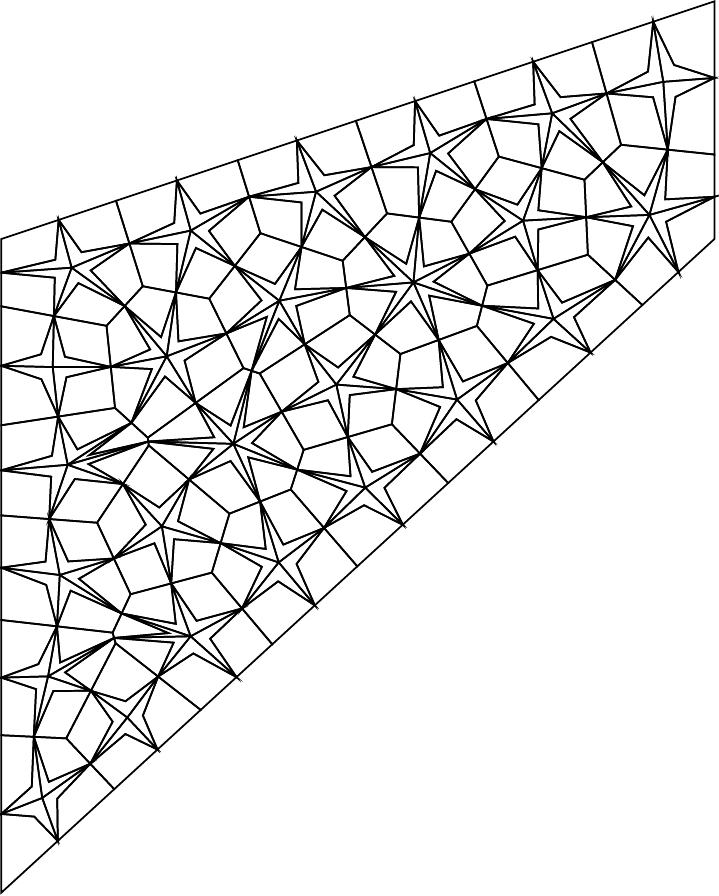}
         \caption{}
     \end{subfigure}
     \hfill
     \begin{subfigure}{.32\textwidth}
         \centering
         \includegraphics[width=\textwidth]{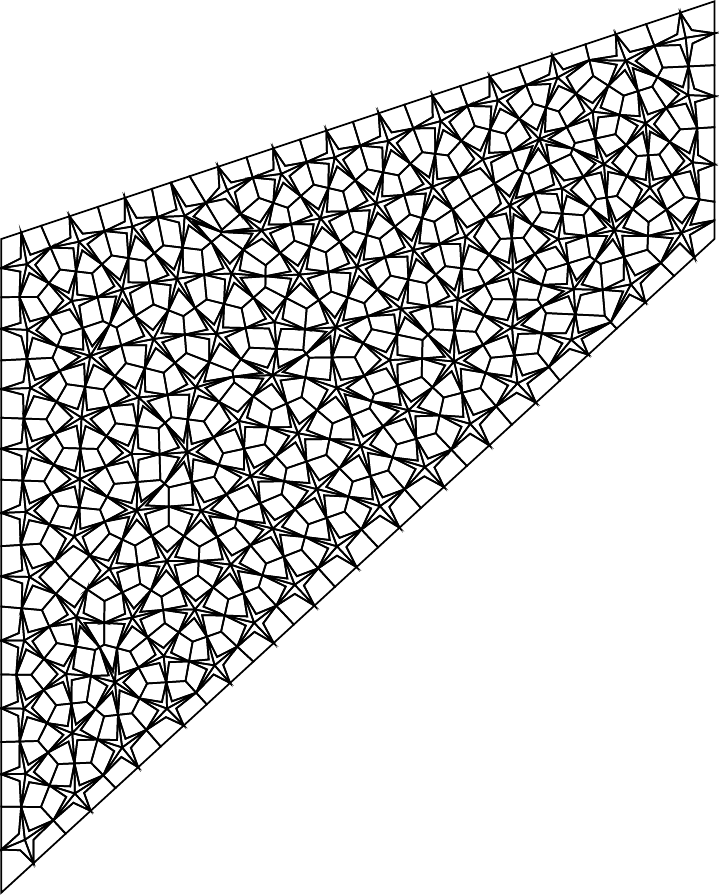}
         \caption{}
     \end{subfigure}
        \caption{Nonconvex quadrilateral meshes for the Cook's membrane problem. (a) 100 elements, (b) 250 elements and (c) 1000 elements.  }
        \label{fig:cooksmesh_nonconvex}
\end{figure}

\begin{figure}[!h]
    \centering
    \begin{subfigure}{0.48\textwidth}
    \centering
    \includegraphics[width=\textwidth]{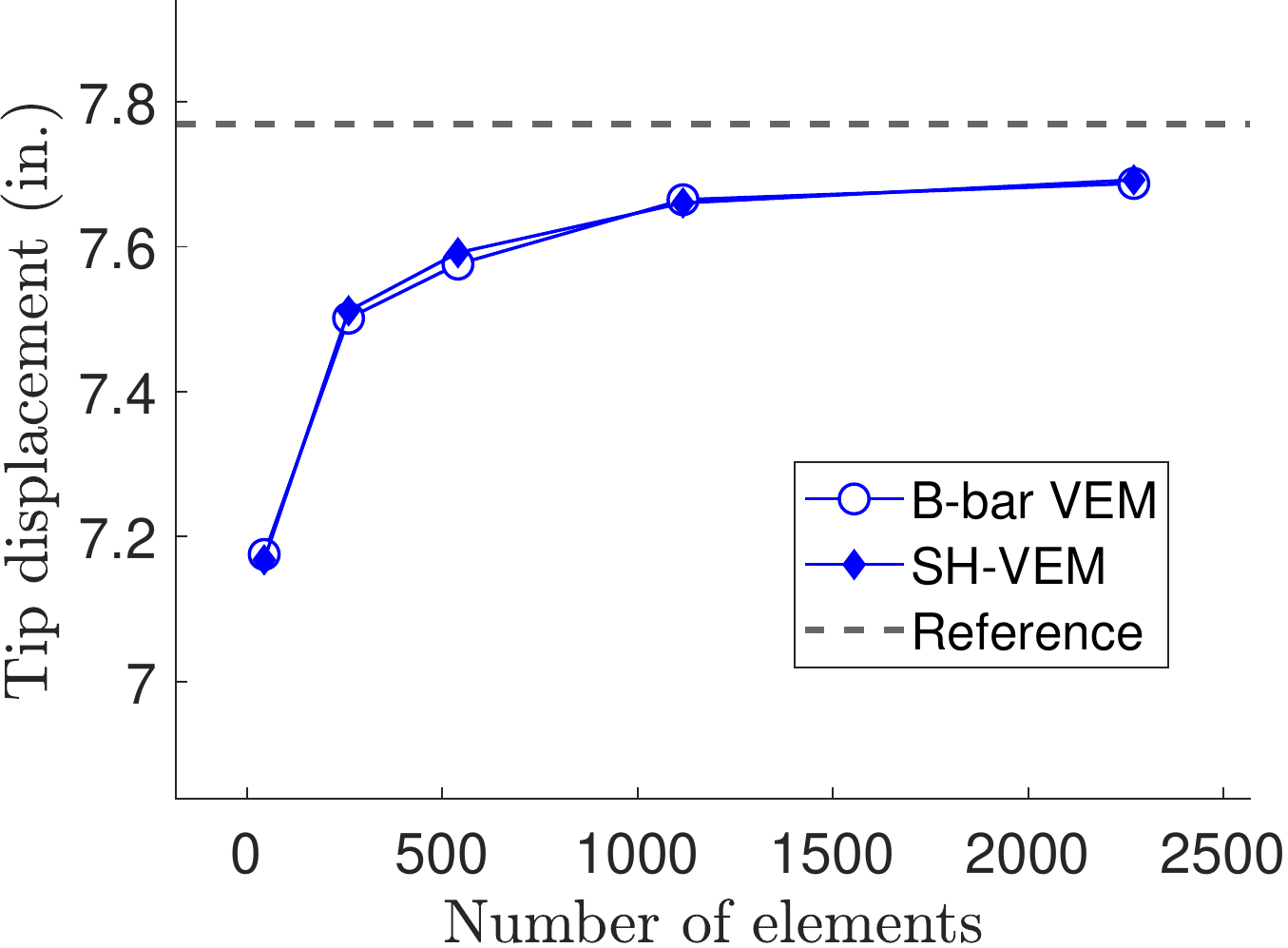}
    \caption{}
    \end{subfigure}
     \hfill
     \begin{subfigure}{0.48\textwidth}
         \centering
         \includegraphics[width=\textwidth]{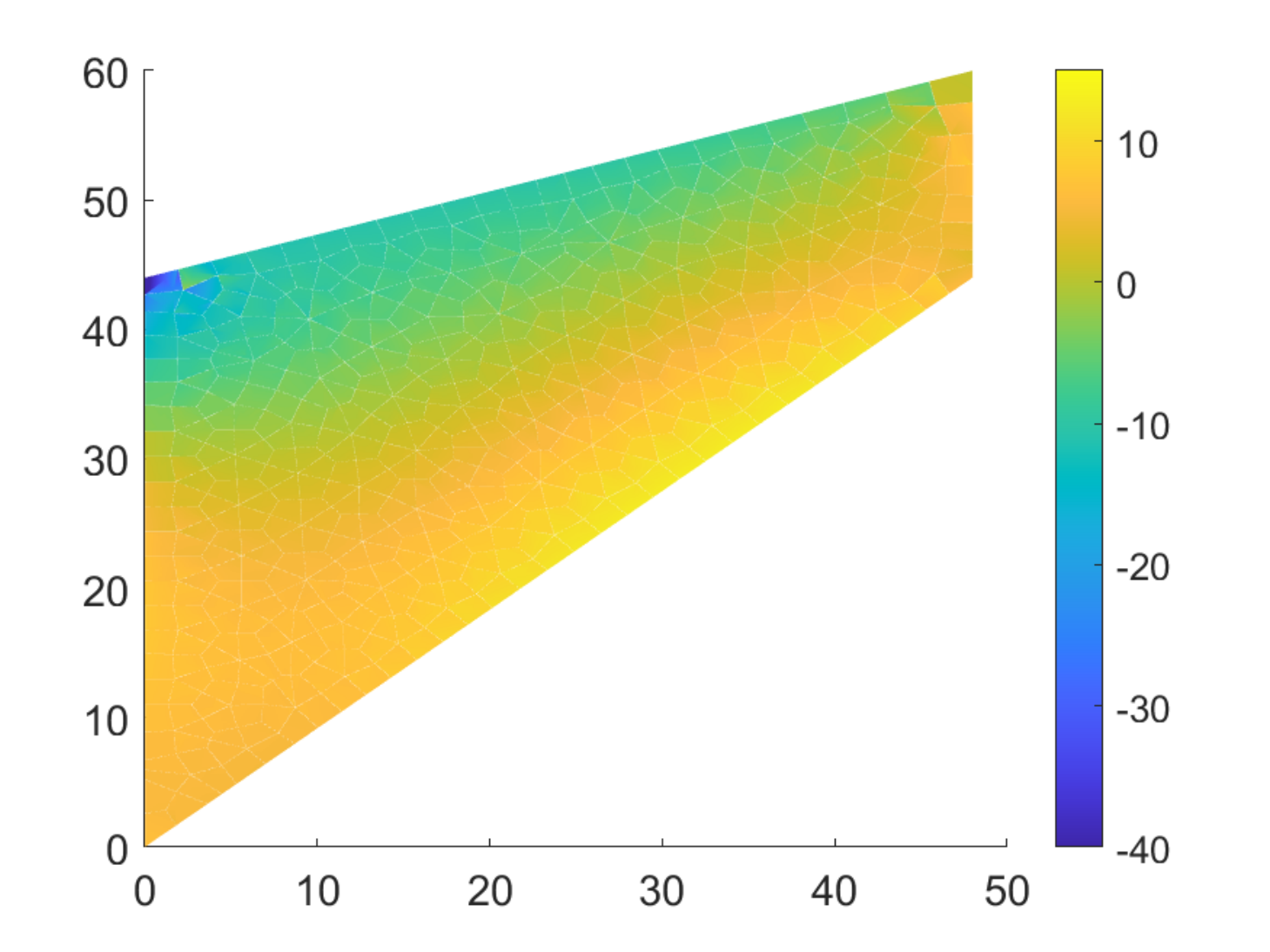}
         \caption{}
     \end{subfigure}
    \caption{(a) Convergence of the tip displacement for Cook's membrane. The mesh consists of nonconvex quadrilaterals \acrev{(see~\fref{fig:cooksmesh_nonconvex})}.
    (b) Contour plot of the hydrostatic stress
    for SH-VEM. }
    \label{fig:cook_tip_dispacement_nonconvex}
\end{figure}

Next, the SH-VEM is now assessed for the Cook's membrane problem on nonconvex meshes. We begin with an unstructured quadrilateral mesh, and then each element is split into a convex and a 
nonconvex quadrilateral. A few representative meshes are shown in Figure~\ref{fig:cooksmesh_nonconvex}. The plots of the convergence of tip displacement and the contour of the hydrostatic stress are
presented in Figure~\ref{fig:cook_tip_dispacement_nonconvex}. The plots show that even on nonconvex meshes, the convergence of the tip displacement of 
B-bar VEM and SH-VEM are proximal, and the 
contours of the 
hydrostatic stress for SH-VEM remains relatively smooth.
\subsection{Plate with a circular hole}
\begin{figure}[!h]
     \centering
     \begin{subfigure}{.32\textwidth}
         \centering
         \includegraphics[width=\textwidth]{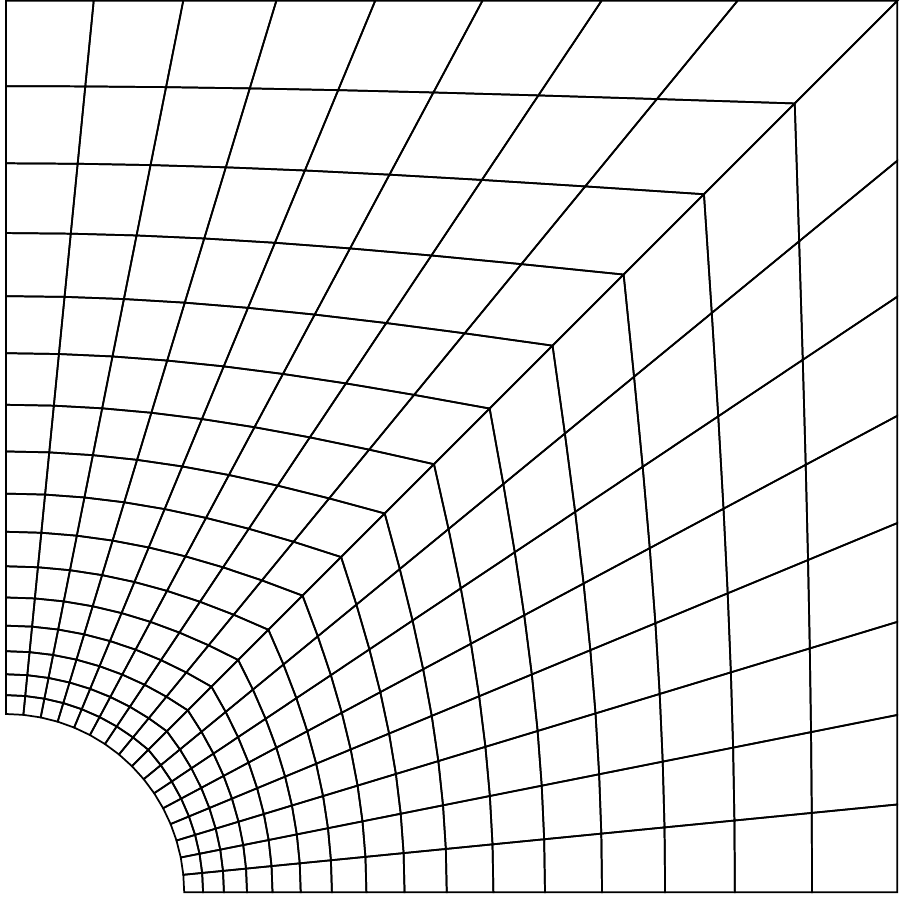}
         \caption{}
     \end{subfigure}
     \hfill
     \begin{subfigure}{.32\textwidth}
         \centering
         \includegraphics[width=\textwidth]{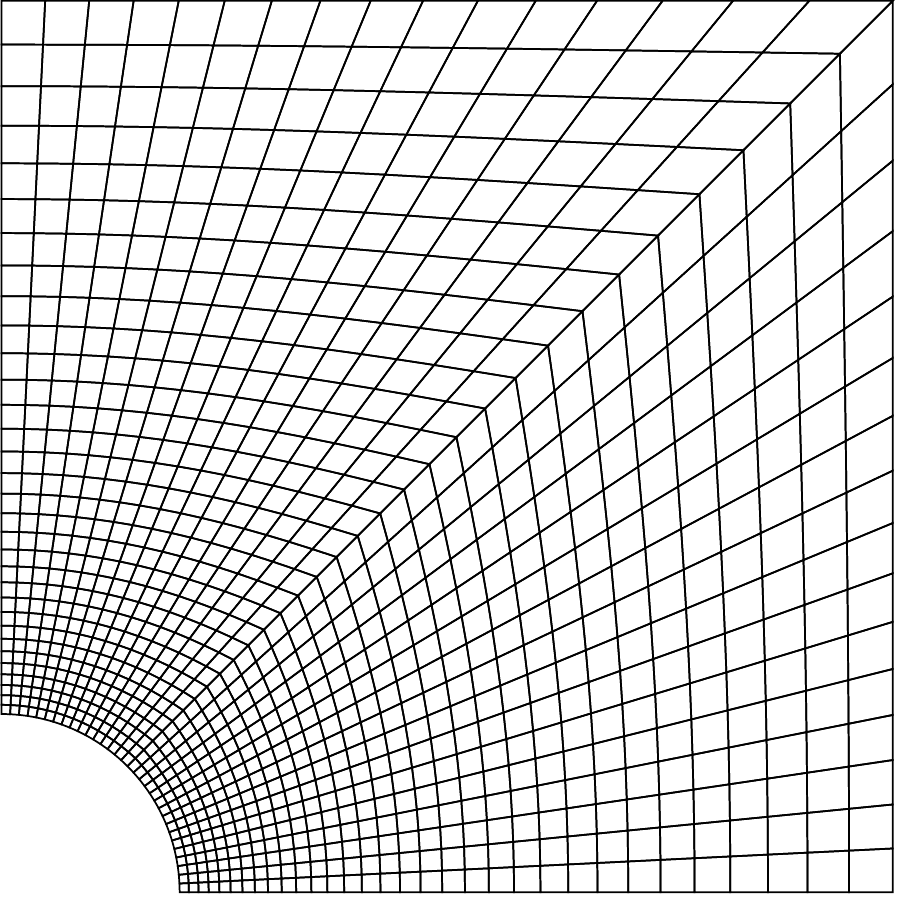}
         \caption{}
     \end{subfigure}
     \hfill
     \begin{subfigure}{.32\textwidth}
         \centering
         \includegraphics[width=\textwidth]{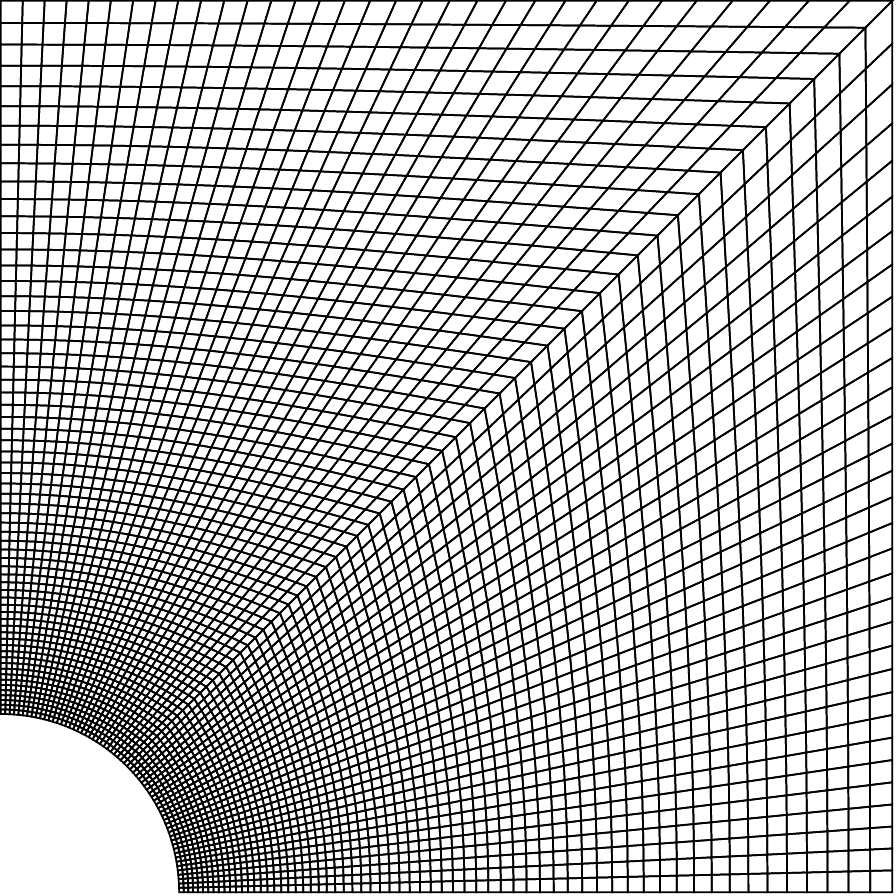}
         \caption{}
     \end{subfigure}
        \caption{Structured quadrilateral meshes for the plate with a hole problem. (a) 256 elements, (b) 1024 elements and (c) 4096 elements.  }
        \label{fig:plate_hole_structured}
\end{figure}

\begin{figure}[!h]
     \centering
     \begin{subfigure}{0.32\textwidth}
         \centering
         \includegraphics[width=\textwidth]{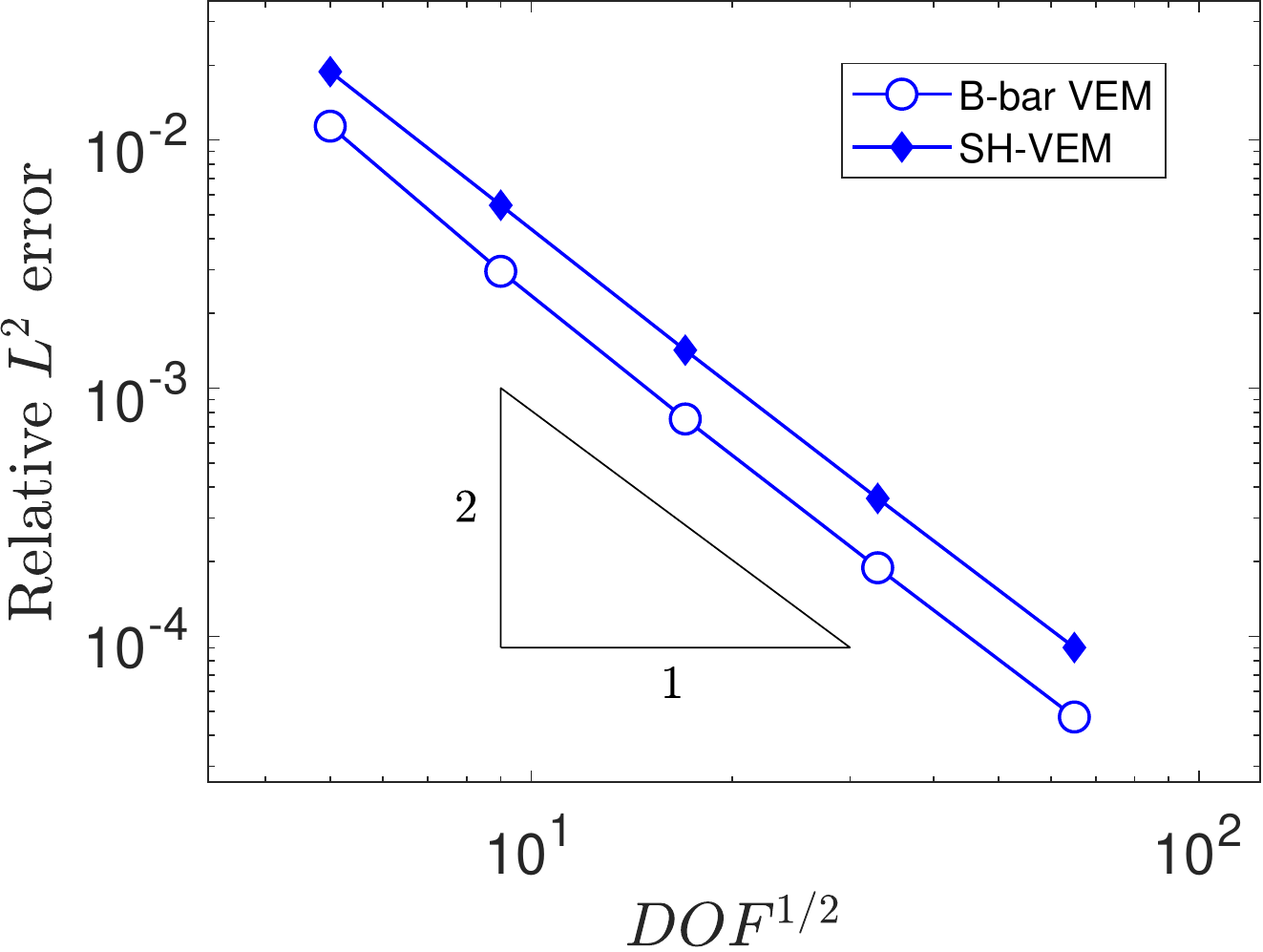}
         \caption{}
     \end{subfigure}
     \hfill
     \begin{subfigure}{0.32\textwidth}
         \centering
         \includegraphics[width=\textwidth]{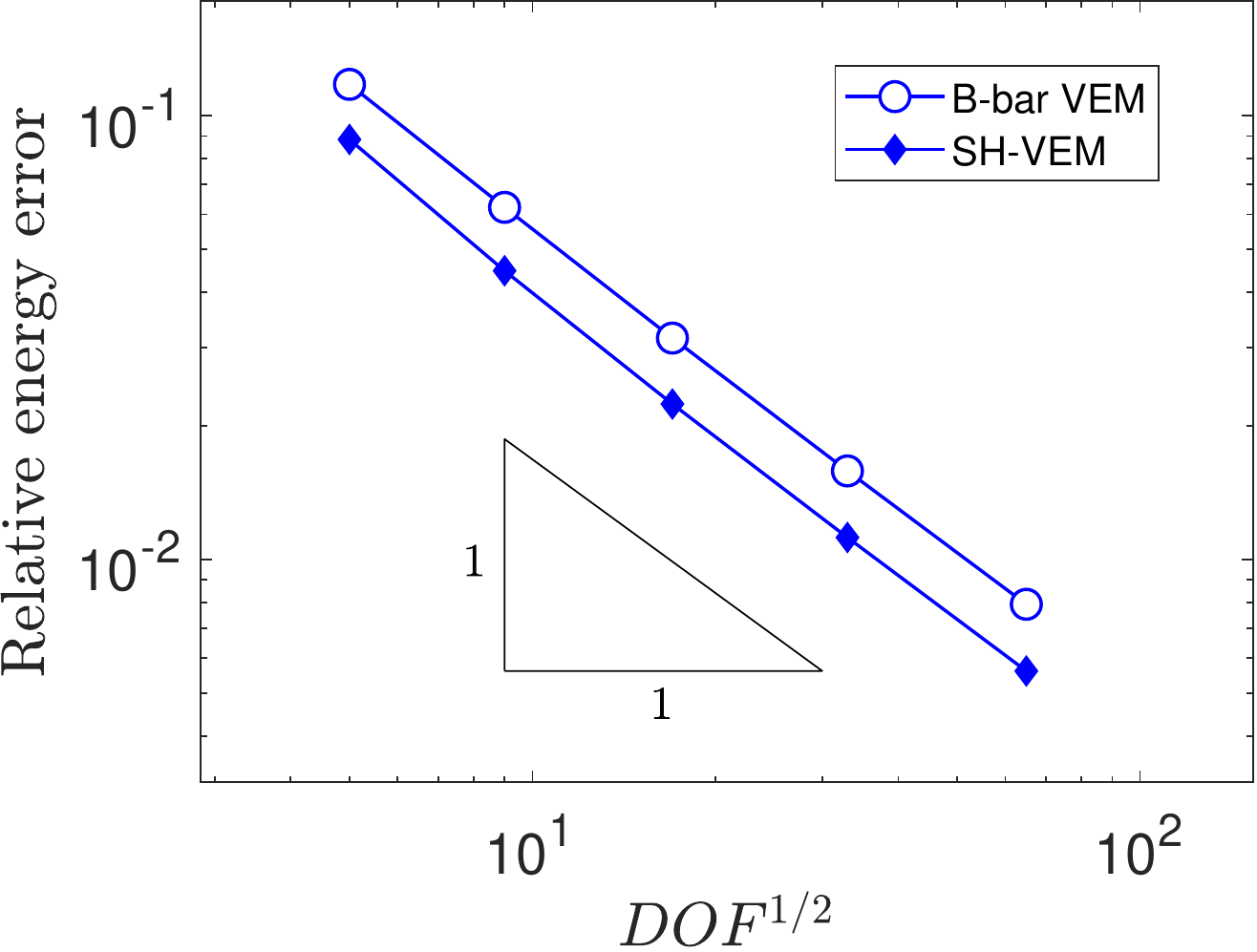}
         \caption{}
     \end{subfigure}
     \hfill
     \begin{subfigure}{0.32\textwidth}
         \centering
         \includegraphics[width=\textwidth]{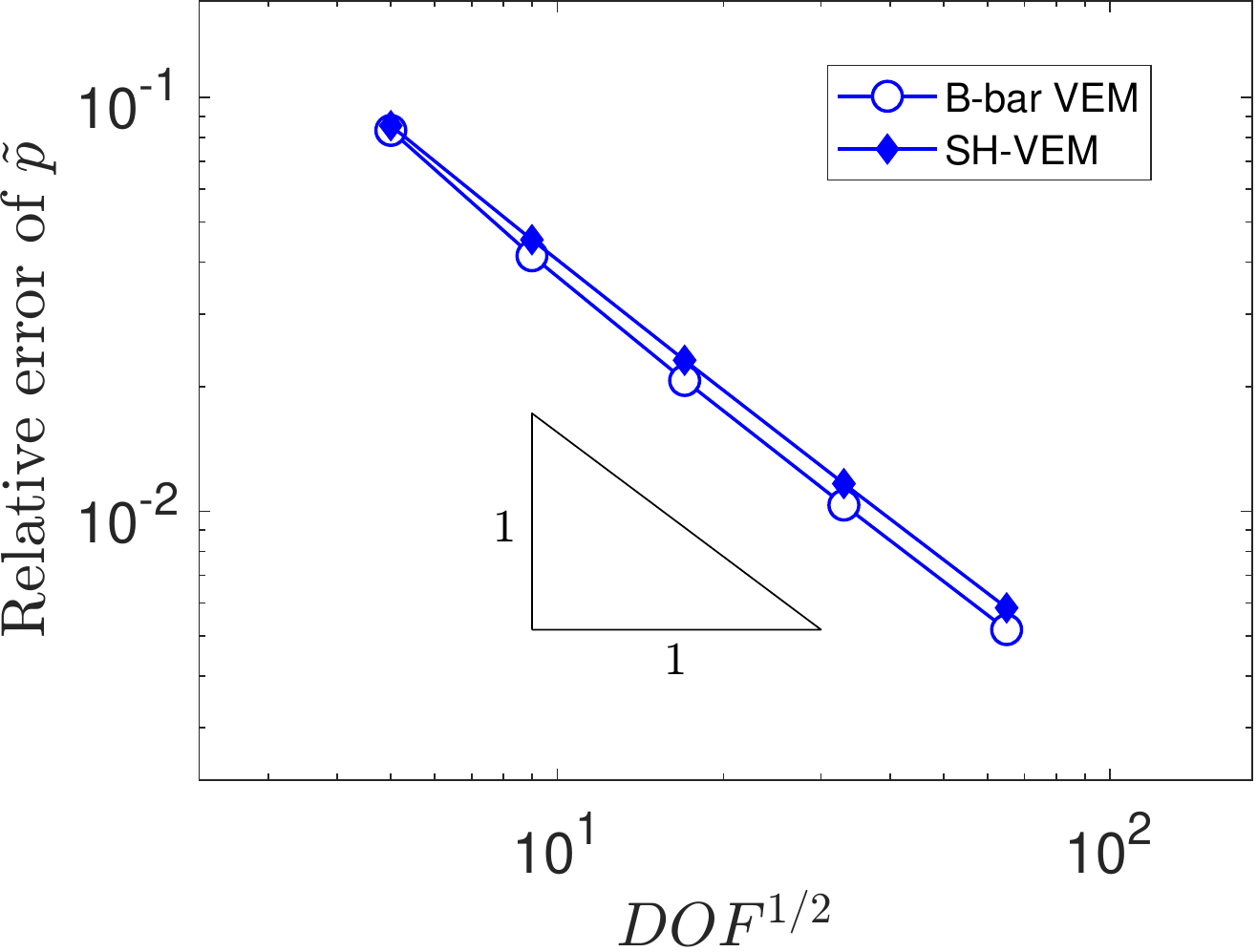}
         \caption{}
     \end{subfigure}
        \caption{Comparison of SH-VEM and B-bar VEM for the plate with a circular hole problem on structured meshes \acrev{(see~\fref{fig:plate_hole_structured})}. (a) 
        $L^2$ error of displacement, (b) energy error and (c) $L^2$ error of the 
        hydrostatic stress. }
        \label{fig:plate_uniform_errors}
\end{figure}

\begin{figure}[!h]
    \centering
    \begin{subfigure}{0.32\textwidth}
        \centering
        \includegraphics[width=\textwidth]{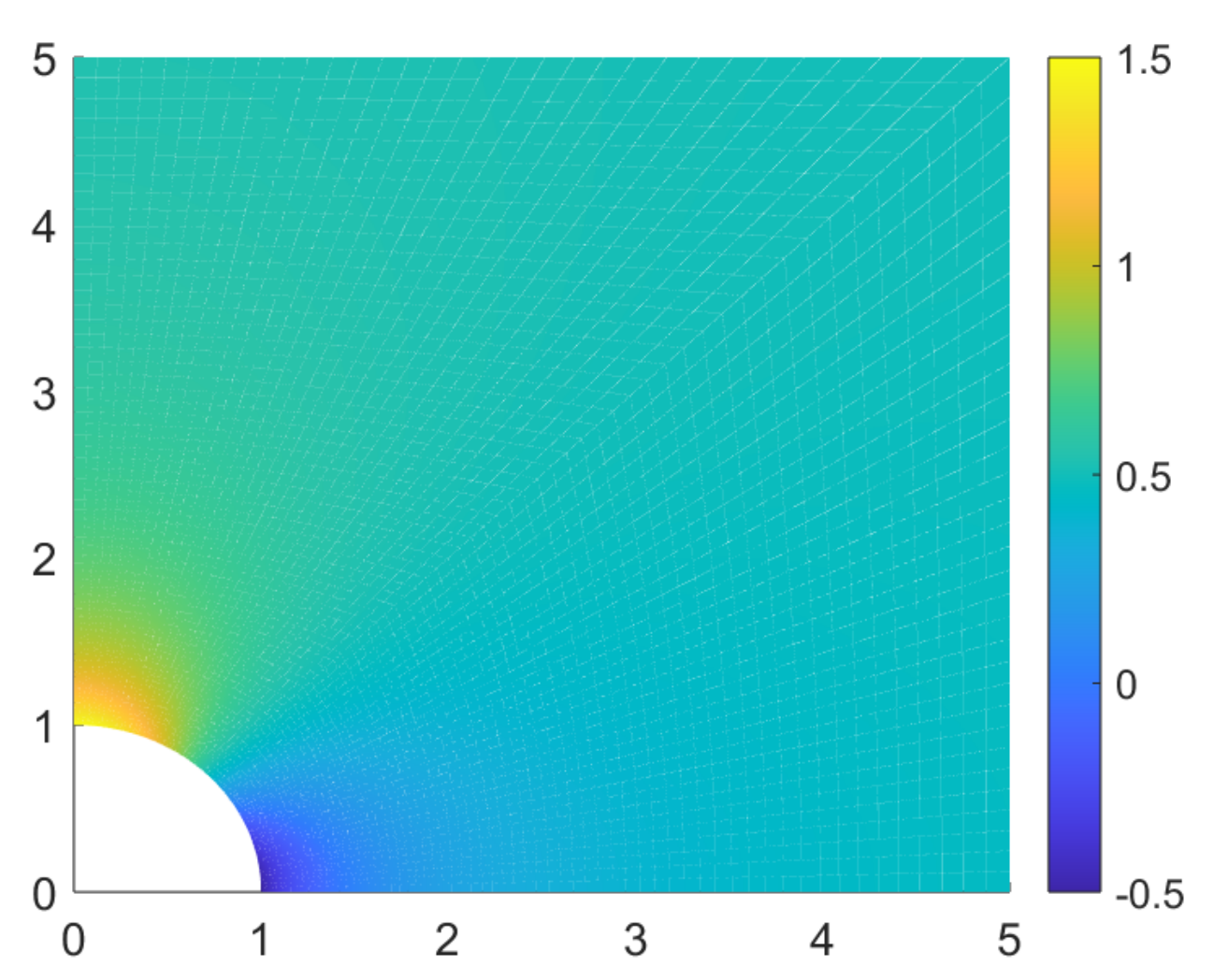}
        \caption{}
    \end{subfigure}
     \hfill
     \begin{subfigure}{0.32\textwidth}
         \centering
         \includegraphics[width=\textwidth]{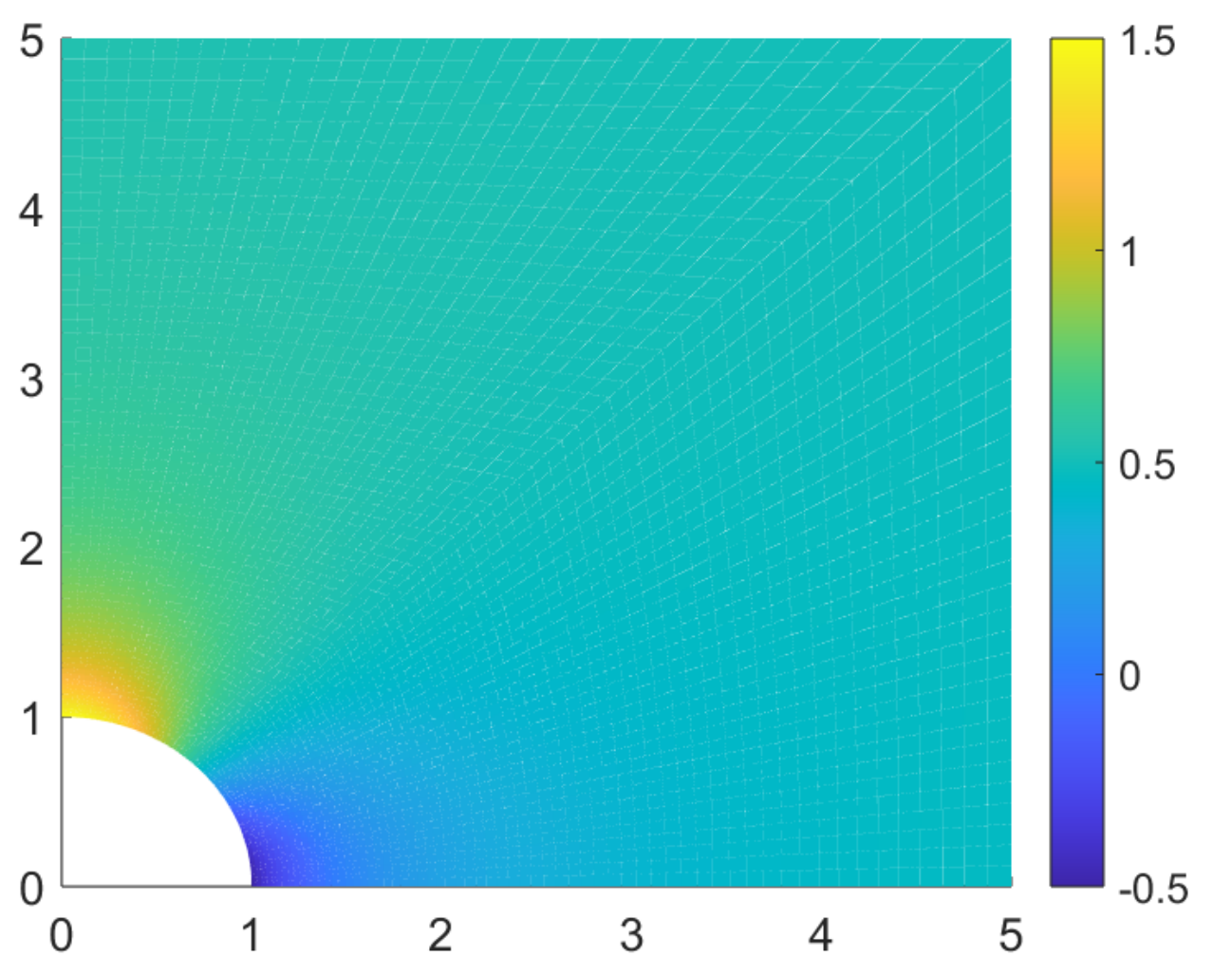}
         \caption{}
     \end{subfigure}
     \hfill
    \begin{subfigure}{0.32\textwidth}
        \centering
        \includegraphics[width=\textwidth]{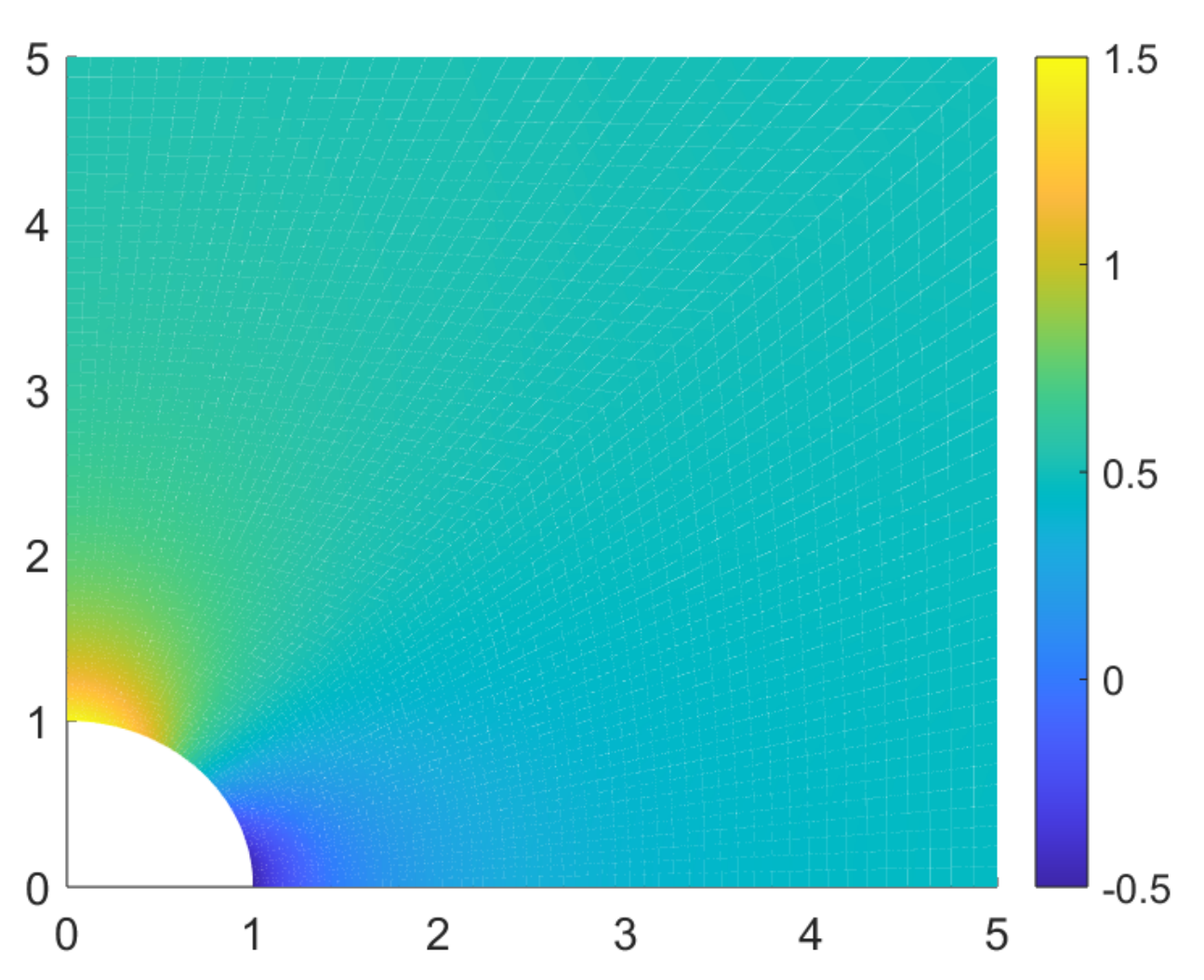}
        \caption{}
    \end{subfigure}
    \caption{Contour plots of the hydrostatic stress on structured meshes \acrev{(see~\fref{fig:plate_hole_structured})} for the plate with a circular hole problem. (a) exact solution, (b) B-bar VEM, (c) SH-VEM.}
    \label{fig:plate_uniform_pressure}
\end{figure}

We consider the problem of an infinite plate with a circular hole under uniaxial tension along the $x$-direction.\cite{timoshenko1951theory} The hole is traction-free, and a far-field tensile load $\sigma_0 = 1$ psi is applied. On using symmetry, we model a quarter of the plate with length $L=5$ inch and a hole of radius $a = 1$ inch. The exact tractions are applied on the traction boundaries. The material has Young's modulus $E_Y = 2\times 10^7$ psi and Poisson's ratio $\nu = 0.4999999$. We first test this problem on structured quadrilateral meshes; a few representative 
meshes are
shown in Figure~\ref{fig:plate_hole_structured}. In Figure~\ref{fig:plate_uniform_errors}, we compare the convergence results of the B-bar formulation and the SH-VEM, and find that both methods deliver optimal convergence rates. In Figure~\ref{fig:plate_uniform_pressure}, we also compare the contours of the hydrostatic stress by the two methods and find that they both are smooth
and have comparable accuracy. 

\begin{figure}[!h]
     \centering
     \begin{subfigure}{.32\textwidth}
         \centering
         \includegraphics[width=\textwidth]{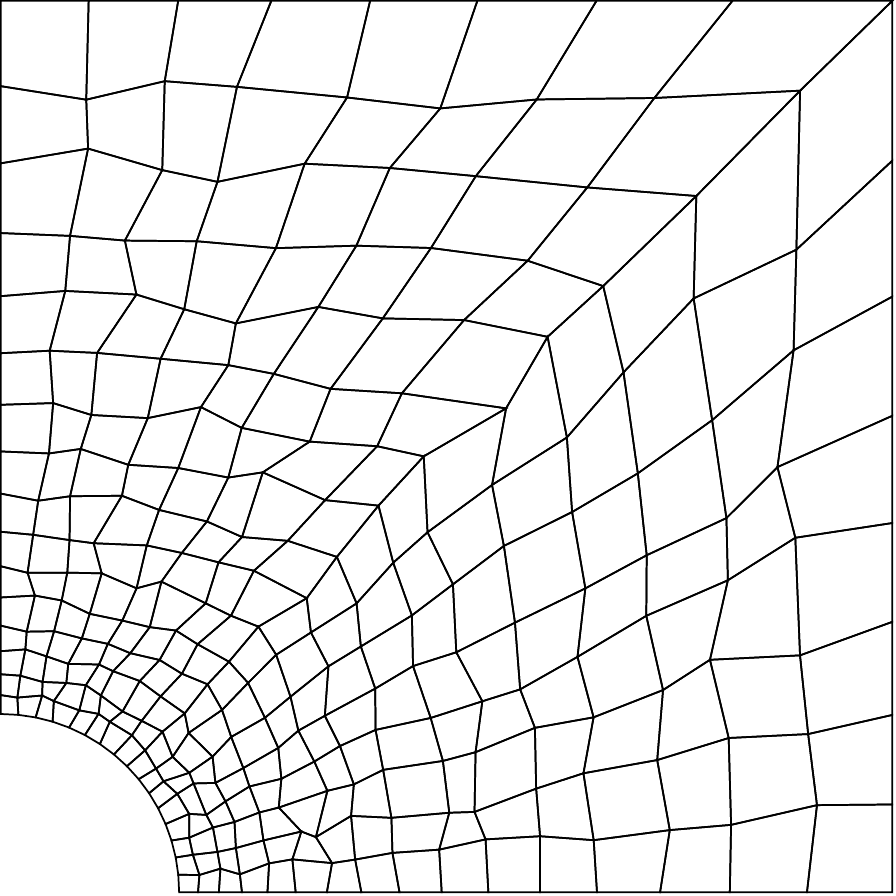}
         \caption{}
     \end{subfigure}
     \hfill
     \begin{subfigure}{.32\textwidth}
         \centering
         \includegraphics[width=\textwidth]{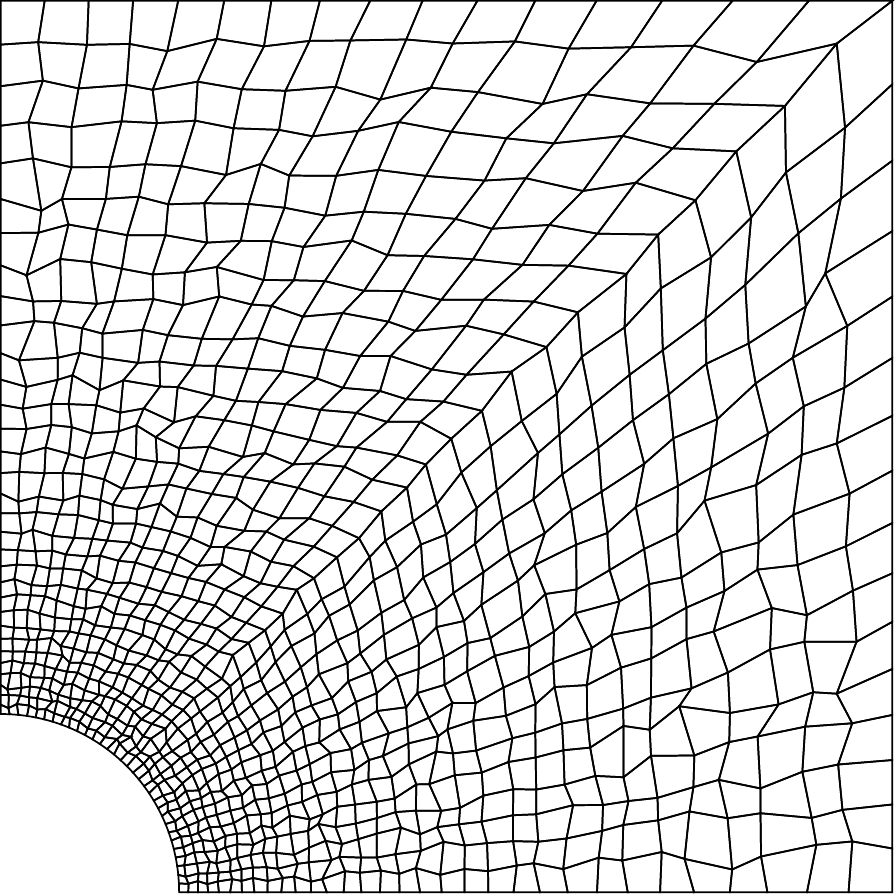}
         \caption{}
     \end{subfigure}
     \hfill
     \begin{subfigure}{.32\textwidth}
         \centering
         \includegraphics[width=\textwidth]{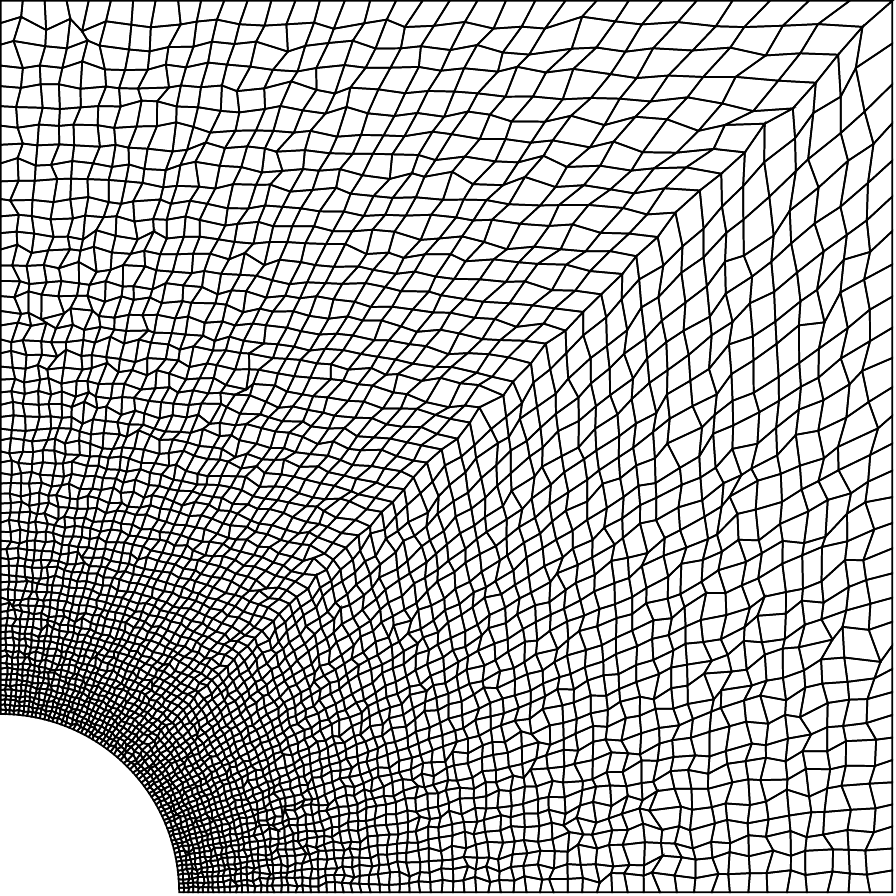}
         \caption{}
     \end{subfigure}
        \caption{Perturbed quadrilateral meshes for the plate with a hole problem. (a) 256 elements, (b) 1024 elements and (c) 4096 elements.  }
        \label{fig:plate_hole_perturbed}
\end{figure}

\begin{figure}[!h]
     \centering
     \begin{subfigure}{0.32\textwidth}
         \centering
         \includegraphics[width=\textwidth]{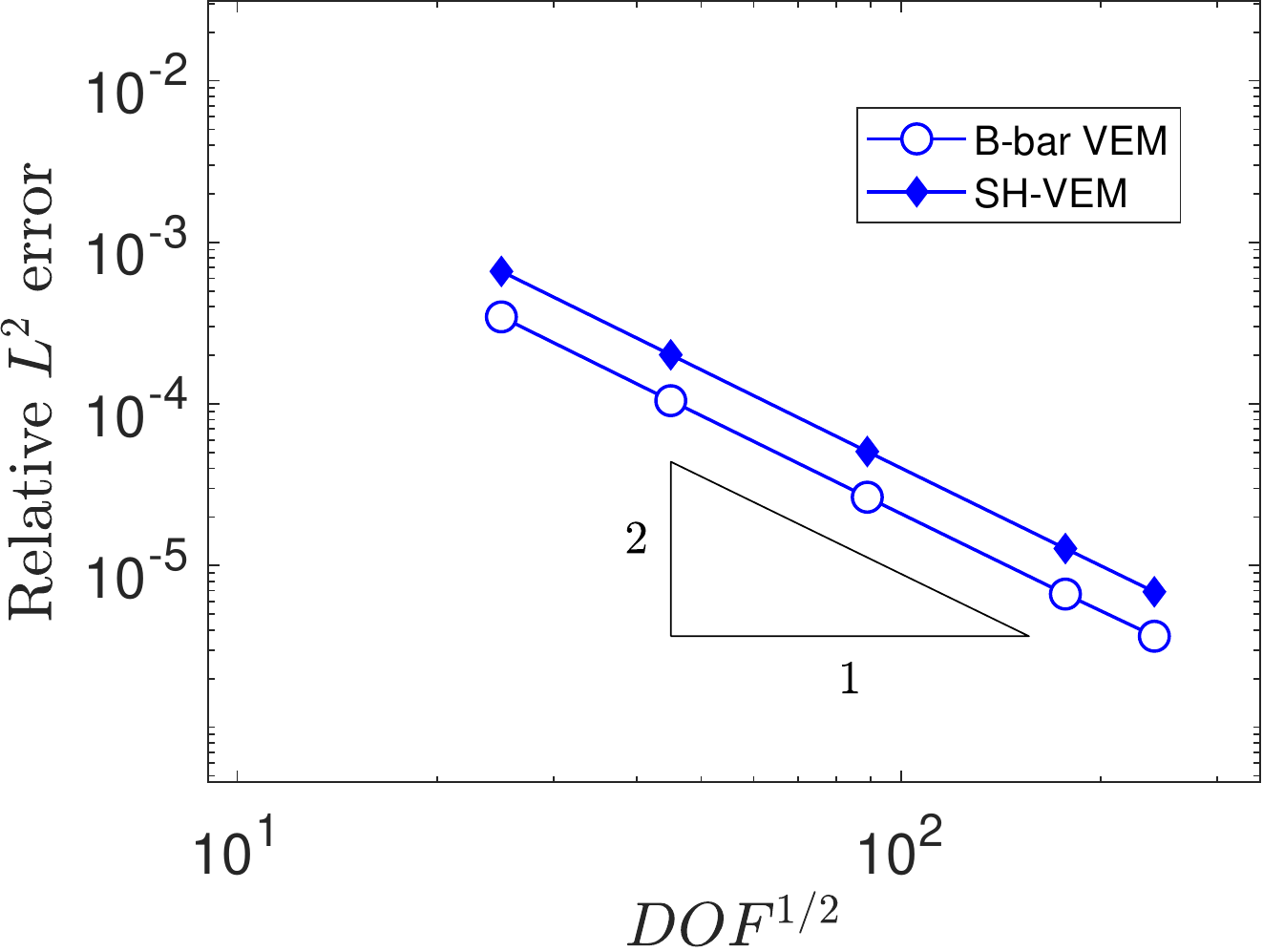}
         \caption{}
     \end{subfigure}
     \hfill
     \begin{subfigure}{0.32\textwidth}
         \centering
         \includegraphics[width=\textwidth]{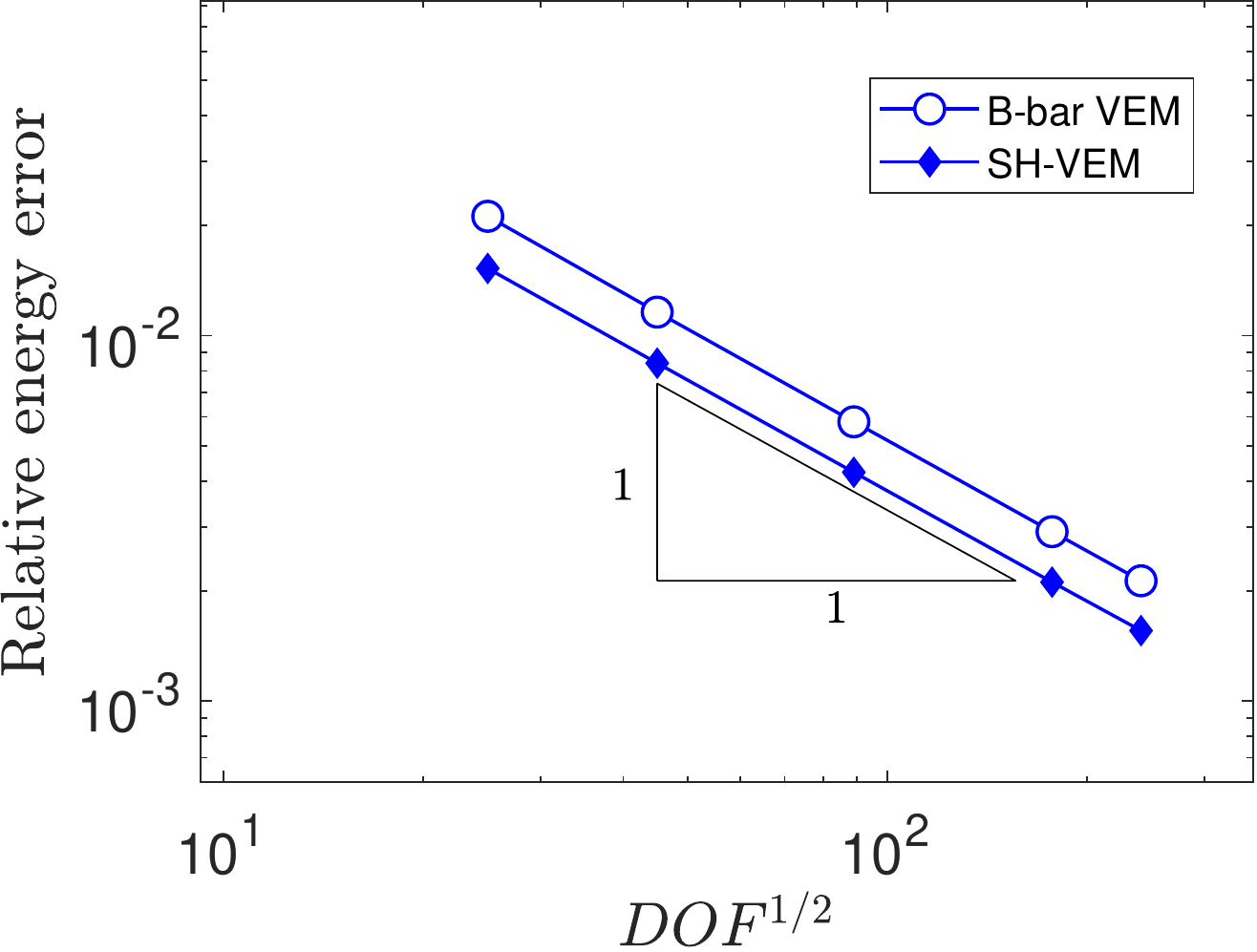}
         \caption{}
     \end{subfigure}
     \hfill
     \begin{subfigure}{0.32\textwidth}
         \centering
         \includegraphics[width=\textwidth]{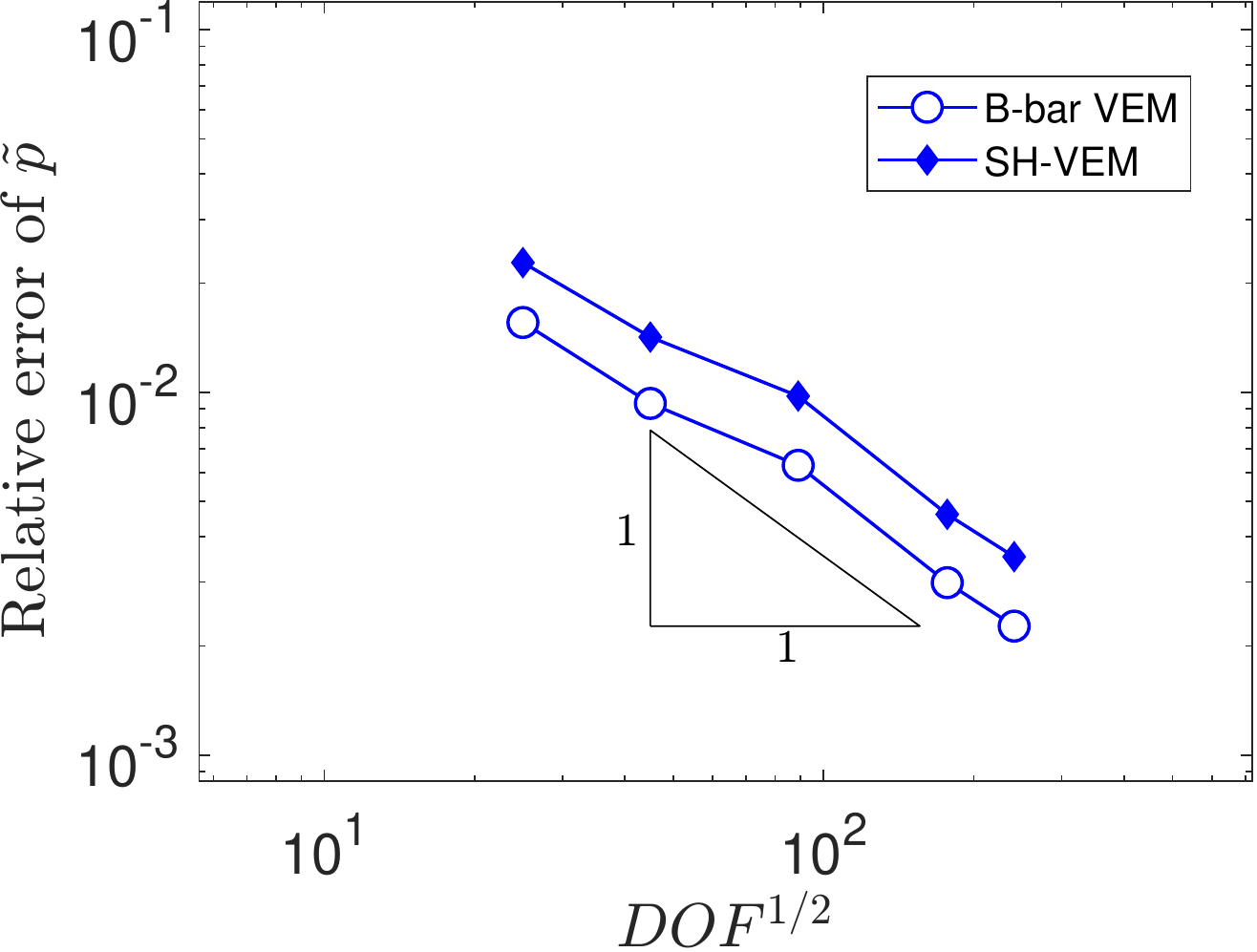}
         \caption{}
     \end{subfigure}
        \caption{Comparison of SH-VEM and B-bar VEM for the plate with a circular hole problem on perturbed meshes \acrev{(see~\fref{fig:plate_hole_perturbed})}. (a) $L^2$ error of displacement, (b) energy error and (c) $L^2$ error of hydrostatic 
        stress. }
        \label{fig:plate_perturbed_errors}
\end{figure}

\begin{figure}[!h]
    \centering
    \begin{subfigure}{0.32\textwidth}
    \centering
    \includegraphics[width=\textwidth]{plate_hole_pressure_exact.pdf}
    \caption{}
    \end{subfigure}
     \hfill
     \begin{subfigure}{0.32\textwidth}
         \centering
         \includegraphics[width=\textwidth]{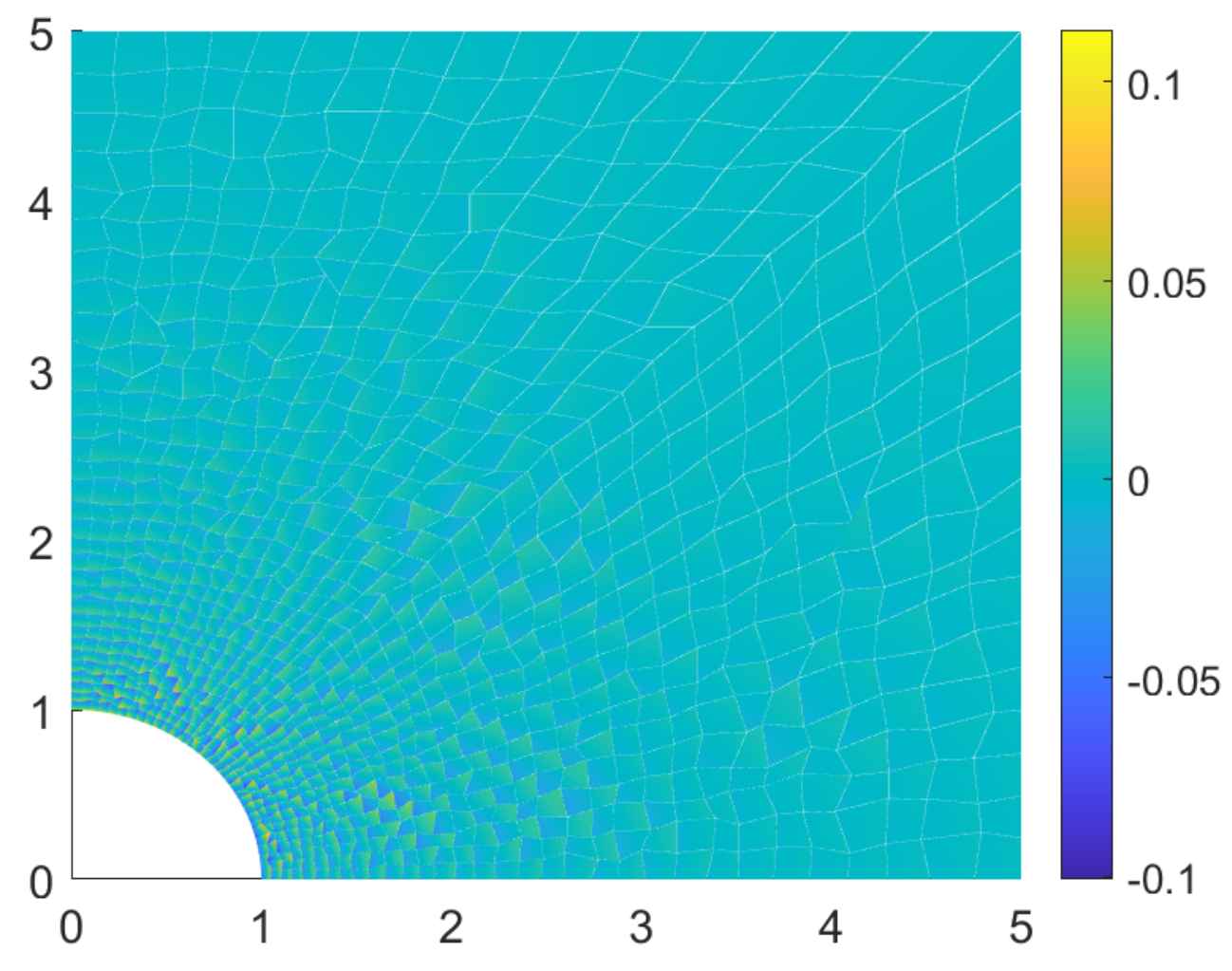}
         \caption{}
     \end{subfigure}
     \hfill
    \begin{subfigure}{0.32\textwidth}
    \centering
    \includegraphics[width=\textwidth]{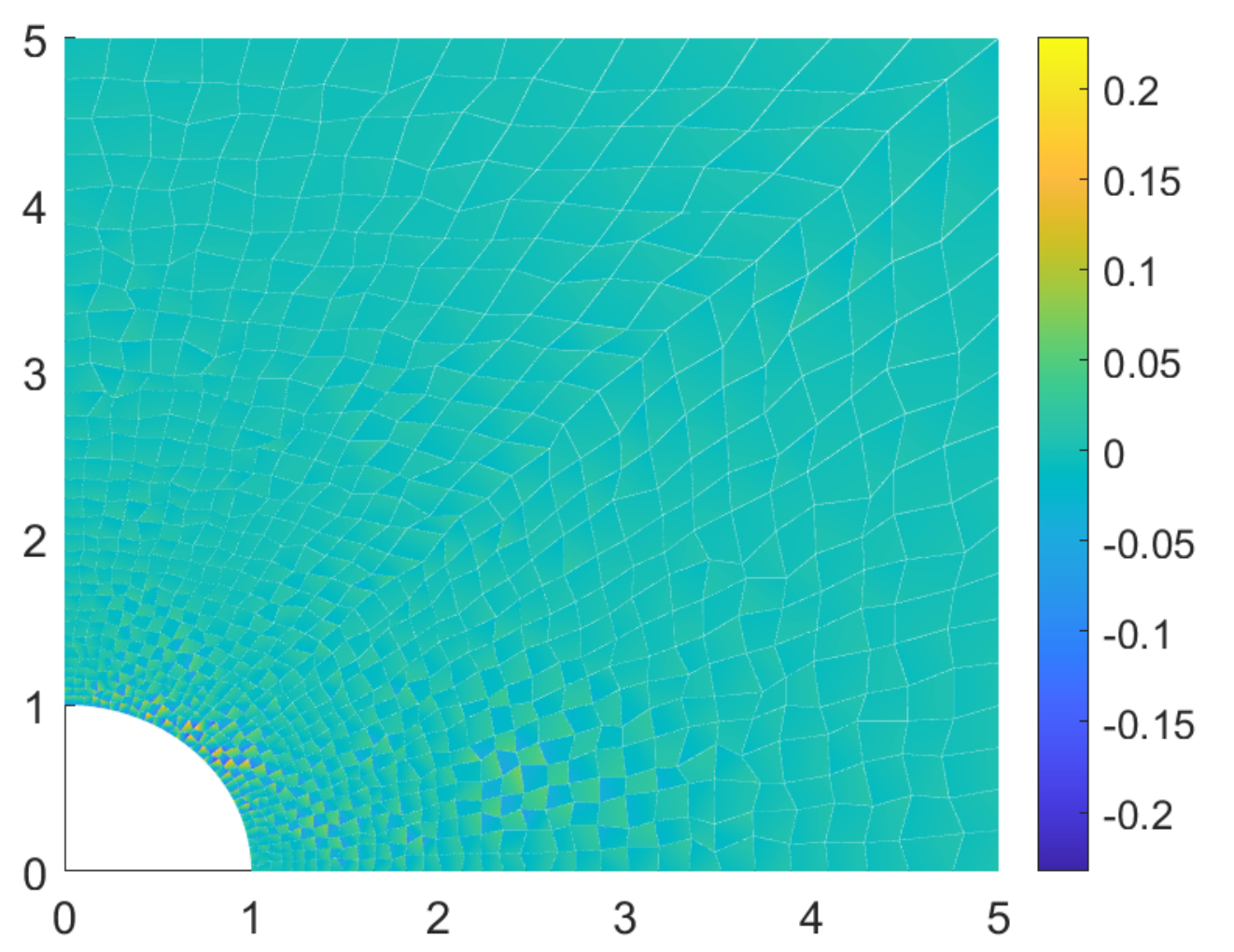}
    \caption{}
    \end{subfigure}
    \caption{Contour plots of the 
    hydrostatic stress for the plate with a circular hole problem \acrev{on perturbed meshes (see~\fref{fig:plate_hole_perturbed})}. (a) exact solution, and error in the hydrostatic stress field, 
    $| \tilde{p} - \tilde{p}_h|$ for (b) B-bar VEM and (c) SH-VEM.} 
    \label{fig:plate_perturbed_pressure}
\end{figure}

We now consider the plate with a circular hole problem on a perturbed mesh. We start with a structured mesh, then for each internal node we perturb its location. Representative meshes are shown in Figure~\ref{fig:plate_hole_perturbed}. In Figure~\ref{fig:plate_perturbed_errors}, we show the convergence rates of the two methods and find that both methods retain optimal convergence on the perturbed mesh. In Figure~\ref{fig:plate_perturbed_pressure}, the exact
hydrostatic stress and contour plots of the error, $| \tilde{p}-\tilde{p}_h|$, are shown. The plots reveal that both methods produce relatively
smooth error distributions of the hydrostatic stress field, with the B-bar VEM having smaller pointwise error than SH-VEM.


\subsection{Hollow cylinder under internal pressure}
\begin{figure}[!h]
     \centering
     \begin{subfigure}{.32\textwidth}
         \centering
         \includegraphics[width=\textwidth]{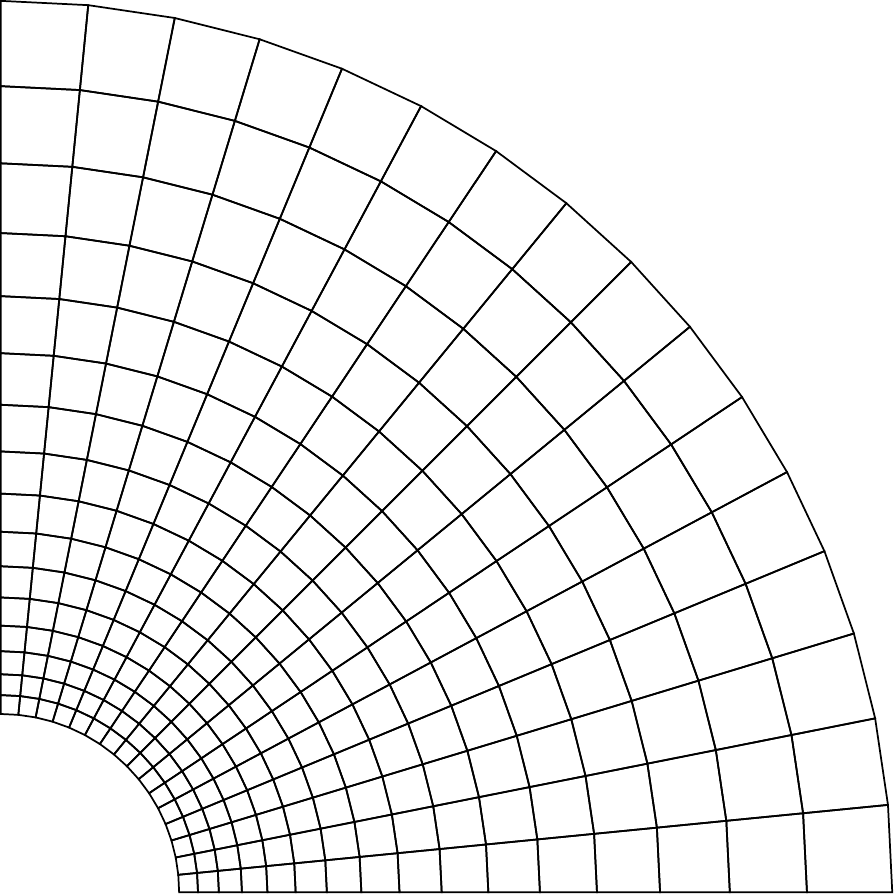}
         \caption{}
     \end{subfigure}
     \hfill
     \begin{subfigure}{.32\textwidth}
         \centering
         \includegraphics[width=\textwidth]{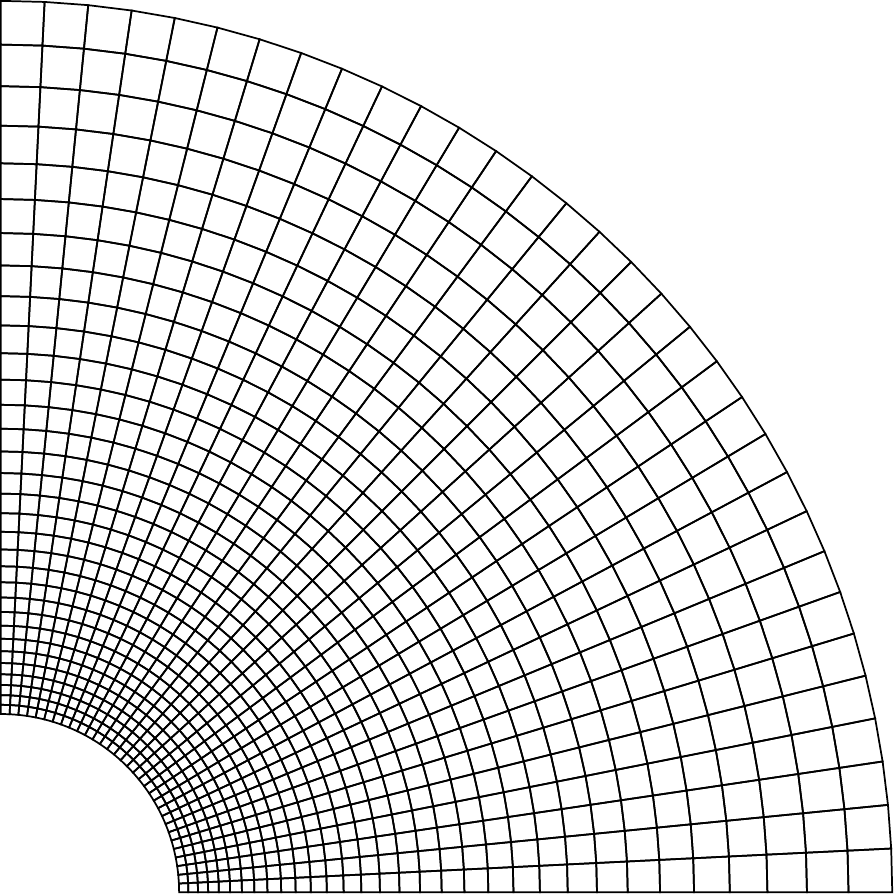}
         \caption{}
     \end{subfigure}
     \hfill
     \begin{subfigure}{.32\textwidth}
         \centering
         \includegraphics[width=\textwidth]{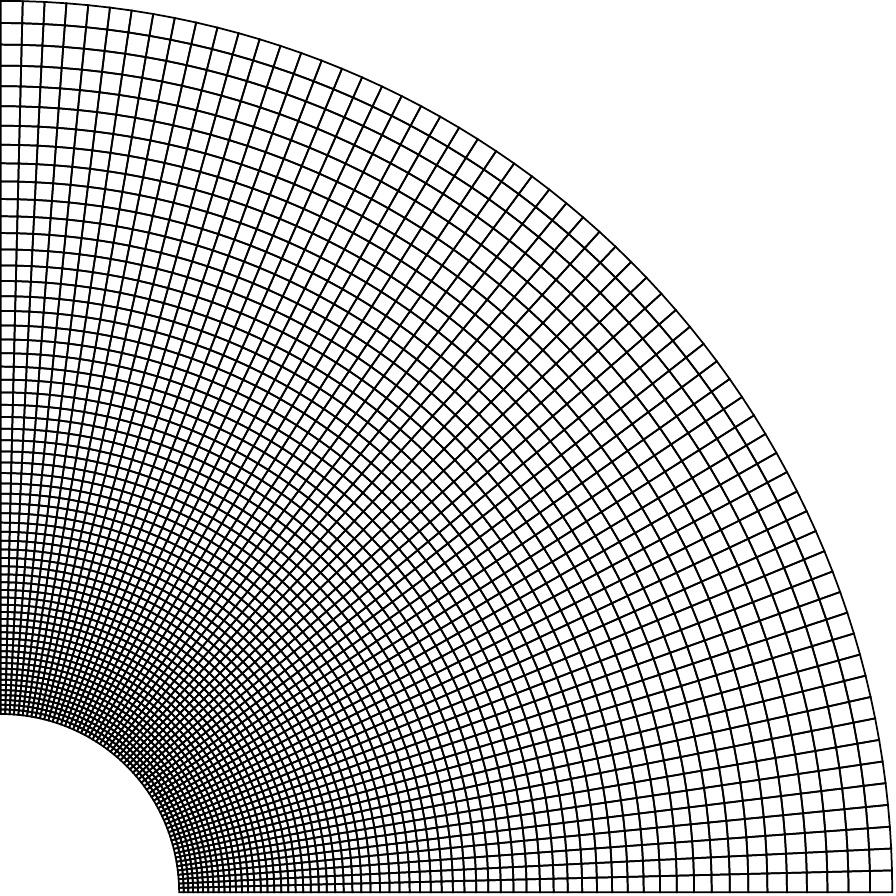}
         \caption{}
     \end{subfigure}
        \caption{Uniform quadrilateral meshes for the hollow cylinder problem. (a) 256 elements, (b) 1024 elements and (c) 4096 elements.  }
        \label{fig:cylinder_uniform}
\end{figure}

\begin{figure}[!h]
     \centering
     \begin{subfigure}{0.32\textwidth}
         \centering
         \includegraphics[width=\textwidth]{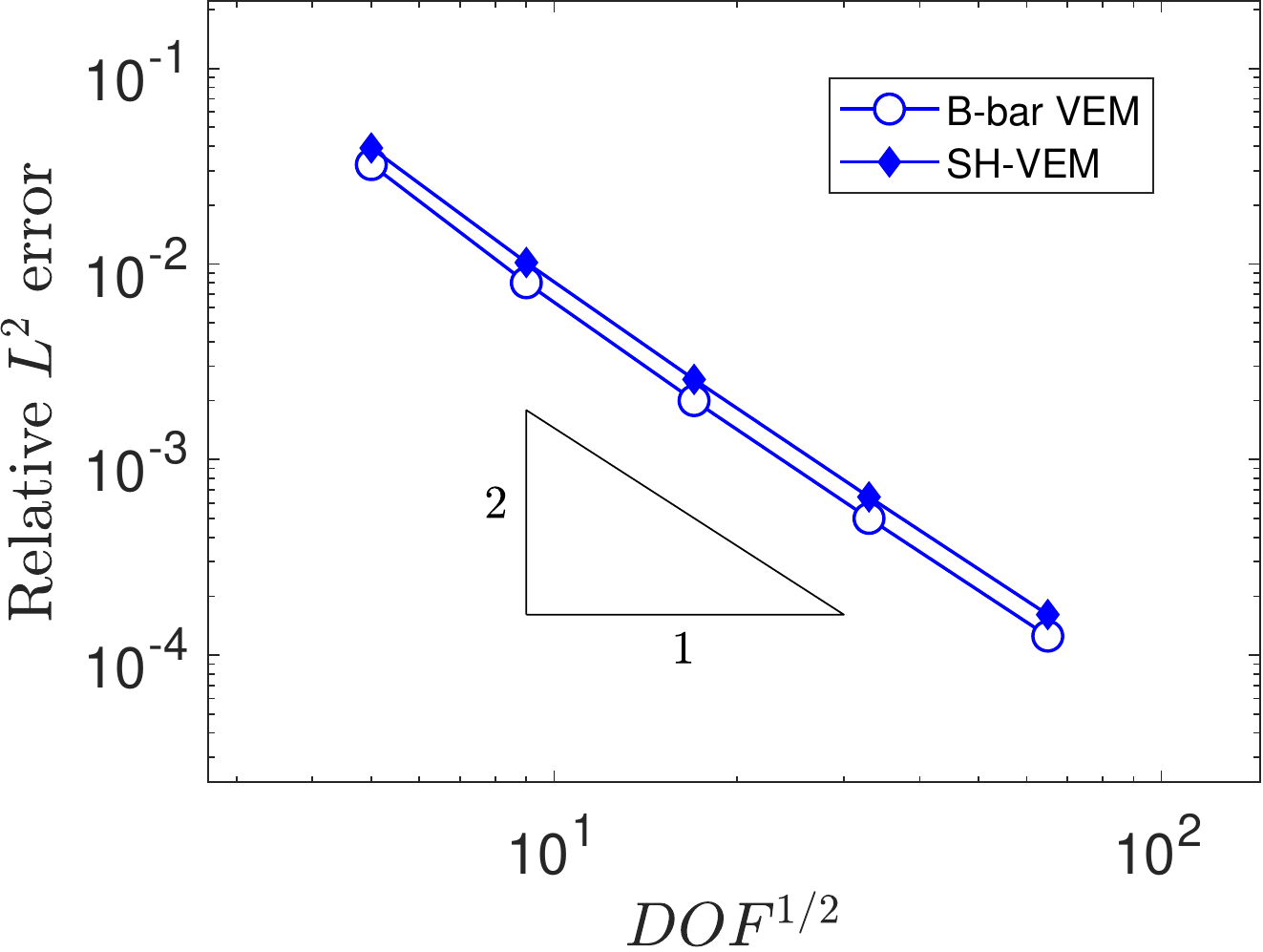}
         \caption{}
     \end{subfigure}
     \hfill
     \begin{subfigure}{0.32\textwidth}
         \centering
         \includegraphics[width=\textwidth]{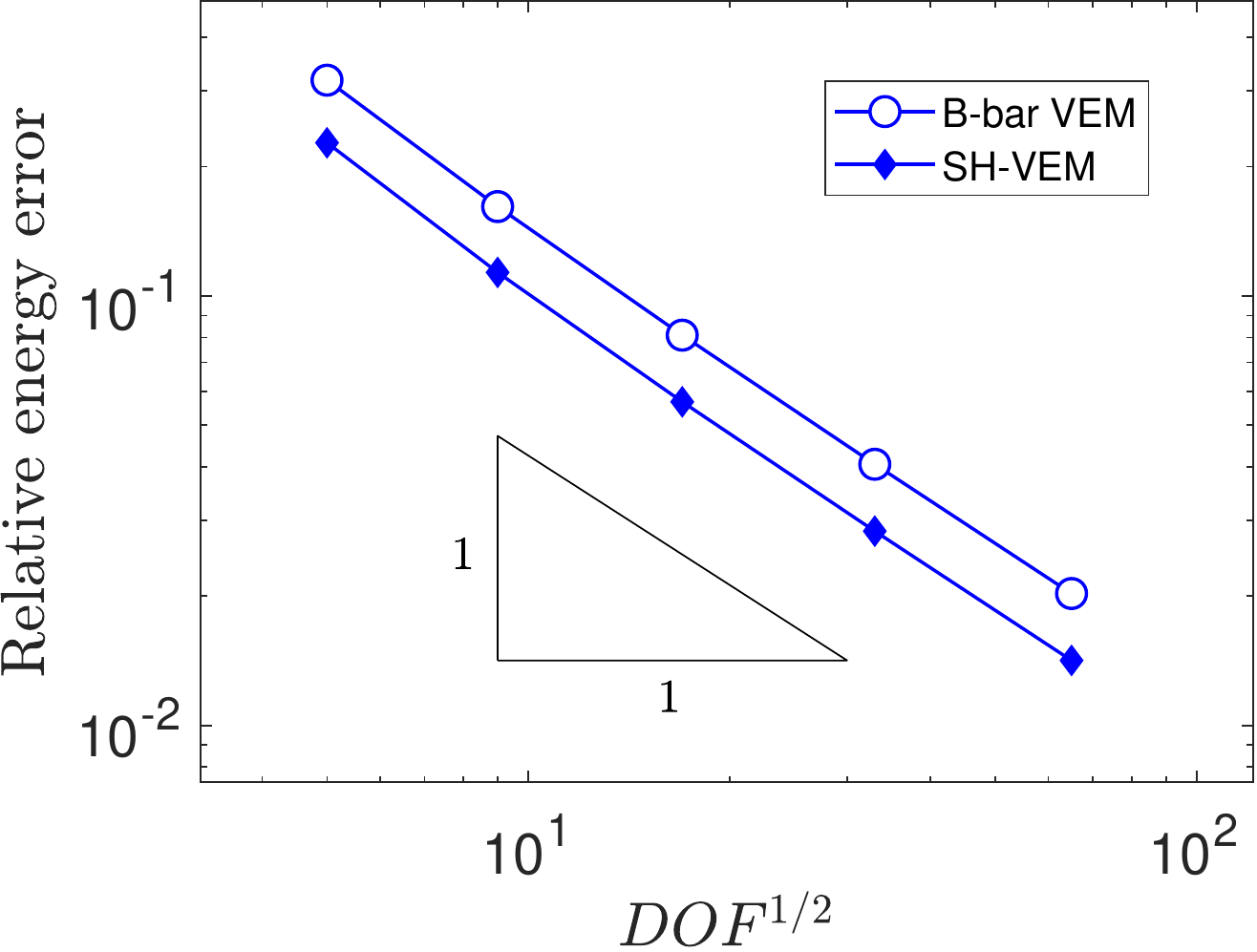}
         \caption{}
     \end{subfigure}
     \hfill
     \begin{subfigure}{0.32\textwidth}
         \centering
         \includegraphics[width=\textwidth]{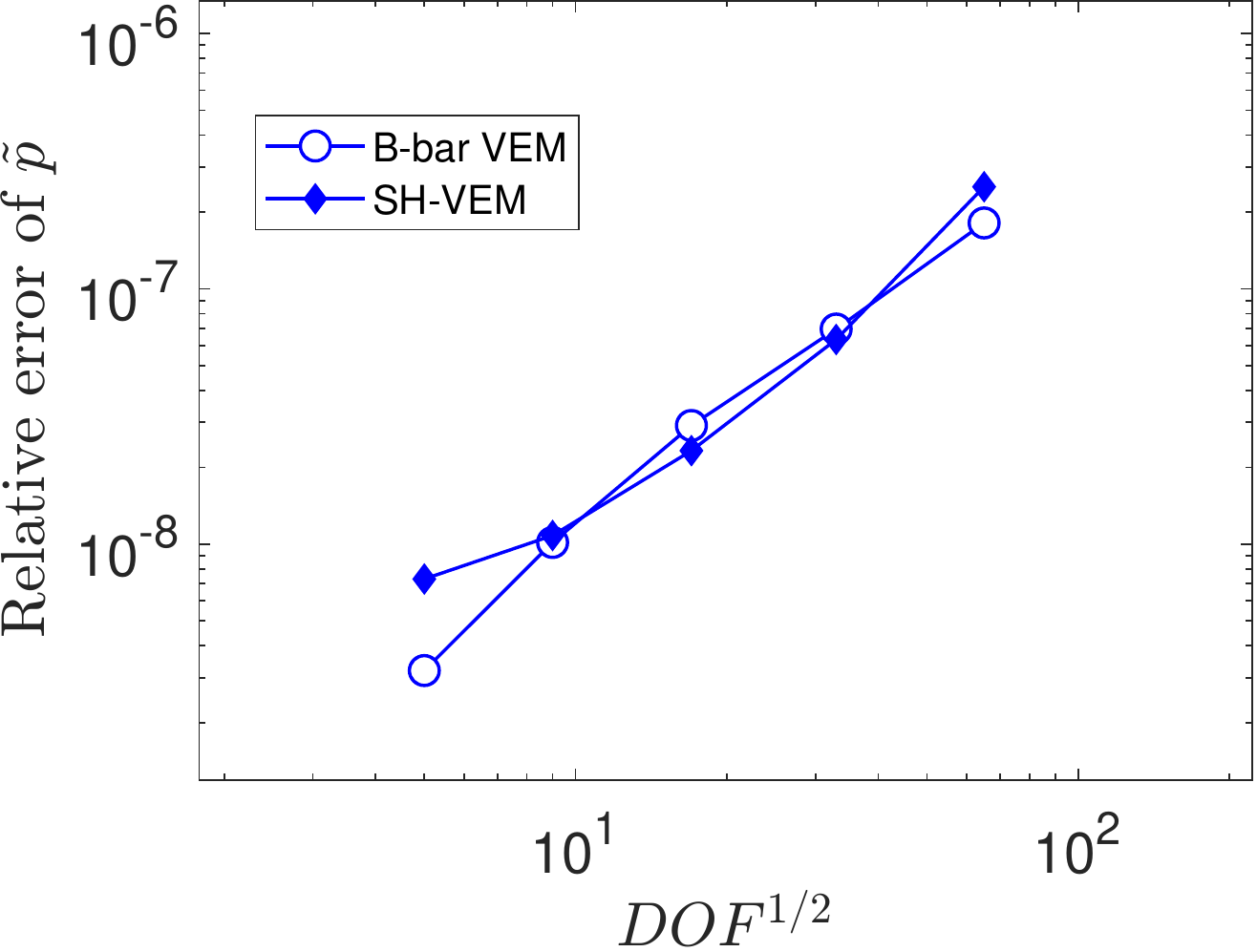}
         \caption{}
     \end{subfigure}
        \caption{Comparison of B-bar VEM and SH-VEM for the hollow cylinder problem on structured meshes \acrev{(see~\fref{fig:cylinder_uniform})}. (a) $L^2$ error of displacement,  
        (b) energy error and (c) $L^2$ error of hydrostatic stress.}
        \label{fig:cylinder_uniform_errors}
\end{figure}

\begin{figure}[!h]
    \centering
     \begin{subfigure}{0.48\textwidth}
         \centering
         \includegraphics[width=\textwidth]{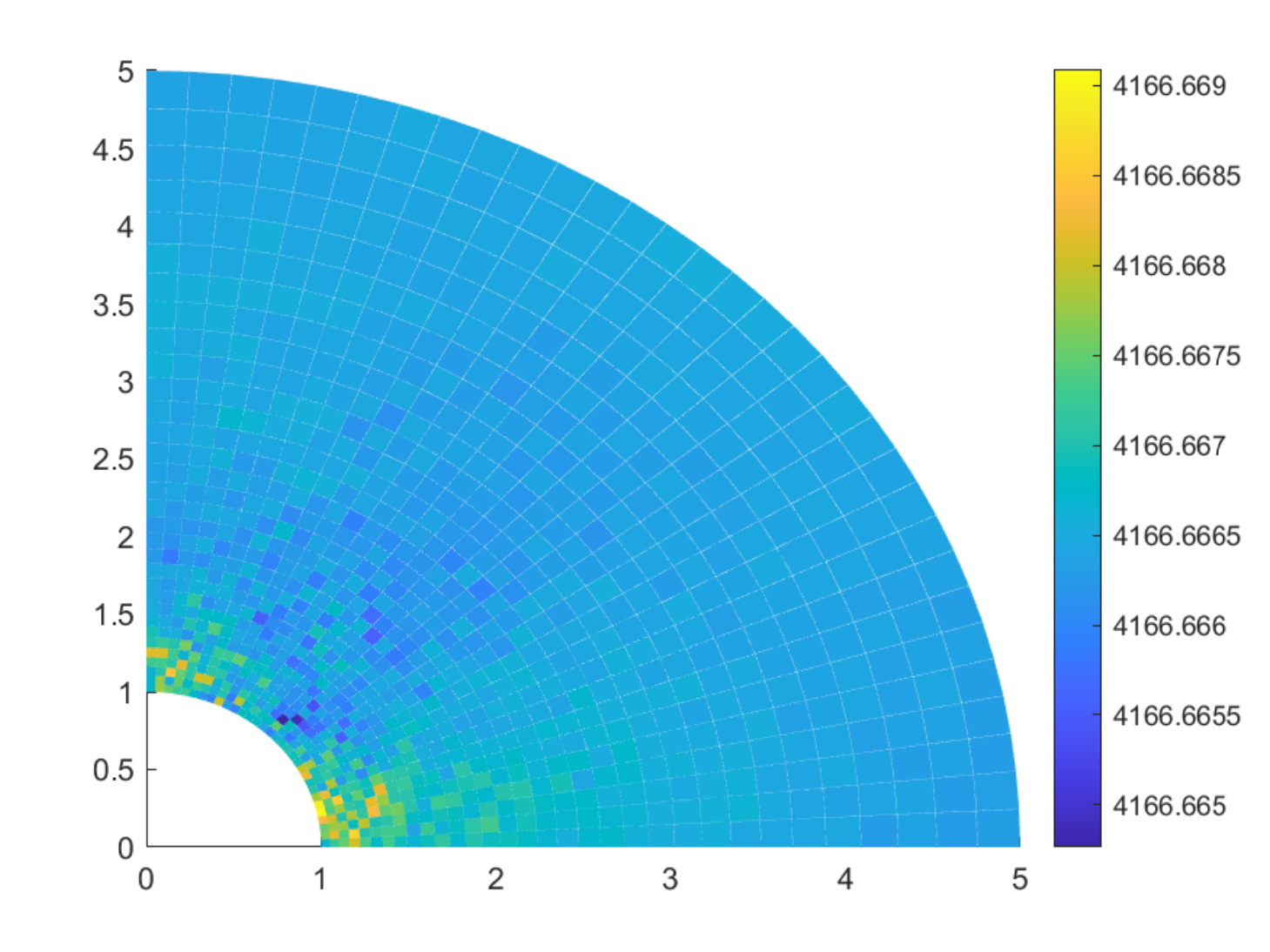}
         \caption{}
     \end{subfigure}
     \hfill
    \begin{subfigure}{0.48\textwidth}
        \centering
        \includegraphics[width=\textwidth]{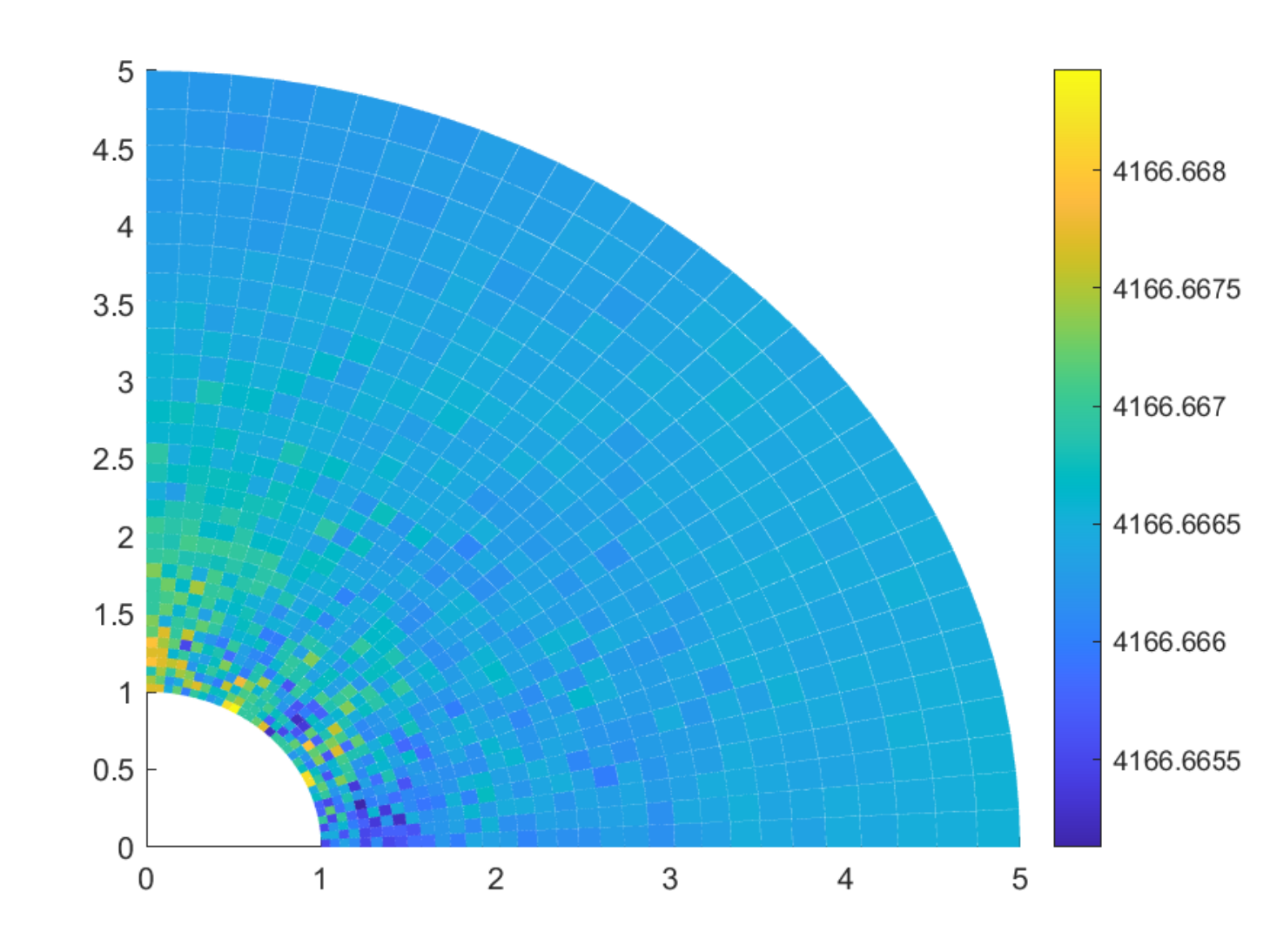}
        \caption{}
    \end{subfigure}
    \caption{Contour plots of the hydrostatic stress field on structured meshes \acrev{(see~\fref{fig:cylinder_uniform})} for the hollow pressurized cylinder problem. 
    The exact hydrostatic stress is 4166.6666 psi.
    (a) B-bar VEM and
    (b) SH-VEM.}
    \label{fig:cylinder_uniform_pressure}
\end{figure}

We consider the problem of a hollow cylinder with inner radius 
$a = 1$ inch and outer radius $b = 5$ inch under internal pressure.~\cite{timoshenko1951theory} Due to symmetry, we model this problem as a quarter cylinder. A uniform pressure of $p=10^5$ psi is applied on the inner radius, while the outer radius is kept traction-free. The material has Young's modulus $E_Y = 2\times 10^5$ psi and Poisson's ratio $\nu=0.4999999$.
\ac{For this example, the hydrostatic stress field is constant; therefore, we use an element averaged approximation to compute the hydrostatic stress $\tilde{p}_h$.} We first examine this problem on structured quadrilateral meshes; a few representative meshes are presented 
in Figure~\ref{fig:cylinder_uniform}. In Figure~\ref{fig:cylinder_uniform_errors}, the convergence rates of 
B-bar VEM and SH-VEM are shown. For both methods, convergence in $L^2$ norm and energy seminorm is optimal. The contour plots in Figure~\ref{fig:cylinder_uniform_pressure} show that both methods are
able to reproduce the constant exact hydrostatic stress field on a uniform mesh. 

\begin{figure}[!h]
     \centering
     \begin{subfigure}{.32\textwidth}
         \centering
         \includegraphics[width=\textwidth]{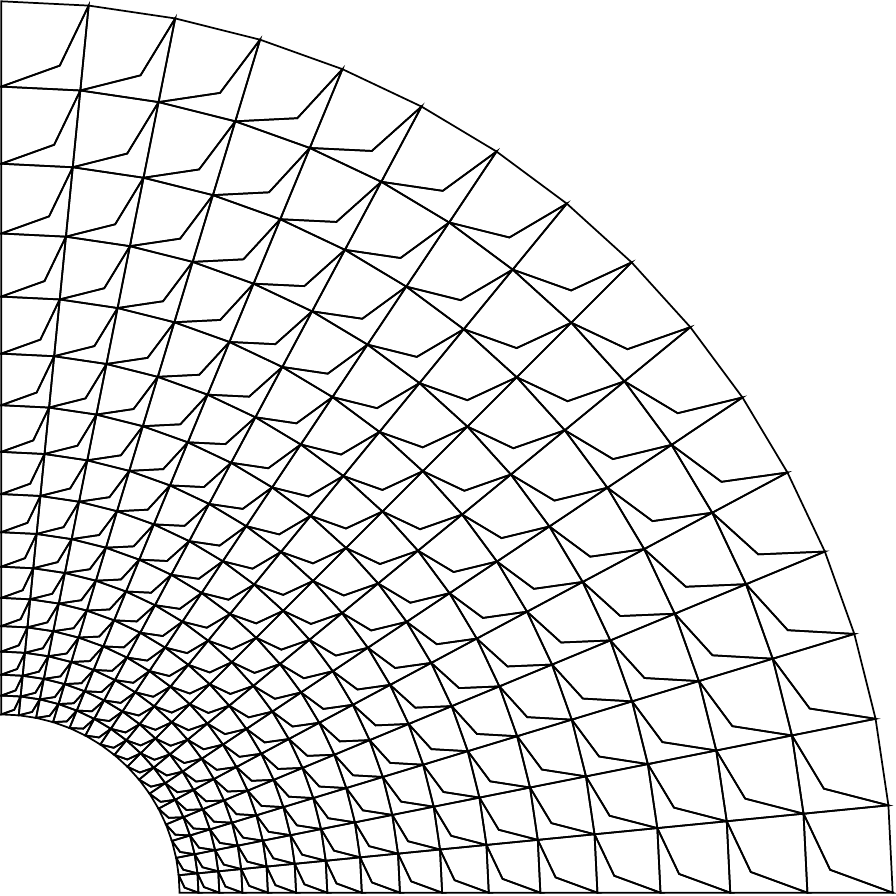}
         \caption{}
     \end{subfigure}
     \hfill
     \begin{subfigure}{.32\textwidth}
         \centering
         \includegraphics[width=\textwidth]{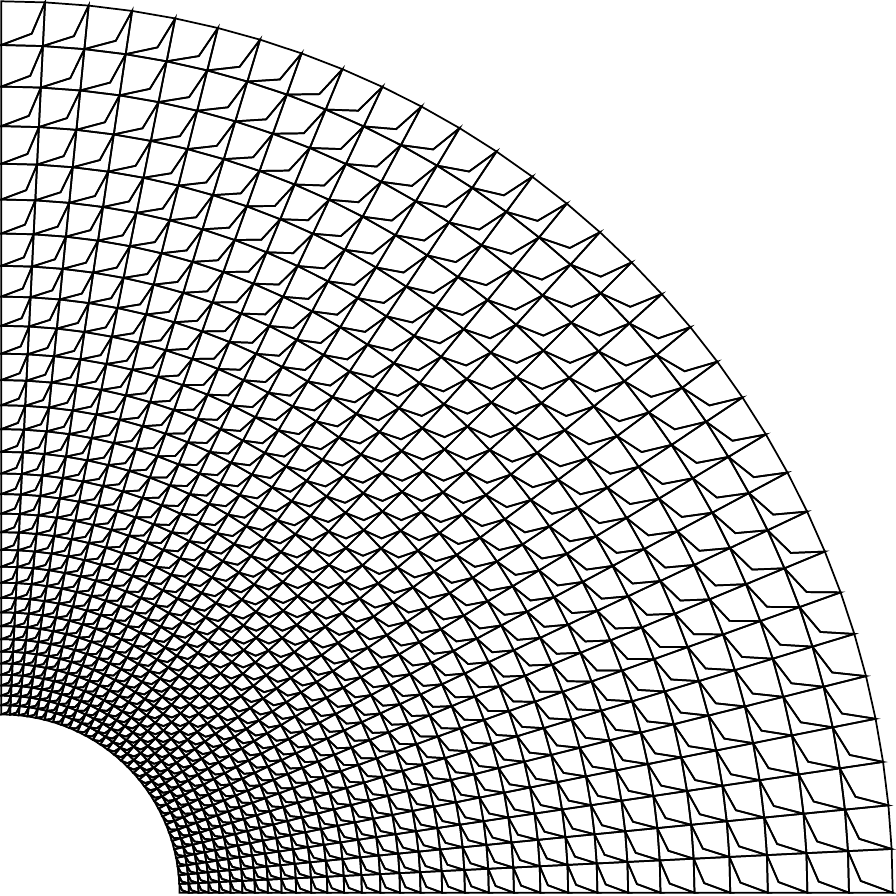}
         \caption{}
     \end{subfigure}
     \hfill
     \begin{subfigure}{.32\textwidth}
         \centering
         \includegraphics[width=\textwidth]{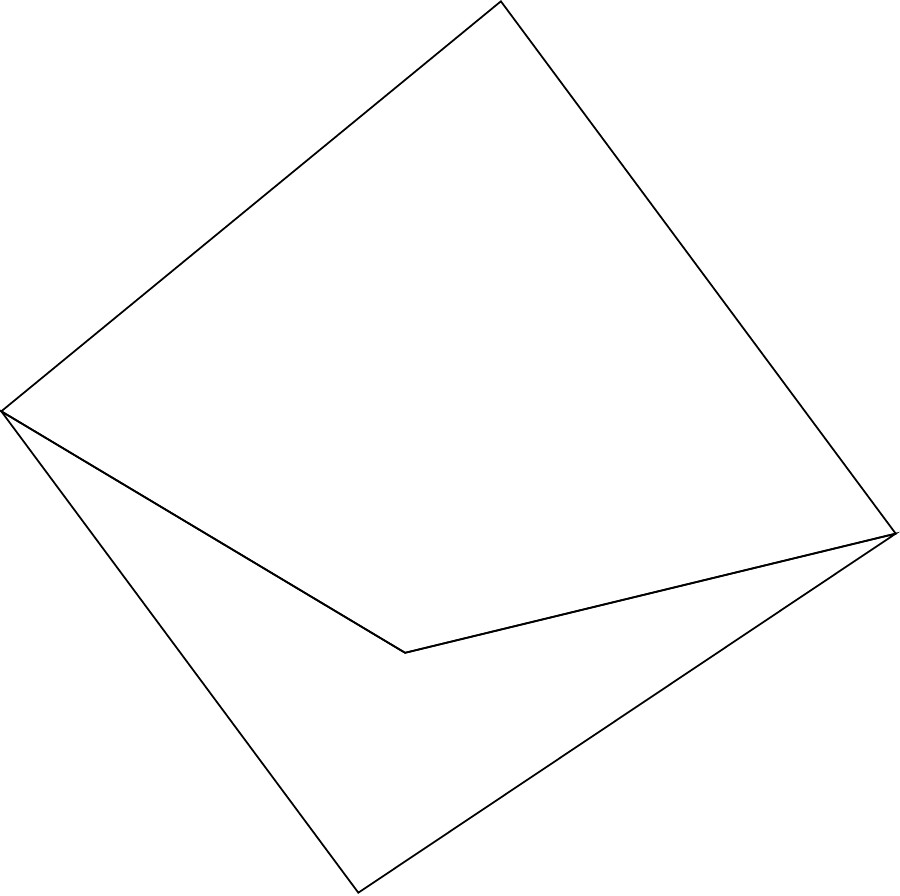}
         \caption{}\label{fig:cylinder_nonconvex-c}
     \end{subfigure}
        \caption{Nonconvex quadrilateral meshes for the hollow cylinder problem. (a) 512 elements, (b) 2048 elements and (c) magnification of a single element split into convex and nonconvex partitions.  }
        \label{fig:cylinder_nonconvex}
\end{figure}

\begin{figure}[!h]
     \centering
     \begin{subfigure}{0.32\textwidth}
         \centering
         \includegraphics[width=\textwidth]{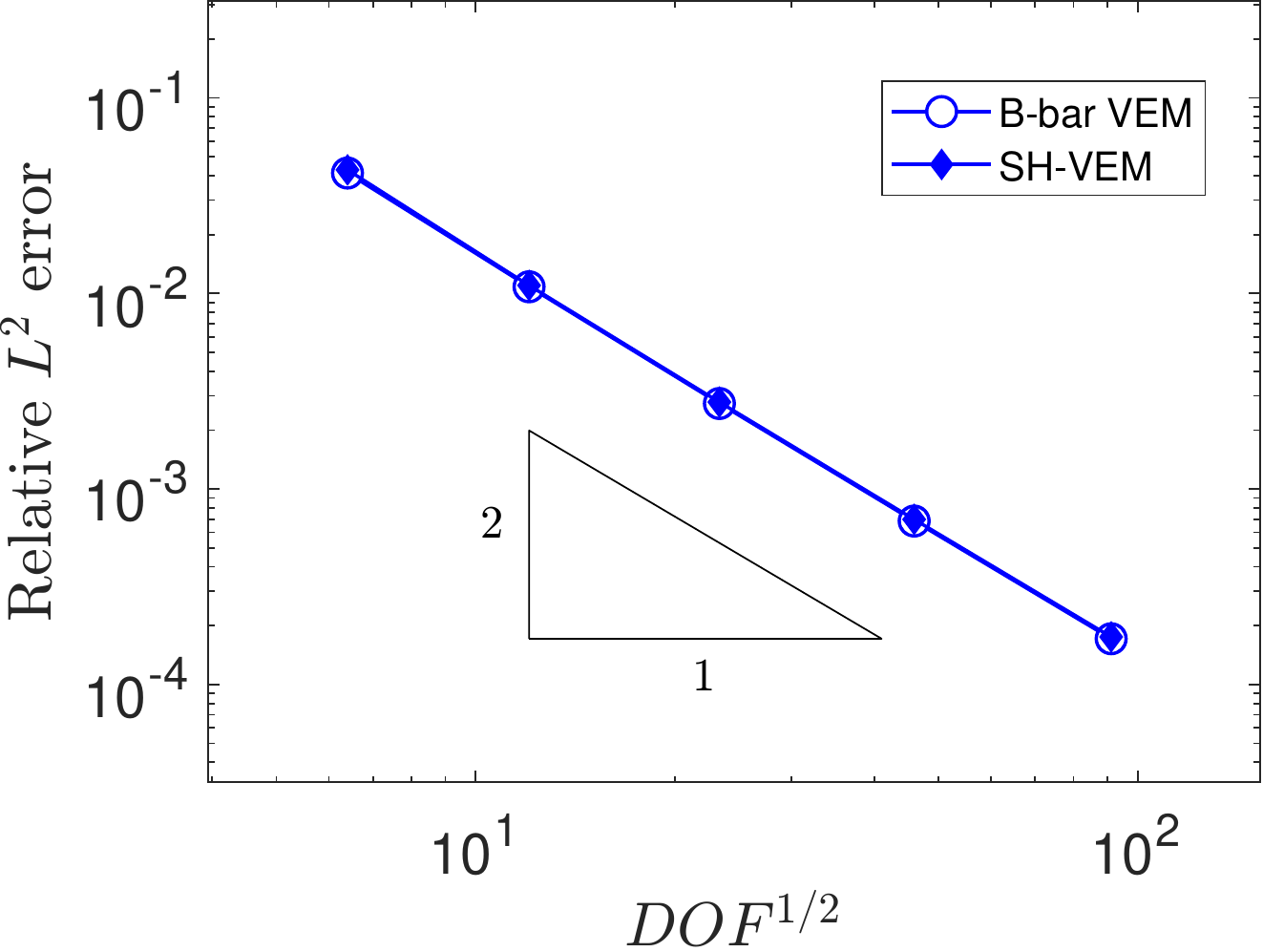}
         \caption{}
     \end{subfigure}
     \hfill
     \begin{subfigure}{0.32\textwidth}
         \centering
         \includegraphics[width=\textwidth]{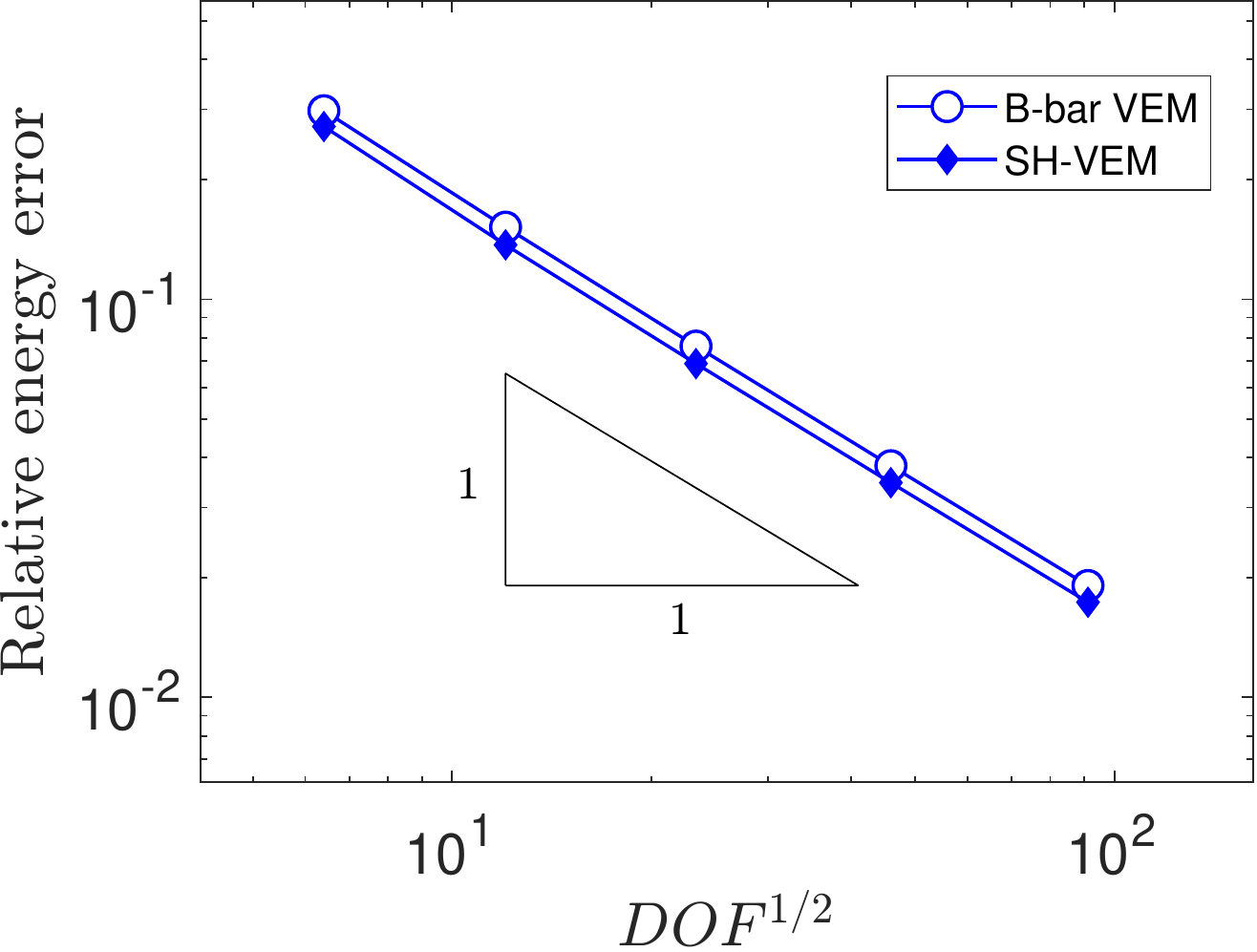}
         \caption{}
     \end{subfigure}
     \hfill
     \begin{subfigure}{0.32\textwidth}
         \centering
         \includegraphics[width=\textwidth]{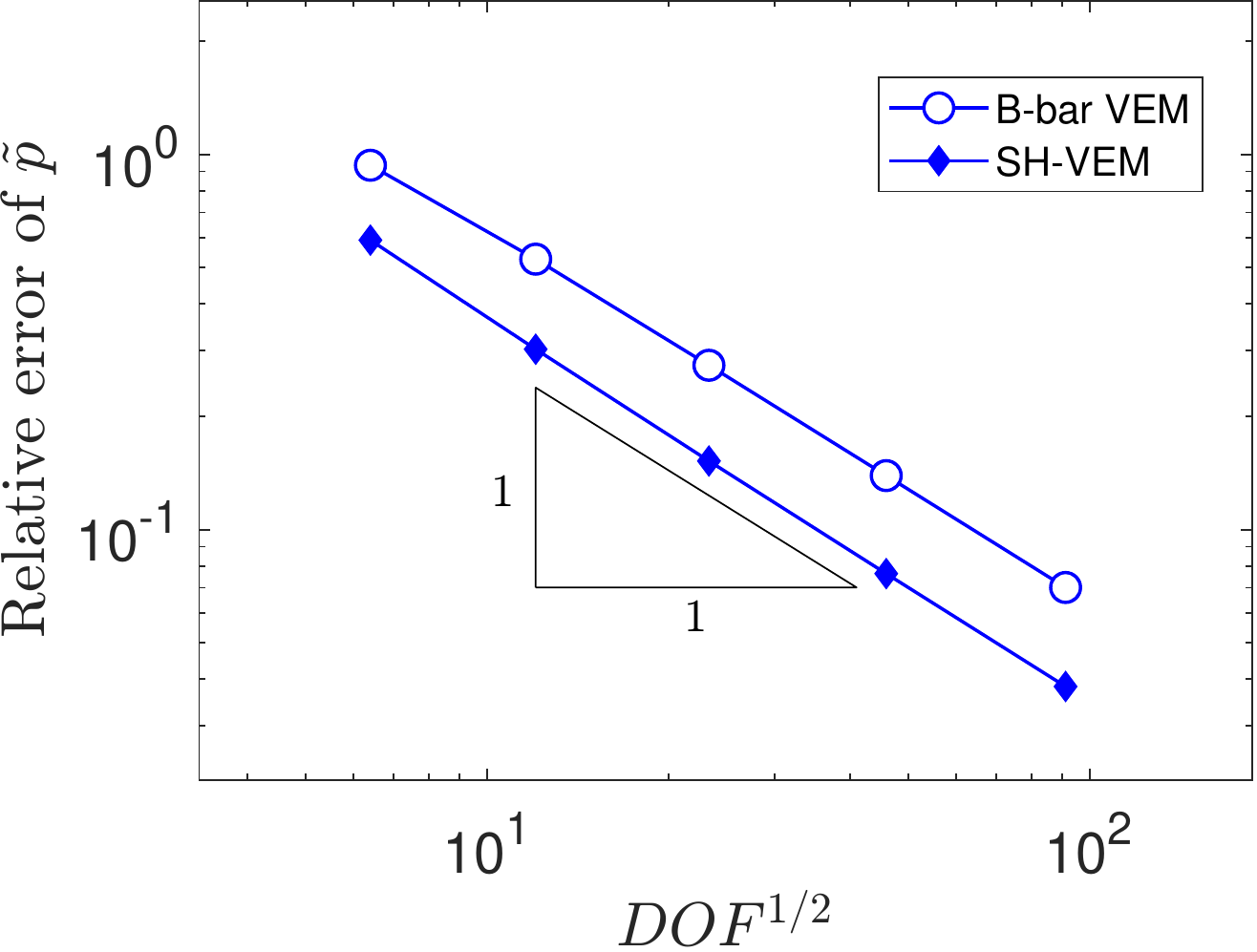}
         \caption{}
     \end{subfigure}
        \caption{Comparison of B-bar VEM and SH-VEM for the hollow cylinder problem on nonconvex meshes \acrev{(see~\fref{fig:cylinder_nonconvex})}. (a) $L^2$ error of
        displacement, (b) energy error and (c) $L^2$ error of hydrostatic stress. }
        \label{fig:cylinder_nonconvex_errors}
\end{figure}

\begin{figure}[!h]
    \centering
     \begin{subfigure}{0.48\textwidth}
         \centering
         \includegraphics[width=\textwidth]{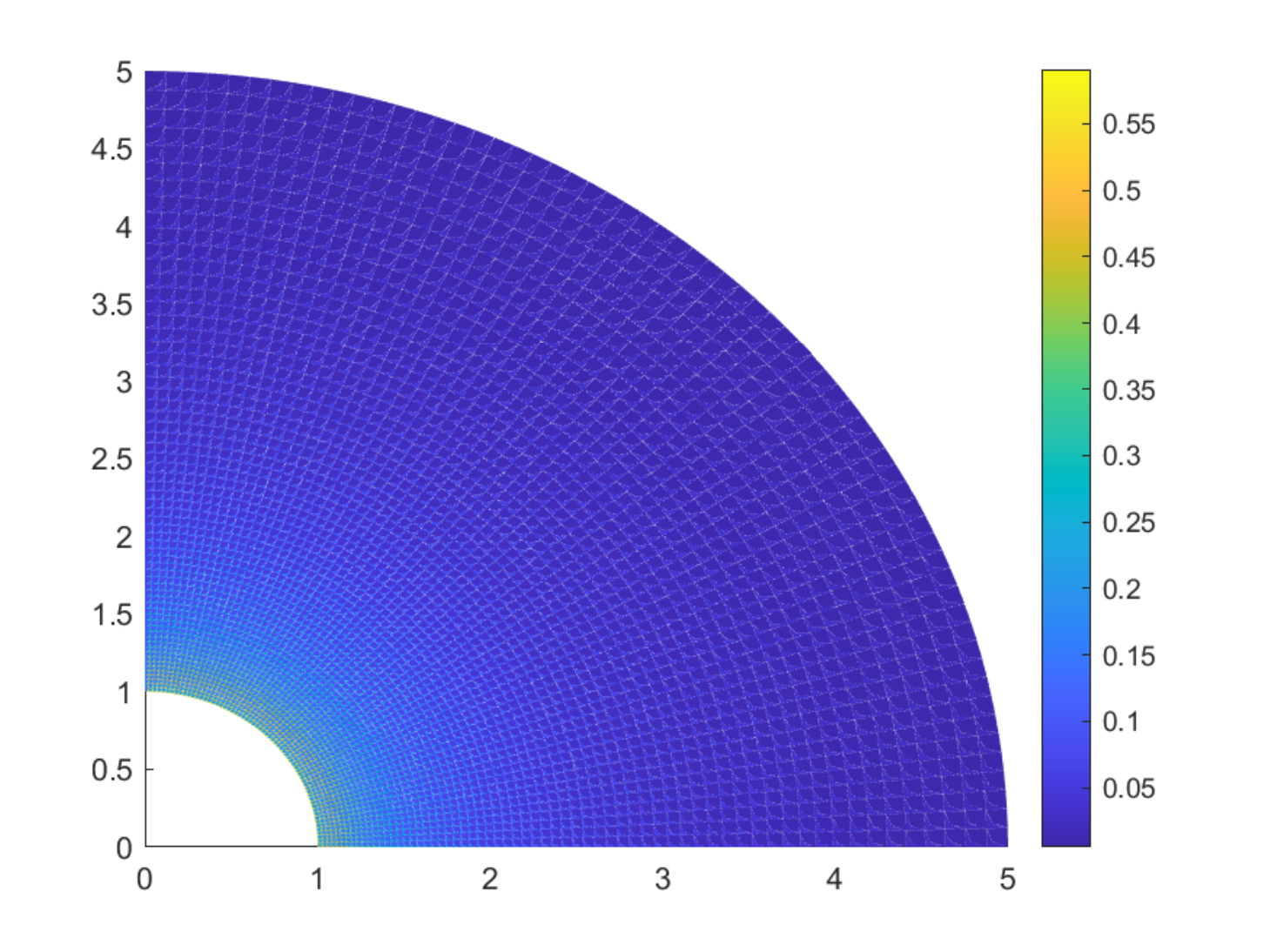}
         \caption{}\label{fig:cylinder_nonconvex_pressure-a}
     \end{subfigure}
     \hfill
    \begin{subfigure}{0.48\textwidth}
        \centering
        \includegraphics[width=\textwidth]{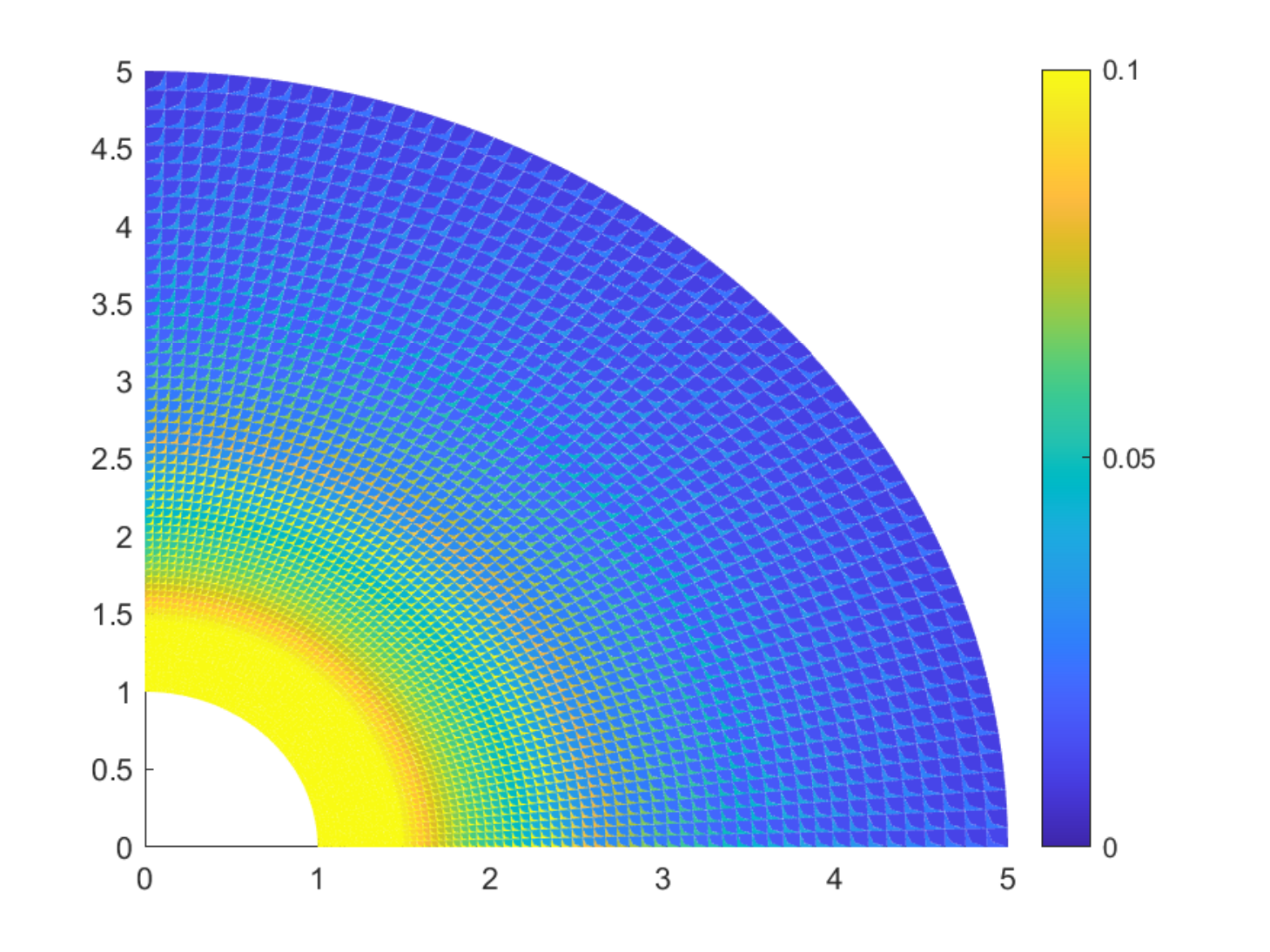}
        \caption{}\label{fig:cylinder_nonconvex_pressure-b}
    \end{subfigure}
    \vfill
    \begin{subfigure}{0.48\textwidth}
         \centering
         \includegraphics[width=\textwidth]{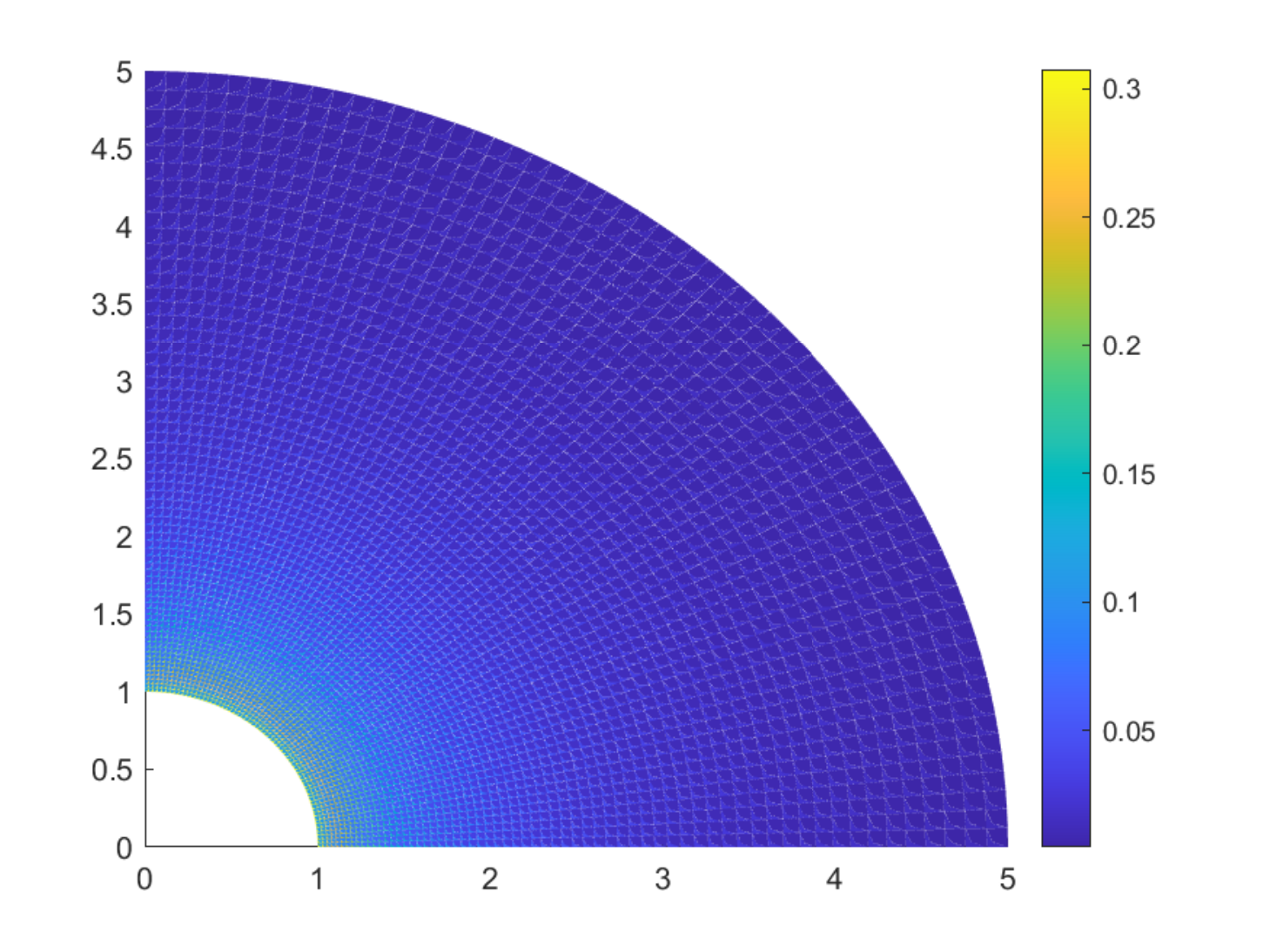}
         \caption{}\label{fig:cylinder_nonconvex_pressure-c}
     \end{subfigure}
     \hfill
    \begin{subfigure}{0.48\textwidth}
        \centering
        \includegraphics[width=\textwidth]{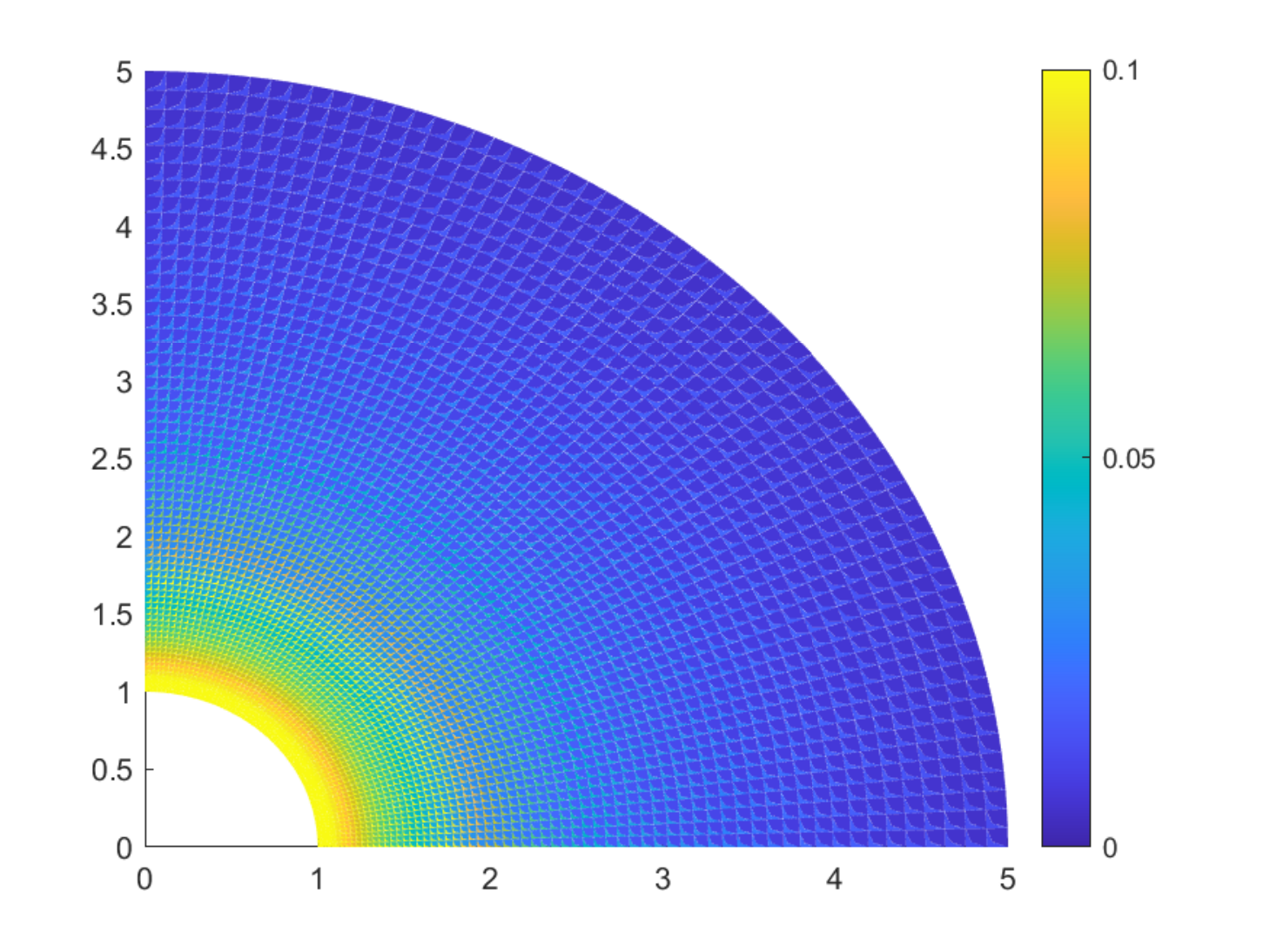}
        \caption{}\label{fig:cylinder_nonconvex_pressure-d}
    \end{subfigure}
    \caption{Contour plots of the relative error in the hydrostatic stress on nonconvex meshes \acrev{(see~\fref{fig:cylinder_nonconvex})} for the hollow pressurized cylinder problem. The exact hydrostatic stress is 4166.6666 psi. 
    (a) B-bar VEM, (b) B-bar VEM (color scale for error is between 0 and 10 percent),
    (c) SH-VEM and (d) SH-VEM 
    (color scale for error is between 0 and 10 percent.}
    \label{fig:cylinder_nonconvex_pressure}
\end{figure}

Now we solve the pressurized cylinder problem on a sequence of nonconvex meshes; a few representation meshes are
shown in Figure~\ref{fig:cylinder_nonconvex}. Figure~\ref{fig:cylinder_nonconvex_errors} shows that both the B-bar VEM and the SH-VEM
deliver optimal convergence rates; however unlike the uniform mesh case, the hydrostatic stress
field is not exactly reproduced by either method. In Figure~\ref{fig:cylinder_nonconvex_pressure}, we compare the contour plots of the error in the 
hydrostatic stress field for the two methods. We observe
that both methods
are very accurate away from the inner circular boundary but produce much larger errors in its vicinity 
(see 
Figures~\ref{fig:cylinder_nonconvex_pressure-b} and~\ref{fig:cylinder_nonconvex_pressure-d}). The maximum
error of the SH-VEM is 30 percent, whereas 
that of B-bar VEM is markedly worse at 55 
percent. Compared to Figure~\ref{fig:cylinder_nonconvex-c}, if
the nonconvex quadrilateral is distorted even more, we find from our simulations that the
maximum error in the hydrostatic stress for SH-VEM
increases to 35 percent, whereas the maximum error using
B-bar VEM has a 10-fold increase. 

\subsection{Flat punch}
\begin{figure}[!h]
    \centering
\begin{tikzpicture}[scale=6]

            %draws the punch main shape
			\filldraw[line width=1.5pt,fill=gray!10] (0,0) -- ++(1,0) -- ++(0,1.05) -- ++(-1, 0) -- cycle;
   
            \draw[line width=1.5pt] (-.05,1.05) -- ++(0,-1.1) -- ++(1.1,0) -- ++(0,1.1);

            %circles on right
            \filldraw[line width=1.5pt,fill=gray]  (-.025,0.035) circle (.025) ++ (0,.1) circle (.025) ++ (0,.1) circle (.025) ++ (0,.1) circle (.025) ++ (0,.1) circle (.025) ++ (0,.1) circle (.025) ++ (0,.1) circle (.025) ++ (0,.1) circle (.025) ++ (0,.1) circle (.025) ++ (0,.1) circle (.025) ++ (0,0.1) circle (.025)  ;

            %circles on bottom
            \filldraw[line width=1.5pt,fill=gray]  (0.005,-.025) circle (.025) ++ (.1,0) circle (.025) ++ (0.1,0) circle (.025) ++ (0.1,0) circle (.025) ++ (0.1,0) circle (.025) ++ (0.1,0) circle (.025) ++ (0.1,0) circle (.025) ++ (0.1,0) circle (.025) ++ (0.1,0) circle (.025) ++ (0.1,0) circle (.025) ++ (0.1,0) circle (.025) ;

            %circles on left
            \filldraw[line width=1.5pt,fill=gray]  (1.025,0.035) circle (.025) ++ (0,0.1) circle (.025) ++ (0,0.1) circle (.025) ++ (0,0.1) circle (.025) ++ (0,0.1) circle (.025) ++ (0,0.1) circle (.025) ++ (0,0.1) circle (.025) ++ (0,0.1) circle (.025) ++ (0,0.1) circle (.025) ++ (0,0.1) circle (.025)  ++ (0,0.1) circle (.025) ;

            %arrows on top
            \draw[line width=1.5pt ] (1/3,1.2) -- ++(1/3,0); 
            \draw[line width=1.5pt , ->] (1/3,1.2) -- ++(0,-.15); 
            \draw[line width=1.5pt , ->] (1/3+1/12,1.2) -- ++(0,-.15); 
            \draw[line width=1.5pt , ->] (1/3+2/12,1.2) -- ++(0,-.15); 
            \draw[line width=1.5pt , ->] (1/3+3/12,1.2) -- ++(0,-.15); 
            \draw[line width=1.5pt , ->] (2/3,1.2) -- ++(0,-.15); 

            %draws the wall on the left hand side
			\draw[] (-.05,1.05)--++(-.04,-.04) ++ (.04,-.04)--++(-.04,-.04) ++ (.04,-.04)--++(-.04,-.04)++ (.04,-.04)--++(-.04,-.04)++ (.04,-.04)--++(-.04,-.04)++ (.04,-.04)--++(-.04,-.04)++ (.04,-.04)--++(-.04,-.04)++ (.04,-.04)--++(-.04,-.04)++ (.04,-.04)--++(-.04,-.04)++ (.04,-.04)--++(-.04,-.04)++ (.04,-.04)--++(-.04,-.04)++ (.04,-.04)--++(-.04,-.04) ++ (.04,-.04)--++(-.04,-.04) ++ (.04,-.04)--++(-.04,-.04)
			; 

            %draws wall on bottom
            \draw[] (0,-.05) --++(-.04,-.04) ++ (.12,.04)--++(-.04,-.04)++ (.12,.04)--++(-.04,-.04)++ (.12,.04)--++(-.04,-.04)++ (.12,.04)--++(-.04,-.04)++ (.12,.04)--++(-.04,-.04)++ (.12,.04)--++(-.04,-.04)++ (.12,.04)--++(-.04,-.04)++ (.12,.04)--++(-.04,-.04)++ (.12,.04)--++(-.04,-.04)++ (.12,.04)--++(-.04,-.04)++ (.12,.04)--++(-.04,-.04)++ (.12,.04)--++(-.04,-.04)++ (.12,.04)--++(-.04,-.04);

            %draws wall on right side 
			\draw[] (1.05,1.05)--++(+.04,-.04) ++(-.04,-.04)--++(+.04,-.04)++(-.04,-.04)--++(+.04,-.04)++(-.04,-.04)--++(+.04,-.04)++(-.04,-.04)--++(+.04,-.04)++(-.04,-.04)--++(+.04,-.04)++(-.04,-.04)--++(+.04,-.04)++(-.04,-.04)--++(+.04,-.04)++(-.04,-.04)--++(+.04,-.04)++(-.04,-.04)--++(+.04,-.04)++(-.04,-.04)--++(+.04,-.04)++(-.04,-.04)--++(+.04,-.04)++(-.04,-.04)--++(+.04,-.04)++(-.04,-.04)--++(+.04,-.04);

            %draw bottom labels
            \draw[] (0,-.15) --++ (0,-.05);
            \draw[] (1,-.15) --++ (0,-.05);
            
            \draw[latex-latex] (0,-.175) --++ (1,0) ;
            \draw[] (.5,-.175) node[anchor = north]{1 };

            %draw side labels
            \draw[] (-.15,1.05) --++ (-.05,0);
            \draw[] (-.15,0) --++ (-.05,0);
            \draw[latex-latex] (-.175,1.05) --++ (0,-1.05) ;
            \draw[] (-.175,.525) node[anchor = east]{1 };

            %draw label for arrows
            \draw[] (1/3,.975) --++ (0,-.05);
            \draw[] (2/3,.975) --++ (0,-.05);
            \draw[latex-latex] (1/3,.95) --++ (1/3,0);
            \draw[] (.5,.95) node[anchor = north]{$\frac{1}{3}$};
            \draw[] (.5,1.2) node[anchor = south]{$v=-0.03$};
        
\end{tikzpicture}
    \caption{ Flat punch problem. }
    \label{fig:punch_bc}
\end{figure}
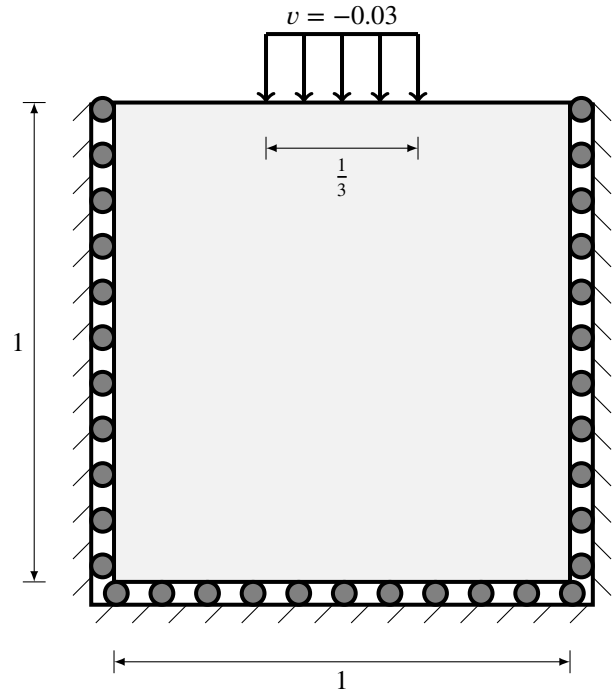

\begin{figure}[!h]
    \centering
    \begin{subfigure}{0.48\textwidth}
    \centering
    \includegraphics[width=\textwidth]{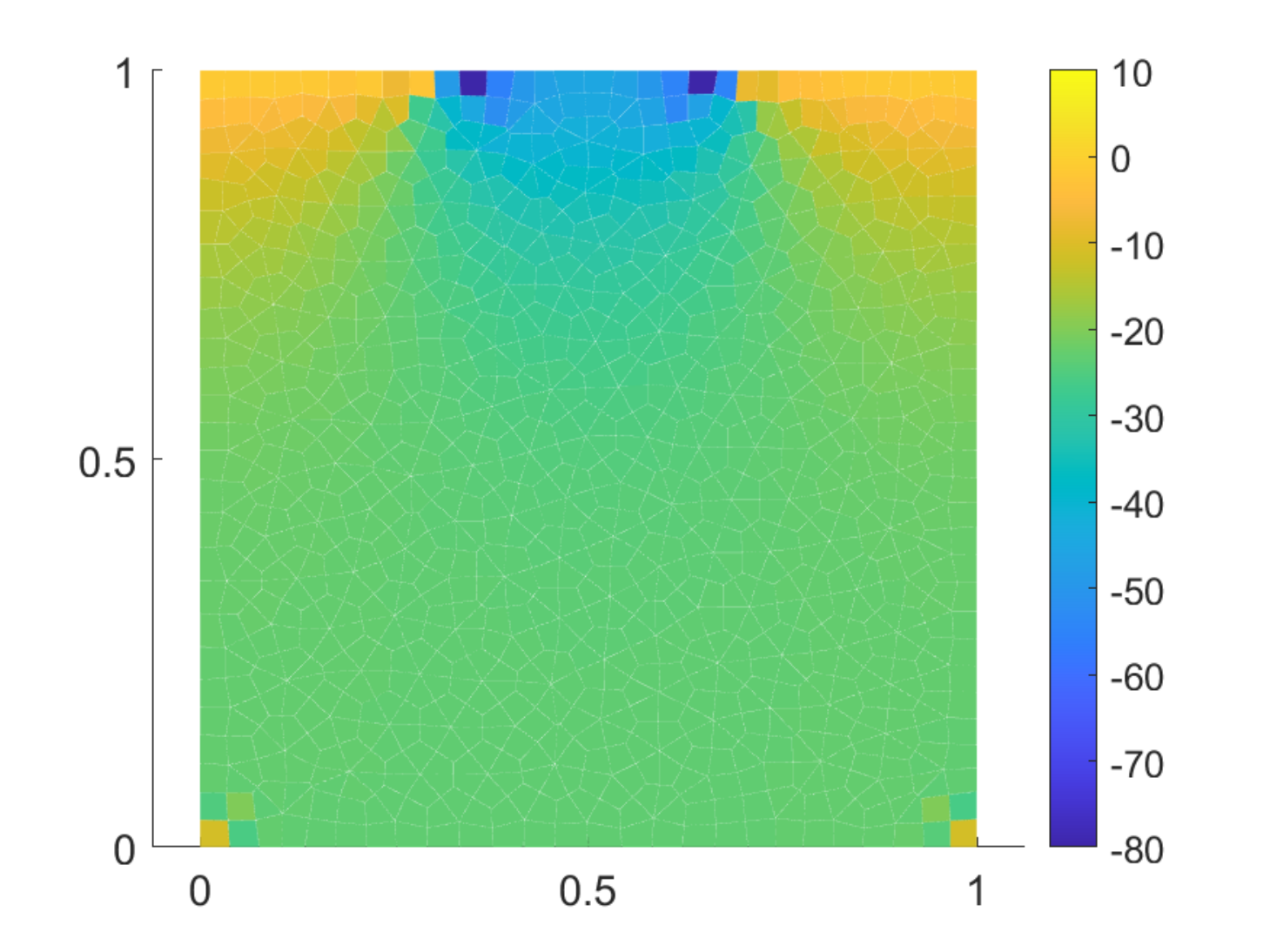}
    \caption{}
    \end{subfigure}
     \hfill
     \begin{subfigure}{0.48\textwidth}
         \centering
         \includegraphics[width=\textwidth]{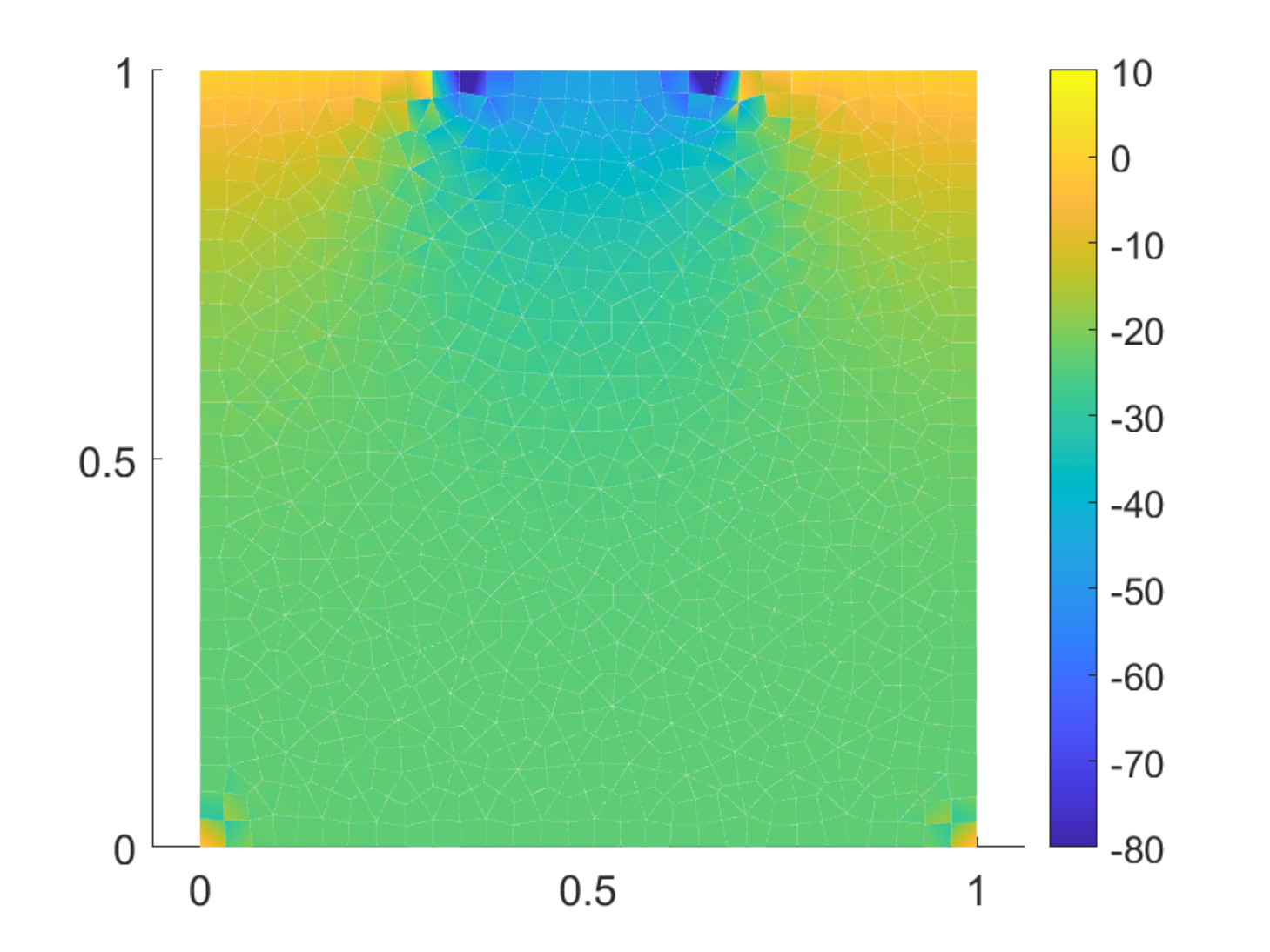}
         \caption{}
     \end{subfigure}
    \caption{Contour plots of the hydrostatic stress on unstructured meshes \acrev{(see~\fref{fig:squaremesh})} for the flat punch problem. (a) B-bar VEM and
    (b) SH-VEM. }
    \label{fig:punch_pressure}
\end{figure}

\begin{figure}[!h]
    \centering
    \begin{subfigure}{0.48\textwidth}
    \centering
    \includegraphics[width=\textwidth]{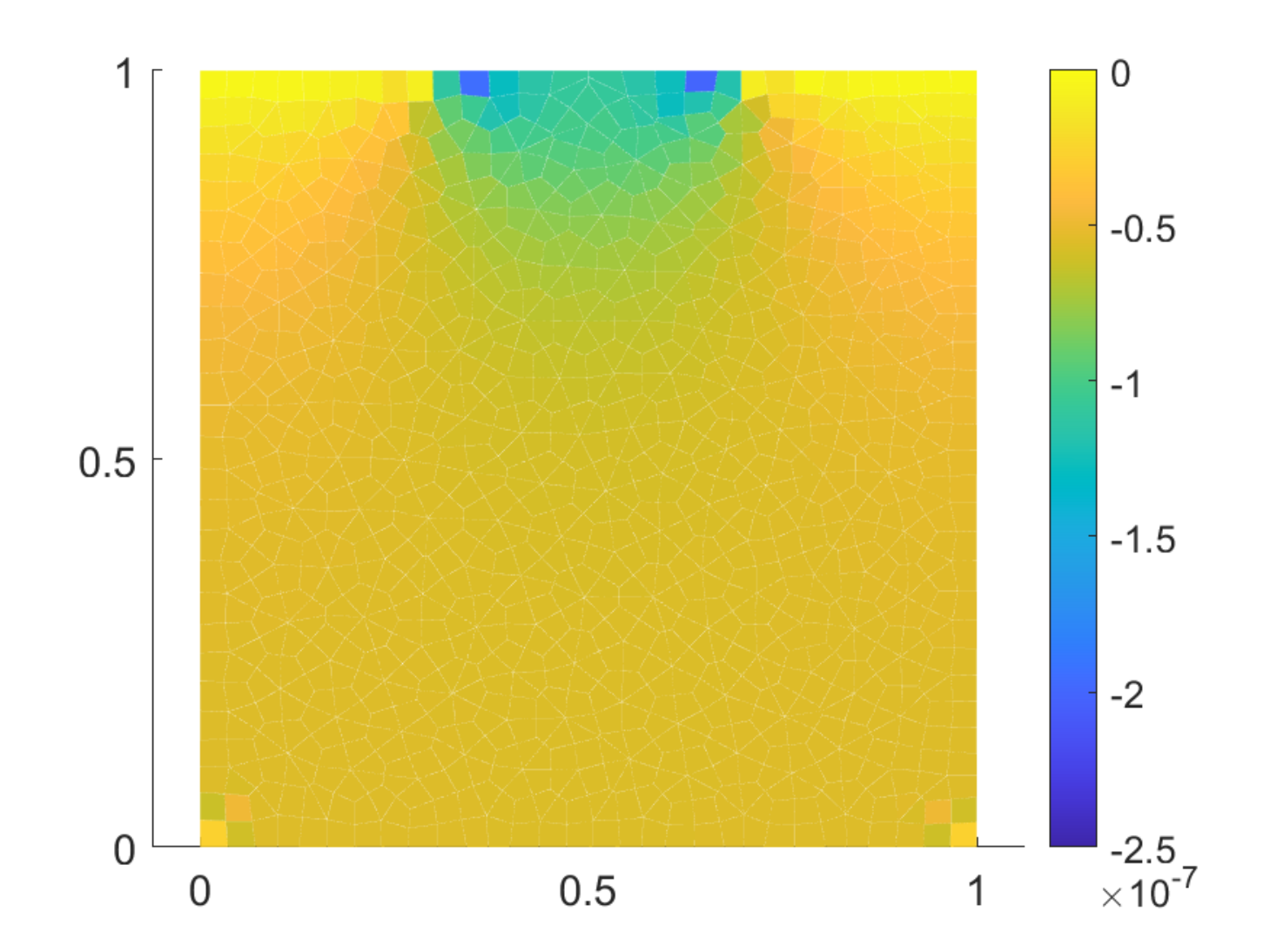}
    \caption{}
    \end{subfigure}
     \hfill
     \begin{subfigure}{0.48\textwidth}
         \centering
         \includegraphics[width=\textwidth]{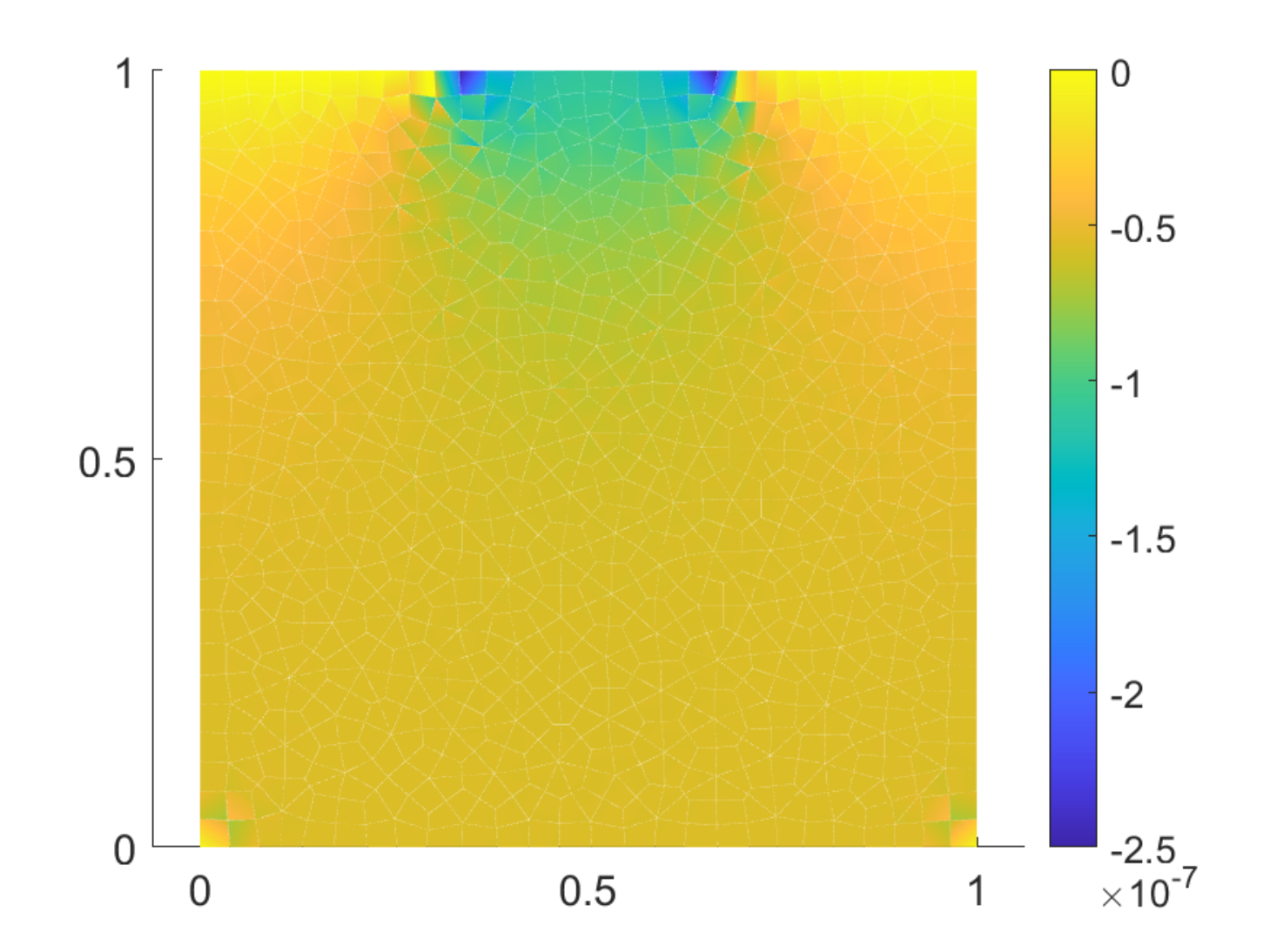}
         \caption{}
     \end{subfigure}
    \caption{Contour plots of the trace of the strain field on unstructured meshes \acrev{(see~\fref{fig:squaremesh})} for the flat punch problem. (a) B-bar VEM and
    (b) SH-VEM. }
    \label{fig:punch_trace}
\end{figure}

Finally, we consider the problem of a flat punch as described in Park et al.~\cite{Park:2020:meccanica} and shown in Figure~\ref{fig:punch_bc}. The domain is the unit square and we choose 
$E_Y=250$ psi and $\nu=0.4999999$. The left, right and bottom edges are constrained in the direction normal to the edges, and the top has a constant vertical displacement of $v = -\,0.03$ applied on the middle third of the edge.  A sequence of unstructured quadrilateral (see Figure~\ref{fig:squaremesh}) is used to solve this problem. The hydrostatic stress field from
both methods are presented in Figure~\ref{fig:punch_pressure}. The plots show that both methods produce relatively smooth hydrostatic stress fields of comparable accuracy. 
In Figure~\ref{fig:punch_trace}, plots of 
the trace of the strain field are shown
for B-bar VEM and SH-VEM, and we
find that consistent with the exact solution the numerically computed strain field is nearly traceless. 

\section{Conclusions}\label{sec:conclusions}
\acrev{In this article, to treat nearly-incompressible materials
in linear
elasticity, we departed from the commonly used
assumed-strain approaches in finite element methods that rely
on the Hu--Washizu three-field variational
principle.\cite{Simo:1986:jam,Simo:1990:ijnme}} Instead, we revisited the assumed stress (or stress-hybrid) formulation that use the two-field Hellinger--Reissner variational principle.\cite{Pian:1984:ijnme} In so doing, we
proposed a stress-hybrid formulation~\cite{Pian:1984:ijnme} of the virtual element method on quadrilateral meshes
for problems in plane linear elasticity. In this approach,  the Hellinger--Reissner functional is used to define  weak imposition of equilibrium equations and the 
strain-displacement relations to determine 
a suitable projection operator for the stress. On each quadrilateral element,
we constructed a local coordinate system~\cite{Cook:1974:jsd,Cook:1975:APS} and used a 5-term divergence-free symmetric
tensor polynomial basis in the local coordinate system. The rotation matrix was then used to 
transform the stress ansatz to
the global Cartesian coordinates so that element stiffness matrix computations could be conducted directly on 
the physical (distorted) element. 
On applying the divergence theorem on each element and using the divergence-free basis, we were able to compute the matrix representation of the stress projection solely from the displacements on the boundary. This resulted in a displacement-based method that was computable using
the virtual element formulation. In the Appendix,
we showed that the proposed approach was equivalent to a stress-hybrid virtual element formulation that follows 
the recipe of Cook~\cite{Cook:1974:jsd} to transform the element stiffness matrix from local to global Cartesian coordinates. 
The SH-formulation was tested for stability, volumetric 
and shear 
locking, and convergence on several benchmark problems. From an element-eigenvalue analysis, we found that the proposed method
was rotationally invariant and remained
stable for a large class of convex and nonconvex elements without needing a stabilization term.  
For a manufactured test problem in the incompressible limit ($\nu \to 0.5$), we showed that the SH-VEM did not suffer from volumetric locking. From the bending of a thin beam and the bending in the Cook's membrane problem, we found that the method was not 
susceptible to shear locking. For a plate with a circular hole, the methods produced optimal convergence rates and
smooth hydrostatic stress fields for both convex and nonconvex meshes. For the pressurized cylinder,
optimal convergence rates in the $L^2$ norm and energy
seminorm of the displacement field were realized, and both the B-bar VEM and
the SH-VEM reproduced close to the exact
hydrostatic stress  
on uniform meshes. However, it was observed that the hydrostatic stress field using the B-bar VEM  and the
SH-VEM
on distorted {\it nonconvex meshes} produced  larger errors, with the latter being more accurate.  In the
problem of a flat punch, the B-bar VEM and 
the SH-VEM produced relatively smooth hydrostatic stress fields that 
were comparable and the strain field was pointwise nearly traceless. \ac{Two of the main advantages of the virtual element method over the finite element method is the ability to handle very general polytopal (polygonal and polyhedral) meshes that provide flexibility in meshing, and to design discrete spaces within each element that are tied to the underlying partial differential equation that is solved.  Since the VEM does not need explicit construction of polytopal basis functions it also allows for a more seamless extension to high-order methods and to higher-dimensional problems.} \ac{Therefore, as future work, 
we plan to extend the stress-hybrid virtual element
method to polytopal meshes for 
compressible and
nearly-incompressible linear elasticity.}

\section*{Acknowledgements}
\ack{The authors acknowledge the research support
of Sandia National Laboratories to the University of 
California at Davis. NS is grateful to 
Professor Chandrashekhar Jog (Department of 
Mechanical Engineering, IISc Bangalore)
for many helpful discussions on
the stress-hybrid finite element method. These discussions in 
July 2022 provided the seeds to pursue the present contribution.
}

\section*{Data Availability Statement}
Data sharing not applicable to this article as no datasets were generated or analysed during the current study.

\appendix
\section{Stress-hybrid formulation based on
Cook's approach}
In this Appendix, we present an alternate formulation of the stress-hybrid virtual element method based on defining the element stiffness matrix on a rotated element as introduced by Cook.\cite{Cook:1974:jsd} Let $\Erot$ be a rotated element, and following~\eqref{eq:H_L_beta}, define the corresponding matrices $\vm{H}^{\prime}$ and $\vm{L}^{\prime}$ by 
\begin{subequations}
     \begin{align}
     \vm{H}^{\prime} &= \int_{\Erot}{(\vm{P}^{\prime})^T \vm{C}^{-1}\vm{P}^{\prime} \, d\vx^{\prime}}, \quad
     \vm{L}^{\prime} = \int_{\partial \Erot}{(\vm{P}^{\prime})^T\vm{N}^{\partial \Erot}\vm{\varphi}^{\prime} \, ds^{\prime}},\\
     \intertext{where $\vm{P}^{\prime}$ is given in~\eqref{eq:Pprime_matrix} and $\vm{\varphi}^{\prime}$ are the virtual element basis functions on $\Erot$. We then solve for the stress coefficients $\vm{\beta}^{\prime}$ in terms of the rotated displacements using}
     \vm{\beta}^{\prime} &= (\vm{H}^{\prime})^{-1}L^{\prime}\vd^\prime := \vm{\Pi}_{\beta}^{\prime}\vd^\prime.\label{eq:beta_displacement_prime}
 \end{align} 
\end{subequations}
The element stiffness matrix on the rotated element $\Erot$ is given as
 \begin{equation}
     \vm{K}^{\prime}_{\Erot} = (\vm{\Pi}^{\prime}_\beta)^T\left(\int_{\Erot}{ (\vm{P}^{\prime})^T\vm{C}^{-1}\vm{P}^{\prime}  \, d\vx^{\prime}}\right)\vm{\Pi}^{\prime}_\beta = (\vm{\Pi}^{\prime}_\beta)^T \vm{H}^{\prime} \vm{\Pi}^{\prime}_\beta, 
 \end{equation}
 and define the rotation matrix $\vm{R}$ as
  \begin{equation}\label{eq:rotation_matrix}
     \vm{R}=\begin{bmatrix}
     \vm{Q} & \vm{0} & \vm{0} & \vm{0}\\
     \vm{0} & \vm{Q} & \vm{0} & \vm{0}\\
     \vm{0} & \vm{0} & \vm{Q} & \vm{0}\\
     \vm{0} & \vm{0} & \vm{0} & \vm{Q}\\
     \end{bmatrix},
 \end{equation}
 where $\vm{Q}$ is given in~\eqref{eq:rotated_coords}, and $\vm{0}$ is the $2\times 2$ zero matrix. Then the element stiffness matrix
 in Cook's formulation on the original element $E$ is 
 recovered by:
 \begin{align}\label{eq:KEC}
     \vm{K}^{\textrm{C}}_E = \vm{R}^T\vm{K}^{\prime}_{\Erot}\vm{R}. 
 \end{align}
Now, on applying~\eqref{eq:beta_displacement_prime} and simplifying, we write the element stiffness matrix as
\begin{align*}
     \vm{K}^{\textrm{C}}_E = \vm{R}^T
     \bigl((\vm{H}^\prime)^{-1}\vm{L}^{\prime}\bigr)^{T}\vm{H}^{\prime}(\vm{H}^\prime)^{-1}\vm{L}^{\prime}\vm{R} = (\vm{L}^{\prime}\vm{R})^T (\vm{H}^{\prime})^{-1} (\vm{L}^{\prime}\vm{R}).
\end{align*}
We now show that the SH-VEM using the basis $\vm{P}^{*}$ in~\eqref{eq:P_matrix_full} is 
identical to Cook's formulation, i.e.,
$\vm{K}_E^* = \vm{K}_E^C$.
\begin{proof}
On expanding the element stiffness matrix of the SH-VEM given in~\eqref{eq:SH-VEM_stiffness} and simplifying, we get
\begin{align*}
    \vm{K}^*_E = \vm{\Pi}_\beta^T \vm{H}^* \vm{\Pi}_\beta = \bigl((\vm{H}^*)^{-1}\vm{L}^*\bigr)^T \vm{H}^* (\vm{H}^{-1}\vm{L}^*) = (\vm{L}^*)^T(\vm{H}^*)^{-1} \vm{L}^* .
\end{align*}
We first examine the matrix $\vm{H}^*$. From~\eqref{eq:HandL_matrix}, we have 
\begin{subequations}\label{eq:KESTAR}
  \begin{align*}
    \vm{H}^* = \int_{E}{(\vm{P}^{*})^T \vm{C}^{-1}\vm{P}^{*} \, d\vx},\\
    \intertext{and after multiplying out the matrices and using an equivalent tensor representation, we write the components of $\vm{H}^*$ as}
    \vm{H}^*_{ij} = \int_{E}{\vm{\mathcal{P}}^{*}_i:\mathbb{C}^{-1}:\vm{\mathcal{P}}^{*}_j \, d\vx},        
\end{align*}  
\end{subequations}
where $\vm{\mathcal{P}^{*}_i}$ is the tensor representation of the $i$-th column of $\vm{P}^{*}$. Using~\eqref{eq:rotated_P_tensor}, we rewrite this integral in terms of the rotated basis $\vm{\mathcal{P}}^{\prime}_i$, that is
\begin{align*}
    \vm{H}^*_{ij} = \int_{E^{\prime}}{\vm{Q}^T\vm{\mathcal{P}}^{\prime}_i\vm{Q}:\mathbb{C}^{-1}:\vm{Q}^T\vm{\mathcal{P}}^{\prime}_j\vm{Q} \, d\vx^{\prime}}.
\end{align*}
It can be shown that for an isotropic material modulli tensor $\mathbb{C}$, that 
\begin{align*}
    \vm{Q}^T\vm{\mathcal{P}}^{\prime}_i\vm{Q}:\mathbb{C}^{-1}:\vm{Q}^T\vm{\mathcal{P}}^{\prime}_j\vm{Q} = \vm{\mathcal{P}}^{\prime}_i:\vm{Q}^T\vm{Q}^T\mathbb{C}^{-1}\vm{Q}\vm{Q}:\vm{\mathcal{P}}^{\prime}_j =\vm{\mathcal{P}}^{\prime}_i:\mathbb{C}^{-1}:\vm{\mathcal{P}}^{\prime}_j.
\end{align*}
Therefore, we now have for all $i,j=1,2,\dots 5$:
\begin{equation}\label{eq:Hstar}
    \vm{H}^*_{ij} = \int_{E^{\prime}}{\vm{\mathcal{P}}^{\prime}_i:\mathbb{C}^{-1}:\vm{\mathcal{P}}^{\prime}_j \, d\vx^{\prime}} = \vm{H}^{\prime}_{ij}.
\end{equation}
Next, we examine the matrix $\vm{L}^*$. From~\eqref{eq:HandL_matrix}, we have 
\begin{align*}
    \vm{L}^* = \int_{\partial E}{(\vm{P}^{*})^T\Nmatrix\vm{\varphi} \, ds},
\end{align*}
After converting to an equivalent tensor representation, we write the components of $\vm{L}^*$ as:
\begin{align*}
    \vm{L}^*_{ij} = \int_{\partial E}{(\vm{\mathcal{P}_i^{*}\cdot \vm{n}})\cdot \vm{\varphi}_{j} \, ds}.
\end{align*}
Since $\vm{\varphi}_j$ and $\vm{\varphi}^{\prime}_j$ are both piecewise
affine functions on $\partial E$ and $\partial E^\prime$, respectively,
it can be shown that the integration of $\vm{\varphi}_j$ along the boundary of an element $E$ is equivalent to integrating $\vm{\varphi}^{\prime}_j$ along the boundary of the rotated element $\Erot$. That is, for any vector field $\vm{f}$, we have
\begin{align*}
    \int_{\partial E} {\vm{f}(\vm{x})\cdot \vm{\varphi}_j\, ds} = \int_{\partial \Erot}{\vm{f}(\vx(\vx^{\prime}))\cdot\vm{\varphi}^{\prime}_j \, ds^{\prime}}.
\end{align*}
With this, we rewrite $\vm{L}^*_{ij}$ in the rotated coordinates as
\begin{align*}
    \vm{L}^*_{ij} =\int_{\partial E}{(\vm{\mathcal{P}_i^{*}\cdot \vm{n}})\cdot \vm{\varphi}_{j} \, ds} = \int_{\partial \Erot}{(\vm{Q}^T\vm{\mathcal{P}}_i^{\prime}\vm{Q}\cdot \vm{Q}^T\vm{n}^{\prime})\cdot \vm{\varphi}^{\prime}_{j} \, ds^{\prime}}= \int_{\partial \Erot}{(\vm{Q}^T\vm{\mathcal{P}}_i^{\prime}\cdot\vm{n}^{\prime})\cdot \vm{\varphi}^{\prime}_{j} \, ds^{\prime}}.
\end{align*}
If we take the basis functions in the standard order $\vm{\varphi}^{\prime}_{2j-1} = (\phi^{\prime}_j, 0)^T $ and $\vm{\varphi}^{\prime}_{2j} = (0,\phi_j^{\prime})^T $, then we can simplify $\vm{L}^*_{ij}$ as:
    \begin{align*}
    \vm{L}^*_{i2j-1} &= c\int_{\partial \Erot}{((\vm{\mathcal{P}}^{\prime}_i)_{11}\vm{n}_1^{\prime}+(\vm{\mathcal{P}}^{\prime}_i)_{12}\vm{n}_2^{\prime})\phi^{\prime}_j \, ds^\prime }-s\int_{\partial \Erot}{((\vm{\mathcal{P}}^{\prime}_i)_{12}\vm{n}_1^{\prime}+(\vm{\mathcal{P}}^{\prime}_i)_{22}\vm{n}_2^{\prime})\phi^{\prime}_j \, ds^\prime } = c\vm{L}^{\prime}_{i2j-1} - s\vm{L}^{\prime}_{i2j}\\ 
    \vm{L}^*_{i2j} &=  s\int_{\partial \Erot}{((\vm{\mathcal{P}}^{\prime}_i)_{11}\vm{n}_1^{\prime}+(\vm{\mathcal{P}}^{\prime}_i)_{12}\vm{n}_2^{\prime})\phi^{\prime}_j \, ds^\prime } +c \int_{\partial \Erot}{((\vm{\mathcal{P}}^{\prime}_i)_{12}\vm{n}_1^{\prime}+(\vm{\mathcal{P}}^{\prime}_i)_{22}\vm{n}_2^{\prime})\phi^{\prime}_j \, ds^\prime } =  s\vm{L}^{\prime}_{i2j-1} +c \vm{L}^{\prime}_{i2j},
\end{align*}
where $c$ and $s$ given in~\eqref{eq:rotated_coords}. On multiplying out the matrix $\vm{L}^{\prime}\vm{R}$, it can be shown that 
\begin{align*}
    (\vm{L}^{\prime}\vm{R})_{i2j-1} &=  c\vm{L}^{\prime}_{i2j-1} - s\vm{L}^{\prime}_{i2j} \\
    (\vm{L}^{\prime}\vm{R})_{i2j} &=s\vm{L}^{\prime}_{i2j-1} +c \vm{L}^{\prime}_{i2j}, 
\end{align*}
and therefore for all $i=1,2,\dots 5$ and $j=1,2,\dots 4$, we have
\begin{subequations}\label{eq:Lstar}
\begin{align}
    \vm{L}^*_{i2j-1} &= c\vm{L}^{\prime}_{i2j-1} - s\vm{L}^{\prime}_{i2j} = (\vm{L}^{\prime}\vm{R})_{i2j-1} \\
    \vm{L}^*_{i2j} &=  s\vm{L}^{\prime}_{i2j-1} +c \vm{L}^{\prime}_{i2j} =  (\vm{L}^{\prime}\vm{R})_{i2j}.  
\end{align}
\end{subequations}
From~\eqref{eq:Hstar} and~\eqref{eq:Lstar}, we 
obtain $\vm{H}^* = \vm{H}^{\prime}$ and $\vm{L}^* = \vm{L}^{\prime}\vm{R}$. On substituting these 
in~\eqref{eq:KESTAR} and 
using~\eqref{eq:KEC} leads us to the desired result:
\begin{equation*}
    \vm{K}^*_E =  (\vm{L}^*)^T(\vm{H}^*)^{-1} \vm{L}^* = (\vm{L}^{\prime}\vm{R})^T (\vm{H}^{\prime})^{-1}(\vm{L}^{\prime}\vm{R}) = \vm{K}^{\textrm{C}}_E.
    \qedhere
\end{equation*}
\end{proof}

\end{document}